\documentclass[11pt,reqno]{amsart}
\usepackage{amsthm,amssymb,amsmath}
\usepackage[T1]{fontenc}
\usepackage{lmodern}
\usepackage[foot]{amsaddr}

\usepackage{xcolor}

\numberwithin{equation}{section}

\newtheorem{theorem}{Theorem}[section]
\newtheorem*{theorem*}{Theorem}
\newtheorem{lemma}{Lemma}[section]

\newtheorem{proposition}{Proposition}[section]

\theoremstyle{definition}

\newtheorem{definition}{\sc Definition}[section]
\newtheorem*{definition*}{\sc Definition}

\newtheorem{remark}{\bf Remark}[section]
\newtheorem*{example*}{\bf Example}
\newtheorem*{examples}{\bf Examples}

\newcommand{\loc}{{\rm loc}}

\newcommand{\Real}{{\rm Re}}

\newcommand{\dist}{\mbox{dist}}

\newcommand{\clos}{{\rm clos}}

\begin{document}

\title{Form-boundedness and SDEs with singular drift}

\author{Damir Kinzebulatov}

\thanks{\today}

\begin{abstract}
We survey and refine recent results on weak and strong well-posedness of stochastic differential equations with singular drift satisfying some minimal assumptions.
\end{abstract}

\email{damir.kinzebulatov@mat.ulaval.ca}

\address{Universit\'{e} Laval, D\'{e}partement de math\'{e}matiques et de statistique, Qu\'{e}bec, QC, Canada}



\thanks{The research of the author is supported by the NSERC (grant RGPIN-2017-05567).}

\maketitle

\tableofcontents

\section{Introduction}

Various applications in physics and technology dictate the need to work with stochastic differential equations 
\begin{equation}
\label{sde1}
dX_t=-b(X_t)dt + \sqrt{2}dW_t, \quad X_0=x \in \mathbb R^d,
\end{equation}
having  an irregular, locally unbounded drift $b:\mathbb R^d \rightarrow \mathbb R^d$.
Here $\{W_t\}_{t \geq 0}$ is a $d$-dimensional Brownian motion in $\mathbb R^d$ defined on some complete filtered probability space $(\Omega,\mathcal F_t,\mathcal F,\mathbf P)$. This naturally leads to the problem of finding the least restrictive assumptions on $b$ that ensure  well-posedness of \eqref{sde1}, in one sense or another. More specifically, one asks: what integral characteristics of $b$ determine whether there exists a unique solution of \eqref{sde1}? The same question arises when one considers
more general SDEs, also dictated by applications:
\begin{equation}
\label{sde2}
dX_t=-b(t,X_t)dt + \sqrt{2}dW_t, 
\end{equation}
with drift $b:\mathbb R \times \mathbb R^d \rightarrow \mathbb R^d$ that can also be singular in time, and
\begin{equation}
\label{sde3}
dX_t=-b(t,X_t)dt + \sigma(t,X_t)dW_t
\end{equation}
with diffusion coefficients $\sigma:\mathbb R_+ \times \mathbb R^d \rightarrow \mathbb R^{d \times d}$ that can be discontinuous. 
Regarding SDE \eqref{sde2}, one illustrative example is the ``passive tracer model'' that 
describes the motion of a small particle in a turbulent flow, i.e.\,\eqref{sde2} with the velocity field $b$ obtained by solving three-dimensional Navier-Stokes equations \cite{MK}. 

\medskip

The paper deals with weak and strong well-posedness  of SDEs \eqref{sde1}-\eqref{sde3}, for every initial point $x \in \mathbb R^d$. Recall that a weak solution of \eqref{sde1}-\eqref{sde3} is a pair of continuous processes $\{(X_t,W_t)\}_{t \geq 0}$ defined on some complete probability space, such that $\{W_t\}$ is a Brownian motion and the identity in \eqref{sde1}-\eqref{sde3} holds a.s.\,for all $t \geq 0$. In turn, a strong solution of \eqref{sde1}-\eqref{sde3} is a continuous process $X_t$ that is adapted to the natural filtration of the Brownian motion $\{W_t\}$, and such that the identity in \eqref{sde1}-\eqref{sde3} holds a.s.\,for $t \geq 0$. That is, $X_t$ is a strong solution if it only depends on time and the driving process $\{W_s\}_{0 \leq s \leq t}$.

\medskip

The question of what local singularities of $b$ are admissible, so that SDEs \eqref{sde1}-\eqref{sde3} are weakly or strongly well-posed, was thoroughly studied in the literature. Below we give a brief outline of the results on multidimensional SDEs with the focus on the singularities of the drift.
We will keep the chronological order of appearance of preprints, where applicable. However, we will be somewhat loose with the terminology by including in ``well-posedness''  uniqueness results of varying strength (in general, within some large classes of solutions).

\medskip

Veretennikov \cite{V} was the first who established, using Zvonkin's method \cite{Zv}, strong well-posedness of \eqref{sde2} when $|b|$ is bounded measurable.
 Portenko \cite{P} considered drift $b$ in the sub-critical Ladyzhenskaya-Prodi-Serrin class
\begin{equation}
\label{L}
|b| \in L^l(\mathbb R,L^p(\mathbb R^d)), \quad p \geq d, l \geq 2, \quad \frac{d}{p}+\frac{2}{l} < 1
\end{equation}
and proved existence of weak solution to SDE \eqref{sde2} and its uniqueness in law.
Krylov-R\"{o}ckner further established, using Yamada-Watanabe theorem, that for such $b$ the SDE \eqref{sde2} is, in fact, strongly well-posed. A number of important results for SDEs with drift satisfying \eqref{L} were established next by X.\,Zhang \cite{Z,Z0,Z4}. 
In the period between the papers of Portenko and Krylov-R\"{o}ckner, Bass-Chen \cite{BC} proved existence and uniqueness in law of weak solutions of \eqref{sde1} for time-homogeneous drift $b=b(x)$ in the Kato class of vector fields, with arbitrarily small $\delta$, cf.\,\eqref{kato}.
The Kato class of vector fields contains $\{|b| \in L^p(\mathbb R^d), p>d\}$ as well as some vector fields with entries not even in $L^{1+\varepsilon}_{\loc}(\mathbb R^d)$, $\varepsilon>0$. However, Kato class does not contain  $\{|b| \in L^d(\mathbb R^d)\}$. 
Speaking of time-homogeneous drifts, the fact that $p=d$ is the optimal exponent on the scale of Lebesgue spaces can be seen already from rescaling the parabolic equation.
In \cite{BFGM}, Beck-Flandoli-Gubinelli-Maurelli  developed an approach to establishing strong well-posedness of \eqref{sde2} with drift $b$ in the critical Ladyzhenskaya-Prodi-Serrin (LPS) class
\begin{equation}
\label{L2}
|b| \in L^l(\mathbb R,L^p(\mathbb R^d)), \quad p \geq d, l \geq 2, \quad \frac{d}{p}+\frac{2}{l} \leq 1.
\end{equation}
 for a.e.\,initial point $x \in \mathbb R^d$,
via stochastic transport and stochastic continuity equations. They also discussed the following attracting drift:
\begin{equation}
\label{b_hardy}
b(x)=\sqrt{\delta} \frac{d-2}{2}\mathbf{1}_{|x|<1}|x|^{-2}x
\end{equation}
(note that $|b|$ misses $L^d(\mathbb R^d)$ by a logarithmic factor) and provided a detailed proof of the fact that for 
\begin{equation}
\label{delta_4}
\delta>4(\frac{d}{d-2})^2,
\end{equation}
 i.e.\,when the attraction to the origin by the drift is strong enough, then SDE \eqref{sde1} with initial point $x=0$ does not have a weak solution.
In \cite{KiS_brownian}, Sem\"{e}nov and the author showed that the constructed in \cite{Ki_a_new_approach} Feller generator $-\Delta + b \cdot \nabla$ for $b$ satisfying weak form-boundedness condition \eqref{f2}, see below,
 determines, for every initial point $x \in \mathbb R^d$, a weak solution to SDE \eqref{sde1} that is unique among weak solutions that can be constructed via approximation. To the best of the author's knowledge, this was the first result on weak well-posedness of \eqref{sde1} that included both $|b| \in L^d(\mathbb R^d)$ and the model vector field \eqref{b_hardy} with $\delta$ small. It also included the elliptic Morrey class $M_{1+\varepsilon}$, see below, and the Kato class considered by Bass-Chen. The construction of the Feller generator with such $b$ used in an essential manner some inequalities for symmetric Markov generators, and hence required time-homogeneity of the drift.
Returning to time-inhomogeneous drifts, we note that almost at the same time Jin \cite{J} proved weak well-posedness of \eqref{sde2} with time-inhomogeneous Kato class drifts, Wei-Lv-Wu \cite{WLW} and Nam \cite{N} obtained results on weak well-posedness of \eqref{sde2}, for every $x \in \mathbb R^d$, for time-inhomogeneous vector fields $b$ that can be more singular than the ones in \eqref{L}. Nevertheless, their results excluded $b=b(x)$ with $|b| \in L^d(\mathbb R^d)$. 
In \cite{XXZZ}, Xia-Xie-Zhang-Zhao established, among other results, weak well-posedness for every initial point of SDE \eqref{sde1} with $b$ having entries in $C_b(\mathbb R,L^d(\mathbb R^d))$. 
In \cite{Kr1,Kr2,Kr3}, Krylov proved weak and strong well-posedness of SDEs \eqref{sde1}-\eqref{sde3} with $|b|$ in $L^d(\mathbb R^d)$ and beyond, e.g.\,in a large Morrey class (in fact, below we use an argument from these papers to prove some gradient bounds). 
In \cite{YZ}, S.\,Yang-T.\,Zhang proved strong well-posedness of \eqref{sde2} for time-inhomogeneous drifts $b$ with $|b|^2$ ``almost'' in the Kato class of potentials, cf.\,\eqref{kato_pt_def} (which, to make the comparison clear at least in the time-homogeneous case, is smaller than the Kato class of vector fields in \cite{BC}; still, the class considered by Yang-Zhang contains some drifts with $|b| \not \in L^{2+\varepsilon}_{\loc}$, $\varepsilon>0$).
In \cite{RZ}, R\"{o}ckner-Zhao  established weak well-posedness of \eqref{sde1}, for any initial point $x \in \mathbb R^d$, for drifts in $L^\infty(\mathbb R,L^{d,w}(\mathbb R^d))$, plus the drifts in the critical LPS class. Here $L^{d,w}(\mathbb R^d)$ denotes the weak $L^d$ class that contains vector fields in $L^d(\mathbb R^d)$, as well as more singular vector fields, such as \eqref{b_hardy}. In \cite{RZ2}, the authors obtained strong well-posedness of \eqref{sde1}, for any initial point, with $b$ in the critical LPS class. In \cite{KiM}, Madou and the author established weak well-posedness of \eqref{sde2}, for every initial point, for $b$ in the class of time-inhomogeneous form-bounded drifts (containing \eqref{f1} below). This class contains $L^\infty(\mathbb R,L^{d,w}(\mathbb R^d))$ and the critical LPS class,  as well as some drifts that are not even in $L^\infty(\mathbb R,L^{2+\varepsilon}(\mathbb R^d))$ for a given $\varepsilon>0$. In \cite{KiS_sharp}, Sem\"{e}nov and the author proved existence of a weak solution to SDE \eqref{sde2} with time-inhomogeneous form-bounded drift, see \eqref{f1}, covering the entire range of admissible form-bounds $\delta<4$, cf.\,\eqref{delta_4}. The critical value $\delta=4$ was recently attained at the PDE level in \cite{Ki_Orlicz}. In \cite{Ki_Morrey}, the author established weak well-posedness of SDE \eqref{sde2} and proved Feller property for time-inhomogeneous drifts in essentially the largest possible parabolic Morrey class, which contains the class of time-inhomogeneous form-bounded drifts, the time-homogeneous Morrey drifts in $M_{1+\varepsilon}$, and allows to include drifts $b$ having critical  behaviour both in time and in space, e.g.
$$
|b(t,x)| \leq \frac{c_1}{|x|}+\frac{c_2}{\sqrt{t}}, \quad t>0, \quad x \in \mathbb R^d.
$$
More recently, Krylov \cite{Kr_weak} established weak well-posedness of SDE \eqref{sde3} with discontinuous diffusion coefficients in the VMO class and time-inhomogeneous drifts in a large parabolic Morrey class; restricted to time-homogeneous drifts his assumptions read as $|b| \in M_{q}$, $q>\frac{d}{2}$. This result was further refined by Krylov in \cite{Kr_Morrey_sdes}, see also \cite{Kr_Morrey_parab} regarding regularity theory of parabolic equations with VMO diffusion coefficients and drift and potential in Morrey classes.

\medskip

Below we survey and refine recent results on weak and strong well-posedness of SDEs \eqref{sde1}-\eqref{sde3} with $|b|$ satisfying some minimal constraints, as in \cite{KiS_brownian,KiM,KiS_sharp, Ki_Morrey} mentioned above. For instance, regarding SDE \eqref{sde1},
our assumption on the order of singularities of $|b|$ is basically that $-\Delta + b \cdot \nabla$ must generate a quasi contraction  strongly continuous semigroup in $L^2$. 
 That is, we will be assuming that $|b|$ is form-bounded:
\begin{equation}
\label{f1}
|b|^2 \leq \delta(-\Delta) + c \quad (\text{in the sense of quadratic forms})
\end{equation}
for some constants $\delta$ and $c$. See Definition \ref{ellip_fb_def} below. This translates, by means of the Cauchy-Schwarz inequality, into  the assumption of smallness  of $b \cdot \nabla$ with respect to $-\Delta$. A broad elementary sufficient condition for \eqref{f1} is the scaling-invariant Morrey class $M_{2+\varepsilon}$, i.e.\,
\begin{equation}
\label{morrey_elliptic_2}
\|b\|_{M_{2+\varepsilon}}:=\sup_{r>0, x \in \mathbb R^d} r\biggl(\frac{1}{|B_r|}\int_{B_r(x)}|b|^{2+\varepsilon}dx \biggr)^{\frac{1}{2+\varepsilon}}<\infty
\end{equation}
where $B_r(x)$ is the ball of radius $r$ centered at $x$. Here $\varepsilon$ can be taken arbitrarily small. One has $\delta=C\|b\|_{M_{2+\varepsilon}}$ for appropriate constant $C$. The class $M_{2+\varepsilon}$ contains all $|b|$  in $L^d$ or, more generally, in the weak $L^d$ class (we recall its definition in Section \ref{notation_sect}). 

\medskip

Regarding the relationship between operator $-\Delta + b \cdot \nabla$ and SDE \eqref{sde1}, one expects that for $b=b(x)$ the function
\begin{equation}
\label{v}
v(t,x):=\mathbf E_{X_0=x}[f(X_t)], 
\end{equation}
solves Cauchy problem $$(\partial_t-\Delta + b \cdot \nabla)v=0, \quad v|_{t=0}=f.$$ 
The intimate relationship between parabolic equations and SDEs is a consequence of the fact that both describe the same physical process of diffusion. 

\medskip

In Section \ref{weak_sde_sect}, we discuss results on weak well-posedness of \eqref{sde1} under more general assumption on the drift than its form-boundedness. Namely, our hypothesis on the order of singularities of $|b|$ will be that $-\Delta + b \cdot \nabla$ generates a quasi contraction strongly continuous semigroup in the Bessel space $\mathcal W^{1/2,2}$, i.e.\,we will require
\begin{equation}
\label{f2}
|b| \leq \delta(\lambda-\Delta)^{\frac{1}{2}} \quad (\text{in the sense of quadratic forms}).
\end{equation}
Such vector fields $b$ are called weakly form-bounded. This class of weakly form-bounded vector fields contains scaling-invariant Morrey class $M_{1+\varepsilon}$. It also contains the Kato class of vector fields.

\medskip

One of our goals will be to bootstrap the semigroup in $L^2$ or in $\mathcal W^{\frac{1}{2},2}$ to a strongly continuous semigroup in $C_\infty$, the space of continuous functions on $\mathbb R^d$ vanishing at infinity, endowed with the $\sup$-norm (that is, to a Feller semigroup). This will come at the cost of restricting admissible values of constant $\delta$ (in terms of the Morrey class, this means that the Morrey norm cannot be too large). We emphasize that while the classes \eqref{f1}, \eqref{f2} determine the order of singularities of drift $b$, the value of $\delta$ measures the magnitude of its singularities. We are particularly interested in the optimal assumptions on $\delta$.

\medskip

The requirement that there should be a properly defined operator behind \eqref{sde1} is reasonable, since it gives a reasonably complete solution theory of the corresponding parabolic equation. 
That being said, there are situations where one does not want to insist on a strong link between parabolic equations and stochastic processes in order to treat, in one sense or another, more singular drifts. See e.g.\,\cite{W} which considers solutions of martingale problem with test functions cutting out the singular set of the drift, or \cite{NU,Z_supercritical} which deal with elliptic or parabolic equations with supercritical (in the sense of scaling) drift where one no longer has H\"{o}lder continuity of solution. We are interested, on the contrary, in finding the maximal singularities of the drift that still give a more or less classical theory of parabolic equations and SDEs, including the possibility to consider SDEs with arbitrarily fixed initial point, e.g.\,in the singular set of the drift.

\medskip

One natural question is: why not restrict attention to the Morrey class of drifts \eqref{morrey_elliptic_2}, a broad subclass of \eqref{f1} defined in elementary terms (and which also allows to use e.g.\,harmonic-analytic arguments that are not available for form-bounded or weakly form-bounded drifts, cf.\,Remark \ref{morrey_vs_fbd})? First, there is a refinement of the Morrey class \eqref{morrey_elliptic_2} due to Chang-Wilson-Wolff, which is still contained in the class of form-bounded vector fields (see Appendix \ref{example_sect}), and there is no reason to believe that their result  itself cannot be refined, in elementary terms, even further. Second, and more importantly, by assuming e.g.\,form-boundedness of $b$ \eqref{f1} we impose a condition on $b$ that is ultimate in some precise sense, i.e.\,the existence of quasi contraction strongly continuous semigroup in $L^2$. The latter is, arguably, an extra hypothesis on the diffusion process, but it is explicitly spelled out. Also, at least in the case of Schr\"{o}dinger operators, see below, form-boundedness is a physical assumption on the potential.

\medskip

We impose similar assumptions on the  time-inhomogeneous drift in SDE \eqref{sde2}. 
We will deal with SDE \eqref{sde3} with diffusion coefficients $\sigma$ that can have some critical discontinuities, but only those that are allowed by the singularities of the drift $b$, which is our main focus.  

\medskip

In Section \ref{stable_sect} we discuss weak well-posedness of SDE
\begin{equation}
\label{sde_s}
X_t=x-\int_0^t b(X_s)ds + Z_t, \quad t \geq 0, \quad x \in \mathbb R^d,
\end{equation} 
driven by rotationally symmetric $\alpha$-stable process $Z_t$, $1<\alpha<2$, with drift $b:\mathbb R^d \rightarrow \mathbb R^d$ satisfying
\begin{equation}
\label{f3}
|b| \leq \delta(\lambda-\Delta)^{\frac{\alpha-1}{2}} \quad (\text{in the sense of quadratic forms}),
\end{equation}
e.g.\,
$$
|b|^{\frac{1}{\alpha-1}} \in M_{1+\varepsilon}
$$
for some $\varepsilon>0$. There is a rich literature devoted to  equation \eqref{sde_s}, which also covers the range $0<\alpha \leq 1$ (see Remark \ref{alpha_small_rem}). In the case $1<\alpha<2$, which allows to deal with locally unbounded drifts, earlier results on weak well-posedness of \eqref{sde_s} include $|b| \in L^\infty$ due to Komatsu \cite{Ko}, $|b| \in L^p$,  $p>\frac{d}{\alpha-1}$ due to Portenko and Podolynny-Portenko \cite{P2,PP}, and, more generally, $|b|$ in appropriate Kato class of vector fields, see Chen-Kim-Song \cite{CKS}, Kim-Song \cite{KSo}, Chen-Wang \cite{CW}. All of these classes are contained in \eqref{f3}.

See also Priola \cite{Pr} and X.\,Zhang \cite{Z5} regarding strong well-posedness of \eqref{sde_s} and its generalizations.

\medskip

We discuss conditions on $b$ stated in terms of $|b|$. 
The latter allows to include measure-valued drifts, see Remark \ref{measure_rem}. However,  distributional drifts lie outside of the scope of this paper (regarding distributional drifts, see \cite{FIR, CM,CJM, PZ, ZZ2} and references therein). We also do not discuss here many interesting results that require additional structure of $b$ such as existence of non-positive divergence ${\rm div\,}b$ or $b$ of the form $b=\nabla V$ for appropriate potential $V$; we only refer to Bresch-Jabin-Wang \cite{BJW}, Fournier-Jourdain \cite{FJ}, R\"{o}ckner-Zhao \cite{RZ} and references therein.

\medskip

The drifts that we consider in this paper in general destroy both the upper and the lower Gaussian bounds on the heat kernel of the corresponding to \eqref{sde1} and \eqref{sde2} parabolic equations. 

\medskip

In this paper we are interested in local singularities of the drift, although our drifts can still be unbounded at infinity along some subsets (e.g.\,$b(x)=\sum_{k=1}^\infty c_k|x-a_k|^{-2}(x-a_k)$, $x \in \mathbb R^d$, where $\sum_{k=1}^\infty |c_k|^{\frac{1}{2}}<\infty$ and $a_k \in \mathbb R^d$, $|a_k| \rightarrow \infty$).

\medskip

Throughout the paper, dimension $d \geq 3$. Dimension $d=2$ does not present an obstacle for our methods, however, in our opinion  it deserves a separate study.
The exposition of the results does not follow the chronological order of their appearance (on arXiv), but proceeds from weaker singularities of the drift  to stronger singularities.

\bigskip

\noindent \textbf{Structure of the paper.} In Section \ref{L2_sect} we introduce the class of time-homogeneous form-bounded vector fields. 

In Section \ref{sharp_sect} we discuss the results on weak solvability of SDE \eqref{sde1} and on solution theory of the corresponding parabolic equation under essentially sharp  assumption on the value of form-bound $\delta$ of the drift. 

In Section \ref{three_sect} we describe two existing approaches to constructing a Feller semigroup associated with $-\Delta + b \cdot \nabla$, and introduce a new one.

In Sections \ref{basic_sde_sect}, \ref{thm1_feller_proof_sect} and \ref{thm1_proof_sect} we prove a detailed weak well-posedness result for \eqref{sde1} with form-bounded drift, however, at expense of requiring smaller $\delta$.

In Sections \ref{parab_sect} and \ref{parab_sde_sect} we discuss results on weak well-posedness of SDE \eqref{sde2} having time-inhomogeneous form-bounded drift. Their proofs use a different (iteration) technique.

In Section \ref{a_sect} we discuss an extension of the previous results to some discontinuous diffusion coefficients.

In Sections \ref{ste_sect} and \ref{rz_sect} we discuss strong well-posedness of SDEs \eqref{sde1} and \eqref{sde2}, via stochastic transport equation and via relative compactness criterion for random fields on the Wiener-Sobolev space (R\"{o}ckner-Zhao's approach).

In Sections \ref{weak_sect} and \ref{weak_sde_sect} we substantially enlarge the class of admissible in SDE \eqref{sde1} time-homogeneous drifts (i.e.\,to weakly form-bounded drifts), but at expense of  $\delta$ that now needs to be smaller than in Section \ref{basic_sde_sect}.

In Section \ref{morrey_sect} we consider again time-inhomogeneous drifts and strengthen and simplify all aspects of the results from Sections \ref{parab_sect} and \ref{parab_sde_sect} except their assumptions on $\delta$. We reach, in particular, critical singularities of  drift in the time variable. Compared to Section \ref{weak_sde_sect} we, however, restrict somewhat the class of the spatial singularities of the drift to essentially the largest possible Morrey class.

In Section \ref{stable_sect} we discuss analogues of the weak well-posedness results from Section \ref{weak_sect}
for the SDE \eqref{sde_s} driven by symmetric $\alpha$-stable process. 

\bigskip

\section{Notations}

\label{notation_sect}

\textbf{1.~}$\mathbb R_+:=[0,\infty[$. In what follows, $B_r(x)$ is the open ball of radius $r$ centered at $x \in \mathbb R^d$. Put
$$
\nabla_if:=\partial_{x_i}f,
$$
where $f=f(x)$ or $f=f(t,x)$, $x=(x_1,\dots,x_d)$.

For $\alpha,\beta \in \mathbb R$, set $$\alpha \vee \beta:=\max\{\alpha,\beta\}, \quad \alpha \wedge \beta:=\min\{\alpha,\beta\}.$$

Let $\mathcal B(X,Y)$ denote the space of bounded linear operators between Banach spaces $X \rightarrow Y$, endowed with the operator norm $\|\cdot\|_{X \rightarrow Y}$. $\mathcal B(X):=\mathcal B(X,X)$.

We write $T=s\mbox{-} X \mbox{-}\lim_n T_n$ for $T$, $T_n \in \mathcal B(X,Y)$ if $$\lim_n\|Tf- T_nf\|_Y=0 \quad \text{ for every $f \in X$}.
$$ 

Put $$L^p=L^p(\mathbb R^d), \quad W^{1,p}=W^{1,p}(\mathbb R^d).$$

Set $$\|\cdot\|_p:=\|\cdot\|_{L^p}$$
and
$$
\|\cdot\|_{p \rightarrow q}:=\|\cdot\|_{L^p \rightarrow L^q}.
$$

Let $\mathcal W^{\alpha,p}$, $\alpha>0$, denote the Bessel potential space on $\mathbb R^d$ endowed with norm $\|u\|_{p,\alpha}:=\|g\|_p$,  
$u=(1-\Delta)^{-\frac{\alpha}{2}}g$, $g \in L^p$. Let $\mathcal W^{-\alpha,p'}$, $p'=p/(p-1)$ denote the anti-dual of $\mathcal W^{\alpha,p}$.

For a given vector field $b$ and $1 \leq p<\infty$, put $$b^{\frac{1}{p}}:=b|b|^{-1+\frac{1}{p}}.$$

Put
$$
\langle f,g\rangle = \langle f \bar{g}\rangle :=\int_{\mathbb R^d}f \bar{g}dx$$ 
(some of our functions will be complex-valued).

$C_c$ ($C_c^\infty$) denotes the space of continuous (smooth) functions on $\mathbb R^d$ having compact support.

$C_\infty:=\{f \in C(\mathbb R^d) \mid \lim_{x \rightarrow \infty}f(x)=0\}$ (with the $\sup$-norm).

$\mathcal S$ is the L.\,Schwartz' space of test functions.

We denote by $\upharpoonright$ the restriction of an operator (or a function) to a subspace (a subset).

Given linear operators $A$, $B$, we write $B \supset A$ if $B$ is an extension of $A$.

Let
$$
\big[A \upharpoonright D(A) \cap L^p \big]^{\rm clos}_{L^p \rightarrow L^p}
$$ denote the closure of operator $A$ as an operator $L^p \rightarrow L^p$ (if it exists).

\medskip

\textbf{2.~}Fix $0<T<\infty$. Let $D([0,T],\mathbb R^d)$, the space of right-continuous functions having left limits, be endowed with the filtration $\mathcal B'_t=\sigma(\omega_r \mid 0\leq r \leq t)$, where $\omega_t$ is the coordinate process on  $D([0,T],\mathbb R^d)$.

We will also need the canonical space $(C([0,T],\mathbb R^d),\mathcal B_t=\sigma(\omega_r \mid 0\leq r \leq t))$, where $\omega_t$ is the coordinate process on $C([0,T],\mathbb R^d)$.

Recall that a probability measure $\mathbb P_x$, $x \in \mathbb R^d$ on $(C([0,T],\mathbb R^d),\mathcal B_t)$ is called a martingale solution to SDE 
\begin{equation}
\label{sde1n}
dX_t=-b(t,X_t)dt + \sqrt{2}dW_t, \quad X_0=x
\end{equation}
if

1) $\mathbb P_x[\omega_0=x]=1$;

2) $$\mathbb E_{x} \int_0^T|b(r,\omega_r)|dr<\infty;$$

3) for every $f \in C_c^2(\mathbb R^d)$ the process
$$
t \mapsto f(\omega_t)-f(x) + \int_0^t (-\Delta f + b \cdot \nabla f)(\omega_r)dt
$$
is a $\mathcal B_t$-martingale under $\mathbb P_x$.

A martingale solution $\mathbb P_x$ of \eqref{sde1n} is called a weak solution if, upon completing filtration $\mathcal B_t$ with respect to $\mathbb P_x$ (to, say, $\hat{\mathcal B}_t$), there exists a Brownian motion $\{W_t\}$ on $\big(C([0,T],\mathbb R^d),\hat{\mathcal B}_t,\mathbb P_x\big)$ such that
$$
\omega_t=x - \int_0^t b(r,\omega_r)dr + \sqrt{2}W_t, \quad 0 \leq t \leq T \quad \mathbb P_x\text{-a.s.}
$$

\medskip

\textbf{3.~}A function $h:\mathbb R^d \rightarrow \mathbb R$ is said to be in the weak $L^d$ class (denoted by $L^{d,w}$) if $$\|h\|_{d,w}:=\sup_{s>0}s|\{x \in \mathbb R^d: |h(x)|>s\}|^{1/d}<\infty.$$ 
Clearly, $L^d \subset L^{d,w}$, but not vice versa, e.g.\,$h(x)=|x|^{-1}$ is in $L^{d,w}$ but not in $L^d$.

\medskip

\textbf{4.~}The De Giorgi mollifier $E_\varepsilon \equiv E_\varepsilon^d$ on $\mathbb R^d$:
$$
E_\varepsilon f(x):=e^{\varepsilon \Delta }f(x), \quad x \in \mathbb R^d, \quad \varepsilon>0, 
$$
where $f \in L^1_{\loc}$.

\medskip

The Friedrichs mollifier $E_\varepsilon \equiv E_\varepsilon^d$ on $\mathbb R^d$:
$$E_\varepsilon f(x):=\eta_\varepsilon \ast f(x),$$
where $\eta_\varepsilon(y):=\frac{1}{\varepsilon^{d}}\eta\left(\frac{y}{\varepsilon}\right)$, $\varepsilon>0$ and
$$
\eta(y):=\left\{
\begin{array}{ll}
c\exp\left(\frac{1}{|y|^2-1}\right)& \text{ if } |y|<1, \\
0, & \text{ if } |y| \geqslant 1,
\end{array}
\right.
$$
with constant $c$ adjusted to $\int_{\mathbb R^d} \eta(x)dx=1$.

\bigskip

\section{Form-bounded drifts. Semigroup in $L^2$}

\label{L2_sect}

First, we discuss sufficient conditions for existence of an operator realization of $-\Delta + b \cdot \nabla$  generating a strongly continuous semigroup in $L^2$.

Assume that a Borel measurable vector field $b:\mathbb R^d \rightarrow \mathbb R^d$ with $|b| \in L^2_{\loc}$ satisfies inequality
\begin{equation}
\label{fb}
\|b\varphi\|_2^2 \leq \delta\|\nabla \varphi\|_2^2+c_\delta\|\varphi\|_2^2 \quad \text{ for all } \varphi \in W^{1,2}
\end{equation}
for finite constants $\delta>0$ and $0 \leq c_\delta<\infty$. 

\begin{definition}
\label{ellip_fb_def}
Such vector fields $b$ are called \textit{form-bounded} (written as $b \in \mathbf{F}_\delta$).
\end{definition}

Inequality \eqref{fb} can be re-stated as an operator norm inequality
\begin{equation}
\label{op_fb}
\|b (\lambda-\Delta)^{-\frac{1}{2}}\|_{2 \rightarrow 2} \leq \sqrt{\delta}
\end{equation}
with $\lambda \equiv \lambda_\delta=c_\delta/\delta$.

\medskip

Using the quadratic (or Cauchy-Schwarz) inequality 
\begin{equation}
\label{quad}
|\langle b \cdot \nabla  \varphi,\varphi\rangle| \leq \varepsilon \|b\varphi\|_2^2+\frac{1}{4\varepsilon}\|\nabla \varphi\|_2^2, \quad \varepsilon>0,
\end{equation}
one can see that the form-boundedness condition \eqref{fb} \textit{with} $\delta<1$ is what is needed to verify conditions of the Lax-Milgram theorem for the bilinear form $\tau[\varphi,\psi]:=\lambda\langle \varphi,\psi\rangle + \langle \nabla \varphi,\nabla \psi\rangle + \langle b \cdot \nabla \varphi,\psi \rangle$ defined on the real space $W^{1,2}$. That is, coercivity for all $\lambda \geq c_\delta/2\sqrt{\delta}$:
\begin{align}
|\tau[\varphi,\varphi]| & \geq \lambda\|\varphi\|_2^2 + \|\nabla \varphi\|_2^2 -\varepsilon \|b\varphi\|_2^2-\frac{1}{4\varepsilon}\|\nabla \varphi\|_2^2 \notag  \\ 
& (\text{we apply \eqref{fb} and select $\varepsilon=\frac{1}{2\sqrt{\delta}}$})  \notag \\
& = (\lambda-\frac{c_\delta}{2\sqrt{\delta}})\|\varphi\|_2^2 + (1-\sqrt{\delta})\|\nabla \varphi\|_2^2 \label{coer}
\end{align}
and boundedness
\begin{equation}
\label{bdd}
|\tau[\varphi,\psi]| \leq C\|\varphi\|_{W^{1,2}}\|\psi\|_{W^{1,2}}.
\end{equation}
So, by the Lax-Milgram theorem, there exists a unique weak solution to elliptic equation $$(\mu-\Delta + b \cdot \nabla)u=f, \quad f \in L^2.$$
Furthermore, form-boundedness \eqref{fb} with $\delta<1$ ensures that the sesquilinear form  of $\lambda-\Delta + b \cdot \nabla$ defined on the complex space $W^{1,2}$ is sectorial, and hence by the KLMN theorem\footnote{Kato-Lions-Lax-Milgam-Nelson theorem, see  \cite[Ch.\,VI]{Ka}, \cite[Ch.1]{O}}
 it determines a unique closed densely defined operator $\Lambda \equiv \Lambda_2(b) $,
\begin{equation}
\label{kolm}
\Lambda\supset (-\Delta + b \cdot \nabla)\upharpoonright C_c^\infty(\mathbb{R}^d),
\end{equation}
 generating a strongly continuous quasi contraction semigroup in $L^2$.

\begin{examples}
Let us mention some elementary sufficient conditions for form-boundedness.

1.~If $|b| \in L^d$ (or $|b| \in L^d+L^\infty$, i.e.\,the sum of two functions, one in $L^d$ and the other one in $L^\infty$), then $b \in \mathbf{F}_\delta$ with $\delta$ that can be chosen arbitrarily small (at expense of increasing $c_\delta$, see Appendix \ref{example_sect}). 

\medskip

2.~There are form-bounded vector fields that have stronger singularities than the ones covered by the class $L^d$. For instance, by Hardy's inequality $$\left(\frac{d-2}{2}\right)^2\||x|^{-1}\varphi\|_2^2 \leq \|\nabla \varphi\|_2^2, \quad \varphi \in W^{1,2},$$ the vector field
\begin{equation}
\label{hardy_vf}
b(x)=\pm \sqrt{\delta} \frac{d-2}{2}\frac{x}{|x|^2},
\end{equation}
is form-bounded with $c_\delta=0$. The constant in Hardy's inequality is sharp, and the last vector field is not in $\mathbf{F}_{\delta'}$ for any $\delta'<\delta$ regardless of the value of $c_{\delta'}$.

As we explain below, the value of constant $\delta$ determines weak solvability of SDEs (see \eqref{sde2_ce} below), and thus is a key characteristics of the vector field $b$. 
However, the dependence of the solution theory of SDEs on $\delta$ is not visible if one deals only with $|b| \in L^d$. In this sense, the vector fields $b$ with entries in $L^d$ are sub-critical.

\medskip

3.~More generally, if vector field $b$ belongs to the scaling-invariant Morrey class $M_{2+\varepsilon}$ for some $\varepsilon>0$ arbitrarily small, i.e.
\begin{equation}
\label{morrey_elliptic}
|b| \in L^{2+\varepsilon}_{\loc} \quad \text{ and } \quad
\|b\|_{M_{2+\varepsilon}}:=\sup_{r>0, x \in \mathbb R^d} r\biggl(\frac{1}{|B_r|}\int_{B_r(x)}|b|^{2+\varepsilon}dx \biggr)^{\frac{1}{2+\varepsilon}}<\infty,
\end{equation}
then $b$ is form-bounded with $\delta = c\|b\|_{M_{2+\varepsilon}}$ for appropriate constant $c=c(d,\varepsilon)$. See \cite{F}, see also  \cite{CF}. Note that $$\|\cdot\|_{M_{q}} \leq \|\cdot\|_{M_{q_1}} \quad \text{ if } q<q_1,$$ so Morrey class becomes larger as $q$ becomes smaller (and so we are interested in the class $M_q$ with $q$ close to $2$). This class contains all $|b| \in L^d$ and $|b| \in L^{d,w}$ (see definition in Section \ref{notation_sect}). It also contains, for every $\epsilon>0$, vector fields $b$ such that $|b| \not\in L^{2+\epsilon}_{\loc}$.

On the other hand, it is easy to show, by considering translates of a bump function, that if $b \in \mathbf{F}_\delta$  (say, with $c_\delta=0$), then $|b| \in M_2$.

\medskip

4. If $|b|^2$ is in the Kato class of potentials $\mathbf{K}^d_\delta$, then vector field $b$ is form-bounded. Recall that $V \in \mathbf{K}^d_\delta$ if 
\begin{equation}
\label{kato_pt_def}
V \in L^1_{\loc} \quad \text{ and } \quad \lim_{\lambda \rightarrow +\infty}\|(\lambda-\Delta)^{-1}|V|\|_\infty \leq \delta.
\end{equation}
The Kato class condition first appeared in a 1961 article by M.\,Birman as an elementary sufficient condition for the form-boundedness of a potential $V$ \cite[Sect.\,2]{Bi}.

Some other  sufficient conditions for the form-boundedness of $b$, including those refining condition $b \in M_{2+\varepsilon}$, are given in Appendix \ref{example_sect}. 

Let us also note that the sum $b_1+b_2$ of form-bounded vector fields $b_1 \in \mathbf{F}_{\delta_1}$, $b_2 \in \mathbf{F}_{\delta_2}$ is form-bounded with form-bound $\delta=(\sqrt{\delta_1}+\sqrt{\delta_2})^2$ (cf.\,\eqref{op_fb}). In particular, vector field
$$b(x):=\sum_{k=1}^\infty c_k\frac{x-a_k}{|x-a_k|^{2}}$$ with $\sum_{k=1}^\infty |c_k|^{\frac{1}{2}}<\infty$ and $\{a_k\}$ constituting e.g.\,a dense subset of $\mathbb R^d$ is form-bounded.
\end{examples}

The form-boundedness condition is well known in the literature on spectral and scattering theory of Schr\"{o}dinger operators, in particular, dealing with the questions of self-adjointness, estimates on the number of bound states, resolvent convergence, see e.g.\,\cite{BS,MV,Si}\footnote{If one is willing to ignore different roles played by the positive and the negative parts of potential $V$ in the theory of Schr\"{o}dinger operator $-\Delta + V$, then $V$ is form-bounded if $\||V|^{\frac{1}{2}}\varphi\|_2^2 \leq \delta\|\nabla \varphi\|_2^2+c_\delta\|\varphi\|_2^2$, $\varphi \in W^{1,2}$.}.  
Regarding Kolmogorov operator \eqref{kolm}, one can show, with little additional effort, that the semigroup $e^{-t\Lambda}$ is a Markov semigroup in $L^2$ (the proof can be found e.g.\,\cite{KiS_theory}). In the probabilistic context, however, one of the basic problems is to construct a Markov process that inherits the essential properties of the Brownian motion. In particular, it is natural to expect that the constructed Markov semigroup would be strongly continuous on the space $C_\infty$ of continuous functions vanishing at infinity endowed with the $\sup$-norm, i.e.\,that it would be a Feller semigroup. However, to show this, one needs to verify the strong continuity in the norm of $C_\infty$, which is, of course, a lot more rigid that the norm of $L^2$ where one defines the form-boundedness of the drift. In this regard, let us note that the strong continuity of this semigroup in $L^p$ for any finite $p \geq 2$, on the contrary, presents no problem. Indeed, since $e^{-t\Lambda}$ is Markov, we have $\|e^{-t\Lambda}f\|_\infty \leq \|f\|_\infty$, $f \in L^2 \cap L^\infty$, so one can use an interpolation argument to define
\begin{equation}
\label{p_semi}
e^{-t\Lambda_p}:=\biggl[e^{-t\Lambda} \upharpoonright L^2 \cap L^p \biggr]_{L^p \rightarrow L^p}^{\rm clos}.
\end{equation}
The left-hand side is a quasi contraction strongly continuous semigroup on $L^p$ (see the proof e.g.\,in \cite{LS}), its generator $\Lambda_p$ is an appropriate operator realization of $-\Delta +b \cdot \nabla $ in $L^p$. So, the difficulty is in the strong continuity in $C_\infty$.
A major advancement came with the fundamental paper \cite{KS} of Kovalenko-Sem\"{e}nov that, among other results, presented a construction of a Feller realization of $-\Delta + b \cdot \nabla$ in $C_\infty$ for form-bounded $b$ using an $L^2 \rightarrow L^\infty$ iteration procedure. To the best of the author's knowledge, surprisingly, the first reaction to this result  was \cite{Ki}, almost three decades later.

Apart from the results described in the present paper, we also refer to \cite{CFKZ,FK} regarding form-boundedness appearing in probabilistic settings.

Form-boundedness and similar conditions, sometimes supplemented with other hypothesis on $b$, also appear in the regularity theory of elliptic and parabolic equations, which includes the Harnack inequality, Gaussian and non-Gaussian heat kernel bounds. See \cite{AD,H,KiS_heat,KiV, LZ,Ph, S} and references therein.

\medskip

Regarding the necessity of the form-boundedness condition \eqref{fb} with $\delta<1$ for the existence of $L^2$ semigroup theory of $-\Delta + b \cdot \nabla$, let us mention the following consequence of the result in \cite{MV}. Let $b$ be a distributional vector field. The sesqulinear form $-\Delta + b\cdot \nabla$ is bounded, i.e.\,
\begin{equation}
\label{MV_bd}
|\langle \nabla \varphi,\nabla \psi\rangle + \langle b\cdot \nabla \varphi,\psi\rangle| \leq C\|\varphi\|_{W^{1,2}}\|\psi\|_{W^{1,2}}
\end{equation}
for some constant $C$, for all $\varphi,\psi \in C_c^\infty$,
if and only if $b$ can be represented as $b=b^{(1)}+b^{(2)}$, where $b^{(1)} \in \mathbf{F}_\delta$ for some $\delta$ and $b^{(2)}$ is divergence-free and is in the class\footnote{i.e.\,the components of $b^{(2)}$ satisfy
$$
b^{(2)}_k=\sum_{i=1}^d \nabla_i F_{ik}, \quad 1 \leq k \leq d,
$$
for a matrix $F$ with entries $F_{ik}=-F_{ki} \in {\rm BMO}$. Recall that a function $f \in {\rm BMO}$ if $f \in L^1_{\loc}$ and
$
\|f\|_{\rm BMO}:=\sup_{Q}\frac{1}{|Q|}\int_{Q}|f-(f)_Q|dx<\infty,
$
where the supremum taken over all cubes $Q \subset \mathbb R^d$ with sides parallel to the axes, and $(f)_Q$ is the average of $f$ over $Q$.} ${\rm BMO}^{-1}$. Thus, since in this paper we are interested in conditions on $|b|$, i.e.\,not assuming any additional structure of $b$ such as zero divergence, our condition \eqref{fb} is essentially necessary for \eqref{MV_bd} to hold. However, \eqref{MV_bd} is not synonymous with the existence of a realization of $-\Delta + b \cdot \nabla$ in $L^2$ generating a strongly continuous semigroup. Should we require or expect \eqref{MV_bd} to hold in order to have $L^2$ semigroup theory? As we explain below, for $-\Delta + b \cdot \nabla$ the answer is ``no'': in the next section we will abandon \eqref{MV_bd} and will go beyond the class $\mathbf{F}_\delta$. However, 
we will have  \eqref{MV_bd} when will be dealing with SDE \eqref{sde3} having discontinuous diffusion coefficients. See Section \ref{a_sect}.

\medskip

 Concerning the difference between a popular condition\footnote{As is well known, on the Lebesgue scale, $|b| \in L^d$ is the best possible condition providing the solvability of \eqref{sde1}. 
} $|b| \in L^d$  and more general condition $b \in \mathbf{F}_\delta$, let us  note the following: if $v$ is a weak solution of the elliptic equation $$(\lambda -\Delta  + b \cdot \nabla )v=f, \quad \lambda>0, \quad f \in C_c^\infty$$ with $|b| \in L^d$ and $v \in W^{1,r}$ for $r$ large (which is valid by a classical result), then, by H\"{o}lder's inequality,
$$
\Delta v \in L^{\frac{rd}{d+r}}_{\loc}. 
$$ 
However, for $b \in \mathbf{F}_\delta$, one can only say that $$\Delta v \in L^{\frac{2d}{d+2}}_{\loc} $$ (one can in fact show that $v \in W^{2,2}$). 
Thus, in case $b \in \mathbf{F}_\delta$, any $W^{2,p}$ estimate on 
the solution $v$ of the elliptic equation for $p$ large is out of question.

\begin{remark}
\label{grad_unb_rem}
That being said, if $\delta<1$, then one has $$v \in \mathcal W^{1+\frac{2}{q},p}, \quad p \in \big[2,\frac{2}{\sqrt{\delta}}[, \quad q>p$$
(see Theorem \ref{thm1_feller} below). After applying the Sobolev embedding theorem, one obtains 
$|\nabla v| \in L^{\gamma}$ for $\gamma<\frac{dp}{d-2}$ arbitrarily close to $\frac{dp}{d-2}$ (depending on how close $q$ is to $p$). Thus, although $p$ can be as large as one wants provided that $\delta$ is chosen sufficiently small, one never arrives at  $|\nabla v| \in L^\infty$. For a form-bounded drift $b$, the gradient of  $v$ is in general unbounded.
\end{remark}

\medskip

\section{Sharp solvability}

\label{sharp_sect}

The constant $c_\delta$ in \eqref{fb} controls the growth of the semigroup $e^{-t\Lambda}$ as $t \rightarrow +\infty$, see \eqref{growth}, and allows to include in the class $\mathbf{F}_\delta$ bounded vector fields. It is of secondary interest to us in this paper.

\medskip

The constant $\delta$ in the assumption $b \in \mathbf{F}_\delta$, on the contrary, is very important to us since it determines weak solvability of SDE \eqref{sde1}. Moreover, there is a quantitative dependence between $\delta$ and the regularity properties of solutions to the corresponding elliptic and parabolic equations, see Theorem \ref{thm1_feller} and other results below. 

\medskip

The following example, analysed in detail in \cite{BFGM}, shows that $\delta$ cannot be too large. 
Consider SDE
\begin{equation}
\label{sde2_ce}
X_t= - \sqrt{\delta} \frac{d-2}{2}\int_0^t |X_r|^{-2}X_r dr + \sqrt{2}W_t,
\end{equation}
which corresponds to the choice of  attracting drift 
\begin{equation}
\label{vf}
b(x)=\sqrt{\delta} \frac{d-2}{2}\frac{x}{|x|^2} \in \mathbf{F}_\delta
\end{equation}
and the initial point $x=0$ in SDE  \eqref{sde1}.
Then, if 
\begin{equation}
\label{delta_counter}
\delta>4\bigg(\frac{d}{d-2}\bigg)^2, 
\end{equation}
the SDE \eqref{sde2_ce} does not have a weak solution. Indeed, suppose  that a weak solution exists.
Then $X(t)=(X^1_t,\dots,X^d_t)$ is a continuous semimartingale with cross-variation $[X^i,X^k]_t=2\delta_{ik}t$.
By It\^{o}'s formula,
$$
|X_t|^2= - 2\int_0^t X_s b(X_s)ds + 2\sqrt{2}\int_0^t X_sdW_s + 2\int_0^t d[W,W]_s,
$$
i.e.~
$$
|X_t|^2=-2  \sqrt{\delta} \frac{d-2}{2} \int_0^t \mathbf{1}_{X_s\neq 0} ds + 2\sqrt{2}\int_0^t X_sdW_s + 2 t d.
$$
One has $\int_0^t \mathbf{1}_{X_s=0}ds=0$ almost surely (see details in \cite{BFGM}), so 
$$
|X_t|^2=2\big(d-\sqrt{\delta} \frac{d-2}{2}\big)  \int_0^t \mathbf{1}_{X_s\neq 0} ds + 2\sqrt{2}\int_0^t X_sdW_s \quad \text{a.s.}
$$
Therefore, $X_t^2 \geq 0$ is a local supermartingale if $\sqrt{\delta} \frac{d-2}{2} > d$.
Then a.s. $X_0=0 \Rightarrow X_t=0$, which contradicts to  $[X^1,X^1]_t=2t$. This argument was used earlier in \cite{CE} in the one-dimensional setting. 
\cite{BFGM} furthermore showed if $\delta>4$, then a trajectory started outside of the origin arrives at $x=0$ in finite time with positive probability; in this regard, see also \cite{W}.

\medskip

Although at the first sight this counterexample can be interpreted (and sometimes was) as a counterexample showing the optimality of the condition $|b| \in L^d$, we argued in \cite{KiS_brownian} that this is a counterexample to admissible values of constant $\delta$. In fact, we have the following theorem.

\begin{theorem}
\label{thm_sharp}
Let $b \in \mathbf{F}_{\delta}$ with
\begin{equation}
\label{c1}
\tag{$C_1$}
\delta<4.
\end{equation}
Then, for every $x \in \mathbb R^d$, the SDE
\begin{equation}
\label{sde_sharp}
X_t=x-\int_0^b b(X_s)ds + \sqrt{2}W_t, \quad t \geq 0,
\end{equation}
has a martingale solution.
\end{theorem}

Theorem \ref{thm_sharp} was proved in \cite{KiS_sharp}. (In fact, it was proved there for time-inhomogeneous form-bounded drifts $b$ with $\delta<4$, see Definition \ref{parab_fb_def}.)

Comparing ``$\delta>4(\frac{d}{d-2})^2$'' in \eqref{delta_counter} and ``$\delta<4$'', one sees that the result is essentially sharp in high dimensions.

Let us note that the well-posedness of SDEs and parabolic equations with $\delta$ reaching the critical value (up to the strict inequality) but for $b$ having additional structure (e.g.\,$b=\nabla V$ for appropriate potential $V$ such as $V(x)=c \log |x|$), was also addressed by Fournier-Jourdan \cite{FJ} (in dimension 2), Bresch-Jabin-Wang \cite{BJW}, see also references therein. A crucial feature of Theorem \ref{thm_sharp} is that it attains the critical threshold $\delta=4$, up to the strict inequality, for the entire class of form-bounded vector fields, i.e.\,without any assumptions on the structure of $b$. 

\medskip

Let us explain where does ``$\delta<4$'' come from (and also how one can handle $1 \leq \delta<4$ given that the KLMN theorem requires $\delta<1$). Multiplying the corresponding to \eqref{sde_sharp} parabolic equation
$(\partial_t-\Delta + b \cdot \nabla)u=0$ by $u|u|^{p-2}$, integrating by parts and using the quadratic inequality and form-boundedness, one obtains that the admissible $p$ that give e.g.\,an energy inequality turn out to be $p>2/(2-\sqrt{\delta})$. In fact, it was proved in \cite{KS} that if $b=b(x)$ has form-bound $\delta<4$, then there exists a realization of $-\Delta + b \cdot \nabla$ in $L^p$, $p>2/(2-\sqrt{\delta})$ generating a quasi contraction strongly continuous semigroup there. The interval of contraction solvability can be extended to $[2/(2-\sqrt{\delta}),\infty[$ and is sharp, see \cite{KiS_theory}.
Now, as $\delta \uparrow 4$, this interval disappears, and with it disappears the theory of operator $-\Delta + b \cdot \nabla$ (see, however, the next section). 

\medskip

The proof of Theorem \ref{thm_sharp} is based on the following analytic result, which allows to use the standard tightness argument (see \cite{RZ}) to construct a martingale solution of \eqref{sde_sharp}. Namely, put 
\begin{equation}
\label{rho}
\rho(x)=(1+\kappa|x|^{-2})^{-\beta}, \quad \beta>\frac{d}{4}, 
\end{equation}
with $\kappa>0$ sufficiently small.
Let 
$u$ be the classical solution to Cauchy problem
\begin{equation}
\label{c_dg}
(\partial_t  - \Delta  + b \cdot \nabla) u=|\mathsf{h}|f, \quad u(0)=0,
\end{equation}
where $f \in C_c$ and
$b \in \mathbf{F}_\delta \cap C^\infty_c$, $\delta<4$ and
$\mathsf{h} \in \mathbf{F}_\nu \cap C^\infty_c$, $\nu<\infty.$
 Fix $T>0$ and $1<\theta<\frac{d}{d-1}$.
For all $p>\frac{2}{2-\sqrt{\delta}}$, $p\geq 2$, there exists a constant $C$ independent of smoothness of $b$ and $\mathsf{h}$ such that 
\begin{align}
\|u\|_{L^\infty([0,T] \times \mathbb R^d)} & \leq C 
\sup_{z \in \mathbb Z^d}\biggl(\int_0^T \bigg\langle \big(\mathbf{1}_{\{|\mathsf{h}| \geq 1\}} + \mathbf{1}_{\{|\mathsf{h}| < 1\}}|\mathsf{h}|^p\big)^{\theta'} |f|^{p\theta' }\rho^2_z\bigg\rangle\biggr)^{\frac{1}{p\theta'}} \label{est_m}
\end{align}
where $\rho_z(x):=\rho(x-z)$.

Now, let $b_n$ be smooth vector fields having compact support, e.g.\,defined by \eqref{b_n1} below, approximating $b$ in the sense of \eqref{b_n2}, \eqref{b_n3}.
Fix $x \in \mathbb R^d$. By a classical result, there exist strong solutions $X^n$ to SDEs
$$
X^n_t=x-\int_0^t b_n(X^n_s)ds + \sqrt{2} W_t, \quad n=1,2,\dots,
$$
where $\{W_t\}_{t \geq 0}$ is a Brownian motion in $\mathbb R^d$ on a fixed complete probability space $(\Omega,\mathcal F,\mathcal F_t,\mathbf P)$. Then \eqref{est_m} yields
\begin{equation}
\label{key_est0}
\left|\mathbf E\int_{t_0}^{t_1} |\mathsf{h}(X^n_s)|f(X^n_s)ds\right| \leq  
C\sup_{z \in \mathbb Z^d}\biggl(\int_{t_0}^{t_1} \bigg\langle \big(\mathbf{1}_{\{|\mathsf{h}| \geq 1\}} + \mathbf{1}_{\{|\mathsf{h}| < 1\}}|\mathsf{h}|^p\big)^{\theta'} |f|^{p\theta' }\rho^2_z\bigg\rangle\biggr)^{\frac{1}{p\theta'}},
\end{equation}
where $0 \leq t_0 <t_1 \leq T$.
We now apply \eqref{key_est0} with $\mathsf{h}=b_n$ and $f \equiv 1$ (here $f \in C_c$ $\Rightarrow$ $f \equiv 1$  using Fatou's Lemma):
\begin{align}
\mathbf E\int_{t_0}^{t_1} |b_n(X^n_s)|ds & \leq   
C\sup_{z \in \mathbb Z^d}\biggl(\int_{t_0}^{t_1} \bigg\langle \big(\mathbf{1}_{\{|b_n| \geq 1\}} + \mathbf{1}_{\{|b_n| < 1\}}|b_n|^p\big)^{\theta'} \rho^2_z\bigg\rangle\biggr)^{\frac{1}{p\theta'}} \notag \\
& \leq C_0 (t_1-t_0)^\mu \quad \text{ for generic }\mu>0 \text{ and } C_0. \label{tight_ineq}
\end{align}
The latter allows to verify tightness of probability measures $$\mathbb P^n_x:=(\mathbf P \circ X^n)^{-1}$$ on $\bigl(C([0,T],\mathbb R^d),\mathcal B_t\bigr)$, so for every $x \in \mathbb R^d$ 
there exists a subsequence $\{\mathbb P_x^{n_k}\}$ and a probability measure $\mathbb P_x$ 
such that
\begin{equation}
\label{w}
\mathbb P_x^{n_k} \rightarrow \mathbb P_x \text{ weakly}.
\end{equation}
Another application of $\eqref{key_est0}$ but with $\mathsf{h}=b_{n_1}-b_{n_2}$ allows to conclude that $\mathbb P_x$ is a martingale solution of \eqref{sde_sharp}, see \cite{KiS_sharp} for details.

\medskip

Estimate \eqref{est_m} is proved using De Giorgi's iterations in $L^p$ for $p>2/(2-\sqrt{\delta})$. Thus, $p>2$  if $1<\delta<4$. 
 In this regard, let us note that passing to $L^p$ right away, using the fact that $u^\frac{p}{2}$ is a sub-solution, and then applying to $u^\frac{p}{2}$ the standard  De Giorgi's iterations in $L^2$, does not allow to treat $1 \leq \delta<4$. We have to 
follow the iteration procedure from the very beginning and adjust it accordingly.

Earlier, De Giorgi's method in $L^2$ was used by Zhang-Zhao \cite{ZZ}, Zhao \cite{Zh}, R\"{o}ckner-Zhao \cite{RZ} to prove, among other results, weak well-posedness of \eqref{sde_sharp} with $b$ having not too singular divergence  and  satisfying $$|b| \in L_{\loc}^q(\mathbb R_+,L^r+L^\infty), \quad \frac{d}{r}+\frac{2}{q} < 2.$$

\begin{remark}
Looking at the counterexample \eqref{sde2_ce}-\eqref{delta_counter} and Theorem \ref{thm_sharp}, one can draw an analogy with the celebrated result of Baras-Goldstein for Schro\"{o}diner operator $$-\Delta - V_0, \quad V_0(x)=\delta \frac{(d-2)^2}{4}|x|^{-2}$$ on $\mathbb R^d$, $d \geq 3$. This $V_0$ is a form-bounded potential, i.e.\,$\langle V_0\varphi,\varphi \rangle \leq \delta \langle \nabla \varphi,\nabla \varphi \rangle$ for all $\varphi \in W^{1,2}$ (this is Hardy's inequality). If $0 <\delta < 1$, then the self-adjoint operator realization $H$ of $-\Delta - V_0$ on $L^2$, defined e.g.\,via the KLMN theorem, satisfies $$e^{-tH}=s{\mbox-}L^2{\mbox-}\lim_{\varepsilon \downarrow 0} e^{-tH(V_{\varepsilon})},$$ where $V_{\varepsilon}(x)=\delta \frac{(d-2)^2}{4}|x|_\varepsilon^{-2}$, $|x|_\varepsilon^{2}:=|x|^2+\varepsilon$, $\varepsilon>0$. 
For $\delta>1$, however, by the result in \cite{BG} (see also \cite{GZa}),
$$
\lim_{\varepsilon \downarrow 0} e^{-tH(V_\varepsilon)}f(x)=\infty, \quad t>0, \quad x \in \mathbb R^d, \quad f \geq 0,\; f \not\equiv 0,
$$
i.e.\,all positive solutions of the corresponding parabolic equation explode instantly at any point. This phenomenon is not observable on any $V_0 \in L^{\frac{d}{2}}$ since any such potential has arbitrarily small form-bound (similarly to how any $b$ with $|b| \in L^d$ has arbitrarily small form-bound $\delta$).
\end{remark}

\medskip

\subsection{Critical magnitude of drift singularities}
It turns out that one still has a strong\footnote{Here ``strong'' refers to the differentiability of the solution in time.} solution (i.e.\,semigroup) theory of parabolic equation $(\partial_t-\Delta + b \cdot \nabla)u=0$ with $b \in \mathbf{F}_\delta$ in the critical borderline case $\delta=4$, although not in $L^p$, as in the previous section, but in an Orlicz space. Orlicz spaces are known to appear in various borderline situations in analysis (e.g.\,Trudiner's theorem).

\medskip

We will work over the $d$-dimensional torus $\Pi^d$ obtained as the quotient of $[-\frac{1}{2},\frac{1}{2}]^d$. This is not a technical assumption since the volume of the torus enters the estimates. Still, since we are interested in the local singularities of $b$, working on a torus is sufficient for our purposes. The definition of form-bounded vector fields does not change, except that now we integrate over the torus. Moreover, the examples of form-bounded vector fields mentioned earlier in the paper remain essentially intact, see \cite{BO,Gu}. In this section, $\langle\cdot,\cdot\rangle$ denotes integration over torus $\Pi^d$.

\medskip

Put 
\begin{align*}
\Phi(y)&:=\cosh(y)-1, \quad y \in \mathbb R
\end{align*}
Expanding $\cosh$ in the Taylor series, one sees right away that $\Phi(y)=\Phi(|y|)$.
 This function  is convex on $\mathbb R$, $\Phi(y)=0$ if and only if $y=0$, $\Phi(y)/y \rightarrow 0$ if $y \rightarrow 0$ and $\Phi(y)/y \rightarrow \infty$ if $y \rightarrow \infty$. Therefore, the space $\mathcal L_\Phi=\mathcal L_\Phi(\mathbb R^d)$ of real-valued measurable functions $f$ on $\Pi^d$ satisfying
\begin{equation}
\label{Phi_norm}
\|f\|_\Phi:=\inf\big\{c>0 \mid \langle \Phi\bigl(\frac{f}{c}\bigr) \rangle \leq 1\big\}<\infty,
\end{equation}
is a Banach space with respect to Orlicz norm  $\|\cdot\|_\Phi$. See e.g.\,\cite[Ch.\,8]{AF}. 

\medskip

Let $L_\Phi$ denote the closure $\mathcal L_\Phi$ of the subspace of smooth functions on $\Pi^d$. This is our Orlicz space.

\medskip

We note that 
\begin{equation}
\label{Lpnorm}
\|\cdot\|_\Phi \geq \frac{1}{(2p)!}\|\cdot\|_{2p}, \quad p=1,2,\dots,
\end{equation}
i.e.\,we are dealing with an Orlicz norm that is stronger that any $L^p$ norm.

\medskip

Let $u_n$ be the classical solution to Cauchy problem 
\begin{equation*}
\left\{
\begin{array}{rr}
(\partial_t - \Delta +b_n \cdot \nabla)u_n=0 \text{ on } [0,\infty[ \times \Pi^d, \\[2mm]
 u_n(0,\cdot)=f(\cdot) \in C^\infty(\Pi^d),
\end{array}
\right.
\end{equation*}
where $b_n$ are  bounded smooth vector fields such that $b_n \in \mathbf{F}_\delta$ with the same $c_\delta$ and
$b_n \rightarrow b$ in $L^2(\Pi^d)$ (e.g.\,obtained upon applying to $b$ the De Giorgi mollifier on $\Pi^d$).
Let $e^{-t\Lambda(b_n)}$, $\Lambda(b_n):=- \Delta +b_n \cdot \nabla$ denote the corresponding semigroup, i.e.
$$e^{-t\Lambda(b_n)}f:=u_n(t).$$
On the smooth initial functions, $[0,\infty[ \ni t \mapsto e^{-t\Lambda(b_n)}f$ is strongly continuous in the norm of $L_\Phi$ since it is strongly continuous in the norm of $L^\infty$.

\begin{theorem} 
\label{thm_Orlicz}
Let $b \in \mathbf{F}_\delta$, $0<\delta \leq 4$. The following are true:

\smallskip

{\rm (\textit{i})} For all $n \geq 1$, $f \in C^\infty(\Pi^d)$, 
$$\|e^{-t\Lambda(b_n)}f\|_\Phi \leq e^{2\frac{c_\delta}{\sqrt{\delta}} t}\|f\|_\Phi, \quad t \geq 0.$$

\smallskip

{\rm (\textit{ii})} There exists a strongly continuous quasi contraction semigroup $e^{-t\Lambda(b)}$ on $L_\Phi$ such that, for every $f \in C^\infty(\Pi^d)$, 
$$
\|e^{-t\Lambda(b)}f-e^{-t\Lambda(b_n)}f\|_{\Phi} \rightarrow 0 \quad \text{ as } n \rightarrow \infty \text{ loc.\,uniformly in $t \geq 0$}.
$$
It follows that $e^{-t\Lambda(b)}$ is a positivity preserving $L^\infty$ contraction. Its generator $\Lambda(b)$ is the appropriate operator realization of the formal operator $-\Delta + b \cdot \nabla$ in $L_\Phi$.

This semigroup is unique in the sense that it does not depend on the choice of smooth vector fields $\{b_n\}$, $b_n \rightarrow b$ in $L^2(\Pi^d)$, as long as they do not increase constants $\delta$, $c_\delta$.

\smallskip

{\rm (\textit{iii})} The following energy inequality holds for $u=e^{-t\Lambda(b_n)}f$:
$$
 \frac{1}{2}\sup_{s \in [0,t]}\langle e^{u^{p}(s)}  \rangle +  4\frac{(p-1)}{p}\langle (\nabla u^{\frac{p}{2}})^2e^{u^{p}}\rangle \leq \langle e^{f^{p}}  \rangle, \quad p=2,4,\dots
$$
provided $\frac{c_\delta}{\sqrt{\delta}}t<\frac{1}{2}$; the last constraint can be removed using the semigroup property.
\end{theorem}

The last assertion is noteworthy, since, at the first sight, it seems like one can reach $\delta=4$ only at the cost of killing off the dispersion term.

\medskip

Theorem \ref{thm_Orlicz} was proved in \cite{Ki_Orlicz}.
The following calculation  illustrates the main observation behind this result.
Below we are looking for integral bounds on $u_n$ that can depend on $\delta$ and $c_\delta$, but not on the smoothness or boundedness of $b_n$. Put $u=u_n$, $b=b_n$.
Replacing $v$ by $v=e^{-\lambda t}u$, $\lambda \geq 0$, we will deal with Cauchy problem
$$
(\lambda+\partial_t-\Delta + b \cdot \nabla)v=0, \quad v(0)=f.
$$
We multiply the equation by $e^v$ and integrate:
$$
\lambda \langle v,e^v\rangle + \langle \partial_t (e^v-1)\rangle + 4\langle (\nabla e^{\frac{v}{2}})^2\rangle + 2\langle b e^{\frac{v}{2}},\nabla e^{\frac{v}{2}}\rangle=0.
$$
By quadratic inequality,
\begin{equation*}
\lambda \langle v,e^v\rangle +  \langle \partial_t (e^v-1)\rangle + 4\langle (\nabla e^{\frac{v}{2}})^2\rangle \leq  \alpha \langle b^2 e^{v}\rangle + \frac{1}{\alpha}\langle (\nabla e^{\frac{v}{2}})^2\rangle.
\end{equation*}
Now, selecting  $\alpha=\frac{1}{\sqrt{\delta}}$ and using $b \in \mathbf{F}_\delta$, we obtain 
\begin{equation*}
\lambda \langle v,e^v\rangle + \langle \partial_t (e^v-1)\rangle + (4-2\sqrt{\delta})\langle (\nabla e^{\frac{v}{2}})^2\rangle \leq \frac{c_\delta}{\sqrt{\delta}} \langle e^v\rangle.
\end{equation*}
Using $\delta \leq 4$, we obtain after integrating in time from $0$ to $t$,
\begin{equation*}
\lambda \int_0^t \langle v,e^{v}\rangle ds + \langle e^{v(t)}-1 \rangle \leq \langle e^{f}-1 \rangle + \frac{c_\delta}{\sqrt{\delta}} \int_0^t \langle e^{v}\rangle ds. 
\end{equation*}
Replacing in the last inequality $u$ by $-u$ and adding up the resulting inequalities, we obtain
\begin{equation*}
\lambda  \int_0^t \langle v\sinh(v) \rangle ds + \langle \cosh (v(t))-1\rangle \leq \langle \cosh (f)-1\rangle + \frac{c_\delta}{\sqrt{\delta}} \int_0^t \langle \cosh(v)\rangle ds.
\end{equation*}
Finally, applying  $v\sinh(v) \geq \cosh(v)-1$, we arrive at 
\begin{equation}
\label{e3}
(\lambda-\frac{c_\delta}{\sqrt{\delta}}) \int_0^t \langle \cosh(v)-1 \rangle ds + \langle \cosh (v(t))-1\rangle \leq \langle \cosh (f)-1\rangle + \frac{c_\delta}{\sqrt{\delta}} t,
\end{equation}
where at the last step we have used the fact that volume $|\Pi^d|=1$. Let $\lambda \geq \frac{c_\delta}{\sqrt{\delta}}$. Estimate \eqref{e3} suggests that one should work in the topology determined by the ``norm'' $\langle  \cosh (v) - 1\rangle$ or, better, in the corresponding Orlicz space $L_\Phi$.

\bigskip

\section{Three approaches to constructing Feller semigroup for $-\Delta + b \cdot \nabla$}

\label{three_sect}

We have at our disposal the following approaches to constructing Feller semigroup associated with the Kolmogorov operator $-\Delta + b \cdot \nabla$ with $b \in \mathbf{F}_\delta$. These approaches require $\delta<1$ ($ \ll 1$ in high dimensions).

\medskip

(1) By using the iteration procedure of \cite{KS} for solutions $u_n$ to elliptic equations $$(\mu-\Delta + b_n\cdot \nabla)u_n=f, \quad f \in C_c,$$ with bounded smooth drifts $b_n$ approximating $b=b(x)$ (in the sense of \eqref{b_n1}, \eqref{b_n2} below). It yields inequality
$$
\|u_n-u_m\|_{\infty} \leq B \|u_n-u_m\|_2^{\gamma} \quad \text{ for some } \gamma>0 \text{ independent of $n$, $m$}.
$$
The latter, in turn,
transfers the verification of the Cauchy criterion in  $C_\infty$ to a much easier to deal with\footnote{Because it has much weaker topology than $C_\infty$ and because form-boundedness is an $L^2$ assumption on $|b|$.} space $L^2$.
The convergence of the iteration procedure depends on the uniform in $n$ gradient bound
$$
\|\nabla u_n\|_\frac{qd}{d-2} \leq C\|f\|_q, \quad q>d-2,
$$
which was also established in \cite{KS}, a pioneer work  on the elliptic regularity theory of $-\Delta + b \cdot \nabla$ with form-bounded $b$, for $q> 2 \vee (d-2)$. This bound requires $\delta<1 \wedge (\frac{2}{d-2})^2$.

We describe this approach, or rather its relatively recent counterpart for time-inhomogeneous form-bounded drift $b=b(t,x)$, in Section \ref{parab_sect}. In the time-inhomogeneous case, one obtains a Feller evolution family (see definition in Section \ref{parab_sect}), and the convergence of the iteration procedure depends on the uniform in $n$ gradient bound for a $q>d$:
\begin{equation}
\label{parab_bd}
\sup_{t \in [s,T]}\|\nabla u_n(t)\|_q^q +c\int_s^T \|\nabla u_n\|_\frac{qd}{d-2}^q dt\leq C\|\nabla f\|_q^q, 
\end{equation}
where $u_n$ is the solution to parabolic equation $(\partial_t-\Delta + b_n(t,x)\cdot\nabla u_n)=0$, $u_n(0)=f \in C_c^1$, and positive constants $c$, $C$ are independent of $n$. See \eqref{bd_2}. The time-inhomogeneous form-bounded vector fields are defined as: $
|b| \in L^2_{\loc}(\mathbb R \times \mathbb R^d)
$
and
there exist a constant $\delta>0$ and a function $0 \leq g\in L^1_{\loc}(\mathbb R)$ such that for a.e.\,$t \in \mathbb R$
\begin{align}
\label{fbb_inhom0}
\|b(t,\cdot)\varphi\|_2^2   \leq \delta \|\nabla \varphi\|_2^2 +g(t)\|\varphi\|_2^2
\end{align}
for all $\varphi \in W^{1,2}$ (Definition \ref{parab_fb_def}).

\bigskip

(2) By constructing the Feller resolvent (Theorem \ref{thm1_feller}):
\begin{align*}
(\mu-\Delta  + b \cdot \nabla)^{-1}f & :=(\mu-\Delta)^{-1}f \\
&- (\mu - \Delta)^{-\frac{1}{2}-\frac{1}{q}} Q_{p}(1 + T_p)^{-1} G_{p} (\mu - \Delta)^{-\frac{1}{2}+\frac{1}{r}}f, 
\end{align*}
where $\mu>0$,
on functions $f \in C_\infty \cap L^p$ for  $r<p<q$.
The operators $Q_p$, $T_p$, $G_p$ are bounded on $L^p$ and, moreover, $\|T_p\|_{p \rightarrow p}<1$ 
under appropriate assumptions on $\delta$, so the geometric series converges. With $p$ chosen larger than $2 \vee (d-2)$, we select $q$ sufficiently close to $p$ so that, by the Sobolev embedding theorem, the free Bessel potential $(\mu - \Delta)^{-1/2-1/q} $ will take us from $L^p$ to $C_\infty$. The difficulty in this approach is ensuring boundedness of $Q_p$, $T_p$, $G_p$ in $L^p$ given that form-boundedness is an $L^2$ assumption on $|b|$. 

We note that $\|T_p\|_{p \rightarrow p}<1$ for $p>2 \vee (d-2)$ provided that $\delta<1 \wedge (\frac{2}{d-2})^2$,
i.e.\,we arrive at the same constraint on $\delta$ as in the previous approach.

This approach was developed later in \cite{Ki_a_new_approach,Ki2}. It is arguably simpler than (1), due to the use of the linear structure of $-\Delta + b \cdot \nabla$. It also gives an explicit representation for the Feller resolvent of $-\Delta + b \cdot \nabla$. 
Throughout most of this paper, we pursue approach (2).

In the case of time-inhomogeneous form-bounded $b=b(t,x)$, one constructs solution to the inhomogeneous parabolic equation $(\mu+\partial_t-\Delta + b \cdot \nabla)u=f$ on $\mathbb R \times \mathbb R^d$ as
\begin{align*}
(\mu+\partial_t-\Delta  + b \cdot \nabla)^{-1}f & :=(\mu+\partial_t-\Delta )^{-1} \\
&- (\mu +\partial_t - \Delta)^{-\frac{1}{2}-\frac{1}{q}} Q_{p}(1 + T_p)^{-1} G_{p} (\mu +\partial_t - \Delta)^{-\frac{1}{2}+\frac{1}{r}}f,
\end{align*}
for appropriately defined \textit{parabolic} operators $Q_p$, $T_p$, $G_p$. Armed with this result, one can obtain the sought Feller evolution family for $\partial_t-\Delta +b \cdot \nabla$, also in an explicit form. This approach is developed in Section \ref{morrey_sect}.

\medskip

The two methods of constructing Feller semigroup, (1) and (2), are quite different, but they give results of more or less the same strength (see Remark \ref{grad_rem}). However, for the larger class of \textit{weakly form-bounded vector fields}, described in Sections \ref{weak_sect}-\ref{weak_sde_sect}, only approach (2) is available. That being said, the iteration procedure in (1) has degrees of flexibility that have not been fully explored yet.

\bigskip

Let us now present the third approach to constructing Feller semigroup (Feller evolution family) associated with $-\Delta + b \cdot \nabla$:

\medskip

(3) Using De Giorgi's method. Let us show that $\{u_n\}$, solutions to parabolic equations $(\partial_t-\Delta + b_n(t,x)\cdot\nabla u_n)=0$, $u_n(0)=f \in C_c^1$, constitute a Cauchy sequence in $L^\infty([0,T] \times \mathbb R^d)$. For simplicity, let us assume that $b$ has compact support (this is not necessary, see Remark \ref{non_comp_rem} below).

Set $g:=u_n-u_m$. We have
$$
\partial_t g-\Delta g + b_n \cdot \nabla g =- (b_n-b_m) \cdot \nabla u_m, \quad g(0)=0.
$$
This is a Cauchy problem for an inhomogeneous parabolic equation of the same form as \eqref{c_dg} (with $\mathsf{h}=b_n-b_m$, $f=\nabla u_m$; the fact that these are vector-valued functions is not an obstacle). Therefore, by \eqref{est_m},
\begin{align}
 \|g\|_{L^\infty([0,T] \times \mathbb R^d)} & \leq C 
\sup_{z \in \mathbb Z^d}\biggl(\int_0^T \bigg\langle \big(\mathbf{1}_{\{|b_n-b_m| \geq 1\}} + \mathbf{1}_{\{|b_n-b_m| < 1\}}|b_n-b_m|^p\big)^{\theta'} |\nabla u_m|^{p\theta' }\rho^2_z\bigg\rangle\biggr)^{\frac{1}{p\theta'}},
\label{g}
\end{align}
where $\theta<\frac{d}{d-1}$ is chosen close to $\frac{d}{d-1}$, so that $\theta'>d$ is close to $d$.
Using H\"{o}lder's inequality, we obtain that
\begin{align*}
& \|g\|_{L^\infty([0,T] \times \mathbb R^d)} \\
& \leq C 
\sup_{z \in \mathbb Z^d}\biggl(\int_0^T \bigg\langle \big(\mathbf{1}_{\{|b_n-b_m| \geq 1\}} + \mathbf{1}_{\{|b_n-b_m| < 1\}}|b_n-b_m|^p\big)^{\frac{s'}{p}} \rho_z^2\bigg\rangle\bigg)^{\frac{1}{s'}}
\biggl(\int_0^T \big\langle|\nabla u_m|^{s}\rho_z^2\big\rangle\biggr)^{\frac{1}{s}}
\end{align*}
for some $s>p\theta'>pd$ is  close to $pd$ (upon selecting $\theta'$ close to $d$).
Since $b$ has compact support, by the Dominated convergence theorem the first multiple in the RHS tends to zero as $n,m \rightarrow \infty$.
Therefore,
\begin{equation}
\label{u_conv}
\|g\|_{L^\infty([0,T] \times \mathbb R^d)} = \|u_n-u_m\|_{L^\infty([0,T] \times \mathbb R^d)} \rightarrow 0 \quad \text{ as } n,m \rightarrow \infty 
\end{equation}
if 
\begin{equation}
\label{nabla_u_m}
\sup_m \|\nabla u_m\|_{L^s([0,T] \times \mathbb R^d)}  <\infty
\end{equation} 
for $s>pd$ is  close to $pd$. By the gradient estimate \eqref{parab_bd}, after applying the interpolation inequality, we have 
\begin{equation}
\label{parab_bd___}
\sup_m \|\nabla u_m\|_{L^\frac{qd}{d-2+\frac{4}{d+2}}([0,T] \times \mathbb R^d)} \leq C \|\nabla f\|_q,
\end{equation}
so \eqref{nabla_u_m} holds with $s=\frac{qd}{d-2+\frac{4}{d+2}}$,
hence ``$s>pd$'' gives us constraint 
\begin{equation}
\label{q_gr}
\frac{qd}{d-2+\frac{4}{d+2}}>dp, \quad p>\frac{2}{2-\sqrt{\delta}}.
\end{equation}
Additionally, the gradient estimate \eqref{parab_bd} imposes its own constraint on $q$, see Remark \ref{grad_q}. 
These two constraints on $q$ (which has to be as small as possible for admissible $\delta$ to be as large as possible) ensure \eqref{u_conv} and hence the existence of the Feller semigroup (evolution family). 

The proof of \eqref{est_m} in \cite{KiS_sharp}, which gives us \eqref{g}, also uses the assumption $p \geq 2$,  so $q$ in \eqref{q_gr} is not as small as one hopes. Hence, we obtain more restrictive assumption on $\delta$ than we have in approaches (1) and (2) (see Sections \ref{basic_sde_sect} and \ref{parab_sect}). It is possible, however, to weaken ``$p \geq 2$'' by modifying \eqref{g}, which we will address elsewhere.

\begin{remark}
\label{non_comp_rem} To exclude the assumption that $b$ has compact support we estimate
\begin{align*}
& \|g\|_{L^\infty([0,T] \times \mathbb R^d)} \\
& \leq C 
\sup_{z \in \mathbb Z^d, |z| \leq R}\biggl(\int_0^T \bigg\langle \big(\mathbf{1}_{\{|b_n-b_m| \geq 1\}} + \mathbf{1}_{\{|b_n-b_m| < 1\}}|b_n-b_m|^p\big)^{\frac{s'}{p}} \rho_z^2\bigg\rangle\bigg)^{\frac{1}{s'}}
\biggl(\int_0^T \big\langle|\nabla u_m|^{s}\rho_z^2\big\rangle\biggr)^{\frac{1}{s}} \\
& + \sup_{z \in \mathbb Z^d, |z| > R}\biggl(\int_0^T \bigg\langle \big(\mathbf{1}_{\{|b_n-b_m| \geq 1\}} + \mathbf{1}_{\{|b_n-b_m| < 1\}}|b_n-b_m|^p\big)^{\frac{s'}{p}} \rho_z^2\bigg\rangle\bigg)^{\frac{1}{s'}}
\biggl(\int_0^T \big\langle|\nabla u_m|^{s}\rho_z^2\big\rangle\biggr)^{\frac{1}{s}},
\end{align*}
where the second sum can be made as small as needed by selecting $R$ sufficiently large. To show this, one needs to use instead of \eqref{parab_bd___} the estimate
$$
\sup_m \int_0^T \langle |\nabla u_m|^s\rho_z^2\rangle^{\frac{1}{s}} \leq C \langle|\nabla f|^q \rho_z^2\rangle^{\frac{1}{q}}, \quad s=\frac{qd}{d-2+\frac{4}{d+2}}
$$
with $C$ independent of $m$ and $z$
(this is a consequence of \eqref{rho_grad_est}). Here $\langle|\nabla f|^q \rho_z^2\rangle$ is small if $|z|>R$ for $R$ large compared to the support of $f$. Now, for $R$ fixed sufficiently large, we can treat the first sum in the same way as in the case of a compact support $b$.
\end{remark}

It should be noted that the methods of constructing a Feller realization of $-\Delta + b \cdot \nabla$
actually yield other regularity results. This is discussed in the rest of this paper.

\bigskip

\section{Basic result on weak well-posedness of SDEs with singular drift} 

\label{basic_sde_sect}

\textbf{1.~}The theory of SDE 
\begin{equation}
\label{sde0}
dX_t=-b(X_t)dt + \sqrt{2}dW_t, \quad X_0=x \in \mathbb R^d,
\end{equation}
with $b \in \mathbf{F}_\delta$, becomes more detailed as form-bound $\delta$ gets smaller. Namely, 
if $\delta < 1 \wedge (\frac{2}{d-2})^2$, then there is a realization of $-\Delta + b \cdot \nabla$ in $C_\infty$ that generates a Feller semigroup. The latter, in turn, determines weak solutions of \eqref{sde1}. For every $x \in \mathbb R^d$, the constructed weak solution is unique in some large classes (e.g.\,in the class of weak solutions satisfying Krylov-type estimates, or in the class of weak solutions that can be obtained via a ``reasonable approximation procedure''). See Theorems \ref{thm1_feller} and \ref{thm1} below.

\medskip

Set $b^{\frac{2}{p}}:=b|b|^{-1+\frac{2}{p}}$. For given $p \in [2,\infty[$ and $1 \leq r<p<q<\infty$, define operators ($\mu>0$)
\begin{align*}
G_p(r)&:=b^{\frac{2}{p}} \cdot \nabla (\mu -\Delta )^{-\frac{1}{2}-\frac{1}{r}},  \\
Q_p(q) \upharpoonright {\mathcal E} & :=(\mu -\Delta )^{-\frac{1}{2}+\frac{1}{q}}|b|^{1-\frac{2}{p}}, \\
T_p \upharpoonright  {\mathcal E} &:=b^{\frac{2}{p}}\cdot \nabla(\mu - \Delta)^{-1}|b|^{1-\frac{2}{p}}.
\end{align*}
We define 
the last two operators on $\mathcal E:=\bigcup_{\varepsilon>0}e^{-\varepsilon|b|}L^p,$ a dense subspace of $L^p$, to remove any question regarding the summability of $|b|^{1-\frac{2}{p}}f$, $f \in L^p$, on which we act with the Bessel potential. 

Set 
$$
c_{\delta,p}:=\bigl(\frac{p}{2}\delta + \frac{p-2}{2}\sqrt{\delta}\bigr)^{\frac{1}{p}}\bigl(p-1-(p-1)\frac{p-2}{2}\sqrt{\delta} - \frac{p(p-2)}{4}\delta\bigr)^{-\frac{1}{p}}.$$

\begin{lemma}[\cite{Ki2}]
\label{lem_GPQ}

Let $b \in \mathbf{F}_\delta$. Then for every $p \in [2,\infty[$, there exists $\mu_0=\mu_0(d,p,q)$ such that the following is true for all $\mu \geq \mu_0$:

\smallskip

{\rm (\textit{i})} $T_p \upharpoonright  {\mathcal E}$  admits extension by continuity to $L^p$, denoted by $T_p$. One has $$\|T_p\|_{p \rightarrow p} \leq c_{\delta,p}.$$
In particular, if $\delta<1$, then $$\|T_p\|_{p \rightarrow p}<1 \quad \text{for every $p \in [2,\frac{2}{\sqrt{\delta}}]$}.$$

\smallskip

{\rm (\textit{ii})} $Q_p(q) \upharpoonright {\mathcal E}$ admits extension by continuity to $L^p$, denoted by $Q_p(q)$.

\smallskip

{\rm (\textit{iii})} $G_p(r)$ is bounded on $L^p$.
\end{lemma}

Lemma \ref{lem_GPQ} is a key result needed to prove Theorem \ref{thm1_feller}. Its proof, which uses only elementary arguments, is included in Appendix \ref{app_lem_GPQ}.

\medskip

Let us fix the following approximation of $b$ by smooth vector fields having compact supports:
\begin{equation}
\label{b_n1}
b_n:= c_n \eta_{\varepsilon_n} \ast (\mathbf{1}_n b), 
\end{equation}
where $\mathbf{1}_n$ is the indicator of $\{x \mid |x| \leq n, |b(x)| \leq n\}$ and $\eta_{\varepsilon_n}$ is the Friedrichs mollifier (Section \ref{notation_sect}), and we choose $\varepsilon_n \downarrow 0$ sufficiently rapidly so that, for appropriate $c_n \uparrow 1$, one has
\begin{equation}
\label{b_n2}
b_n \rightarrow b \quad \text{ in } L^2_{\loc}(\mathbb R^d,\mathbb R^d)
\end{equation}
and
\begin{equation}
\label{b_n3}
b_n \in \mathbf{F}_{\delta} \quad \text{ with some $c_\delta$ independent of $n=1,2,\dots$}  
\end{equation}
see Appendix \ref{approx_app}.1. 
Actually, in the next theorem any bounded smooth $b_n$ satisfying \eqref{b_n2}, \eqref{b_n3} will do, not necessarily the ones defined by \eqref{b_n1}, but at the moment we save ourselves some 
efforts by considering \eqref{b_n1} since these are, essentially, cutoffs of $b$. Later, however, we will consider any bounded smooth $\{b_n\}$ satisfying \eqref{b_n2}, \eqref{b_n3}. This is important because it will give us a uniqueness result on its own: the constructed  semigroups or weak solutions to SDEs will not depend on the choice of approximating $\{b_n\}$, as long as they satisfy \eqref{b_n2}, \eqref{b_n3}. In this regard, it is worth introducing the following definition:

\begin{definition}
We say that vector fields $b_n$ satisfying \eqref{b_n3} are uniformly (in $n$) form-bounded.
\end{definition}

\begin{theorem}
\label{thm1_feller}
Let $b \in \mathbf{F}_\delta$, $\delta<1$. 
There exists $\mu_0>0$ such that the following is true:

\medskip

{\rm (\textit{i})} For every $p \in [2,\frac{2}{\sqrt{\delta}}[$, for all $\mu \geq \mu_0$ the function\footnote{Note that after expanding $(1 + T_p)^{-1}$ in \eqref{form_neu} in the geometric series, we obtain  the formal Neumann series representation for $u$.}
\begin{equation}
\label{form_neu}
u=(\mu - \Delta)^{-1}f - (\mu - \Delta)^{-\frac{1}{2}-\frac{1}{q}} Q_{p}(q) (1 + T_p)^{-1} G_{p}(r) (\mu - \Delta)^{-\frac{1}{2}+\frac{1}{r}}f, \quad f \in L^p
\end{equation}
is a weak solution to equation 
\begin{equation}
\label{eq_e}
(\mu-\Delta + b \cdot \nabla)u=f
\end{equation}
i.e.
$$
\mu \langle u,\psi\rangle + \langle \nabla u,\nabla \psi \rangle + \langle b \cdot \nabla u,\psi\rangle=\langle f,\psi\rangle \quad \text{ for all }\psi \in C_c^\infty.
$$
Moreover, if $f \in L^p \cap L^2$, then $u$ is the unique in $W^{1,2}$ weak solution to \eqref{eq_e}.

\smallskip

{\rm (\textit{ii})}  It follows from \eqref{form_neu} that
\begin{equation}
\label{grad}
u \in \mathcal W^{1+\frac{2}{q},p} \;\;(\text{Bessel potential space}), \quad q>p.
\end{equation}
In particular,  if $\delta<1 \wedge \big(\frac{2}{d-2}\bigr)^2$, then in the interval $p \in [2,\frac{2}{\sqrt{\delta}}]$ we can select $p>d-2$, and then select $q$ sufficiently close to $p$, so that by the Sobolev embedding theorem $u$ is H\"{o}lder continuous (possibly after a modification on a measure zero set), with the H\"{o}lder continuity exponent less than but arbitrarily close to $1-\frac{d-2}{p}$.

\smallskip

{\rm (\textit{iii})} The operator-valued function in \eqref{form_neu} $$\Theta_p(\mu,b):=(\mu - \Delta)^{-1}- (\mu - \Delta)^{-\frac{1}{2}-\frac{1}{q}} Q_{p}(q) (1 + T_p)^{-1} G_{p}(r) (\mu - \Delta)^{-\frac{1}{2}+\frac{1}{r}},$$ defined on $\{\mu \geq \mu_0\}$, takes values in $\mathcal B(\mathcal W^{-1+\frac{2}{r},p},\mathcal W^{1+\frac{2}{q},p})$.

\smallskip

{\rm (\textit{iv})} Let
$\delta<1 \wedge \big(\frac{2}{d-2}\bigr)^2$. 
Fix $p \in ]d-2,\frac{2}{\sqrt{\delta}}[$ if $d \geq 4$, or $p \in [2,\frac{2}{\sqrt{\delta}}[$ if $d=3$. Then
\begin{equation}
\label{feller_resolvent}
(\mu+\Lambda_{C_\infty}(b))^{-1}:=\bigl(\Theta_p(\mu,b) \upharpoonright L^p \cap C_\infty \bigr)_{C_\infty \rightarrow C_\infty}^{\clos}, \quad \mu \geq \mu_0,
\end{equation}
determines the resolvent of a Feller generator on $C_\infty$. This semigroup satisfies
\begin{equation*}
e^{-t\Lambda_{C_\infty}(b)}=\text{\small $s\text{-}C_\infty\text{-}$}\lim_n e^{-t\Lambda_{C_\infty}(b_n)} \quad \text{locally uniformly in $t \geq 0$},
\end{equation*} 
where bounded smooth $b_n$ are defined by \eqref{b_n1} and the approximating operators
$\Lambda_{C_\infty}(b_n):=-\Delta + b_n \cdot \nabla$ with domain $D(\Lambda_{C_\infty}(b_n)):=(1-\Delta)^{-1}C_\infty$ are, by a classical result, Feller generators.

\smallskip

{\rm (\textit{v})} Feller semigroup $e^{-t\Lambda_{C_\infty}(b)}$ is conservative, i.e.\,its integral kernel $e^{-t\Lambda_{C_\infty}}(x,\cdot)$ satisfies 
\begin{equation}
\label{cons}
\langle e^{-t\Lambda_{C_\infty}(b)}(x,\cdot) \rangle=1 \quad \text{ for all } x \in \mathbb R^d,t>0.
\end{equation}
\end{theorem}

Excluding the cases where $b$ is sufficiently regular, one has little information about the domain $D(\Lambda_{C_\infty}(b))$ of the Feller generator $\Lambda_{C_\infty}(b)$. But, it is easily seen, one can be certain that $C_c^2 \not \subset D(\Lambda_{C_\infty}(b))$ already if $|b| \in L^\infty - C_b$.

\begin{remark}
\label{grad_rem}
The construction of the Feller semigroup  via the iteration procedure of \cite{KS} requires gradient bounds on solutions $u_n$ of $(\mu-\Delta + b_n \cdot \nabla )u_n=f$, $f \in C_c^\infty$. Namely, it is established in \cite[proof of Lemma 5]{KS} that if $\delta<1 \wedge (\frac{2}{d-2})^2$, then for $p>2 \vee (d-2)$ close to $2 \vee (d-2)$
\begin{equation}
\label{grad_KS}
\|\nabla|\nabla u_n|^{\frac{p}{2}}\|_2^2 \leq K \|f\|_p^p,
\end{equation}
for a constant $K$ independent of $n$. (We discuss a parabolic analogue of \eqref{grad_KS} below.) Then, applying the Sobolev embedding theorem twice, one obtains that $u_n$ has H\"{o}lder continuity exponent $1-\frac{d-2}{p}$ (independent of $n$). Taking into account that $p$ satisfies a strict inequality, one thus obtains the same H\"{o}lder continuity result as in Theorem \ref{thm1_feller}(\textit{ii}). However, there is substantial difference between gradient bounds \eqref{grad} and \eqref{grad_KS}: \eqref{grad} allows to control order $>1$ derivatives of $u_n$, while \eqref{grad_KS} does not need additional strict inequalities such as ``$q>p$'' in \eqref{grad} (clearly, if in \eqref{grad} one could take $q=p$, then it would make the regularity result stronger, but one cannot do this).
\end{remark}

Above we mentioned that the approach of Theorem \ref{thm1} is somewhat simpler than \cite{KS}. One reason is that it uses to the full extent the linear structure of the equation. There is another reason: in the iteration procedure of \cite{KS} one shows that the solutions of the approximating equations constitute a Cauchy sequence, while in the proof of Theorem \ref{thm1} one already has a candidate for the limit of this sequence, cf.\,\eqref{feller_resolvent}.

\medskip 

The free Bessel potential $(\mu - \Delta)^{-\frac{1}{2}-\frac{1}{q}}$ in \eqref{form_neu} provides a ``trampoline'' from $L^p$ to $C_\infty$. One advantage of working in $L^p$ with $p$ large is that it simplifies the proof of approximation results as e.g.\,in assertion (\textit{iv}) of Theorem \ref{thm1_feller}, since it is easier to prove convergence in $L^p$ than in $C_\infty$. In fact, this ``trampoline'' was used in \cite{Ko} who considered bounded drifts of $\alpha$-stable process. There, however, working in $L^p$ is a matter of convenience (one can also stay entirely in $C_\infty$ while dealing with unbounded drifts, see \cite{BC}). Our goal, however, is to transition from an $L^2$ assumptions on the drift (i.e.\,the form-boundedness), via $L^p$, to a semigroup in $C_\infty$. It is the transition from $L^2$ to $L^p$ with $p$ large that is the most difficult one.

\medskip

The operator-valued function $\mu \mapsto \Theta_p(\mu,b) \in \mathcal B(L^p)$ in Theorem \ref{thm1_feller} is itself the resolvent of the generator of a quasi contraction semigroup in $L^p$. 
In fact, this semigroup coincides with the semigroup $e^{-t\Lambda_p(b)}$ constructed via \eqref{p_semi}, and satisfies
\begin{equation}
\label{p_semi_lim}
e^{-t\Lambda_{p}(b)}=\text{\small $s\text{-}C_\infty\text{-}$}\lim_n e^{-t\Lambda_{p}(b_n)} \quad \text{locally uniformly in $t \geq 0$},
\end{equation} 
where $\Lambda_{p}(b_n)=-\Delta + b_n\cdot \nabla$ having domain $\mathcal W^{2,p}$.
See \cite{Ki2} for details.

By the way, it is easy to write a similar to $\Theta_p(\mu,b)$ operator-valued function representation for the resolvent of $-\Delta - \nabla \cdot b$ with $b \in \mathbf{F}_\delta$, and to modify the proofs in \cite{Ki2} to work for this operator.

Combining \eqref{p_semi_lim} with Theorem \ref{thm1_feller}, one of course obtains
\begin{equation}
\label{consist}
e^{-t\Lambda_{C_\infty}(b)} \upharpoonright L^2 \cap C_\infty = e^{-t\Lambda_{p}(b)} \upharpoonright L^2 \cap C_\infty = e^{-t\Lambda(b)} \upharpoonright L^2 \cap C_\infty,
\end{equation}
where the last semigroup is provided by the KLMN theorem.

\medskip

The resolvent representations of the type $\Theta_p(\mu,b)$ were considered earlier for Schr\"{o}dinger operators with form-bounded potentials, to obtain information about Sobolev regularity of their domains and hence regularity of their eigenfunctions, see \cite{BS,LS}. Let us note, however, that a direct comparison between Feller theories of Kolmogorov and Schr\"{o}dinger operators is not possible, see Remark \ref{kato_rem}.

\medskip

The semigroup  $e^{-t\Lambda_p(b)}$ can, in fact, be defined via the limit \eqref{p_semi_lim} for all $\delta<4$, $p>\frac{2}{2-\sqrt{\delta}}$ \cite{KS}, and satisfies
\begin{equation}
\label{growth}
\|e^{-t\Lambda_p(b)}f\|_{q} \leq ce^{t\omega_p}t^{-\frac{d}{2}(\frac{1}{p}-\frac{1}{q})}\|f\|_p, \quad t>0, \quad \omega_p:=\frac{c_\delta}{2(p-1)},
\end{equation}
for all $f \in L^q \cap L^p$.
In view of \eqref{consist}, the same estimate on $e^{-t\Lambda_{C_\infty}(b)}f$, $f \in L^p \cap C_\infty$ is valid (of course, under the assumptions on $\delta$ and $p$ of Theorem \ref{thm1_feller}(\textit{iv})).
Regarding the properties of $e^{-t\Lambda_p(b)}$ (in particular, regarding extending the interval $p \in ]\frac{2}{2-\sqrt{\delta}},\infty[$ to a larger interval), see \cite{KiS_theory}.

\medskip

Nevertheless, both the upper and the lower Gaussian bounds on the heat kernel of $-\Delta + b \cdot \nabla$ are easily destroyed by form-bounded drifts. (The heat kernel is defined, up to a modification on a measure zero set, as the integral kernel $e^{-t\Lambda(b)}(x,y)$ of $e^{-t\Lambda_p(b)}$. The heat kernel does not depend on $p$.) In fact, it suffices to consider drifts 
\begin{equation}
\label{hardy2}
b(x)=\pm\sqrt{\delta}\frac{d-2}{2}|x|^{-2}x,
\end{equation}
which are form-bounded.
The singularity of \eqref{hardy2} is so strong that it introduces an extra factor $\varphi_t(y)$ in the Gaussian bounds, 
$$
c_1t^{-\frac{d}{2}}e^{-\frac{|x-y|^2}{c_2t}}\varphi_t(y) \leq e^{-t\Lambda(b)}(x,y) \leq c_3t^{-\frac{d}{2}}e^{-\frac{|x-y|^2}{c_4t}}\varphi_t(y), 
$$
where $\varphi_t(t)$ either explodes or vanishes at the origin depending on the sign in front of $\sqrt{\delta}$ in \eqref{hardy2} (moreover, the rate of explosion or vanishing is an explicit function of $\delta$) \cite{MNS}. Of course, if one considers a sum (or a series) of such drifts with singularities at different points (which is still form-bounded), the situation at the level of heat kernel bounds becomes even more complicated.
A detailed discussion on this subject can be found in \cite{KiS_heat} where the authors prove Gaussian lower and/or upper bound on the heat kernel of $-\nabla \cdot a \cdot \nabla + b \cdot \nabla$ with measurable uniformly elliptic matrix $a$ and  drift $b$ that is form-bounded or even more singular, under additional constraints on ${\rm div\,}b$. 

\medskip

Let us also note that the two-sided Gaussian bound on the heat kernel of $-\Delta + b \cdot \nabla$ or of $-\nabla \cdot a \cdot \nabla + b \cdot \nabla$ hold, without any assumptions on ${\rm div\,}b$, when $b$ is in the Kato class of vector fields\footnote{See definition in Section \ref{weak_sect}.} \cite{Z_Gaussian},
 or in the Nash class\footnote{$b=b(x)$ is in the Nash class $\mathbf{N}_\delta$ if $|b| \in L^2_{\loc}$ and
$$
\inf_{h>0}\sup_{x \in \mathbb{R}^d} \int_0^h \sqrt{e^{t\Delta}|b|^2(x) } \; \frac{dt}{\sqrt{t}} \leq \delta 
$$. It contains $b$ with $|b| \in L^p$, $p>d$, and it also contains vector fields with $|b| \not \in L^{2+\varepsilon}_{\loc}$, $\varepsilon>0$.} \cite{S2,KiS_Nash}, respectively. Moreover, one can go beyond the Kato class and prove Guassian bounds for distributional drifts, see \cite{PZ,ZZ2}. Both the Kato class and the Nash class neither contain the class $\mathbf{F}_\delta$ nor are contained in it (but the Kato class of vector fields is contained in the class of \textit{weakly form-bounded vector fields}, considered in the next section).

\begin{remark}
\label{second_deriv_rem}
It is clear from \eqref{form_neu} that we do not have (and cannot have) information about $L^p$ summability of the second derivatives of $u$ for $p>2$ large (in this regard, see discussion before Theorem \ref{thm_sharp}). However, we have weighted estimates on the second derivatives. Assume for simplicity that $b$ is bounded and smooth, so we are looking for estimates with constants that do not depend on smoothness of $b$ or its boundedness. It follows from \eqref{form_neu} that 
$$
(|b|+1)^{-1+\frac{2}{p}}(\mu-\Delta)u=(|b|+1)^{-1+\frac{2}{p}}f-\bigg(\frac{|b|}{|b|+1}\bigg)^{1-\frac{2}{p}}(1 + T_p)^{-1} T_p'\,(|b|+1)^{-1+\frac{2}{p}}f
$$
where $T_p'=b^{\frac{2}{p}} \cdot \nabla (\mu-\Delta)^{-1} (|b|+1)^{1-\frac{2}{p}}$ is bounded on $L^p$, just like $T_p$,
and so 
\begin{align*}
\|(|b|+1)^{-1+\frac{2}{p}}(\mu-\Delta)u\|_p  \leq 
(1+C)\|(|b|+1)^{-1+\frac{2}{p}}f\|_p, \quad C<\infty.
\end{align*}
Thus, if $p>2$, then at the points where $|b|$ is infinite the factor $(|b|+1)^{-1+\frac{2}{p}}$ vanishes, and so around the singular set of $b$ the information about $L^p$ summability of $(\mu-\Delta)u$ disappears, but in a controlled way.
\end{remark}

Finally, let us add that in Theorem \ref{thm1_feller}
we could also write
$$
u=(\mu - \Delta)^{-1}f - (\mu-\Delta)^{-\frac{1}{2}-\frac{1}{q}}\hat{Q}_p(q) (1 + \hat{T}_p)^{-1} \hat{G}_{p}(\mu-\Delta)^{-\frac{1}{2}+\frac{1}{r}} f,
$$
where operators $\hat{Q}_p$, $\hat{T}_p$, $\hat{R}_p(r)$ have ``classical'' form but are bounded in the weighted $L^p$ space with weight $|b|^2_1:=(|b|+1)^2$: $$\hat{Q}_p(q)=(\mu - \Delta)^{-\frac{1}{2}+\frac{1}{q}} b \cdot \quad \text{ is in } \mathcal B\big(\big[L^p(\mathbb R^d,|b|_1^2dx)\big]^d,L^p(\mathbb R^d)\big),$$
$$\hat{T}_p=\nabla(\mu-\Delta)^{-1} b \cdot \quad \text{ is in } \mathcal B\big(\big[L^p(\mathbb R^d,|b|_1^2dx)\big]^d\big),$$
$$\hat{G}_p(r)=\nabla (\mu-\Delta)^{-\frac{1}{2}-\frac{1}{r}} \quad \text{ is in } \mathcal B\big(L^p(\mathbb R^d),\big[L^p(\mathbb R^d,|b|_1^2dx)\big]^d\big)$$
where $[\cdot]^d$ denotes vector with $d$ components.
Their boundedness follows from the boundedness of $Q_p$, $T_p$, $G_p$ on $L^p=L^p(\mathbb R^d,dx)$.

\medskip

\textbf{2.~}We are now in position to prove the following result on weak well-posedness of SDE
\begin{equation}
\label{sde1_}
X_t=x-\int_0^t b(X_r)dr + \sqrt{2}W_t, \quad 0 \leq t \leq T
\end{equation}
with $x \in \mathbb R^d$ fixed.

\begin{theorem}
\label{thm1} Let $b \in \mathbf{F}_\delta$ with $\delta<1 \wedge \big(\frac{2}{d-2}\bigr)^2$.
Let $e^{-t\Lambda_{C_\infty}(b)}$ be the Feller semigroup constructed in Theorem \ref{thm1_feller}. Fix $T>0$.
The following is true:

\smallskip

{\rm (\textit{i})} There exist probability measures $\{\mathbb P_x\}_{x \in \mathbb R^d}$ on the canonical space $(C[0,T],\mathcal B_t)$ such that 
\begin{equation*}
\mathbb E_{\mathbb P_x}[f(X_t)]=(e^{-t\Lambda_{C_\infty}(b)} f)(x), \quad f \in C_\infty, \quad x \in \mathbb R^d.
\end{equation*}
For every $x \in \mathbb R^d$ the measure $\mathbb P_x$ is a weak solution to SDE \eqref{sde1_}.

\medskip

{\rm (\textit{ii})} If $\delta$ is sufficiently small and, additionally, $|b| \in L^{\frac{d}{2}+\varepsilon}$ for some $\varepsilon>0$, then the  constructed in {\rm(\textit{i})} weak solution $\mathbb P_x$ belongs to and is unique in the class of weak solutions satisfying the following Krylov-type estimate for $q>\frac{d}{2}$ sufficiently close to $\frac{d}{2}$ (depending on how small $\varepsilon$ is):
\begin{align}
\label{krylov_class}
\mathbb E_{\mathbb P_x}\int_0^T |h(t,\omega_t)|dt  \leq c\|h\|_{L^q([0,T] \times \mathbb R^d)}.
\end{align}
for all $h \in C_c(\mathbb R^{d+1})$, for generic constant $c$.
\end{theorem}

\begin{remark}
\label{approx_uniq}
Arguing as in \cite{KiS_brownian}, one can also prove the following ``approximation uniqueness'' result. If $\{\mathbb Q_x\}_{x \in \mathbb R^d}$ is another weak solution to \eqref{sde1_} such that
$$
\mathbb Q_x=w{\mbox-}\lim_n \mathbb P_x(\tilde{b}_n) \quad \text{for every $x \in \mathbb R^d$},
$$
for some $\{\tilde{b}_n\} \subset \mathbf{F}_{\delta_1}$ with $\delta<1 \wedge \big(\frac{2}{d-2}\bigr)^2$ and $c_{\delta}$ independent of $n$, then $\{\mathbb Q_x\}_{x \in \mathbb R^d}=\{\mathbb P_x\}_{x \in \mathbb R^d}.$  In other words, the constructed weak solutions $\mathbb P_x$ of \eqref{sde1_} are unique among those that can be obtained via an approximation procedure. Note that we do not require here any convergence of $\tilde{b}_n$ to $b$.
\end{remark}

In Section \ref{weak_sde_sect} we will discuss analogues of Theorems \ref{thm1_feller}, \ref{thm1} for drifts that can be more singular than the form-bounded drifts. See Theorems \ref{thm1_feller_weak}, \ref{thm1_weak}. This, however, will come at the cost of imposing more restrictive assumption on $\delta$, and losing the possibility to include discontinuous diffusion coefficients, as we do for the form-bounded drifts in end of this section. Also, while the proof of the analogue of Lemma \ref{lem_GPQ} in Section \ref{weak_sde_sect} (i.e.\,Lemma \ref{lem_GPQ_weak} there) relies on some operator inequalities for fractional powers of the Laplacian, the proof of Lemma \ref{lem_GPQ} uses only elementary arguments.

\begin{remark}
\label{elliptic_rem}
The proof of Theorem \ref{thm1} and the construction of the Feller semigroup in Theorem \ref{thm1_feller}(\textit{iv})-(\textit{v})  rely on the elliptic regularity result of Theorem \ref{thm1_feller}(\textit{i})-(\textit{iii}). In the next sections we will be working directly in the parabolic setting, thus avoiding the use of the Trotter approximation theorem and, generally speaking, arriving at shorter proofs. However, by working with resolvents and using the semigroup theory (i.e.\,Trotter's theorem),  we can construct a  Feller semigroup departing from $L^p$ for a smaller $p$, which leads to less restrictive assumptions on $\delta$. For instance, one can compare Theorem \ref{thm1_feller}(\textit{iv}) where $p$ is chosen to be strictly greater than $2 \vee (d-2)$, and we require $\delta<1 \wedge (\frac{2}{d-2})^2$, and Theorem \ref{feller_parab} where $p$ has to be strictly greater than $d$, and $\delta<\frac{1}{d^2}$; or, better, $\delta$ satisfies \eqref{c5} below, which is still more restrictive than $\delta<1 \wedge (\frac{2}{d-2})^2$.
\end{remark}

\medskip

\section{Proof of Theorem \ref{thm1_feller}}

\label{thm1_feller_proof_sect}

Assertion (\textit{iii}) follows right away from Lemma \ref{lem_GPQ}. The first part of assertion (\textit{i}), and assertion (\textit{ii}), follow from (\textit{iii}). The proof of the second part of (\textit{i}), i.e.\,the characterization of $u$ as the unique weak solution to the elliptic equation, 
is standard and we will attend to it in the end.

(\textit{iv}) For every $n=1,2,\dots$, the operator-valued function $\Theta_p(\mu,b_n)$ is a pseudo-resolvent, i.e.\,it satisfies
\begin{equation}
\label{theta0}
\Theta_p(\mu,b_n) - \Theta_p(\eta,b_n) = (\nu - \mu) \Theta_p(\mu,b_n)\Theta_p(\nu,b_n), \quad \mu,\nu \geq \mu_0,
\end{equation}
where $\mu_0$ is from Lemma \ref{lem_GPQ}.
Identity \eqref{theta0} is verified via direct calculation. See \cite[proof of Prop.\,2]{Ki_a_new_approach} for details\footnote{Alternatively, one can verify, using the KLMN theorem, that $\Theta_2(\mu,b_n)$, $\mu \geq \mu_0$, for a $\mu_0>0$ independent of $n$, is the resolvent of $-\Delta + b_n \cdot \nabla$ in $L^2$, and so \eqref{theta0} holds on $L^2 \cap L^p$ and hence on $L^p$. But we do not need the $L^2$ theory of $-\Delta + b \cdot \nabla$ in the proof of Theorem \ref{thm1_feller}}.

By the classical theory, for every $n=1,2,\dots$, the resolvent of the approximating operator $(\mu+\Lambda_{C_\infty}(b_n))^{-1}$
is defined on $\{\mu \geq \mu_n\}$, where $\mu_n$ depends e.g.\,on $\|b_n\|_\infty$.
Our first observation is that we can replace $\mu_n$ by a $\mu_0$ independent of $n$ by establishing a link between $(\mu+\Lambda_{C_\infty}(b_n))^{-1}$ and the operator-valued function $\Theta_p(\mu,b_n)$. 
That is,
\begin{equation}
\label{a1}
(\mu+\Lambda_{C_\infty}(b_n))^{-1} \upharpoonright \mathcal S = \Theta_p(\mu,b_n) \upharpoonright \mathcal S \quad \text {for all $\mu \geq \mu_0$,}
\end{equation} 
for some $\mu_0$ independent of $n$. Indeed, we have 
\begin{equation*}
\Theta_p(\mu_n,b_n) \upharpoonright \mathcal S =(\mu_n+\Lambda_{C_\infty}(b_n))^{-1} \upharpoonright \mathcal S
\end{equation*}
for all sufficiently large $\mu_n$.
By $\Theta_p(\mu,b_n)\mathcal S\subset \mathcal S$, the previous identity and the resolvent identity \eqref{theta0},
$$
\Theta_p(\mu,b_n) \upharpoonright \mathcal S= (\mu_n+\Lambda_{C_\infty}(b_n))^{-1}\bigl(1 + (\mu_n - \mu) \Theta_p(\mu,b_n)\bigr) \upharpoonright \mathcal S, \quad \mu \geq \mu_0,
$$
so $\Theta_p(\mu,b_n) \upharpoonright \mathcal S$ is the right inverse of $\mu + \Lambda_{C_\infty}(b_n) \upharpoonright \mathcal S$ on $\mu \geq \mu_0$. Similarly, it is seen that
$\Theta_p(\mu,b_n)|_{\mathcal S}$ is the left inverse of $\mu + \Lambda_{C_\infty}(b_n) \upharpoonright \mathcal S$ on $\mu \geq \mu_0$. This gives \eqref{a1}.

Second, let us show that for every $\mu \geq \mu_0$
\begin{equation}
\label{a2}
\Theta_p(\mu,b)\mathcal S \subset L^p \cap C_\infty  \quad \text{(after a modification on a measure zero set)}, 
\end{equation}
and
\begin{equation}
\label{a2_}
\Theta_p(\mu,b_n) \overset{s}{\rightarrow} \Theta_p(\mu,b) \text{ in $L^p \cap C_\infty$} \quad (n \rightarrow \infty).
\end{equation}
The inclusion into $C_\infty$ in \eqref{a2} is immediate due to the factor $(\mu - \Delta)^{-\frac{1}{2}-\frac{1}{q}}$ in the definition of $\Theta_p(\mu,b_n)$, upon applying the Sobolev embedding theorem (here we use the assumption $p>2 \vee (d-2)$ and the choice of $q>p$ close to $p$). The second assertion \eqref{a2} is proved using again the Sobolev embedding theorem and the convergence
\begin{equation*}
Q_p(q,b_n) \overset{s}{\rightarrow} Q_p(q,b), \quad T_p(b_n) \overset{s}{\rightarrow} T_p(b), \quad G_p(r,b_n) \overset{s}{\rightarrow} G_p(r,b).
\end{equation*}
The latter, in turn, follow from the Dominated convergence theorem since $b_n$ defined by \eqref{b_n1} are, essentially, cutoffs of $b$ (for details, if needed, see the proof of \cite[Prop.\,7]{Ki_a_new_approach}).

Third, we have
\begin{equation}
\label{a3}
\sup_n\|\mu(\mu+\Lambda_{C_\infty}(b_n))^{-1}\|_{\infty \rightarrow \infty} \leq 1 \quad \text{ for all $\mu \geq \mu_0$.}
\end{equation}
Indeed, for every $n=1,2,\dots$, the semigroup $e^{-t\Lambda_{C_\infty}(b_n)}$ is an $L^\infty$ contraction, so, integrating $\|e^{-\mu t}e^{-t\Lambda_{C_\infty}(b_n)}f\|_\infty \leq e^{-\mu t}\|f\|_\infty$ in $t$ from $0$ to $\infty$ we arrive at \eqref{a3}.

Fourth, we note that
\begin{equation}
\label{a4}
\mu \Theta_{p}(\mu,b_n) \overset{s}{\rightarrow} 1 \text{  in $C_\infty$ as $\mu \uparrow \infty$ uniformly in $n$}.
\end{equation}
Indeed, in view of \eqref{a3}, since $\mathcal S$ is dense in $C_\infty$, it suffices to prove that $\mu \Theta_{p}(\mu,b_n)f \rightarrow f$ in $C_\infty$ as $\mu\uparrow \infty$ for all $f \in \mathcal S$. In turn, since $\lim_{\mu\to\infty}\|\mu(\mu-\Delta)^{-1}f-f\|_\infty=0$,
it suffices to show that $\sup_n\|\mu \Theta_p(\mu,b_n)f- \mu(\mu-\Delta)^{-1}f\|_\infty \rightarrow 0$ as $\mu\uparrow \infty$. We have
\begin{align*}
\Theta_{p}(\mu,b_n)f &- (\mu-\Delta)^{-1}f =\\
& -(\mu-\Delta)^{-\frac{1}{2}-\frac{1}{q}} Q_{p}(q)   
\bigl(1+T_{p} \bigr)^{-1} b_n^{\frac{2}{p}} \cdot \nabla (\lambda-\Delta)^{-1} (\mu-\Delta)^{-1}  (\lambda-\Delta) f,
\end{align*}
with $q>p$
where $\lambda$ is sufficiently large but fixed. 
Note, that $$\|(\mu-\Delta)^{-\frac{1}{2}-\frac{1}{q}}\|_{p \rightarrow \infty} \leq c \mu^{-\frac{1}{2}+\frac{d}{2p}-\frac{1}{q}},$$ and $$
\|b_n^{\frac{2}{p}} \cdot \nabla (\lambda-\Delta)^{-1}\|_{p \rightarrow p} \leq c_1 
$$
with $c_1$ independent of $n$ and $\mu$
(since $\|G_p(r)\|_{p \rightarrow p} \equiv \|b_n^{\frac{2}{p}} \cdot \nabla (\lambda-\Delta)^{-\frac{1}{2}-\frac{1}{r}}\|_{p \rightarrow p} $ is uniformly bounded in $n$). Thus
$$
\|\Theta_p(\mu,b_n)f - (\mu-\Delta)^{-1}f\|_\infty 
\leq C \mu^{-\frac{1}{2}+\frac{d}{2p} -\frac{1}{q}} \mu^{ -1}  \|(\lambda-\Delta) f\|_p.
$$
Since $p>2 \vee (d-2)$, we select $q$ sufficiently close to $p$ so that $-\frac{1}{2}+\frac{d}{2p}-\frac{1}{q}-1<-1,$ and hence $\sup_n\|\mu \Theta_p(\mu,b_n)f - \mu(\mu-\Delta)^{-1}f\|_\infty \rightarrow 0$ as $\mu\uparrow \infty$, which yields \eqref{a4}.

We now prove assertion (\textit{iv}) of the theorem using the Trotter approximation theorem (Appendix \ref{trotter_sect}). Its conditions 
$$\text{$\sup_n\|(\mu+\Lambda_{C_\infty}(b_n))^{-1}\|_{\infty \rightarrow \infty} \leq \frac{1}{\mu}$ \quad for all $\mu \geq \mu_0$,}$$
$$\text{there exists $s\text{-}C_\infty\text{-}\lim_n (\mu+\Lambda_{C_\infty}(b_n))^{-1}$ \quad for some $\mu \geq \mu_0$},$$
$$\text{$\mu (\mu+\Lambda_{C_\infty}(b_n))^{-1} \rightarrow 1$ \quad in $C_\infty$ as $\mu \uparrow \infty$ uniformly in $n$}$$

\noindent are verified in \eqref{a1}-\eqref{a4}. So, Trotter's theorem yields (\textit{iv}) including the strong convergence of semigroups in $C_\infty$.

\medskip

(\textit{v}) The proof of \eqref{cons} uses the localized estimate
\begin{equation}
\label{j_1_w}
\bigl\|\rho (\mu+\Lambda_{C_\infty}(b_n))^{-1} h \bigr\|_{\infty} \leq C_3\|\rho h\|_p,
\end{equation}
where $$\rho(x):=(1+\kappa |x|^2)^{-\nu}$$
with $\nu > \frac{d}{2p}$ fixed (then $\rho \in L^p$) and $\kappa>0$ to be chosen sufficiently small. We comment on the proof of \eqref{j_1_w} below. Estimate \eqref{j_1_w} yields: for every fixed $x \in \mathbb R^d$  there is $C$ such that $$|(\mu+\Lambda_{C_\infty}(b_n))^{-1} h(x)| \leq C\|\rho h\|_p.$$ By considering an increasing sequence $h \uparrow 1-\mathbf{1}_{B_R(0)}$ we obtain
$$
\langle (\mu+\Lambda_{C_\infty}(b_n))^{-1}(x,\cdot) (1-\mathbf{1}_{B_R(0)}(\cdot)) \rangle \leq C\|\rho(1-\mathbf{1}_{B_R(0)}\|_p,
$$
where, as is evident from the definition of $\rho$, the right-hand side can be made smaller than any $\varepsilon$ uniformly in $n$ by selecting radius $R>0$ sufficiently large.

 Since $\langle (\mu+\Lambda_{C_\infty}(b_n))^{-1}(x,\cdot) \rangle=\mu^{-1}$, $n=1,2,\dots$, we have
$$\mu^{-1}-\varepsilon \leq \langle (\mu+\Lambda_{C_\infty}(b_n))^{-1}(x,\cdot) \mathbf{1}_{B_R(0)}(\cdot) \rangle \leq \mu^{-1}.$$
By passing to the limit in $n$ we obtain $ \mu^{-1}-\varepsilon \leq \langle (\mu+\Lambda_{C_\infty}(b))^{-1}(x,\cdot) \mathbf{1}_{B_R(0)}(\cdot) \rangle \leq \mu$. Finally, sending $R \rightarrow \infty$, and then $\varepsilon \downarrow 0$, we arrive at $$\langle (\mu+\Lambda_{C_\infty}(b))^{-1}(x,\cdot)  \rangle = \mu^{-1}$$ which gives us \eqref{cons}.

Regarding the proof of \eqref{j_1_w}, we can either commute $\rho$ with the operators that constitute $\Theta_p(\mu,b_n) \equiv (\mu+\Lambda_{C_\infty}(b_n))^{-1}$ (see \cite{KiS_brownian} for details), or we cannote that the equation for $u_n$, i.e.
$
(\mu - \Delta + b_n \cdot \nabla)u_n =h,
$
can be rewritten as 
$$
\mu \rho u_n  - \Delta (\rho u_n) + \tilde{b}_n \cdot \nabla (\rho u_n) = \rho h + K,
$$
where $$\tilde{b}_n:=b_n+2\frac{\nabla \rho}{\rho}, \qquad
K=2 \frac{(\nabla \rho)^2}{\rho} u_n + (-\Delta \rho)u_n + b_n u_n \cdot \nabla \rho.
$$
Now, we apply bounds
$$
|\nabla \rho(x)| \leq 2\nu\sqrt{\kappa}\rho(x), \quad \frac{|\nabla \rho(x)|^2}{\rho(x)} \leq 4\nu^2\kappa
\rho(x), \quad |\Delta \rho(x)| \leq \big(4\nu^2 + (4+2d)\nu\big) \kappa \rho(x), \quad   x \in \mathbb R^d.
$$
By selecting $\kappa$ sufficiently small, one can make the form-bound of $\tilde{b}$ as close to the form-bound $\delta$ of $b$ as needed. Furthermore, the first two terms in $K$ can be absorbed by $\mu \rho u_n$ (at the expense of replacing $\mu$ with $\mu-\mu_1$ for appropriate $\mu_1=\mu_1(\kappa)>0$). We have $$\rho u_n=(\mu-\mu_1+\Lambda_{C_\infty}(\tilde{b}_n))^{-1} (\rho h + b_n u_n \cdot \nabla \rho)=\Theta_p(\mu-\mu_1,\tilde{b}_n)(\rho h + b_n u_n \cdot \nabla \rho), \quad \mu > \mu_1,$$
so we can apply Lemma \ref{lem_GPQ} to obtain \eqref{j_1_w}. See details in \cite{Ki_Morrey}.

\noindent\textit{Returning to the proof of {\rm (\textit{i})}}, we note that 
$u_n$ coincides with the classical solution to 
$
(\mu-\Delta + b_n\cdot \nabla)u_n=f.
$
Now, using convergence
\begin{equation*}
Q_p(q,b_n) \overset{s}{\rightarrow} Q_p(q,b), \quad T_p(b_n) \overset{s}{\rightarrow} T_p(b), \quad G_p(r,b_n) \overset{s}{\rightarrow} G_p(r,b),
\end{equation*}
(see the prof of \eqref{a2_}) it is easy to pass to the limit in $n$ in
$$
\mu \langle u_n,\psi\rangle + \langle \nabla u_n,\nabla \psi \rangle + \langle b_n \cdot \nabla u_n,\psi\rangle=\langle f,\psi\rangle \quad \text{ for all }\psi \in C_c^\infty,
$$
which shows that $u$ is a weak solution.
A standard argument (i.e.\,the Lax-Milgram theorem) yields that $u$ is the unique weak solution.
\hfill \qed

\section{Proof of Theorem \ref{thm1}}
\label{thm1_proof_sect}

(\textit{i}) The following estimate will be needed: for all $h \in C_c$ and all $\mu \geq \mu_0$
\begin{equation}
\label{j_2}
\bigl\|(\mu+\Lambda_{C_\infty}(b))^{-1}|b_m| h \bigr\|_{\infty} \leq C_1\||b_m|^{\frac{2}{p}}h\|_p,  
\end{equation}
\begin{equation}
\label{j_3}
\|(\mu+\Lambda_{C_\infty}(b))^{-1}|b_m-b_n|h\|_\infty  \leq C_2 \bigl\| |b_m-b_n|^{\frac{2}{p}} h\bigr\|_p.
\end{equation}
for appropriate constants $C_i=C_i(\delta,p)$, $i=1,2$, where $p>2 \vee (d-2)$ is fixed.
These estimates follow right away from the construction of $(\mu+\Lambda_{C_\infty}(b))^{-1}$ via $\Theta_p(\mu,b)$ in Theorem \ref{thm1_feller}(\textit{iv}). Namely, since $|b_m|h \in C_c$, 
\begin{align}
(\mu+\Lambda_{C_\infty}(b))^{-1}|b_m|h & = \Theta_p(\mu,b)|b_m|h \notag \\
& = (\mu - \Delta)^{-\frac{1}{2}-\frac{1}{q}} G_p(q,b_m)|b_m|^{\frac{2}{p}}h \notag \\
& - (\mu - \Delta)^{-\frac{1}{2}-\frac{1}{q}} Q_{p}(q,b) (1 + T_p(b))^{-1} T_p(b,b_m)|b_m|^{\frac{2}{p}}h. \label{j_2_repr}
\end{align}
The operator $T_p(b, b_m):=b^{\frac{2}{p}} \cdot \nabla (\mu-\Delta)|b_m|^{1-\frac{2}{p}}$ is ``almost $T_p(b)$''. In fact, repeating the proof of Lemma \ref{lem_GPQ}(\textit{i}), we obtain $\|T_p(b, b_m)\|_{p \rightarrow p} \leq c_\delta'$ with constant $c_\delta'$ independent of $m$. Now, applying Lemma \ref{lem_GPQ} in \eqref{j_2_repr} and using the Sobolev embedding theorem (recall $p>2 \vee (d-2)$ and $q>p$ is close to $p$), we obtain \eqref{j_2}. The same argument gives \eqref{j_3}.

By a standard result (see e.g.\,\cite[Sect.\,I.9]{BGe}), given a conservative Feller semigroup $e^{-t\Lambda_{C_\infty}(b)}$, there exist probability measures $\mathbb P_x$ $(x \in \mathbb R^d)$ on $(D(\mathbb R_+,\mathbb R^d),\mathcal B'_t=\sigma(\omega_r, 0 \leq r \leq t))$, where $D([0,T],\mathbb R^d)$ is the space of right-continuous functions having left limits, and $\omega_t$ is the coordinate process on $D(\mathbb R_+,\mathbb R^d)$, such that
$$
\mathbb E_{x}[f(\omega_t)]=e^{-t\Lambda_{C_\infty}(b)}f(x), \quad f \in C_\infty, \quad t>0.
$$
Here and below, $\mathbb E_x:=\mathbb E_{\mathbb P_x}$. We will show that $\mathbb P_x$ are actually concentrated on $(C(\mathbb R_+,\mathbb R^d),\mathcal B_t)$.

For every $n=1,2,\dots$, let $X^n_t=X^n_{t,x}$ denote the strong solution to the approximating SDE
$$
X^n_t=x-\int_0^t b_n(s,X^n_s)ds + \sqrt{2}W_t, \quad x \in \mathbb R^d,
$$
on a complete probability space ($\Omega$, $\mathcal F$, $\mathcal F_t$, $\mathbf{P}$). Put $\mathbb P_x^n:=(\mathbf PX^n)^{-1}$, $n=1,2,\dots$, and set $\mathbb E_x^n:=\mathbb E_{\mathbb P_x^n}$.

\smallskip

Fix $\mu \geq \mu_0$. In what follows, $0<t\leq T<\infty$.
For every $g \in C_c^2$, the following is true:

\medskip

$(\mathbf a)$ $\mathbb E_x \int_0^t \bigl|b\cdot\nabla g \bigr|(\omega_s)ds<\infty$.

\smallskip

\noindent Indeed, since $b_n \rightarrow b$ everywhere outside of a measure zero set, we have by Fatou's lemma\footnote{When applying Fatou's lemma, we use
$$
\mathbb E_x\int_0^t \bigl|b\cdot\nabla g \bigr|(\omega_s)ds=\mathbb E_x\int_0^t \liminf_n\bigl|b_n\cdot\nabla g \bigr|(\omega_s)ds.
$$
Indeed, $\xi:=\bigl|b\cdot\nabla g \bigr| - \liminf_n\bigl|b_n\cdot\nabla g \bigr|=0$ a.e.\,on $\mathbb R^d$, 
but $|\mathbb E_x\int_0^t \xi(\omega_s)ds |=0$ (as follows e.g.\,by representing $\{\xi \neq 0\}=\cap_{k}U_k$ for a decreasing sequence of open sets $U_k$ such that $|U_k| \downarrow 0$, smoothing out $\mathbf{1}_{U_k}$ by replacing it with $e^{\varepsilon_k\Delta}\mathbf{1}_{U_k}$ with $\varepsilon_k \downarrow 0$ rapidly, and then applying $\mathbb E_x\int_0^t e^{\varepsilon_k\Delta}\mathbf{1}_{U_k}(\omega_s)ds \leq e^{\mu T}(\mu+\Lambda_{C_\infty}(b))^{-1}e^{\varepsilon_k\Delta}\mathbf{1}_{U_k}(x) \leq Ce^{\mu T}\|e^{\varepsilon_k\Delta}\mathbf{1}_{U_k}\|_p \downarrow 0$ as $k \rightarrow \infty$. The last inequality follows from the construction of $(\mu+\Lambda_{C_\infty}(b))^{-1}$ via the operator-valued function $\Theta_p(\mu,b)$.}
\begin{align*}
& \notag \mathbb E_x \int_0^t \bigl|b\cdot\nabla g \bigr|(\omega_s)ds  \\
& \notag \leq \liminf_n \mathbb E_x \int_0^t \bigl|b_n\cdot\nabla g \bigr|(\omega_s)ds  = \liminf_n  \int_0^t e^{-s\Lambda_{C_\infty}(b)}\bigl|b_n\cdot\nabla g \bigr|(x)ds   \\
& \notag \leq e^{\mu T} \liminf_n (\mu+\Lambda_{C_\infty}(b))^{-1}|b_n| |\nabla g|(x)
\end{align*}
Now, applying \eqref{j_2} with $h=|\nabla g|$,
we obtain
\begin{align*}
\mathbb E_x \int_0^t \bigl|b\cdot\nabla g \bigr|(\omega_s)ds & \leq C_1 e^{\mu T}\liminf_n  \langle |b_n|^2 |\nabla g|^p \rangle^{\frac{2}{p}} \\
&  = C_1 e^{\mu T}\langle |b|^2 |\nabla g|^p \rangle^\frac{2}{p} <\infty \quad (\text{by $|b| \in L^2_{\loc}$}).
\end{align*}

\medskip

$(\mathbf b)$ $$\mathbb E_x\int_0^t (b_n\cdot\nabla g)(\omega_s)ds - \mathbb E^n_x\int_0^t (b_n\cdot\nabla g)(\omega_s)ds\rightarrow 0 \quad \text{ as } n \rightarrow \infty.$$
Indeed, we have:
\begin{align*}
& \mathbb E_x\int_0^t (b_n\cdot\nabla g)(\omega_s)ds - \mathbb E^n_x\int_0^t (b_n\cdot\nabla g)(\omega_s)ds \\
& = \int_0^t \left( e^{-s\Lambda_{C_\infty}(b)} - e^{-s\Lambda_{C_\infty}(b_n)} \right) (b_n\cdot\nabla g)(x)ds \\
& =\int_0^t \left( e^{-s\Lambda_{C_\infty}(b)} - e^{-s\Lambda_{C_\infty}(b_n)} \right) ((b_n-b_m)\cdot\nabla g)(x)ds  \\
& + \int_0^t \left( e^{-s\Lambda_{C_\infty}(b)} - e^{-s\Lambda_{C_\infty}(b_n)} \right) (b_m\cdot\nabla g)(x)ds \\
& =: S_1 + S_2,
\end{align*}
where $m$ is to be chosen. 
Reducing the estimates on the expectations of time integrals to the estimates on resolvents as in the proof of $(\mathbf a)$, we obtain:
\begin{align*}
S_1(x) & \leq   
e^{\mu T} (\mu+\Lambda_{C_\infty}(b))^{-1}|(b_n-b_m)\cdot\nabla g|(x)  + e^{\mu T} (\mu+\Lambda_{C_\infty}(b_n))^{-1}|(b_n-b_m)\cdot\nabla g|(x).
\end{align*}
Using \eqref{j_3} and the convergence $b_n-b_m \rightarrow 0$ in $L^2_{\loc}$ as $n,m \uparrow \infty$, we obtain $S_1 \rightarrow 0$ as $n,m \uparrow \infty$.
Now, let us fix a sufficiently large $m$. Since $e^{-s\Lambda_{C_\infty}(b)}=s\text{-}C_\infty\text{-}\lim_n e^{-s\Lambda_{C_\infty}(b_n)}$ uniformly in $0 \leq s \leq T$ (i.e.\,assertion (\textit{iv}) of Theorem \ref{thm1_feller}),
we have $S_2 \rightarrow 0$ as $n \uparrow \infty$. 
The proof of $(\mathbf b)$ is completed.

\medskip

$(\mathbf c)$ 
$$
\mathbb E^n_x[g(\omega_t)] \rightarrow \mathbb E_x[g(\omega_t)] 
$$
and
$$\mathbb E_x\int_0^t (b\cdot\nabla g)(\omega_s)ds - \mathbb E^n_x\int_0^t (b_n\cdot\nabla g)(\omega_s)ds\rightarrow 0 \quad \text{ as } n \rightarrow \infty.$$
The first convergence is direct a consequence of  $e^{-s\Lambda_{C_\infty}(b)}=s\text{-}C_\infty\text{-}\lim_n e^{-s\Lambda_{C_\infty}(b_n)}$ uniformly in $0 \leq s \leq T$. The second convergence is a consequence of (\textbf{b}) and $\mathbb E_x\int_0^t ((b_n-b) \cdot\nabla g)(\omega_s)ds \rightarrow 0$ as $n \rightarrow \infty$, as follows from \eqref{j_3} upon applying Fatou's lemma in $m$ there.

\medskip

Now, since
$$
M_{r,m}^{g}:=g(\omega_r)-g(x) + \int_0^r (-\Delta g + b_m \cdot \nabla g)(\omega_t)dt
$$
is a $\mathcal B'_r$-martingale under $\mathbb P^m_x$, 
$$
x \mapsto \mathbb E^m_x[g(\omega_r)] - g(x) +\mathbb E^m_x\int_0^r (-\Delta g + b_m\cdot\nabla g)(\omega_t)dt \quad \text{ is identically zero on } \mathbb R^d,
$$
and so by (\textbf{c})
$$
x \mapsto \mathbb E_x[g(\omega_r)] - g(x) +\mathbb E_x\int_0^r (-\Delta g + b\cdot\nabla g)(\omega_t)dt \quad  \text{ is identically zero in } \mathbb R^d.
$$
Since $\{\mathbb P_x\}_{x \in \mathbb R^d}$ are determined by a Feller semigroup, and thus constitute a Markov process, we can conclude (see e.g.\,the proof of \cite[Lemma 2.2]{Kr1}) that 
$$
M^g_r:=g(\omega_r) - g(x) + \int_0^r (-\Delta g + b\cdot\nabla g)(\omega_t)dt,
$$
 is a $\mathcal B'_r$-martingale under $\mathbb P_x$.

\medskip

Let us show now that
$\{\mathbb P_x\}_{x \in \mathbb R^d}$ are concentrated on 
$(C([0,T],\mathbb R^d),\mathcal B_t)$.
Since $\omega_t$ is a semimartingale under $\mathbb P_x$, It\^{o}'s formula yields, for every $g \in C_c^\infty(\mathbb R^d)$, that
\begin{equation}
\label{g_}
g(\omega_t)-g(x)=\sum_{s \leq t}\bigl(g(\omega_s)-g(\omega_{s-})\bigr) + S_t,
\end{equation}
where $S_t$ is defined in terms of some  integrals and sums of $(\partial_{x_i}g)(\omega_{s-})$ and $(\partial_{x_i}\partial_{x_j}g)(\omega_{s-})$ in $s$, see \cite[Sect.\,2]{CKS} for details.
Now,
let $A$, $B$ be arbitrary compact sets in $\mathbb R^d$ such that $\dist(A,B)>0$. 
Fix $g \in C_c^\infty(\mathbb R^d)$ that separates $A$, $B$, say, $g = 0$ on $A$, $g = 1$ on $B$. Set
$$
K^g_t:=\int_0^t \mathbf{1}_A(\omega_{s-})dM_s.
$$
In view of \eqref{g_}, when evaluating $K^g_t$ one needs to integrate $\mathbf{1}_A(\omega_{s-})$ with respect to $S_t$, however, one obtains zero since $(\partial_{x_i}g)(\omega_{s-})=(\partial_{x_i}\partial_{x_j}g)(\omega_{s-})=0$ if $\omega_{s-} \in A$. Thus,
\begin{align*}
K^g_t &=\sum_{s \leq t} \mathbf{1}_A \left(\omega_{s-}\right)g(\omega_{s}) +
\int_0^t \mathbf{1}_A(\omega_{s-})\bigl(-\Delta g + b\cdot\nabla g \bigr)(\omega_s)ds \\
&=\sum_{s \leq t} \mathbf{1}_A \left(\omega_{s-}\right)g(\omega_s).
\end{align*}
Since $M^g_t$ is a martingale, so is $K^g_t$. Thus, $\mathbb{E}_x\bigl[\sum_{s \leq t} \mathbf{1}_A (\omega_{s-})g(\omega_s)\bigr]=0.$ Using the Dominated convergence theorem, we further obtain
$\mathbb{E}_x\bigl[\sum_{s \leq t} \mathbf{1}_A (\omega_{s-})\mathbf{1}_B(\omega_s)\bigr]=0$, which yields the required. By the way, this construction, in a more general form, was used to control the jumps of stable process perturbed by a drift, see \cite{CKS}.

\medskip

We denote the restriction of $\mathbb P_x$ from $(D([0,T],\mathbb R^d), \mathcal B_t')$  to $(C([0,T],\mathbb R^d),\mathcal B_t)$ again by $\mathbb P_x$, and thus obtain that
for every $x \in \mathbb R^d$ and all $g \in C_c^2(\mathbb R^d)$ 
$$
M_r^g=g(\omega_r) - g(x) + \int_0^r (-\Delta g + b\cdot\nabla g)(\omega_t)dt, \quad \omega \in C([0,T],\mathbb R^d),$$
is a $\mathcal B_r$-martingale under $\mathbb P_x$. Thus, $\mathbb P_x$ is a martingale solution to \eqref{sde1_}.

To show that $\mathbb P_x$ is a weak solution it suffices to show that $M_r^g$ is also a martingale for $g(x)=x_i$ and $g(x)=x_ix_j$ (proving along the way that $\mathbb E_x \int_0^t |b|(X(s))ds<\infty$), which can be done by following closely \cite[proof of Lemma 6]{KiS_brownian}. 

\smallskip

(\textit{ii})  is obtained via a simple modification of the proof of the uniqueness result of Theorem \ref{thm2_morrey}(\textit{iv}) below. \hfill \qed

\bigskip

\section{Time-inhomogeneous form-bounded drifts and Feller theory via iterations} 

\label{parab_sect}

\textbf{1.~}The following is the time-inhomogeneous counterpart of Definition \ref{ellip_fb_def}.

\begin{definition}
\label{def1}
\label{parab_fb_def}
A Borel measurable vector field $b:\mathbb R_+ \times \mathbb R^d \rightarrow \mathbb R^d$  is said to be form-bounded  if 
$$
|b| \in L^2_{\loc}(\mathbb R_+ \times \mathbb R^d)
$$
and
there exist a constant $\delta>0$ and a function $0 \leq g\in L^1_{\loc}(\mathbb R_+)$ such that for a.e.\,$t \in \mathbb R_+$
\begin{align}
\label{fbb_inhom}
\|b(t,\cdot)\varphi\|_2^2   \leq \delta \|\nabla \varphi\|_2^2 +g(t)\|\varphi\|_2^2
\end{align}
for all $\varphi \in W^{1,2}$. 
\end{definition}

This will be written as $b \in L^\infty\mathbf{F}_\delta + L^2_{\loc}(\mathbb R_+).$

\medskip

An equivalent form of the a.e.\,inequality \eqref{fbb_inhom} is: for every $0<T<\infty$,
$$
\int_{0}^T \|b(t)\psi(t)\|_2^2 dt \leq \delta \int_{0}^T \|\nabla \psi(t)\|_2^2 dt + \int_{0}^T g(t) \|\psi(t)\|_2^2 dt
$$
for all $\psi \in L^\infty(\mathbb R_+,W^{1,2})$.

\begin{examples}
The class of time-inhomogeneous form-bounded vector fields includes e.g.\,the critical Ladyzhenskaya-Prodi-Serrin class
$$
|b| \in L_{\loc}^q(\mathbb R_+,L^r+L^\infty), \quad \frac{d}{r}+\frac{2}{q} \leq 1, \quad 2 \leq q \leq \infty,
$$
as well as vector fields having stronger spatial singularities, see Appendix \ref{example_sect}.
\end{examples}

We fix an approximation of $b \in  L^\infty\mathbf{F}_\delta + L^2_{\loc}$ by smooth bounded vector fields $b_m$ that preserve the form-bound $\delta$ and have functions $g_n$ locally uniformly bounded in $L^1(\mathbb R_+)$, i.e.
\begin{equation}
\label{b_m_cond1}
b_n \rightarrow b \quad \text{ in $L^2_{\loc}(\mathbb R_+ \times \mathbb R^d, \mathbb R^d$}) 
\end{equation}
and for all $t \geq 0$
\begin{equation}
\label{b_m_cond2}
\|b_n(t)\varphi\|_2^2   \leq \delta \|\nabla \varphi\|_2^2 +g_n(t)\|\varphi\|_2^2
\end{equation}
 with $g_n$ such that 
\begin{equation}
\label{b_n_cond3}
\sup_n\int_0^Tg_n(s)ds < \infty \quad \text{ for any }0< T<\infty.
\end{equation}

\begin{examples}
It is easy to show that the following $b_n$, with $\varepsilon_n \downarrow 0$ sufficiently rapidly and $c_n \uparrow 1$ sufficiently slow, satisfy \eqref{b_m_cond1}-\eqref{b_n_cond3}. 
$$
b_n:=c_n E^{1+d}_{\varepsilon_n} (\mathbf 1_n b),
$$
where $\mathbf 1_n$ is the indicator of $\{(t,x) \mid |b(t,x)| \leq n, |x| \leq n, |t| \leq n\}$ (say, $b$ is extended by $0$ to $t<0$),  $E^{1+d}_\varepsilon$ is the De Giorgi or the Friedrichs mollifier on $\mathbb R \times \mathbb R^d$. See details in Appendix \ref{approx_app}.1. Note that, by selecting $\varepsilon_n \downarrow 0$ rapidly, one can treat $b_n$ as basically a cutoff of $b$ times constant $c_n$.

Moreover, with some additional effort, one can simplify this approximation to 
$$
b_\varepsilon:=E_{\varepsilon}^1 E_{\varepsilon}^d b, \quad \varepsilon \downarrow 0,
$$
where $E_{\varepsilon}^1$ in the Friedrichs mollifier on $\mathbb R$, and $E_{\varepsilon}^d$ is the De Giorgi or the Friedrichs mollifier on $\mathbb R^d$. 
See Appendix \ref{approx_app}.3. 

The last approximation is important if one needs e.g.\,to transfer the form-boundedness assumption on ``potential'' ${\rm div\,}b$ to the uniform form-boundedness of ${\rm div\,}b_\varepsilon$ since then one commute ${\rm div\,}$ and the mollifiers, although we are not concerned with this here.
\end{examples}

\medskip

\textbf{2.~}Our first goal is to construct the corresponding to $(\partial_t-\Delta + b \cdot \nabla)u=0$,  $b  \in L^\infty\mathbf{F}_\delta + L^2_{\loc}$, Feller evolution family on $D_T=\{(s,t) \mid 0 \leq s \leq t \leq T\}$ for $T>0$ fixed, i.e.\,a family of operators
$\{U^{t,s}\}_{(s,t) \in D_T}$ that are bounded on $C_\infty$, and

1) $U^{t,r}U^{r,s}=U^{t,s}$, $r \in [s,t]$, and $U^{s,s}={\rm I},$

2) $\|U^{t,s}f\|_\infty \leq \|f\|_\infty$, $U^{t,s}[C_\infty^+] \subset C_\infty^+$,

3) $$U^{r,s}=s\mbox{-}C_\infty\mbox{-}\lim_{t \downarrow r}U^{t,s}, \quad r \geq s,$$

\noindent and which will, additionally, satisfy: $u(t):=U^{t,s}f$ is the unique weak solution of $(\partial_t-\Delta + b \cdot \nabla)u=0$, $u(s)=f \in C_\infty \cap L^2$.

\medskip

The sought Feller evolution family is produced as the limit $L^\infty(D_T,L^\infty)$ of  $$
U_n^{t,s}f(\cdot):=u_n(t,\cdot), \quad (s,t) \in D_T
$$
and $u_n$ is the classical solution to initial problem
\begin{equation}
\label{cp}
\bigl(\partial_t - \Delta + b_n\cdot \nabla \bigr)u_n=0, \quad
 u_n(s,\cdot)=f(\cdot) \in C_c^\infty.
\end{equation}
We will prove the uniform convergence of the functions $\{(t,s,x) \mapsto U^{t,s}_n f(x)\}$ on $D_T \times \mathbb R^d$ by showing that they constitute a Cauchy sequence in $L^\infty(D_T,L^\infty)$. 
To that end, we will employ a parabolic variant of the iteration procedure of  \cite{KS}. This parabolic variant first appeared in \cite{Ki} and was recently refined in \cite{KiS_note}.

Namely, subtracting the equations for $u_m$, $u_n$ and setting $$h:=u_m-u_n,$$ one obtains
\begin{equation}
\label{eq_h}
\partial_t h - \Delta h + b_m \cdot \nabla h + (b_m-b_n) \cdot \nabla u_n=0, \quad h(s,\cdot)=0
\end{equation}
Multiplying the last equation by $h|h|^{r-2}$, $r>\frac{2}{2-\sqrt{\delta}}$, integrating over $[s,T] \times \mathbb R^d$ and applying the Sobolev embedding theorem, we arrive at the inequality
\[
c_dr^k \|h\|_{L^\infty([s,T],L^r)}^r + \|h\|^r_{L^r([s,T],L^{\frac{rd}{d-2}})}\leq C r^{2k}e^{C_T} \sup_{\tau \in [s,T]}\|\nabla u_n(\tau)\|_q^2 \int_s^T\|h(\tau)\|^{r-2}_{\frac{q}{q-2}(r-2)} d\tau, 
\]
where for a fixed $q>d$, for constants $C$, $C_T$ that are independent of $m$, $n$. Hence, applying the interpolation inequality in the left-hand side and setting $K=Ce^{C_T}$, one obtains
$$
\|h\|_{\frac{r}{1-\beta},\frac{rd}{d-2+2\beta}}\leq K^\frac{1}{r} (r^\frac{1}{r})^{2k} \biggl(\sup_{\tau \in [s,T]}\|\nabla u_n(\tau)\|_q\biggr)^\frac{2}{r} \|h\|^{1-\frac{2}{r}}_{L^{r-2}([s,T], L^{\frac{q}{q-2}(r-2)})}. 
$$
(we only need $\delta<4$ to prove this inequality). 
Now, with appropriate choice of $\beta$, one can iterate this inequality in essentially the same way as it was done \cite{KS} \textit{provided that one has uniform in $n$ bound on $\sup_{\tau \in [s,T]}\|\nabla u(\tau)\|_q$} (see below), arriving at
$$
\|u_m-u_n\|_{L^\infty([s,T],L^\infty)} \leq C_1\|u_m-u_n\|^\gamma_{L^{r_0}([s,T],L^{r_0})}
$$
for $r_0>\frac{2}{2-\sqrt{\delta}}$ and some $\gamma>0$ (this strict inequality is the main concern of the iteration procedure). Now, a standard argument yields $$\|u_m-u_n\|^\gamma_{L^{r_0}([s,T],L^{r_0})} \rightarrow 0 \quad \text{as $m,n \rightarrow \infty$},$$ see e.g.\,\cite{Ki},\cite{KiS_note},
and so we have our Cauchy sequence:
$$\|u_m-u_n\|_{L^\infty([s,T],L^\infty)} \rightarrow 0 \quad \text{ as $m,n\rightarrow \infty$},
$$
moreover, the convergence is uniform in $s \in [0,T]$. So, we can define the sought Feller evolution family by
$$
U^{t,s}f:=s\mbox{-}C_\infty\mbox{-}\lim_{n}u_n(t), \quad (s,t) \in D_T,
$$
for $f \in C_c^\infty$, as was assumed above, and then extend operators $U^{t,s}$ to all $f \in C_\infty$ by continuity using the fact that $U^{t,s}$ inherits the $L^\infty$ contraction property from $U_n^{t,s}$. Let us emphasize that the a priori assumption $f \in C_c^\infty$ is needed for the uniform in $n$ bound on $\sup_{\tau \in [s,T]}\|\nabla u(\tau)\|_p$, $p>d$.

\begin{remark}
Working in the elliptic setting (i.e.\,as in the proof of Theorem \ref{thm1_feller} or in \cite{KS}), after showing that solutions to the approximating elliptic equation converge in $C_\infty$, one needs to verify the other conditions of the Trotter approximation theorem. This is not needed when one is working directly with the parabolic equation, so we arrive at shorter proofs even if $b=b(x)$, however, at expense of requiring smaller $\delta$. We discussed this effect in Remark \ref{elliptic_rem}. It is fair to say that  there is a fundamental difference between time-homogeneous and time-inhomogeneous cases when one is dealing with singular drifts. 
\end{remark}

\textbf{3.~}To make the iteration procedure converge, one needs gradient bound 
\begin{equation}
\label{grad_bd}
\sup_n\sup_{\tau \in [s,T]}\|\nabla u_n(\tau)\|_p<\infty, \quad \text{ for some } p>d 
\end{equation}
for $f$ in a dense subset of $C_\infty$ (e.g.\,for $f \in C_c^\infty$). To obtain such a bound, one can differentiate the initial problem \eqref{cp}. Namely, writing for brevity $$u:=u_n, \quad b:=b_n$$ and
$$ w:=\nabla u, \quad w_i:=\nabla_{i} u,$$ we obtain
\begin{equation}
\label{w_j1}
\partial_t w_j -\Delta w_j + b \cdot \nabla w_j + (\nabla_j b)\cdot w=0, \quad w_j(0)=\nabla_j f, \quad 1 \leq j \leq d.
\end{equation}
Now, one needs to ``wrap up'' this system and, additionally, get rid of the derivative $\nabla_j b$. For instance, one can consider all products $w_{i_1}\dots w_{i_m}$ for $m$ fixed and then sum them up, as was done in \cite{BFGM} for solutions of stochastic transport equation and, after them, in \cite{KiSS_transport}. However, we are interested in arguments that give less restrictive assumptions on $\delta$.  Another argument was used in \cite{Kr2}, although in a different situation dealing with a more sophisticated system of parabolic equations. This argument still imposes more restrictive assumption  on $\delta$ than the arguments employed \cite{Ki} and \cite{KiS_note} (compare \eqref{c3} with \eqref{c4}, \eqref{c5}). However, it is quite nice and simple (and, again, works for other equations), so we describe it here.
Assume that that the form-bound $\delta$ of $b$ satisfies
\begin{equation}
\label{c3}
\tag{$C_3$}
\sqrt{\delta}<\frac{d-1}{d(d+1)}. 
\end{equation}
Put
$$
w_\eta:=\eta \cdot w, \quad \eta=(\eta_j)_{j=1}^d \in \mathbb R^d.
$$
We differentiate \eqref{cp} in the direction $\eta$, i.e.\,multiply \eqref{w_j1} by $\eta_j$ and add the resulting identities in $j=1,\dots,d$ to obtain
$$
\partial_t w_\eta - \Delta_x w_\eta - b \cdot \nabla_x w_\eta - \sum_{i=1}^d (\eta \cdot \nabla b_i)\nabla_{\eta_i}w_\eta=0,
$$
where in the last term we have used $w_i=\nabla_{\eta_i}w_\eta$. Given a function $g(x,\eta)$, we denote by $\langle g\rangle_x$, $\langle g\rangle_x$ the integral of $g$ over $\mathbb R^d$ in $x$ and in $\eta$, respectively. Let  $\langle g \rangle_{x,\eta}$ denote the corresponding repeated integral over $\mathbb R^d\times \mathbb R^d$.
Set 
$$
h(\eta):=(1+\kappa|\eta|^2)^{-\theta},
$$
where $\kappa>0$ is fixed arbitrarily, and $\theta>\frac{d+q}{2}$ so that $\langle |\eta|^qh\rangle_\eta<\infty$.
Let $q \geq 2$ (in the iteration procedure we need $q>d$). Also, without loss of generality, $q$ is rational with odd denominator, so we, if needed, we can raise negative numbers of power $q$. Multiply the previous identity by $h w_\eta^{q-1}$ and 
integrate in $(x,\eta) \in \mathbb R^{2d}$ to obtain 
\begin{align}
\frac{1}{q}\partial_t \langle hw_\eta^q \rangle_{x,\eta} & + \frac{4(q-1)}{q^2}\langle h |\nabla w^{\frac{q}{2}}_\eta|^2\rangle_{x,\eta} \notag \\
& - \frac{2}{q}\langle b \cdot \nabla w_\eta^{\frac{q}{2}},h w_\eta^{\frac{q}{2}}\rangle_{x,\eta}  - \sum_{i=1}^d \langle(\eta \cdot \nabla b_i)\partial_{\eta_i} w_\eta, hw_\eta^{q-1}\rangle_{x,\eta}=0, \label{id_b}
\end{align}
where $b_i$ are the components of $b$.
The last term in the left-hand side is dealt with as follows:
\begin{align*}
-\sum_{i=1}^d \langle (\eta \cdot \nabla b_i)\partial_{\eta_i} w_\eta, hw_\eta^{q-1}\rangle_{x,\eta} & =  \sum_{i=1}^d \langle \eta b_i \nabla \partial_{\eta_i}w_\eta,hw_\eta^{q-1}\rangle_{x,\eta} \\
& + (q-1)\sum_{i=1}^d \langle \eta b_i \partial_{\eta_i}w_\eta,hw_\eta^{q-2}\nabla w_\eta\rangle_{x,\eta} \\
& (\text{now we integrate by parts in $\eta_i$ in the first term}) \\
& = - \langle b \cdot \nabla w_\eta,hw_\eta^{q-1}\rangle_{x,\eta} - \sum_{i=1}^d \langle \eta b_i \nabla w_\eta, (\partial_{\eta_i}h)w_\eta^{q-1}\rangle_{x,\eta}.
\end{align*}
Hence, \eqref{id_b} becomes
\begin{align*}
\frac{1}{q}\partial_t \langle hw_\eta^q \rangle_{x,\eta} & + \frac{4(q-1)}{q^2}\langle h |\nabla w^{\frac{q}{2}}_\eta|^2\rangle_{x,\eta} \notag \\
& - \frac{4}{q}\langle b\cdot \nabla w_\eta^{\frac{q}{2}},h w_\eta^{\frac{q}{2}}\rangle_{x,\eta} - \frac{2}{q}\sum_{i=1}^d \langle \eta b_i \nabla w_\eta^{\frac{q}{2}}, (\partial_{\eta_i}h)w_\eta^{\frac{q}{2}}\rangle_{x,\eta}=0,
\end{align*}
so
\begin{align*}
\partial_t \langle hw_\eta^q \rangle_{x,\eta} & + \frac{4(q-1)}{q}\langle h |\nabla w^{\frac{q}{2}}_\eta|^2\rangle_{x,\eta} \notag \\
& - 4\langle b\cdot \nabla w_\eta^{\frac{q}{2}},h w_\eta^{\frac{q}{2}}\rangle_{x,\eta} + 4\theta\sum_{i=1}^d \langle  b_i \nabla w_\eta^{\frac{q}{2}}, \frac{\kappa \eta_i \eta}{1+\kappa |\eta|^2} h w_\eta^{\frac{q}{2}}\rangle_{x,\eta}=0.
\end{align*}
Hence
\begin{align}
\label{ineq_kr}
\partial_t \langle hw_\eta^q \rangle_{x,\eta} & + \frac{4(q-1)}{q}\langle h |\nabla w^{\frac{q}{2}}_\eta|^2\rangle_{x,\eta} - (4+4\theta)\langle |b| |\nabla w_\eta|^{\frac{q}{2}},h w_\eta^{\frac{q}{2}}\rangle_{x,\eta} \leq 0.
\end{align}
We estimate
\begin{align*}
\langle |b| |\nabla w_\eta|^{\frac{q}{2}},h w_\eta^{\frac{q}{2}}\rangle_{x,\eta} & \leq \varepsilon \langle |b|^2,h w_\eta^{q}\rangle_{x,\eta}+\frac{1}{4\varepsilon}\langle h|\nabla w_\eta|^{q}\rangle_{x,\eta} \\
& \leq \varepsilon\biggl(\delta \langle h|\nabla w^{\frac{q}{2}}_\eta|^2\rangle + g(t)\langle w_\eta^q\rangle\biggr) +\frac{1}{4\varepsilon}\langle h|\nabla w_\eta|^{q}\rangle_{x,\eta} \qquad \varepsilon:=\frac{1}{2\sqrt{\delta}}.
\end{align*}
Thus, integrating \eqref{ineq_kr} from $s$ to $t$, one obtains
$$
\langle hw_\eta^q(t)\rangle_{x,\eta} + 4\biggl(\frac{q-1}{q} - \sqrt{\delta}(1+\theta) \biggr)\int_s^t \langle h |\nabla w^{\frac{q}{2}}_\eta(\tau)|^2\rangle_{x,\eta} d\tau \leq \frac{1}{2\sqrt{\delta}}\int_s^t g(\tau)\langle w_\eta^q(\tau)\rangle d\tau + \langle h(\nabla f \cdot \eta)^q\rangle,
$$
which allows to conclude
\begin{equation}
\label{bd_h}
\sup_{t \in [s,T]}\langle hw_\eta^q(t)\rangle_{x,\eta} + \int_s^t \langle h |\nabla w^{\frac{q}{2}}_\eta(\tau)|^2\rangle_{x,\eta} d\tau \leq C\|\nabla f\|_q^q
\end{equation}
for some for some $q>d$, for constant $C$ independent of $n$,
provided that
$$
\frac{q-1}{q} - \sqrt{\delta}(1+\theta) > 0 \text{ for some $\theta>\frac{d+q}{2}$} \quad \Leftrightarrow \quad \sqrt{\delta}<\frac{d-1}{d(d+1)}.
$$

It remains to derive \eqref{grad_bd} from \eqref{bd_h}.
Put $$A_{t,x}:=\{\eta \in \mathbb R^d \mid |\eta-\frac{w(t,x)}{|w(t,x)|}|<\frac{1}{2}\}, \quad x \in \mathbb R^d$$ (if $w(t,x)=0$, fix $\eta=(1,0,\dots,0)$). Thus, $A_{t,x}$ is a ball of radius $\frac{1}{2}$ with centre placed at distance $1$ from the origin. The angle between $w(t,x)$ and $\eta \in A_{t,x}$ is bounded from above by a generic constant, hence $|\eta \cdot w(t,x)| \geq c|w(t,x)|$ for some $c>0$ independent of $(t,x)$, for all $\eta \in A_{t,x}$. Therefore, for all $t \in [s,T]$
\begin{align*}
\langle hw_\eta^q\rangle_{x,\eta} & =\langle \langle h(\eta)|\eta \cdot w(x)|^q\rangle_\eta\rangle_x \\
& \geq \langle \langle h(\eta)|\eta \cdot w(x)|^q\mathbf{1}_{A_{t,x}}(\eta)\rangle_\eta\rangle_x \\
& \geq \langle c^q |w(x)|^q \langle h(\eta)\mathbf{1}_{A_{t,x}}(\eta)\rangle_\eta\rangle_x =C\langle |w|^q\rangle_x, \quad C>0.
\end{align*}
Thus, we obtain from \eqref{bd_h}
\begin{equation}
\label{bd_2}
\sup_{t \in [s,T]}\|w(t)\|_q^q +c\int_s^T\|\nabla |w|^\frac{q}{2}\|_2^2dt\leq C_T\|\nabla f\|_q^q, \quad c>0.
\end{equation}
Hence, one can run the iteration procedure under the assumption \eqref{c3}.

\medskip

\textbf{4.~}The proofs of \eqref{bd_2} in \cite{Ki,KiS_note} choose a specific direction of the differentiation (following \cite{KS} which, by the way, appeared earlier than the other papers cited above): $$\eta=\frac{w}{|w|},$$ which maximizes the directional derivative $w_\eta=w \cdot \eta$. Put differently, one multiplies the parabolic equation in \eqref{cp} by the test function
\begin{equation}
\label{phi_test}
\varphi=-\nabla \cdot (\frac{w}{|w|}|w|^{q-1})
\end{equation}
and then integrates by parts (the same test function is used in the proof of Lemma \ref{lem_GPQ}). This choice of  the test function (or direction) leads to better assumptions on $\delta$ than \eqref{c3}. Indeed, the identity
\[
\langle \partial_tu,\varphi\rangle + \langle -\Delta u,\varphi\rangle + \langle b_m \cdot w,\varphi\rangle =0 ,
\]
yields
\begin{equation}
\label{I_q_id}
\frac{1}{q}\partial_t\|w\|_q^q + I_q + (q-2)J_q=\langle b_m\cdot w,\nabla\cdot(w|w|^{q-2})\rangle
\end{equation}
where 
\[
I_q=\sum_{i=1}^d\langle |\nabla w_i|^2,|w|^{q-2}\rangle, \quad J_q=\langle |\nabla |w||^2,|w|^{q-2}\rangle
\]
are the ``good'' terms, i.e.\,the right-hand side of \eqref{I_q_id}  will be estimated in terms of $I_q$, $J_q$ multiplied by coefficients that, thus, cannot be too large, hence our assumptions on $\delta$. Namely, we represent
\begin{align*}
\langle b_m\cdot w,\nabla\cdot(w|w|^{q-2})\rangle &= \langle b_m\cdot w,\Delta u |w|^{q-2}\rangle+(q-2)\langle b_m\cdot w,|w|^{q-3}w\cdot\nabla|w|\rangle\\
&:= S_1+S_2.
\end{align*}
Put $B_q=\langle(b_m\cdot w)^2,|w|^{q-2}\rangle$, then
\begin{equation}
\label{S2_bd}
S_2 \leq(q-2)B_q^\frac{1}{2}J_q^\frac{1}{2},
\end{equation}
where $B_q$  is estimated using \eqref{b_m_cond2}, i.e.\,the form-boundedness of $b_m$:
\begin{equation}
\label{bd_B}
B_q\leq  \frac{q\sqrt{\delta}}{2}J_q+g_m\|w\|_q^q.
\end{equation}
Thus, $S_2$ is estimated in terms of $J_q$, and one can apply the resulting bound on $S_2$ in \eqref{I_q_id}.

In \cite{Ki}, after estimating $S_1$ as
\begin{equation}
\label{S1_est}
|S_1| \leq B_q^\frac{1}{2}\langle|\Delta u|^2,|w|^{q-2}\rangle^\frac{1}{2},
\end{equation}
the factor $\langle|\Delta u|^2,|w|^{q-2}\rangle$ was bounded by $I_q$ and $J_q$  without appealing to the equation, by representing $|\Delta u|^2=(\nabla \cdot w)^2$
and integrating by parts twice: 
$$
\langle |w|^{q-2}|\Delta u|^2 \rangle = -\langle  w \cdot \nabla |w|^{q-2}, \Delta u \rangle 
+\sum_{r=1}^d  \left\langle w \cdot \nabla w_r, \nabla_r |w|^{q-2}\right\rangle  
+ I_q,
$$
where 
$$
|\langle  w \cdot \nabla |w|^{q-2}, \Delta u \rangle| \leqslant (q-2)\left( \frac{1}{4\varkappa}  \langle |w|^{q-2}|\Delta u|^2\rangle + \varkappa J_q\right), \quad \kappa>0,
$$
and
$$
\quad |\sum_{r=1}^d  \left\langle w \cdot \nabla w_r, \nabla_r |w|^{q-2}\right\rangle | \leqslant (q-2)\left( \frac{1}{2} I_q+ \frac{1}{2} J_q\right).
$$
Hence
\begin{equation}
\label{bd_laplace}
\biggl(1-\frac{q-2}{4\varkappa} \biggr)\langle |w|^{q-2}|\Delta u|^2 \rangle \leqslant I_q+(q-2)\left(\varkappa J_q + \frac{1}{2} I_q + \frac{1}{2} J_q \right), \quad \varkappa>\frac{q-2}{4}.
\end{equation}
The resulting from \eqref{S1_est}, \eqref{bd_laplace} bound on $|S_1|$, combined with \eqref{S2_bd}, \eqref{bd_B}, led in \cite{Ki} to the gradient estimate \eqref{bd_2}
for a $q>d$ close to $d$ provided that
\begin{equation}
\label{c4}
\tag{$C_4$}
\sqrt{\delta}<\frac{1}{d}.
\end{equation}

One important advantage of working with the test function \eqref{phi_test} is that one can ``evaluate'' it by representing  $\Delta u=\partial_tu + b_m \cdot \nabla u$, thus using the equation one more time. This is what \cite{KS} did. The same can be done in the parabolic setting, and it leads to better assumptions on $\delta$ than \eqref{c4} \cite{KiS_note}. Specifically, in dimensions $3\leq d\leq 6$, one abandons estimate \eqref{S1_est} and represents $\Delta u=\partial_tu + b_m \cdot \nabla u$ to evaluate
\[
S_1=-\frac{1}{q}\frac{d}{dt}\|w\|_q^q +B_q-\langle|\partial_tu|^2,|w|^{q-2}\rangle-(q-2)\langle|w|^{q-3}w\cdot\nabla |w|,\partial_tu\rangle,
\]
which, upon applying the quadratic inequality, gives $$S_1 \leq -\frac{1}{q}\partial_t\|w\|_q^q +B_q +\frac{(q-2)^2}{4}J_q.$$ 
This bound, \eqref{S2_bd} and \eqref{bd_B} applied in \eqref{I_q_id} give the desired gradient bound \eqref{bd_2} for $\sqrt{\delta}\leq \big(\sqrt{q-1}-\frac{q-2}{2}\big)\frac{2}{q}$. In dimensions $d=3$ and $d=4$ this gives significantly less restrictive assumption on $\delta$ than \eqref{c4}, see \eqref{c5} below.

In dimensions $d\geq 5$, one starts with two representations for $S_1$:
\begin{align*}
S_1 & = \langle - \partial_t u + \Delta u, |w|^{q-2} \Delta u \rangle \\
& = \langle |\Delta u|^2, |w|^{q-2} \rangle - \Real\langle \partial_t u , |w|^{q-2} \Delta u \rangle,
\end{align*}
\begin{align*}
S_1 & =\langle b_n \cdot w, |w|^{q-2} (\partial_t u + b_n \cdot w) \rangle \\
& = B_q + \langle b_n \cdot w, |w|^{q-2} \partial_t u \rangle.
\end{align*}
Equating the right-hand sides, one obtains
\begin{align*} 
\langle |\Delta u|^2, |w|^{q-2} \rangle & = B_q + \langle \partial_t u, |w|^{q-2} ( b_n \cdot w + \Delta u  )\rangle \\[3mm]
& = B_q + \langle \partial_t u, |w|^{q-2} ( - \partial_t u + 2 \Delta u  )\rangle\\[3mm]
& = B_q - \langle |\partial_t u|^2, |w|^{q-2} \rangle + 2 \langle \partial_t u,|w|^{q-2} \Delta u \rangle \\[2mm]
& = B_q - \langle |\partial_t u|^2, |w|^{q-2} \rangle - \frac{2}{q}\frac{d}{dt}\|w\|_q^q -2 (q-2) \langle \partial_t u,|w|^{q-3} w \cdot \nabla |w| \rangle \\
& \leq B_q - \langle |\partial_t u|^2, |w|^{q-2} \rangle - \frac{2}{q}\frac{d}{dt}\|w\|_q^q + (q-2)^2 J_q + \langle |\partial_t u|^2,|w|^{q-2} \rangle \\
& = B_q - \frac{2}{q}\partial_t\|w\|_q^q + (q-2)^2 J_q.
\end{align*}
This  estimate on $\langle |\Delta u|^2, |w|^{q-2} \rangle$, which is more efficient than \eqref{bd_laplace}, when applied in \eqref{S1_est} leads, together with \eqref{S2_bd}, \eqref{bd_B}, to  the following.
\textit{If form-bound $\delta$ satisfies}
\begin{equation}
\label{c5}
\tag{$C_5$}
\begin{array}{ll}
d\geq 5 & \sqrt{\delta} \text{ satisfies } \frac{d \sqrt{\delta}}{2} \bigg( \sqrt{\frac{d^2 \delta}{4}+ (d-2)^2} + d-2 \bigg) < d-1, \\
d=4 & \sqrt{\delta} < \frac{2(\sqrt{3}-1)}{d} \approx 0.36602,\\[2mm]
d=3 & \sqrt{\delta}< \frac{2\sqrt{2}-1}{d} \approx 0.60947,
\end{array}
\end{equation}
\textit{then gradient estimate \eqref{bd_2} holds for a $q>d$ close to $d$.} See \cite{KiS_note} for the proof.

 The assumption \eqref{c5} is less restrictive than \eqref{c4}. In fact, if one assumes $\sqrt{\delta}=\frac{1}{d}$, then \eqref{bd_2} holds even for $q=d+1$.

\begin{remark} 
\label{grad_q}
The gradient bounds in \cite{KS,Ki, KiS_note} are proved not only for $q$ close to $d$, but for the entire range of $(\delta,q)$ satisfying some algebraic inequalities. In particular, in \cite{KiS_note}, 
$$
 q-1 - \frac{q \sqrt{\delta}}{2} \bigg( \sqrt{\frac{q^2 \delta}{4}+ (q-2)^2} + q-2 \bigg) > 0 \quad \text{ if $d \geq 5$}
$$
and, as was mentioned above, for $d=3,4$, $\sqrt{\delta}\leq \big(\sqrt{q-1}-\frac{q-2}{2}\big)\frac{2}{q}$.
\end{remark}

Thus, we have the following result.

\begin{theorem}[\cite{Ki}, \cite{KiS_note}]
\label{feller_parab} 
Assume that  $b \in L^\infty\mathbf{F}_\delta+L^2_{\loc}(\mathbb R_+)$ with $\delta$ that satisfying \eqref{c4} or, better, \eqref{c5}. Then the following is true:

{\rm (\textit{i})} The limit
$$
U^{t,s}f:=s\mbox{-}C_\infty\mbox{-}\lim_n U_n^{t,s}f \quad \text{ uniformly in $(s,t) \in D_T$}
$$
exists for all $f \in C_c^\infty$ and satisfies $\|U^{t,s}f\|_\infty \leq \|f\|_\infty$. Upon extending operators $U^{t,s}$ by continuity to all $f \in C_\infty$, one obtains a Feller evolution family.

\smallskip

{\rm (\textit{ii})} The Feller evolution family $\{U^{t,s}\}_{(s,t) \in D_T}$ is unique in the sense that 
it does not depend on the choice of the approximation vector fields $\{b_n\}$, as long as they satisfy \eqref{b_m_cond1}, \eqref{b_m_cond2}, \eqref{b_n_cond3}.

\smallskip

{\rm (\textit{iii})} If $f \in C_\infty \cap L^2$, then $u(t):=U^{t,s}f$ is the unique weak solution of $(\partial_t-\Delta + b \cdot \nabla)u=0$, $u(s)=f$, in $L^2$.
\end{theorem}

\bigskip

\section{SDEs with time-inhomogeneous form-bounded drifts}

\label{parab_sde_sect}

We return to the discussion of weak well-posedness of SDE
\begin{equation}
\label{parab_SDE}
X_t=x - \int_0^b b(s,X_s)ds + \sqrt{2}W_t, \quad t \geq 0,
\end{equation}
where $x \in \mathbb R^d$ is fixed and $b \in L^\infty\mathbf{F}_\delta+L^2_{\loc}(\mathbb R_+)$, where, we assume for simplicity, \eqref{c4} holds.

\medskip

We  need to supplement Theorem \ref{feller_parab} with a localized analogue of \eqref{bd_2} for inhomogeneous parabolic equations, proved in \cite{KiM}.
Let $\mathsf{f} \in L^\infty\mathbf{F}_\nu+L^2_{\loc}(\mathbb R_+)$, $\nu<\infty$,  define $\mathsf{f}_k$ similarly to $b_m$ in \eqref{b_m_cond1}-\eqref{b_n_cond3}.
 Let $h \in C([s,T], \mathcal S)$, $ g \in C_c^\infty(\mathbb R^d)$. Fix $T>s$. Let $u=u_{m,k}$ be the solution to Cauchy problem on $[s,T]$
\begin{equation}
\label{Cauchy3}
(\partial_t - \Delta + b_m \cdot \nabla)u=|\mathsf{f}_k|h, \quad u(s,\cdot)=g.
\end{equation}
Then, for every $q \in ]d,\delta^{-\frac{1}{2}}[$, there exist constants $C$ and $\kappa$ such that, for all $0 \leq s \leq r \leq T$, 
\begin{align}
\|u\|_{L^\infty([s,r],L_\rho^q)}^q  + \|\nabla u\|^q_{L^\infty([s,r],L_\rho^q)} & + \|\nabla|\nabla u|^{\frac{q}{2}}\|_{L^2([s,r],L_\rho^2)}^2 \notag \\
& \leq C \bigl(\|\mathsf{f} |h|^{\frac{q}{2}}\|^2_{L^2([s,r],L_\rho^2)} + \|\nabla g\|^q_{L_\rho^q} + \|g\|^q_{L^q_\rho}\big). \label{rho_grad_est}
\end{align}
Here $\rho(x)
:=(1+\kappa |x|^2)^{-\theta}$ ($x \in \mathbb R^d$),
where $\theta>\frac{d}{2}$ is fixed, and $L^2_\rho:=L^2(\mathbb R^d,\rho dx)$.

\medskip

Define backward Feller evolution family ($0 \leq t \leq r \leq T$)
\begin{equation*}
P^{t,r}(b)=U^{T-t,T-r}(\tilde{b}), \quad \tilde{b}(t,x)=b(T-t,x),
\end{equation*}
where $U^{t,s}$ is the Feller evolution family from Theorem \ref{feller_parab}.
Using \eqref{rho_grad_est} with $\mathsf{h}=0$ and arguing essentially as in the proof of Theorem \ref{thm1_feller}(\textit{v}), one obtains that $\{P^{t,r}(b)\}_{0 \leq t \leq r \leq T}$ is conservative, i.e.\,for all $x \in \mathbb R^d$
$
\langle P^{t,r}(x,\cdot)\rangle=1.
$
Now, by a standard result (see e.g.\,\cite[Ch.\,2]{GvC}), given a conservative backward Feller evolution family, there exist probability measures $\mathbb P_x$ $(x \in \mathbb R^d)$ on 
$(D([0,T],\mathbb R^d),\mathcal B'_t)$, such that
$$
\mathbb E_{x}[f(\omega_r)]=P^{0,r}f(x), \quad 0 \leq r \leq T. 
$$
Here and below, $\mathbb E_x:=\mathbb E_{\mathbb P_x}$.

Let $X_t^m$ ($m=1,2,\dots$) be the strong solution of
$$
X^m_t=x-\int_0^t b_m(r,X^m_r)dr + \sqrt{2}W_t, \quad x \in \mathbb R^d
$$
defined on some complete probability space $(\Omega,\mathcal F,\mathcal F_t,\mathbf{P})$.

We will require the following estimate:
there exists a constant $C>0$ independent of $m$, $k$ such that
\begin{equation}
\label{est_F}
\sup_m \sup_{x \in \mathbb R^d}\mathbf E\int_{s}^r |b_k(t,X^m_{t})|dt \leq CF(r-s)
\end{equation}
for $0 \leq s \leq r \leq T$, where $F(h):=h + \sup_{s \in [0,T-h]}\int_s^{s+h}g(t)dt$. Here we assume, without loss of generality, that $b_k \in L^\infty\mathbf{F}_\delta + L^2_{\loc}(\mathbb R_+)$ with the same function $g$ as $b$ (if not, then we can increase $g$, cf.\,\eqref{b_n_cond3}). 

Indeed, let $v=v_{m,k}$ be the solution to the terminal-value problem
\begin{equation}
\label{veq}
(\partial_t + \Delta - b_m \cdot \nabla)v=-|b_k|, \quad v(r,\cdot)=0, \quad t \leq r.
\end{equation}
By It\^{o}'s formula,
$$
v(r,X_r^m)=v(s,X_s^m) + \int_s^r (\partial_t v + \Delta v - b_m \cdot \nabla v)(t,X^m_t) dt + \sqrt{2}\int_s^r \nabla v(t,X_t^m)dW_t.
$$
Taking expectation, we obtain
$$
\mathbb E\int_s^r |b_k(t,X_t^m)|dt = \mathbb E v(s,X_s^m).
$$
Now,  \eqref{rho_grad_est} applied to equation \eqref{veq} (so, we reverse the direction of time and take $|\mathsf{f}|:=|b_k|$, $h \equiv 1$ and $g=0$) yields, upon applying the Sobolev embedding theorem,
$$
\|v\|_{L^\infty([s,r] \times B_{1}(0))} \leq C_1 \|b_k \sqrt{\rho}\|^2_{L^2([s,r],L^2)},
$$
where $0$ is the ``centre'' of the weight $\rho$. Thus, considering its translates $\rho_z:=\rho(x-z)$, we obtain
$$
\|v\|_{L^\infty([s,r] \times \mathbb R^d)} \leq C_2\sup_{z \in \mathbb Z^d}\|b_k \sqrt{\rho_{z}}\|_{L^2([s,r],L^2)}.
$$
Since $\mathbb E v(s,X_s^m) \leq \|v(s,\cdot)\|_\infty$, we obtain 
$$
\mathbb E\int_s^r |b_k(t,X_t^m)|dt \leq C_2\sup_{z \in \mathbb Z^d}\|b_k \sqrt{\rho_{z}}\|_{L^2([s,r],L^2)}.
$$
Since $b_k \in L^\infty\mathbf{F}_\delta+L^2_{\loc}$, we have
\begin{align*}
\|b_k \sqrt{\rho_{z}}\|^2_{L^2([s,r],L^2)} & \leq \frac{\delta}{4}\int_s^r \langle \frac{|\nabla \rho_z|^2}{\rho_z}\rangle dt + \int_s^r g(t)\langle \rho_z \rangle dt \\
& (\text{we are using $|\nabla \rho| \leq \theta\sqrt{\kappa}\rho$ and $\|\sqrt{\rho}\|_2<\infty$}) \\
&  \leq CF(r-s)
\end{align*}
for $0 \leq s \leq r \leq T$, so \eqref{est_F} follows.

Now, define probability measures $\mathbb P^n_x:=(\mathbf P \circ X^n)^{-1}$ on $\bigl(C([0,T],\mathbb R^d),\mathcal B_t\bigr)$, so \eqref{est_F} takes form $$\sup_m \sup_{x \in \mathbb R^d}\mathbb E_{\mathbb P_x^m}\int_{s}^r |b_k(t,\omega_t)|dt \leq CF(r-s),$$ where $\omega_t$ is the coordinate process. We apply in \eqref{est_F} the convergence result of Theorem \ref{feller_parab} (in $m$) and then Fatou's lemma (in $k$) to obtain 
$$
\mathbb E\int_{s}^r |b(t,\omega_t)|dt \leq CF(r-s)<\infty
$$
(which is one of the requirements in the definition of a martingale solution).
Arguing similarly, we obtain, for every $f \in C_c^2$,
\begin{align*}
\mathbb E_{\mathbb P_x^m}\bigg|\int_0^r \big((b_m-b_n)\cdot \nabla f\big)(t,\omega_t)dr \bigg| & \leq C\|(b_m-b_n)|\nabla f|^{\frac{q}{2}}\|_{L^2([0,r],L^2)} \\
& \rightarrow 0 \quad (m,n \rightarrow \infty)
\end{align*}
since $b_m \rightarrow b$ in $L^2_{\loc}(\mathbb R_+ \times \mathbb R^d)$ and $f$ has compact support.
This, and the convergence result of Theorem \ref{feller_parab}, allow to pass to the limit in the martingale problem for $b_m$ in essentially the same way as in the proof of Theorem \ref{thm1} to show that $\mathbb P_x$ is a martingale solution of \eqref{parab_SDE} but on $(D([0,T],\mathbb R^d),\mathcal B_t')$. The latter allows to prove, arguing again as in the proof of Theorem \ref{thm1}, that $\mathbb P_x$ are actually concentrated on continuous trajectories.
We arrive at the following result.

\begin{theorem}
\label{sde_parab}
Under the assumptions of Theorem \ref{feller_parab}, let us also assume\footnote{We can assume \eqref{c5}, but then we need to adjust interval $q \in ]d,\delta^{-\frac{1}{2}}[$. For simplicitiy, we will not do this here.}  \eqref{c4}. The following is true:

\smallskip

{\rm (\textit{i})} For every $x \in \mathbb R^d$, the probability measure $\mathbb P_x$ is
a weak solution
to SDE \eqref{parab_SDE}. 

\smallskip

{\rm (\textit{ii})} $\mathbb P_x$ satisfies, for all $\mathsf{f} \in L^\infty\mathbf{F}_\nu+L^2_{\loc}(\mathbb R_+)$, $\nu<\infty$,  $ h \in C([0,T], \mathcal S)$, for all $q \in ]d,\delta^{-\frac{1}{2}}[$,  the estimate
\begin{equation}
\label{kr_est0}
\mathbb E_{\mathbb P_x}\int_0^T |\mathsf{f}(r,\omega_r)h(r,\omega_r)|dr \leq c\|\mathsf{f}|h|^{\frac{q}{2}}\|^{\frac{2}{q}}_{L^2([0,T] \times \mathbb R^d)}.
\end{equation}
On the other hand, if, for some $x \in \mathbb R^d$, $\mathbb P'_x$ is a martingale solution of \eqref{parab_SDE} that satisfies \eqref{kr_est0} for some $q \in ]d,\delta^{-\frac{1}{2}}[$ with $\mathsf{f}=b$,
then it coincides with $\mathbb P_x$.

\smallskip

{\rm (\textit{iii})} $\mathbb P_x$ satisfies, for a given $\nu>\frac{d+2}{2}$, for all $ h \in C([0,T], \mathcal S)$ the following Krylov-type bound:
\begin{equation}
\label{krylov_type}
\mathbb E_{\mathbb P_x}\int_0^T |h(r,\omega_r)|dr \leq c\|h\|_{L^\nu([0,T] \times \mathbb R^d)}.
\end{equation} On the other hand, if additionally $|b| \in L_{\loc}^{\frac{d+2}{2}+\varepsilon}(\mathbb R_+ \times \mathbb R^d)$ for some $\varepsilon>0$ and $\delta$ is sufficiently small, then any martingale solution $\mathbb P'_x$ of \eqref{parab_SDE} that satisfies \eqref{krylov_type} for some $\nu>\frac{d+2}{2}$ sufficiently close to $\frac{d+2}{2}$ (depending on how small $\varepsilon$ is)
coincides with $\mathbb P_x$.

\end{theorem}

The first two assertions of Theorem \ref{sde_parab} were proved in \cite{KiM}, the last assertion will be proved in the next section, in fact, for a substantially larger than $L^\infty\mathbf{F}_\delta + L^2_{\loc}(\mathbb R_+)$ class of drifts.

\begin{remark}
\label{unique_comp_rem}
One advantage of the uniqueness class in (\textit{ii}), i.e.
\begin{equation}
\label{unique_cl}
\mathbb E_{\mathbb P_x}\int_0^T |b(r,\omega_r)h(r,\omega_r)|dr \leq c\|b|h|^{\frac{q}{2}}\|^{\frac{2}{q}}_{L^2([0,T] \times \mathbb R^d)}
\end{equation}
\textit{for some} $q \in ]d,\delta^{-\frac{1}{2}}[$
is that it senses the value of $\delta$. Namely, as $\delta$ becomes smaller, one can take $q$ larger, and so the verification of \eqref{unique_cl}, in principle, becomes easier (e.g.\,$|b|$ is bounded, then $\delta$ can be arbitrarily small and hence $q$ can be arbitrarily large).
\end{remark}

Regarding the proof of Theorem \ref{sde_parab}(\textit{i}), let us note that we can alternatively use the tightness argument, also employed in the proof of Theorem \ref{thm_sharp}, and then apply the convergence result of Theorem \ref{feller_parab}. See details in \cite{KiM}. The uniqueness results in assertions (\textit{ii}), (\textit{iii}), however, require gradient bounds \eqref{rho_grad_est}.

\bigskip

\section{``Form-bounded'' diffusion coefficients} \

\label{a_sect}

\medskip

\textbf{1.~}The results of the previous two sections can be extended to It\^{o} and Stratonovich SDEs
\begin{equation}
\label{sde3_}
X_t=x-\int_0^t b(s,X_s)ds + \sqrt{2}\int_0^t\sigma(s,X_s)dW_s, \quad \in \mathbb R^d
\end{equation}
\begin{equation}
\label{sde3_str}
X_t=x-\int_0^t b(s,X_s)ds + \sqrt{2}\int_0^t\sigma(s,X_s) \circ dW_s,
\end{equation}
where the drift $b:\mathbb R^d \rightarrow \mathbb R^d$ is form-bounded and the diffusion coefficient $\sigma:\mathbb R^d \rightarrow \mathbb R^{d \times d}$ are bounded, uniformly elliptic and can be discontinuous. For  time-homogeneous $b$ and $\sigma$ such extension was carried out in \cite{KiS_Osaka}.
Namely, let $a:=\sigma \sigma^{\scriptscriptstyle \top}$ satisfy,
\begin{equation}
\label{H_}
\tag{$H_{\sigma,\xi}$}
\sigma I \leq a \leq \xi I \quad \text{ a.e.\,on } \mathbb R^d
\end{equation}
for some $0<\sigma \leq \xi<\infty$, and assume that the entries $a_{ij}$ of $a$ have form-bounded derivatives, that is,  
\begin{equation}
\label{nabla_a}
(\nabla_{r}a_{ij})_{i=1}^d \in \mathbf{F}_{\delta_{rj}}
\end{equation}
for some $\delta_{rj}>0$. Equivalently, since the entries of $\sigma$ are bounded, we can replace \eqref{nabla_a} with $(\nabla_{r}\sigma_{ij})_{i=1}^d \in \mathbf{F}_{\delta_{rj}'}$ for appropriate $\delta'_{rj}$.

\begin{examples}
1.~If $a \in W^{1,d}(\mathbb R^d,\mathbb R^{d \times d})$, then  \eqref{nabla_a} holds with $\delta_{rj}$ that can be chosen arbitrarily small (Appendix \ref{example_sect}).
More generally, if the derivatives of  $a_{ij}$ are in the Morrey class $M_{2+\varepsilon}$, then \eqref{nabla_a} holds. 

2.~Here is a concrete example of matrix $a$ satisfying \eqref{nabla_a} and having a critical discontinuity at the origin:
\begin{equation}
\label{gs}
a(x)=I+c\frac{x\otimes x}{|x|^2}, \quad \text{ the constant } c>-1.
\end{equation}
Indeed, $\nabla_r a_{ij}=c\mathbf{1}_{r=i}\frac{x_j}{|x|^2} + c\mathbf{1}_{r=j}\frac{x_i}{|x|^2} + cx_ix_j \frac{2x_r}{|x|^4}$, so  $$|(\nabla_r a_{ij})_{i=1}^d| \leq 2|c||x|^{-1} \quad \Rightarrow \quad (\nabla_r a_{ij})_{i=1}^d \in \mathbf{F}_{\delta_{rj}}, \quad \delta_{rj}=(4c)^2/(d-2)^2$$ by the Hardy inequality.
Another example is $$a(x)=I+c (\sin \log(|x|))^2 e \otimes e, \quad e \in \mathbb R^d, |e|=1$$
(indeed, $\nabla_r a_{ij}=2c (\sin \log |x|)(\cos \log |x|)|x|^{-2}x_r\, e_{i}e_j$, so using that the Hardy vector field \eqref{hardy_vf} is form-bounded one obtains the required).

More generally, \eqref{nabla_a} holds for $a$ that is an infinite sum of such matrices (properly normalized so that the series converges) with their points of discontinuity constituting e.g.\,a dense subset of $\mathbb R^d$. 
\end{examples}

Without loss of generality, in \eqref{H_} $\sigma=1$. In \cite{KiS_Osaka}, assuming that $b \in \mathbf{F}_\delta$ and
\begin{equation}
\label{osaka_cond}
\text{$a$ satisfies \eqref{H_} and \eqref{nabla_a},} \quad \nabla a \in \mathbf{F}_{\delta_a},
\end{equation}
 where $(\nabla a)_k:=\sum_{i=1}^d (\nabla_i a_{ik})$, with
 $\delta$, $\delta_a$ and $\delta_{rn}$ satisfying, for some $q>2 \vee (d-2)$,
\begin{equation}
\label{nabla_a_ineq}
1-\frac{q}{4}(\sqrt{\gamma} + \|a-I\|_\infty \sqrt{\delta+\delta_a}) > 0,
\end{equation}
where $\gamma:=\sum_{r,n=1}^d \delta_{rn}$, and
\begin{align}
(q-1)\big(1-\frac{q\sqrt{\gamma}}{2}\big) & - (\sqrt{\delta+\delta_a}\sqrt{\delta_a} \notag + \delta+\delta_a)\frac{q^2}{4} \\
& - (q-2)\frac{q\sqrt{\delta+\delta_a}}{2} -\|a-I\|_\infty \frac{q\sqrt{\delta+\delta_a}}{2}>0, \label{nabla_a_ineq2}
\end{align}
the authors
constructed a Feller semigroup and proved an analogue of Theorem \ref{thm1}(\textit{i}), including the ``approximation uniqueness'' result in Remark \ref{approx_uniq}, for the It\^{o} SDE
\begin{equation}
\label{sde3_n}
X_t=x-\int_0^t b(X_s)ds + \sqrt{2}\int_0^t\sigma(X_s)dW_s.
\end{equation}
The result for the Stratonovich SDE
\begin{equation}
\label{sde3_str_n}
X_t=x-\int_0^t b(X_s)ds + \sqrt{2}\int_0^t\sigma(X_s) \circ dW_s,
\end{equation}
in  \cite{KiS_Osaka} is valid under assumption \eqref{nabla_a_ineq}, \eqref{nabla_a_ineq2} but with $\delta$ replaced by $\delta+\delta_a+\delta_c$, where $\delta_c$ is the form-bound of $$c:=(c^i)_{i=1}^d, \quad \text{ where } c^i:=\frac{1}{\sqrt{2}}\sum_{r,j=1}^d  (\nabla_r \sigma_{ij})\sigma_{rj}.$$

\medskip

The assumptions \eqref{nabla_a_ineq}, \eqref{nabla_a_ineq2} imply that $\delta$, $\delta_a$ and $\delta_{rj}$ cannot be too large. It is also easily seen that if $a=I$, then these assumptions reduce to $\delta<1 \wedge (\frac{d}{d-2})^2$, i.e.\,then there exists $q>2 \vee (d-2)$ such that \eqref{nabla_a_ineq}, \eqref{nabla_a_ineq2} hold.

\medskip
The assumptions on diffusion coefficients $\sigma$ of the form \eqref{nabla_a} go back to  Veretennikov \cite{V} who proved strong well-posedness of \eqref{sde3_} for bounded measurable  $b$ and $\nabla_r \sigma_{ij} \in L^{2d}_{\loc}$. There are many other papers that consider assumptions of this type, see e.g.\,\cite{Z4} and \cite{Kr2,Kr3} who considered $\nabla_r \sigma_{ij} \in L^{p}_{\loc}$ ($p>d$) and $\nabla_r \sigma_{ij} \in L^{d}_{\loc}$.

\medskip

The construction of the Feller semigroup in \cite{KiS_Osaka} is based on an extension of the iteration procedure described in Section \ref{parab_sect} (in the elliptic setting) to solutions $u_n$ of divergence-form equations
\begin{equation}
\label{div_eq}
(\mu+\Lambda(a_n,b_n))u_n=f, \quad f \in C_c^\infty, \quad \mu \geq \mu_0,
\end{equation}
where $\Lambda(a_n,b_n)=-\nabla \cdot a_n \cdot \nabla + b_n \cdot \nabla$,
and uses the gradient bound 
\begin{equation}
\label{grad_a}
\|\nabla u_n\|_{\frac{qd}{d-2}} \leq K \|f\|_q, \quad \text{$K$ is independent of $n$}.
\end{equation}
Here $a_n$, $b_n$ are bounded and smooth, $\{b_n\}$ is uniformly (in $n$) form-bounded, and $\{a_n\}$ satisfy the same assumptions as $a$ above (thus, with constants independent on $n$).
This iteration procedure and \eqref{grad_a} yield Feller semigroup $e^{-t\Lambda_{C_\infty}(a,b)}$. 
Then $e^{-t\Lambda_{C_\infty}(a,\nabla a + b)}$ is the sought Feller semigroup that produces weak solution to It\^{o} SDE \eqref{sde3_n}, where we used the identity
\begin{equation}
\label{div_nondiv_id}
-a_n \cdot \nabla^2 + b_n \cdot \nabla=-\nabla \cdot a_n \cdot \nabla + (\nabla a_n + b_n )\cdot \nabla.
\end{equation}
For Stratonovich SDE \eqref{sde3_str_n} one needs Feller semigroup $e^{-t\Lambda_{C_\infty}(a,\nabla a + b-c)}$.

\begin{remark}For instance, the approximating vector fields $b_n$ can be defined via \eqref{b_n1}, so they are uniformly (in $n$) in $\mathbf{F}_\delta$. The approximating matrices $a_n$ can be defined via
$$
a_n=\eta_{\varepsilon_n} \ast a,
$$
where $\eta_{\varepsilon_n}$ is the Friedrichs mollifier, and $\varepsilon_n \downarrow 0$.
To see that $a_n$ are such that $\nabla a_n$ are indeed uniformly in $\mathbf{F}_{\delta_a}$ and satisfy \eqref{nabla_a} with the same $\delta_{rn}$, also uniformly in $n$, one can apply the result of Appendix \ref{example_sect}.3. In \cite{KiS_Osaka}, there was an additional cutoff function under the mollifier, which is not necessary. 
\end{remark}

Let us add that already the task of proving the uniform in $n$ gradient bound \eqref{grad_a} for 
solutions $u_n$ to the divergence-form equation \eqref{div_eq} (which was the original interest of the authors of \cite{KiS_Osaka}) leads to condition \eqref{nabla_a}.

\medskip

Theorem \ref{feller_parab} and the existence part of Theorem \ref{sde_parab} can be extended to SDEs \eqref{sde3_} and \eqref{sde3_str}
with time-inhomogeneous $b \in L^\infty\mathbf{F}_\delta + L^2_{\loc}(\mathbb R_+)$
and time-inhomogeneous bounded $\sigma$ such that $a=\sigma\sigma^{\scriptscriptstyle \top}$ is uniformly elliptic  and satisfies
$$
\nabla a \in L^\infty\mathbf{F}_{\delta_a}+L^2_{\loc}(\mathbb R_+), \quad (\nabla_{r}a_{ij})_{i=1}^d \in L^\infty\mathbf{F}_{\delta_{rj}}+L^2_{\loc}(\mathbb R_+)
$$
for all $1 \leq r,j \leq d$, for appropriate $\delta_a$ and $\delta_{rj}$. A direct extension of the weak uniqueness part of Theorem  \ref{sde_parab} to  $a_{ij}$ as above is problematic: one has to control second-order derivatives of $u_n$, but these are destroyed by form-bounded $b$. Thus, one needs extra assumptions both on $b$ and $a$, such as Morrey class in \cite{Kr_weak}, see below.

\medskip

\textbf{2.~}Below we assume, for simplicity, that $a$, $b$ are time-homogeneous, although most of the results cited below are valid for time-inhomogeneous coefficients. 

The condition \eqref{nabla_a} puts $a$ in the class ${\rm VMO}$\footnote{In \cite{KiS_Osaka} there is an incorrect statement that there are matrices $a$ satisfying \eqref{nabla_a}  and not contained in the {\rm VMO} class.}, see \cite{Kr3}. Recall that matrix $a$ is in ${\rm VMO}$ if
$$
\sup_{B_r}\frac{1}{|B_r|}\int_{B_r}\big|a-(a)_{B_r}\big|dx \rightarrow 0 \quad \text{ as } \rho \downarrow 0,
$$
where the supremum is taken over all balls of radius $\leq \rho$. Here $(a)_{B_r}$ denotes the average of $a$ on $B_r$.

There is a very rich literature on well-posedness of parabolic equations and SDEs  with ${\rm VMO}$  diffusion matrix and singular drift $b$ satisfying more restrictive assumptions than the form-boundedness.
The strongest result on weak well-posedness of SDE \eqref{sde3_} with ${\rm VMO}$ diffusion coefficients is a very recent result of  Krylov \cite{Kr_weak} who proved that there exist positive constants $\theta$ (sufficiently small) and $\rho_a$ such that if $a$ in the BMO class with norm $\leq \theta$, i.e.
\begin{equation}
\label{a_kr}
\sup_{B_r,\,r \leq \rho}\frac{1}{|B_r|}\int_{B_r}\big|a-(a)_{B_r}\big|dx \leq \theta \quad \text{ for all } \rho \leq \rho_a,
\end{equation}
(e.g.\,if $a \in {\rm VMO}$) and 
\begin{equation}
\label{M}
|b| \in M_{\frac{d}{2}+\varepsilon}, \quad \varepsilon>0,
\end{equation}
with sufficiently small norm,
then \eqref{sde3_} is weakly well-posed, i.e.\,the solution exists and is unique in a class similar to the one in Theorem \ref{sde_parab}.

In the case of constant diffusion coefficients $\sigma$  one can prove weak well-posedness of \eqref{sde3_} for substantially larger than $\mathbf{F}_\delta$ class of \textit{weakly form-bounded drifts}, discussed in Sections \ref{weak_sect}-\ref{weak_sde_sect}. This class contains e.g.
\begin{equation}
\label{M3}
|b| \in M_{1+\varepsilon}
\end{equation}
with arbitrarily small $\varepsilon>0$.

If we were to exploit the relationship between non-divergence and divergence form operators in the case the matrix $a$ is sufficiently discontinuous, cf.\,\eqref{div_nondiv_id}, then the sesquilinear form of the divergence-form operator will have to satisfy
\begin{equation}
\label{MV_bd2}
|\langle a \cdot \nabla \varphi,\nabla \psi\rangle + \langle b\cdot \nabla \varphi,\psi\rangle| \leq C\|\varphi\|_{W^{1,2}}\|\psi\|_{W^{1,2}}, \quad \varphi,\psi \in W^{1,2},
\end{equation}
which, by \cite{MV}, would make $\nabla a + b$ form-bounded (modulo a divergence-free component in ${\rm BMO}^{-1}$, see the discussion around \eqref{MV_bd}). See also Remarks \ref{a_weak} and \ref{a_rem} below.

\bigskip

\section{Stochastic transport equation and strong solutions to SDEs}

\label{ste_sect}

\medskip
In \cite{BFGM}, Beck-Flandoli-Gubinelli-Maurelli presented, among many results, an approach to proving strong well-posedness of SDE 
\begin{equation}
\label{sde_ste}
X_t=x-\int_0^t b(r,X_r)dr+\sigma W_t,
\end{equation}
for a.e.\,$x \in \mathbb R^d$, where $\{W_t\}_{t \geq 0}$ is a $d$-dimensional Brownian motion in $\mathbb R^d$ defined on a complete filtered probability space $(\Omega,\mathcal F_t,\mathcal F,\mathbf P)$, with $b$ in the critical Ladyzhenskaya-Prodi-Serrin class. Below we will discuss the time-homogeneous case $b:\mathbb R^d \rightarrow \mathbb R^d$, so their assumption on $b$ reads as $|b| \in L^d + L^\infty$.
Their approach is based on a detailed regularity theory of the stochastic transport equation (STE)
\begin{equation}
\label{eq1}
\begin{array}{c}
du+b \cdot\nabla u dt + \sigma \nabla u \circ dW_t=0 \quad \text{ on } (0,\infty) \times \mathbb R^d, 
\\[3mm]
u|_{t=0} = f ,
\end{array}
\end{equation}
where $u(t,x)$ is a scalar random field, $\sigma \neq 0$ is a constant, 
$f$ is in $L^p$ or $W^{1, p}$, $\circ$ is the Statonovich multiplication.

Speaking of the STE \eqref{eq1}, let us mention that the Cauchy problem for the deterministic transport equation $\partial_t u + b \cdot \nabla u=0$ is in general not well posed already for a bounded but discontinuous $b$. Moreover, in that case, even if the initial function $f$ is regular, one cannot hope that 
the corresponding solution $u$ will be regular 
immediately after $t=0$.
This, however, changes if one adds the noise term $\sigma \nabla u \circ dW_t$, $\sigma>0$. 
For the stochastic STE \eqref{eq1}, a unique weak solution exists and is regular for
some discontinuous $b$. 
This effect of regularization and well-posedness by noise, demonstrated by the STE, attracted considerable interest in the past few years, as a part of the more general program of establishing well-posedness by noise for SPDEs whose deterministic counterparts arising in fluid dynamics are not well-posed, see Flandoli-Gubinelli-Priola \cite{FGP}, Gess-Maurielli \cite{GM} for detailed discussions and references.

\medskip

Let us make a few preliminary remarks regarding STE \eqref{eq1} in the case the drift is smooth. 

\smallskip

1.~Let $b \in C_c^\infty(\mathbb R^d,\mathbb R^d)$ and 
$f \in C_c^\infty$.
Then there exists (see 
\cite[Theorem 6.1.9]{Ku})  a unique adapted strong solution to
\begin{equation}
\label{eq1_}
u(t) - f+\int_0^t b \cdot \nabla uds + \sigma\int_0^t \nabla u \circ dW_s=0 \text{ a.s.},
 \quad t \in [0,T],
\end{equation}
given by
\begin{equation}
\label{u_Psi}
u(t):=f(\Psi_{t}^{-1}), \quad t \geqslant 0,
\end{equation}
where $\Psi_{t}:\mathbb R^d \times \Omega \rightarrow \mathbb R^d$ is the stochastic flow for the SDE \eqref{sde_ste}.
The latter means that there exists $\Omega_0 \subset \Omega$, $\mathbb P(\Omega_0)=1$, such that, for all $\omega \in \Omega_0$,
$$\Psi_{t}(\cdot,\omega) \Psi_{s}(\cdot,\omega) = \Psi_{t+s}(\cdot,\omega), \quad \Psi_{0}(x,\omega)=x,$$ for every $x \in \mathbb R^d$, the process $t \mapsto \Psi_{t}(x,\omega)$ is a strong solution to \eqref{sde_ste},
and
$\Psi_{t}(x,\omega)$ is continuous in $(t,x)$, $\Psi_{t}(\cdot,\omega):\mathbb R^d \rightarrow \mathbb R^d$ are homeomorphisms  and $\Psi_{t}(\cdot,\omega)$, $\Psi_{t}^{-1}(\cdot,\omega) \in C^\infty(\mathbb R^d,\mathbb R^d)$.

\medskip

2.~Applying It\^{o}'s formula, one easily obtain that for every $\mu \geq 0$, 
\begin{equation*}
u(t)=e^{\mu t } f(\Psi_{t}^{-1}), \quad t \geqslant 0
\end{equation*}
solves
\begin{equation}
\label{eq2_}
u(t) - f+ \mu \int_0^t u ds +\int_0^t b \cdot \nabla uds + \sigma\int_0^t \nabla u \circ dW_s=0 \text{ a.s.},
 \quad t \in [0,T].
\end{equation}
Thus, solutions of 
the Cauchy problems \eqref{eq1_} and \eqref{eq2_} differ by the factor $e^{-\mu t}$.

\medskip

3.~One can rewrite the equation in \eqref{eq2_}, using the identity relating Stratonovich and It\^{o} integrals
\begin{equation}
\label{strat_id}
\int_0^t \nabla u \circ dW_s=\int_0^t \nabla u dW_s-\frac{1}{2}\sum_{k=1}^d[\partial_{x_k}u,W^k]_t, \qquad W_t=(W_t^k)_{k=1}^d,
\end{equation}
as
\begin{equation}
\label{eq_ito}
du + \mu u dt+b \cdot\nabla u dt + \sigma \nabla u  dW_t-\frac{\sigma^2}{2}\Delta u dt =0
\end{equation}
(the It\^{o} form of the STE). Now, taking expectation, one obtains that $v:=\mathbf{E}[u]$ solves Cauchy problem for the deterministic parabolic equation
\begin{equation*}
\partial_t v + \mu v + b \cdot \nabla v - \frac{\sigma^2}{2}\Delta v  =0, \quad
v|_{t=0} = f.
\end{equation*}

\medskip

Let now $b$ be discontinuous. The authors of \cite{BFGM}, in a sense, reversed \eqref{u_Psi}, i.e.\,given a $|b| \in L^d$ (in the time-homogeneous case) they used their Sobolev regularity theory of \eqref{eq1} to prove strong well-posedness of SDE \eqref{sde_ste} for a.e.\,initial point $x \in \mathbb R^d$. 

\medskip

In \cite{KiSS_transport}, the authors extended the approach of \cite{BFGM} to (time-homogeneous) form-bounded drifts $b$. We describe these results below.

\medskip

Set
\begin{align*}
\rho(x)  \equiv \rho_{\kappa, \theta}(x) 
:=(1+\kappa |x|^2)^{-\theta}, \quad \kappa>0,  \quad \theta>\frac{d}{2}, \quad x \in \mathbb R^d.
\end{align*}
Let $L^p_\rho \equiv L^p(\mathbb R^d,\rho dx)$. 
Denote by $\|\cdot\|_{p,\rho}$ the norm in $L^p_\rho$, and by $\langle \cdot ,\cdot \rangle_\rho$ the inner product in $L^2_\rho$.
Set $$W^{1,2}_\rho:=
\{g \in W^{1,2}_{\loc} \mid 
\|g\|_{W^{1,2}_\rho}:=\|g\|_{2,\rho}+\|\nabla g\|_{2,\rho}<\infty\}.$$
Fix $T>0$ and put $$\beta_{2q}:=1+4qd, \quad q=1,2,\dots$$

\begin{theorem}
\label{thm_ste}Let $b \in \mathbf{F}_\delta$ with
$\sqrt{\delta}<\frac{\sigma^2}{2\beta_2}$.
Let $p \geq 2$. 
Provided that $\kappa$ is chosen sufficiently small, there are generic constants
$\mu_1 \geq 0$, $C_1>0, C_2>0$ {\rm (}i.e.\,they depend only on $\delta$, $c_\delta$, $p$ and $T${\rm)} such that
for any $\mu \geq \mu_1$,  
for every $f \in L^{2p}$
there exists a function $u \in L^\infty([0,T], L^2(\Omega,L^2_\rho))$ 
for which the following are true:

{\rm(\textit{i})} For a.e.~$\omega \in \Omega$, $$\nabla \int_0^T u(s,\cdot,\omega)ds \in L_{\loc}^2(\mathbb R^d, \mathbb R^d),$$ so $$b \cdot \nabla \int_0^T u(s,\cdot,\omega)ds \in L_{\loc}^1,$$ and for every test function 
$\varphi \in C_c^\infty$, 
we have a.s.\,for all $t \in [0,T]$,
\begin{align}
& \langle u(t),\varphi \rangle - \langle f,\varphi\rangle  \notag \\
& + \mu\langle \int_0^t u ds,\varphi \rangle + \bigl\langle b \cdot \nabla \int_0^t  uds,\varphi \bigr\rangle - \sigma \bigl\langle \int_0^t u dW_s,\nabla \varphi \bigr\rangle + \frac{\sigma^2}{2}\bigl\langle \nabla \int_0^t u ds, \nabla \varphi\bigr\rangle=0.   \label{eq__}
\end{align}

{\rm(\textit{ii})} For any sequence of smooth vector fields $b_m \in C_c^\infty(\mathbb R^d,\mathbb R^d)$, $m=1,2,\dots,$ that are uniformly form-bounded in the sense that $b_m \in \mathbf{F}_\delta$ with $c_\delta$ independent of $m$, and are such that 
$$b_m \rightarrow b \text{ in $L^2_{\loc}(\mathbb R^d,\mathbb R^d)$ as $m \rightarrow \infty$},$$
we have for 
any initial function $f \in C_c^\infty$,
\begin{equation*}
u_m(t) \rightarrow u(t) \quad \text{ in } L^2(\Omega,L_\rho^2) \quad \text{ uniformly in $t \in [0,T]$},
\end{equation*}
where $u_m$ is the unique strong solution to \eqref{eq_ito} {\rm(}with $b=b_m${\rm)} with initial condition $u_m|_{t=0}=f$.
\end{theorem}

The last property implies that $u$ does not depend on the choice of the approximating sequence $\{b_m\}$ as long as it preserves the class of form-bounded vector fields. This can be viewed as a uniqueness result on its own.

\medskip

The next theorem establishes Sobolev regularity of $u$  
up to the initial time $t=0$.

\begin{theorem}
\label{thm2_ste}
Let $b \in \mathbf{F}_\delta$ with 
$\sqrt{\delta}<\frac{\sigma^2}{2\beta_2}$ and $f \in W^{1,4}$.  
Let $\kappa$ be sufficiently small and $\mu_1$ be the constant in Theorem \ref{thm_ste} with $p=2$. For $\mu\ge 
\mu_1$, let $u$ be the process constructed in Theorem \ref{thm_ste}. There exists generic constant $\mu_2 \ge \mu_1$
such that for $\mu\ge \mu_2$,
the following are true:

\medskip

{\rm(a)} $\mathbf Eu^2$, $\mathbf E|\nabla u|^2 \in
L^\infty([0,T],L^2)$, so $u \in L^\infty([0,T], 
L^2(\Omega,W^{1,2}_\rho))$.

\medskip

{\rm(b)} 
For any test function $\varphi \in C_c^\infty$, 
the process
$
t \mapsto \langle u(t), \varphi\rangle
$
is $(\mathcal F_t)$-progressively measurable and has a continuous $(\mathcal F_t)$-semi-martingale modification that satisfies a.s.\,for every $t \in [0,T]$,
\begin{align}
& \langle u(t),\varphi \rangle - \langle f,\varphi\rangle \notag    \\
& + \mu\int_0^t \langle u,\varphi \rangle ds + \int_0^t \bigl \langle b \cdot \nabla  u,\varphi \bigr\rangle ds  - \sigma  \int_0^t \langle u,\nabla \varphi \rangle dW_s + \frac{\sigma^2}{2}\int_0^t \bigl\langle  u, \Delta \varphi\bigr\rangle ds =0 \label{eq___}.
\end{align}
Moreover, if $\sqrt{\delta}<\frac{\sigma^2}{2\beta_{2q}}$ for some  $q=1,2,\dots$, 
then there exist generic constants 
$\mu_2(q)\ge \mu_1$ (with $\mu_2(1)$ equal to the $\mu_2$ above) and $C_1>0$ such that when $\mu\ge \mu_2(q)$ and $f \in W^{1,4q}$, we have
\begin{equation}
\label{grad_reg2}
\sup_{0 \leq \alpha \leq 1}\bigl\|\mathbf E|\nabla u|^{2q}\bigr\|_{L^{\frac{2}{1-\alpha}}([0,T],L^{\frac{2d}{d-2+2\alpha}})} \leq 
C_1\|\nabla f\|^{2q}_{4q}.
\end{equation}
In particular,  there exists generic $C_2>0$ such that
\begin{equation*}
\sup_{t \in [0,T]}\mathbf E\langle \rho |\nabla u|^{2q}\rangle \leq  C_2\|\nabla f\|^{2q}_{4q}.
\end{equation*}
If $2q>d$, then for a.e. $\omega \in \Omega$, $t \in [0,T]$, the function $x \mapsto u(t,x,\omega)$ is H\"{o}lder continuous, possibly after 
modification on a set of measure zero in $\mathbb R^d$ (in general, depending on $\omega$).
\end{theorem}

The estimate \eqref{grad_reg2} can be viewed as a counterpart of \eqref{bd_2}.

\medskip

A function satisfying (a), (b) of Theorem \ref{thm2_ste}
will be called a weak solution of Cauchy problem
\begin{equation}
\label{eq2}
\begin{array}{c}
du  + \mu\, u dt  +b \cdot\nabla u dt+ \sigma \nabla u \circ dW_t=0 \quad \text{ on } (0,\infty) \times \mathbb R^d, \\[3mm]
u|_{t=0} = f \in L^p, \quad p \geq 2.
\end{array}
\end{equation}
 This definition of weak solution is  close to \cite[Definition 2.13]{BFGM}. 

\begin{theorem}
\label{thm3_ste}
Let $b \in \mathbf{F}_\delta$ with 
$\sqrt{\delta}<\frac{\sigma^2}{2\beta_2}$ and $f\in W^{1, 4}$. 
Provided $\kappa$ is sufficiently small, there exists generic $\mu_3 \ge 0$ such that for $\mu\ge \mu_3$, 
the Cauchy problem \eqref{eq2}
has a unique weak solution in the class of functions satisfying 
{\rm(a)}, {\rm(b)} of Theorem \ref{thm2_ste}.
\end{theorem}

Theorems \ref{thm_ste}-\ref{thm3_ste} were proved in \cite{KiSS_transport}. Theorem \ref{thm2_ste} extends a similar result in \cite{BFGM} for (in the time-homogeneous case) $|b| \in L^d$. The proof of the uniqueness result in Theorem \ref{thm3_ste}  adopts the method of \cite[Sect.\,3]{BFGM}. 

\medskip

It should be noted that the authors in \cite{BFGM} prove their uniqueness result in a larger class of weak solutions (not requiring any differentiability, see \cite[Definition 3.3]{BFGM}) but under additional assumptions on $b$. 
Specialized to the time-dependent case, they assume that $b$ satisfies
\begin{equation}
\label{div}
{\rm div\,}b \in L^d + L^\infty 
\end{equation}
in addition to $b \in L^d + L^\infty$.
The latter is needed to establish \eqref{grad_reg2} for solutions of the adjoint equation to the STE, i.e.\,the stochastic continuity equation (which allows to prove an even stronger 
result: the uniqueness of 
weak solution to the corresponding  random transport equation),
see \cite[Sect.\,3]{BFGM}. 

\medskip

Armed with Theorems \ref{thm_ste}, \ref{thm2_ste}, one can repeat the argument in \cite[Sect.\,4]{BFGM} to prove the following result.
\textit{Assuming that $b \in \mathbf{F}_\delta$ with $\delta$ sufficiently small,  
there exists a stochastic Lagrangian flow for SDE \eqref{sde_ste}, i.e.\,a measurable map $\Phi:[0,T] \times \mathbb R^d \times \Omega \rightarrow \mathbb R^d$ such that, for a.e.~$x \in \mathbb R^d$, the process $t \mapsto \Phi_t(x,\omega)$ is a strong solution of the SDE \eqref{sde_ste}:
\begin{equation*}
\Phi_t(x,\omega)=x-\int_0^t b(s,\Phi_r(x,\omega))dr+\sigma W_t(\omega), \quad \text{a.s.}, \quad t \in [0,T],
\end{equation*}
and $\Phi_t(x,\cdot)$ is $\mathcal F_t$-progressively measurable.
If also $\sqrt{\delta}<\frac{\sigma^2}{2\beta_{2q}}$, $q=1,2,\dots$, then $\Phi_t(\cdot,\omega) \in W_{\loc}^{1,2q}$ ($t \in [0,T]$) for a.e. $\omega \in \Omega$. Moreover, $\Phi_t$ is unique, i.e.\,any two such stochastic flows coincide a.s. for every $t>0$ for a.e.\,$x$. }

\medskip

The restriction ``for a.e.\,initial point'' in the above strong existence result for  SDE \eqref{sde_ste} can be removed using a different method discussed in the next section.

\bigskip

\section{Strong well-posedness via R\"{o}ckner-Zhao's approach}

\label{rz_sect}

In \cite{RZ2}, R\"{o}ckner and Zhao proved strong existence and uniqueness in a large class of strong solutions satisfying Krylov estimate (cf.\,\eqref{krylov_rz} below) for SDE
\begin{equation}
\label{sde_rz}
X_t^x=x+\int_0^t b(s,X_s^x)ds + W_t, \quad 0 \leq t \leq T,
\end{equation}
for every initial point $x \in \mathbb R^d$, provided that drift $b:\mathbb R^{d+1} \rightarrow \mathbb R^d$ satisfies the critical Ladyzhenskaya-Prodi-Serrin condition
\begin{equation}
\label{LPS_rz}
b \in L^p(\mathbb R,L^q(\mathbb R^d)), \quad \frac{d}{q}+\frac{2}{p} \leq 1, \quad p > 2, \quad q \geq d.
\end{equation}
Above $\{W_t\}_{0 \leq t \leq T}$ denotes, as before, a Brownian motion on a complete filtered probability space $(\Omega,\{\mathcal F_t\}_{0 \leq t \leq T},\mathcal F,\mathbf{P})$.

\medskip

The method of R\"{o}ckner-Zhao is different from the other methods used in the literature on strong well-posedness of \eqref{sde_rz} (cf.\,\cite{BFGM} and \cite{Kr1}-\cite{Kr3}). Their proof of strong existence uses a relative compactness criterion for random fields on the Wiener-Sobolev space. Their proof of uniqueness uses their weak uniqueness result from \cite{RZ} and Cherny's theorem \cite{C} (briefly, strong existence and weak uniqueness $\Rightarrow$ strong uniqueness). The method of \cite{RZ2} is a far-reaching strengthening of the methods of Meyer-Brandis and Proske \cite{MP}, Mohammed-Nilsen-Proske \cite{MNP} (for $b \in L^\infty(\mathbb R \times \mathbb R^d)$) and Rezakhanlou \cite{R} (for $b$ in the sub-critical Ladyzhenskaya-Prodi-Serrin class).

\medskip

Let us give a brief outline of their method. For a given vector field $Y=(Y_i)_{i=1}^d:\mathbb R^{k} \rightarrow \mathbb R^m$, denote 
\begin{equation}
\label{nabla_Y}
\nabla Y=\nabla_x Y(x):=\left( 
\begin{array}{cccc}
\nabla_1 Y_1 & \nabla_2 Y_1 & \dots & \nabla_k Y_1 \\
& &  \dots & \\
\nabla_1 Y_m & \nabla_2 Y_m & \dots & \nabla_k Y_m
\end{array}
\right).
\end{equation}
The proof of strong existence in \cite{RZ2} is based on the following estimates. Let $b$ be additionally bounded and smooth.
For every $r \geq 1$, there exist constants $K_1$, $K_2$ independent of smoothness or boundedness of $b$ such that

\medskip

{\rm (\textit{i})} $\|\nabla X_t^x - I\|_{L^{2r}(\mathbb R^d,L^r(\Omega))} \leq K_1 t^{\frac{1}{2r}}$ for all $0 \leq t \leq T$;
\medskip

{\rm (\textit{ii})} $\|D_sX_t^x - I\|_{L^{2r}(\mathbb R^d,L^r(\Omega))} \leq K_1(t-s)^{\frac{1}{4r}}$ for a.e.\,$s \in [0,T]$ and $0 \leq s \leq t \leq T$,
where $D_sX_t^x$ denotes the Malliavin derivative; 
\medskip

{\rm (\textit{iii})} $\|D_s X_t^x - D_{s'}X_t^x\|_{L^{2r}(\mathbb R^d,L^r(\Omega))} \leq K_2|s-s'|^{\frac{1}{4r}}$ for a.e.\,$s,s' \in [0,T]$ and $0 \leq s,s' \leq t \leq T$,
\medskip

These estimates allow \cite{RZ2} to apply the relative compactness criterion  on the Wiener-Sobolev space in order to construct a strong solution to \eqref{sde_rz}. 

\medskip

The first step in the proof of (\textit{i})-(\textit{iii}) is  to differentiate SDE \eqref{sde_rz}, e.g.\,for (\textit{i})
$$
\nabla X_t^x-I=\int_0^t \nabla b(s,X_s^x)\nabla X_s^x ds,
$$
with the goal of iterating this identity and, in the end, obtaining an expression for the left-hand side that one can control. This goal is achieved using the following bound (for (\textit{i})): there exist positive generic constants $C_0$, $K$ such that, for every $n \geq 1$,
\begin{equation}
\label{est_rz}
\int_{\mathbb R^d }\left|\mathbf{E} \int_{\Delta_{n}(T_0,T_1)} \prod_{i=1}^n \nabla_{\alpha_i} f_i(t_i, X_{t_i}^x)dt_1 \dots dt_n \right|^2 dx \leq C_0K^n(T_1-T_0),
\end{equation}
where $1 \leq \alpha_i \leq d$ ($i \geq 1$),  $\nabla_{i}:=\partial_{x_i}$, $\Delta_n(T_0,T_1):=\{(t_1,\dots,t_n) \mid T_0 \leq t_1 \leq \dots \leq t_n\ \leq T_1\}$, and functions $f_i$ are taken to be the components of drift $b$. In \cite{RZ2}, \eqref{est_rz} is proved using Sobolev regularity estimates for solutions of parabolic equations with distributional right-hand side. These estimates are quite strong and are interesting on their own.

\medskip

In \cite{KiM_strong}, the authors noticed that \eqref{est_rz} can be proved for time-inhomogeneous form-bounded drifts $$b \in L^\infty\mathbf{F}_\delta+L^{2+\varepsilon}_{\loc}(\mathbb R), \quad \text{(i.e.\,with function $g_\delta \in L_{\loc}^{1+\varepsilon/2}(\mathbb R)$)}, \quad \varepsilon>0, 
$$
that can have stronger spatial singularities than than the drifts in \eqref{LPS_rz},
using a different argument which applies repeatedly integration by parts, quadratic inequality and the form-boundedness of $b$. The rest of the proof of strong existence essentially repeats \cite{RZ2}. This, combined with the weak uniqueness results from \cite{KiM} and \cite{Ki_Morrey}  via Cherny's theorem (see Sections \ref{parab_sde_sect} and \ref{morrey_sect}),  yields the following result.

\begin{theorem}[\cite{KiM_strong}]
\label{thm_strong}
Assume that $b \in L^\infty\mathbf{F}_\delta+L^{2+\varepsilon}_{\loc}(\mathbb R)$ for some $\varepsilon>0$. Also, assume that $b$ has compact support. Then, provided that form-bound $\delta$ is sufficiently small, for every $x \in \mathbb R^d$, SDE \eqref{sde_rz} has a strong solution $X_t^x$. This strong solution satisfies the following Krylov-type bounds:

{\rm 1)} For a given $q \in ]d,\delta^{-\frac{1}{2}}[$ and any vector field $\mathsf{g} \in L^\infty\mathbf{F}_{\delta_1}+L^{2+\varepsilon}_{\loc}(\mathbb R)$, $\delta_1<\infty$,
\begin{equation}
\label{krylov0_rz}
\mathbf E \int_0^T |\mathsf{g}h|(\tau,X_{0,\tau}^x)d\tau \leq c\|\mathsf{g}|h|^{\frac{q}{2}}\|^{\frac{2}{q}}_{L^2([0,T] \times \mathbb R^d)} \quad \text{ for all }h \in C_c([0,T] \times \mathbb R^d).
\end{equation}

{\rm 2)} For a given $\mu>\frac{d+2}{2}$, there exists constant $C$ such that
\begin{equation}
\label{krylov_rz}
\mathbf{E}\bigg[\int_{0}^{T}|h(\tau,X_{0,\tau}^{x})|d\tau\bigg] \leq C\|h\|_{L^{\mu}([0,T] \times \mathbb R^d)} \quad \text{ for all }h \in C_c([0,T] \times \mathbb R^d).
\end{equation}

Solution $X_t^x$ is unique among strong solutions to \eqref{sde_rz} that satisfy \eqref{krylov0_rz} for some $q \in ]d,\delta^{-\frac{1}{2}}[$ with $\mathsf{g}=1$ and with $\mathsf{g}=b$. 
If, in addition to the hypothesis on $b$, one has $|b| \in L^{\frac{d+2}{2}+\epsilon}$ for some $\epsilon>0$, then $X_t^x$ is unique among strong solutions to \eqref{sde_rz} that satisfy \eqref{krylov_rz}.
\end{theorem}

The assumption that $b$ has  compact support can be removed with an additional effort (using weight \eqref{rho}).
\medskip

Speaking of ``$\varepsilon>0$'' in $b \in L^\infty\mathbf{F}_\delta+L^{2+\varepsilon}_{\loc}(\mathbb R)$, this assumption does not allow us to include completely the critical Ladyzhenskaya-Prodi-Serrin class even with $p>2$, see above, as is assumed in \cite{RZ2}. It does include, however, the case that interests us the most: $p=\infty$, $q=d$. It also includes with case $p>2$, $q=\infty$.

\medskip

The weak well-posedness of SDE \eqref{sde_rz} is known to hold for larger classes of singular drifts than class $L^\infty \mathbf{F}_\delta + L^2_{\loc}(\mathbb R)$ discussed in Theorem \ref{thm_strong}. This is the subject of the next three sections.

\bigskip

\section{More singular than form-bounded. Semigroup in $\mathcal W^{{\scriptscriptstyle \frac{1}{2},2}}$}

\label{weak_sect}

In Sections \ref{weak_sect}, \ref{weak_sde_sect} and \ref{morrey_sect} we expand the classes of singular drifts considered in Theorems \ref{thm1_feller}, \ref{thm1} and \ref{feller_parab}, \ref{sde_parab}, although at expense of requiring smaller $\delta$.

\begin{definition}
\label{wfb}
A Borel measurable vector field $b:\mathbb R^d \rightarrow \mathbb R^d$ with $|b|\in L^1_{\loc}$ is said to be weakly form-bounded if there exists a constant $\delta>0$ such that
$$
\||b|^{\frac{1}{2}}(\lambda-\Delta)^{-\frac{1}{4}}\|_{2 \rightarrow 2} \leq \sqrt{\delta}
$$
for some $\lambda=\lambda_\delta \geq 0$. This is written as $b \in \mathbf{F}_\delta^{\scriptscriptstyle \frac{1}{2}}$.
\end{definition}

There is an important difference between form-bounded vector fields and weakly form-bounded vector fields. 
Namely, when we dealt with $b \in \mathbf{F}_\delta$, we controlled the gradient term $b \cdot \nabla$ in the Kolmogorov operator $-\Delta + b \cdot \nabla$ using the quadratic (Cauchy-Schwarz) inequality
\begin{equation}
\label{quad2}
|\langle b \cdot \nabla  \varphi,\varphi\rangle| \leq \varepsilon \|b\varphi\|_2^2+\frac{1}{4\varepsilon}\|\nabla \varphi\|_2^2, \quad \varepsilon>0
\end{equation}
(e.g.\,\eqref{quad} in the verification of conditions of the Lax-Milgram theorem and the KLMN theorem in $L^2$, in the proof of Lemma \ref{lem_GPQ}, in the iteration procedure and gradient bound \eqref{grad_bd}, etc). \textit{We can no longer do this  when dealing with $b \in \mathbf{F}_\delta^{\scriptscriptstyle 1/2}$ if only because $|b|$ is in general no longer locally in $L^2$.}

The form-bounded vector fields are weakly form-bounded. To show this, let us recall that the condition $b \in \mathbf{F}_\delta$ can be stated as an operator-norm inequality
$
\||b|(\lambda-\Delta)^{-\frac{1}{2}}\|_{2 \rightarrow 2} \leq \sqrt{\delta}.
$
The Heinz-Kato inequality \cite{He} allows us to take square roots in the operators that constitute the left-hand side, so we arrive at
$
\||b|^{\frac{1}{2}}(\lambda-\Delta)^{-\frac{1}{4}}\|_{2 \rightarrow 2} \leq \delta^{\frac{1}{4}},
$
and hence
$$
b \in \mathbf{F}_\delta \quad \Rightarrow \quad b \in \mathbf{F}_{\sqrt{\delta}}^{\scriptscriptstyle \frac{1}{2}}.
$$
The opposite inclusion is invalid: there are weakly form-bounded vector fields that are not form-bounded:

\begin{examples}

1.~If $|b|$ belongs to scaling-invariant Morrey class $M_{1+\varepsilon}$, for arbitrarily fixed small $\varepsilon>0$, i.e.
\begin{equation}
\label{elliptic_morrey2}
|b| \in L^{1+\varepsilon}_{\loc} \quad \text{ and } \quad
\|b\|_{M_{1+\varepsilon}}:=\sup_{r>0, x \in \mathbb R^d} r\biggl(\frac{1}{|B_r|}\int_{B_r(x)}|b|^{1+\varepsilon}dx \biggr)^{\frac{1}{1+\varepsilon}}<\infty,
\end{equation}
then $b \in \mathbf{F}_\delta^{\scriptscriptstyle 1/2}$. The proof of this inclusion follows right away from
\cite[Theorem 7.3]{A}. 

Recalling that the class of form-bounded vector fields $\mathbf{F}_\delta$ satisfies $M_{2+\varepsilon} \subset \mathbf{F}_\delta \subset M_2$ (say, with $c_\delta=0$), one can see that we gain quite a lot in admissible singularities of $b$  by working with $\mathbf{F}_\delta^{\scriptscriptstyle 1/2}$. In particular, we gain all $b$ with $|b| \in M_{1+\varepsilon} - M_2$.

\smallskip

2.~Recall that a vector field $b:\mathbb R^d \rightarrow \mathbb R^d$ is said to belong to the Kato class if $|b| \in L^1_\loc$ and
\begin{equation}
\label{kato}
\|(\lambda-\Delta)^{-\frac{1}{2}}|b|\|_{\infty}  \leq \sqrt{\delta}
\end{equation}
for some $\delta>0$ and $\lambda=\lambda_\delta \geq 0$. This is written as $b \in \mathbf{K}_\delta^{d+1}$. The Kato class vector fields 
are weakly form-bounded. 
Indeed, if $b \in \mathbf{K}_\delta^{d+1}$, then by duality one has 
\begin{equation}
\label{kato_dual}
\||b|(\lambda-\Delta)^{-\frac{1}{2}}\|_{1 \rightarrow 1}  \leq \sqrt{\delta}. 
\end{equation}
Applying Stein's interpolation between \eqref{kato} and \eqref{kato_dual}, one has
$\||b|^\frac{1}{2}(\lambda-\Delta)^{-\frac{1}{2}}|b|^{\frac{1}{2}}\|_{2 \rightarrow 2} \leq \sqrt{\delta}$; in the left-hand side of the last inequality one has the norm of the product of $|b|^\frac{1}{2}(\lambda-\Delta)^{-\frac{1}{4}}$ and its adjoint. By a standard result this yields $\||b|^\frac{1}{2}(\lambda-\Delta)^{-\frac{1}{4}}\|_{2 \rightarrow 2} \leq \delta^{\frac{1}{4}}$. Thus, 
$$b \in \mathbf{K}^{d+1}_\delta \quad \Rightarrow \quad b \in \mathbf{F}_{\sqrt{\delta}}^{\scriptscriptstyle 1/2}.$$

As was mentioned earlier, the Kato class $b \in \mathbf{K}_\delta^{d+1}$ with $\delta$ sufficiently small provides two-sided Gaussian bounds on the heat kernel of Kolmogorov operator $-\Delta + b \cdot \nabla$ \cite{Z_Gaussian}. The Kato class also provides  uniqueness in law for SDE \eqref{sde_w}, see \cite{BC}. There the authors required $\delta$ to be arbitrarily small. In fact, they show that under the Kato class assumption on $b$ the gradient of solution of elliptic equation $(\mu-\Delta + b \cdot \nabla)v=f$ is bounded. The reader can compare this with Remark \ref{grad_unb_rem} concerning the gradient of $v$ for a form-bounded $b$.

The Kato class $\mathbf{K}_\delta^{d+1}$ does not contain a popular class $|b| \in L^d$ (if only because there are vector fields $b$ with $|b| \in L^d$ that destroy two-sided Gaussian bounds on $-\Delta + b \cdot \nabla$), and so it does not contain $\mathbf{F}_\delta$. On the other hand, the Kato class is not contained in $\mathbf{F}_\delta$. However, both form-bounded and Kato class vector fields are contained in $\mathbf{F}_\delta^{\scriptscriptstyle 1/2}$.

\end{examples}

Let us demonstrate one way to arrive at the condition $b \in \mathbf{F}_\delta^{\scriptscriptstyle 1/2}$, $\delta<1$. First, let $b \in \mathbf{F}_\delta$, $\delta<1$, and, for brevity, assume that $c_\delta=0$. Also, let $b$ be bounded and smooth so that all manipulations with the equation are justified, but the constants in the estimates below will not depend on the smoothness or boundedness of $b$. One can prove the following two  $L^2$ regularity results 
for Cauchy problem $(\partial_t - \Delta + b \cdot\nabla)u = 0$, $u(0)=f$ with such $b$:

-- Multiplying $(\partial_t - \Delta + b \cdot\nabla)u = 0$ by $u$, integrating over $[0,t] \times \mathbb R^d$ and applying $b \in \mathbf{F}_\delta$ and quadratic inequality \eqref{quad2}, we obtain energy inequality
\begin{equation}
\label{u_r}
\frac{1}{2}\|u(t)\|_2^2 + (1-\sqrt{\delta})\int_0^t \|\nabla u\|_2^2dr \leq \frac{1}{2}\|f\|_2^2.
\end{equation}

-- Multiplying $(\partial_t - \Delta + b \cdot\nabla)u = 0$ by $-\Delta u$, integrating over $[0,t] \times \mathbb R^d$, we obtain
$$
\frac{1}{2}\partial_t \|\nabla u\|_2^2 + \|\Delta u\|_2^2 + \langle b \cdot \nabla u,-\Delta u\rangle=0,
$$
where we further estimate, using $b \in \mathbf{F}_\delta$ with $c_\delta=0$,
\begin{align*}
|\langle b \cdot \nabla u,-\Delta u\rangle| & \leq \varepsilon \|b|\nabla u|\|_2 + \frac{1}{4\varepsilon}\|\Delta u\|_2^2  \\
& \leq \varepsilon \delta \|(-\Delta)^{\frac{1}{2}}|\nabla u|\|_2^2 + \frac{1}{4\varepsilon}\|\Delta u\|_2^2 \\
& (\text{use Beurling-Deny-type inequality $ \|(-\Delta)^{\frac{1}{2}}|\nabla u|\|^2_2 \leq \|(-\Delta)^{\frac{1}{2}}\nabla u\|_2^2 \equiv \|\Delta u\|_2^2$}\\
& \text{and select $\varepsilon=\frac{1}{2\sqrt{\delta}}$}) \\
& \leq  \sqrt{\delta}\|\Delta u\|_2^2.
\end{align*}
Thus, we obtain another ``energy inequality'':
\begin{equation}
\label{u2_r}
\frac{1}{2}\|\nabla u(t)\|_2^2 + (1-\sqrt{\delta})\int_0^t \|\Delta u\|_2^2dr \leq \frac{1}{2}\|\nabla f\|_2^2.
\end{equation}

-- Now, one can ask what happens if we multiply $(\partial_t - \Delta + b \cdot\nabla)u = 0$ by an intermediate test function $(-\Delta)^s u$, $0<s<1$. One obtains an intermediate result between \eqref{u_r} and \eqref{u2_r}, but for a larger class of vector fields $b$, which becomes maximal if one multiplies by $(-\Delta)^{\frac{1}{2}}u$:
\begin{align*}
|\langle b \cdot \nabla u, (-\Delta)^\frac{1}{2} u\rangle| & =|\langle b^{\frac{1}{2}}(-\Delta)^{-\frac{1}{4}}\nabla (-\Delta)^{\frac{1}{4}}u,|b|^{\frac{1}{2}}(-\Delta)^{-\frac{1}{4}}(-\Delta)^{\frac{3}{4}}u\rangle| \qquad b^{\frac{1}{2}}:=b|b|^{-\frac{1}{2}}\\
& \leq \||b|^{\frac{1}{2}}(-\Delta)^{-\frac{1}{4}}\|_{2 \rightarrow 2} \|\nabla (-\Delta)^{\frac{1}{4}}u\|_2 \, \||b|^{\frac{1}{2}}(-\Delta)^{-\frac{1}{4}}\|_{2 \rightarrow 2} \|(-\Delta)^{\frac{3}{4}}u\|_2 \\
& = \||b|^{\frac{1}{2}}(-\Delta)^{-\frac{1}{4}}\|^2_{2 \rightarrow 2} \|(-\Delta)^{\frac{3}{4}}u\|^2_2.
\end{align*}
Thus, requiring $$\||b|^{\frac{1}{2}}(-\Delta)^{-\frac{1}{4}}\|_{2 \rightarrow 2} \leq \sqrt{\delta} \quad \text{(i.e.\,$b \in \mathbf{F}_\delta^{\scriptscriptstyle 1/2}$ with $\lambda=0$)}, \quad \delta<1,$$ one obtains the following ``energy inequality'':
\begin{equation}
\label{e_2}
\frac{1}{2}\|(-\Delta)^{\frac{1}{4}} u(t)\|_2^2 + (1-\delta)\int_0^t \|(-\Delta)^{\frac{3}{4}} u\|_2^2dr \leq \frac{1}{2}\|(-\Delta)^{\frac{1}{4}} f\|_2^2.
\end{equation}
From the look of \eqref{e_2}, it is seen that one needs to work with the chain of Bessel spaces \begin{equation}
\label{scale}
\mathcal W^{\frac{3}{2},2} \subset \mathcal W^{\frac{1}{2},2} \subset \mathcal W^{-\frac{1}{2},2}
\end{equation}
rather than the standard Sobolev triple $W^{1,2} \subset L^2  \subset W^{-1,2}$ (above $\mathcal W^{-\frac{1}{2},2}$ is the dual of $\mathcal W^{\frac{3}{2},2}$ with respect to the inner product in $\mathcal W^{\frac{1}{2},2}$) . Of course, by doing that, one sacrifices \eqref{MV_bd} and loses the possibility to consider  general operator $-\nabla \cdot a \cdot \nabla + b \cdot \nabla$ unless uniformly elliptic matrix $a$ satisfies additional regularity assumptions that make $-\nabla \cdot a \cdot \nabla$ a bounded operator from $\mathcal W^{\frac{3}{2},2}$ to $\mathcal W^{-\frac{1}{2},2}$, see Remark \ref{a_weak}. 
In fact, the following result is true:

\begin{proposition} \label{prop_weak}Let $b \in \mathbf{F}_\delta^{\scriptscriptstyle \frac{1}{2}}$, $\delta<1$.
Then for every $f \in \mathcal W^{\frac{1}{2},2}$ there exists a unique weak solution to Cauchy problem 
\begin{equation}
\label{cp2}
(\partial_t + \lambda- \Delta + b \cdot \nabla)u=0, \quad u(0+)=f
\end{equation}
where $\lambda$ is from the condition $b \in \mathbf{F}_\delta^{\scriptscriptstyle 1/2}$, i.e.\,a unique
 in $L^\infty_{\loc}(]0,\infty[,\mathcal W^{\frac{1}{2},2}) \cap L^2_{\loc}(]0,\infty[,\mathcal W^{\frac{3}{2},2})$ function $u$  satisfying 
\begin{align*}
 \int_{0}^{\infty} \langle (\lambda-\Delta)^{\frac{1}{4}}u,\partial_t (\lambda-\Delta)^{\frac{1}{4}}\varphi \rangle dt & = \int_{0}^{\infty} \langle (\lambda-\Delta)^{\frac{3}{4}}u,(\lambda-\Delta)^{\frac{3}{4}}\varphi \rangle dt  \notag
\\
& + \int_{0}^{\infty} \langle b \cdot \nabla u,(\lambda-\Delta)^{\frac{1}{2}}\varphi\rangle
\end{align*}
for all $\varphi \in C_c^\infty(]0,\infty[,\mathcal S)$ and
$w\mbox{-}{\mathcal W^{\frac{1}{2},2}}\mbox{-}\lim_{t \downarrow 0}u(t)=f$.
One has $u \in C(\mathbb R_+,\mathcal W^{\frac{1}{2},2})$, the following energy inequality holds:
\begin{equation*}
\|u(t)\|_{\mathcal W^{\frac{1}{2},2}}^2 + (1-\delta)\int_{0}^{t} \|u\|^2_{\mathcal W^{\frac{3}{2},2}}dr \leq  \|f\|_{\mathcal W^{\frac{1}{2},2}}^2, \quad t \geq 0,
\end{equation*}
and $T^{t}f(\cdot):=u(t,\cdot)$ is a contraction strongly continuous in $\mathcal W^{\frac{1}{2},2}$ Markov semigroup. 
If $\{b_\varepsilon\}_{\varepsilon>0}$ is a family of bounded smooth vector fields such that $b_\varepsilon \in \mathbf{F}_\delta^{\scriptscriptstyle 1/2}$ with the same $\lambda$ as $b$, 
$
b_\varepsilon \rightarrow b$ in $[L^1_{\loc}]^d$ as $\varepsilon \to 0
$, and
if $u_\varepsilon$ denotes the  solution to Cauchy problem \eqref{cp2} with the vector field $b_\varepsilon$, then
$$
u_\varepsilon \rightarrow u \quad \text{weakly in } L^2_{\loc}(\mathbb R_+,\mathcal W^{\frac{3}{2},2}) \text{ as } \varepsilon \to 0.
$$

\end{proposition}

The proof uses the standard J.-L.\,Lions approach in the scale \eqref{scale}, and was carried out in \cite{KiS_BMO}, in fact, in greater generality: for time-inhomogeneous $b \in L^\infty\mathbf{F}_\delta^{\scriptscriptstyle 1/2}$, that is, satisfying for a.e.\,$t \in \mathbb R_+$ the operator inequality $$\||b(t)|^{\frac{1}{2}}(\lambda-\Delta)^{-\frac{1}{4}}\|_{2 \rightarrow 2} \leq \sqrt{\delta}$$ for some fixed $\lambda=\lambda_\delta$.

\medskip

The above argument leading to the energy inequality \eqref{e_2} is not how the class $\mathbf{F}_\delta^{\scriptscriptstyle 1/2}$ first appeared in the literature. The ($L^p,L^q$) estimate
\begin{equation}
\label{growth2}
\|e^{-t\Lambda_p(b)}f\|_{q} \leq ce^{t\omega_p}t^{-\frac{d}{2}(\frac{1}{p}-\frac{1}{q})}\|f\|_p, \quad t>0, \quad \omega_p:=\frac{c_\delta}{2(p-1)},
\end{equation}
can be proved separately for $b \in \mathbf{F}_\delta$, see \eqref{growth}, and for $b \in \mathbf{K}^{d+1}_\delta$ (as was mentioned above, for the Kato class one even has two-sided Gaussian bounds on the integral kernel $e^{-t\Lambda(b)}(x,y)$ of $e^{-t\Lambda_p(b)}$). 
It was noticed in \cite{S} that the validity of estimate \eqref{growth2} depends, in fact, only on the weaker
condition $\||b|^{\frac{1}{2}}(\lambda-\Delta)^{-\frac{1}{4}}\|_{2 \rightarrow 2} \leq \sqrt{\delta}$ with $\delta<1$, which led to the introduction of the class $\mathbf{F}_\delta^{\scriptscriptstyle 1/2}$. 
Also, \cite{S} proposed a way to construct a \textit{quasi bounded} semigroup in $L^2$ associated with $-\Delta + b \cdot \nabla$ with weakly form-bounded $b \in \mathbf{F}_\delta^{\scriptscriptstyle 1/2}$ by constructing its resolvent as the operator-valued function
$$
\Phi(\zeta,b):=(\zeta-\Delta)^{-\frac{3}{4}}(1+H S)^{-1}(\zeta-\Delta)^{-\frac{1}{4}}, \qquad \Real\,\zeta \geq \lambda_\delta,
$$
where, by $b \in \mathbf{F}_\delta^{\scriptscriptstyle 1/2}$, operators 
\begin{equation}
\label{HS}
H:= (\zeta-\Delta)^{-\frac{1}{4}}|b|^{\frac{1}{2}}, \quad S:=b^{\frac{1}{2}} \cdot \nabla(\zeta-\Delta)^{-\frac{3}{4}}
\end{equation}
 are bounded on $L^2$ with norm $\sqrt{\delta}$ each (for operator $S$, taking into account that $\nabla (\zeta-\Delta)^{-\frac{1}{2}}$ is bounded on $L^2$ with norm $1$), and so $\Phi(\zeta,b)$ is bounded on $L^2$. Then, it is easily seen, $\Phi(\zeta,b)$ satisfies 
\begin{equation}
\label{hol}
\|\Phi(\zeta,b)\|_{2 \rightarrow 2} \leq (1-\delta)^{-1}|\zeta|^{-1}, \quad \text{ on } \{\Real\,\zeta \geq \lambda_\delta\}.
\end{equation}
 The proof that $\Phi(\zeta,b)$ is indeed the resolvent of the generator $\Lambda$ of a quasi bounded strongly continuous semigroup\footnote{But not a quasi contraction semigroup, as in Section \ref{basic_sde_sect}.}  $e^{-t\Lambda}$ on $L^2$ uses the general Trotter approximation theorem. The latter, in practice, requires the uniform in $n$ estimate \eqref{hol} for $\Phi(\zeta,b_n)$, where $b_n$ are approximating vector fields for $b$ (cf.\,\eqref{b_n2_weak}, \eqref{b_n3_weak}). In other words, it is essential for the construction that one is working with a holomorphic semigroup.
We refer to  \cite{KiS_theory} for detailed discussion.

\medskip

Note also that one no longer has \eqref{kolm}, i.e.\,
$$
\Lambda \not\supset -\Delta + b \cdot \nabla \upharpoonright C_c^\infty.
$$
The reason is that for a weakly form-bounded $b$ its norm $|b|$ is in general not in $L^2_{\loc}$.

\medskip

If $b$ is form-bounded or in the Kato class of vector fields -- two standard assumptions -- then one can construct a realization of $-\Delta +b \cdot \nabla$ as the generator of a strongly continuous semigroup \textit{in some }$L^p$ by invoking the KLMN theorem in $L^2$ or the Miyadera theorem in $L^1$, respectively (see e.g.\,\cite{KiS_theory}). However these two theorems (and, generally speaking, the standard perturbation-theoretic tools) are inapplicable to $-\Delta +b \cdot \nabla$  in any $L^p$ if $b$ is in the class of weakly form-bounded drifts $\mathbf{F}_\delta^{\scriptscriptstyle 1/2}$.

\begin{remark}
\label{a_weak}
Let us comment on the assumptions on a measurable uniformly elliptic matrix $a$ (i.e.\,$\sigma I \leq a \leq \xi I$ a.e.\,on $\mathbb R^d$ for $0<\sigma \leq \xi<\infty$) that would allow to extend Proposition \ref{prop_weak} to operator $-\nabla \cdot a \cdot \nabla + b \cdot \nabla$. If $b \in \mathbf{F}_\delta^{\scriptscriptstyle 1/2}$, then it is easily seen that $$b \cdot \nabla \in \mathcal B(\mathcal W^{\frac{3}{2},2},\mathcal W^{-\frac{1}{2},2}).$$ The matrix $a$ has to be such that 
\begin{equation}
\label{a_emb}
-\nabla \cdot a \cdot \nabla \in \mathcal B(\mathcal W^{\frac{3}{2},2},\mathcal W^{-\frac{1}{2},2})
\end{equation}
(of course, if $a$ is only measurable uniformly elliptic, then one only has $-\nabla \cdot a \cdot \nabla \in \mathcal B(W^{1,2},W^{-1,2})$). Let us mention one elementary sufficient condition for \eqref{a_emb}. For simplicity we will stay at the a priori level (i.e.\,the matrix is smooth but the norm of the operator in \eqref{a_emb} does not depend on smoothness of $a$). Also, assume that $a=I+a^0$ where $a^0$ has entries $a^0_{ij}$ in $\mathcal S$. For given $\varphi \in \mathcal W^{\frac{3}{2},2}$, $\psi \in \mathcal W^{\frac{1}{2},2}$, we have
\begin{align*}
|\langle -\nabla \cdot a^0 \cdot \nabla \varphi,\psi \rangle| & = |\langle (1-\Delta)^{\frac{1}{4}} a^0  (1-\Delta)^{-\frac{1}{4}}\cdot(1-\Delta)^{\frac{1}{4}} \nabla \varphi,(1-\Delta)^{-\frac{1}{4}}\nabla \psi\rangle| \\
& \leq \| (1-\Delta)^{\frac{1}{4}} a^0 (1-\Delta)^{-\frac{1}{4}}\|_{2 \rightarrow 2}\|\varphi\|_{\mathcal W^{\frac{3}{2},2}}\|\psi\|_{\mathcal W^{\frac{1}{2},2}},
\end{align*}
where, in turn, by the Kato-Ponce inequality ($\equiv$\,\,fractional Leibnitz rule) \cite{GO}, for all $1 \leq i,j \leq d$,
\begin{align*}
\|(1-\Delta)^{\frac{1}{4}} a^0_{ij} (1-\Delta)^{-\frac{1}{4}}f\|_2 & \leq \|(1-\Delta)^{\frac{1}{4}}a^0_{ij}\|_{2d}\|(1-\Delta)^{-\frac{1}{4}}f\|_{\frac{2d}{d-1}} + \|a_{ij}\|_{\infty}\|f\|_{2}.
\end{align*}
Thus, if 
$$\|a^0_{ij}\|_{\mathcal W^{\frac{1}{2},2d}} \leq c < \infty$$
for all $i, j$, for a generic constant $c$ (i.e.\,a constant that does not depend on the smoothness of $a^0_{ij}$),
then, using $\|(1-\Delta)^{-\frac{1}{4}}f\|_{\frac{2d}{d-1}} \leq C\|f\|_2$, we obtain $\|(1-\Delta)^{\frac{1}{4}} a^0 (1-\Delta)^{-\frac{1}{4}}\|_{2 \rightarrow 2} \leq c'$ for a generic $c'$, and hence \eqref{a_emb} with the operator norm bounded by a generic constant.
\end{remark}

\begin{remark}
\label{kato_rem}
At the level of Feller semigroups, one cannot really draw a parallel between the Kato class of potentials $\mathbf{K}^{d}_\delta$ (see \eqref{kato_pt_def}) and the Kato class of drifts $\mathbf{K}^{d+1}_\delta$. Indeed, in view of the results \cite{OSSV}, for the Schr\"{o}dinger operator $-\Delta + V$ condition $V \in \mathbf{K}^{d}_\delta$ is, basically, necessary and sufficient for the Feller property ($\equiv$\,\,strong continuity of the semigroup on $C_\infty$) to hold. For the Kolmogorov operator $-\Delta + b \cdot \nabla$,  condition $b \in \mathbf{K}^{d+1}_{\delta}$ is only sufficient for the Feller property; as Theorem \ref{thm1_feller_weak} below shows, one can go much farther, to weakly form-bounded drifts.
\end{remark}

\bigskip

\section{Weakly form-bounded drifts and SDEs} 

\label{weak_sde_sect}

In this section we construct the Feller semigroup for $-\Delta + b \cdot \nabla$ and prove weak well-posedness of SDE
\begin{equation}
\label{sde_w}
X_t=x-\int_0^t b(X_r)dr + \sqrt{2}W_t, \quad t \geq 0, 
\end{equation}
for a fixed $x \in \mathbb R^d$,
with $b:\mathbb R^d \rightarrow \mathbb R^d$ in the class of \textit{weakly form-bounded vector fields} $\mathbf{F}_\delta^{\scriptscriptstyle 1/2}$.

\medskip

\textbf{1.~}We construct the sought Feller generator by arguing essentially as in the proof of Theorem \ref{thm1_feller} (so, in particular, we do not use $L^2$ theory of $-\Delta + b \cdot \nabla$, $b \in \mathbf{F}_\delta^{\scriptscriptstyle 1/2}$). First, we prove an analogue of Lemma \ref{lem_GPQ} for weakly form-bounded drifts.
Namely, for given $p \in ]1,\infty[$, $1 \leq r<p<q<\infty$ and $\mu>0$, define operators
\begin{align*}
G_p(r)&:=b^{\frac{1}{p}} \cdot \nabla (\mu -\Delta )^{-\frac{1}{2}-\frac{1}{2r}}, \\
Q_p(q) \upharpoonright {\mathcal E}&:=(\mu -\Delta )^{-\frac{1}{2}+\frac{1}{2q}}|b|^{1-\frac{1}{p}}, \\
T_p \upharpoonright  {\mathcal E}&:=b^{\frac{1}{p}}\cdot \nabla(\mu - \Delta)^{-1}|b|^{1-\frac{1}{p}}.
\end{align*}
where, recall, $\mathcal E:=\bigcup_{\varepsilon>0}e^{-\varepsilon|b|}L^p$ is a dense subspace of $L^p$.
Notice that the power $\frac{2}{p}$ in the definition of operators $G_p(r)$, $T_p$, $Q_p(q)$ in Lemma \ref{lem_GPQ} is now replaced with $\frac{1}{p}$.
Set $$
m_d := 
\pi^{\frac{1}{2}} (2e)^{-\frac{1}{2}} d^\frac{d}{2} (d-1)^{\frac{1-d}{2}}
, \qquad c_p:= pp'/4. 
$$

\begin{lemma}
\label{lem_GPQ_weak} Let $b \in \mathbf{F}_\delta^{\scriptscriptstyle 1/2}$.
For every $p \in ]1,\infty[$, the following is true for all $\mu \geq \kappa_d\lambda_\delta$,

\smallskip

{\rm (\textit{i})} $T_p \upharpoonright  {\mathcal E}$  admits extension by continuity to $L^p$, denoted by $T_p$. One has $$\|T\|_{p \rightarrow p} \leq m_d c_p\delta.$$
In particular, if $\delta$ satisfies $m_d\delta<1$, then for every $p \in I_\delta:=]\frac{2}{1+\sqrt{1-m_d\delta}},\frac{2}{1-\sqrt{1-m_d\delta}}[$ one has $\|T_p\|_{p \rightarrow p}<1.$

\smallskip

{\rm (\textit{ii})} $Q_p(q) \upharpoonright {\mathcal E}$ admits extension by continuity to $L^p$, denoted by $Q_p(q)$.

\smallskip

{\rm (\textit{iii})} $G_p(r)$ is bounded on $L^p$.

\end{lemma}

Lemma \ref{lem_GPQ_weak} was proved in \cite{Ki_a_new_approach}. Let us demonstrate the proof of (\textit{i}) to make it easier to compare Lemma \ref{lem_GPQ_weak} with Lemma \ref{lem_GPQ} (proved in Appendix \ref{app_lem_GPQ}). Define in $L^2$ operator $A = (\mu-\Delta)^{\frac{1}{2}}$, $D(A)= W^{1,2}.$ This is a \textit{symmetric} Markov generator. Therefore, we have for $p \in ]1,\infty[$:
\begin{equation*}
0 \leqslant u \in D(A_p) \quad \Rightarrow \quad  u^\frac{p}{2} \in D(A^\frac{1}{2}) 
\end{equation*}
and the following inequality (sometimes called the Stroock-Varopoulos inequality) is valid:
\begin{equation}
\label{markov_ineq}
c_p^{-1} \|A^\frac{1}{2} u^\frac{p}{2} \|_2^2 \leqslant \langle A_p u, u^{p-1} \rangle, \quad c_p:=\frac{pp'}{4}, \quad p'=\frac{p}{p-1}
\end{equation}
(see \cite[Theorem 2.1]{LS}, see also \cite{KiS_RIMS} for a useful vector-valued analogue of these inequalities).
Here $A_p$ is the generator of strongly continuous semigroup $e^{-tA_p}:=[e^{-tA} \upharpoonright L^2 \cap L^p]^{\clos}_{L^p \rightarrow L^p}$, cf.\,discussion in the beginning of the previous section.
Now, 
let $u$ be the solution of equation\footnote{We can also carry out the proof of Lemma \ref{lem_GPQ_weak} for bounded $b_n$ as e.g.\,defined below, and then pass to the limit in $n$.} $A_p u = |b|^{1/p'}|f|$, $f \in \mathcal{E}.$
The condition $b \in \mathbf{F}_\delta^{\scriptscriptstyle 1/2}$ yields, provided  $\mu \geq \lambda_\delta$,
$$
\| |b|^\frac{1}{2} u^\frac{p}{2} \|_2^2 \leq \delta \|A^\frac{1}{2} u^\frac{p}{2} \|_2^2.
$$
Hence, by \eqref{markov_ineq},
$$
(c_p \delta)^{-1} \| |b|^\frac{1}{2} u^\frac{p}{2} \|_2^2 \leqslant \langle A_p u, u^{p-1} \rangle.
$$
Now, noting that $\| |b|^\frac{1}{2} u^\frac{p}{2} \|_2^2=\| |b|^\frac{1}{p} u\|_p^p$ and using $A_p u = |b|^\frac{1}{p'}|f|$, we obtain
\[
\| |b|^\frac{1}{p} u\|_p^p \leqslant c_p\delta \|f\|_p \| |b|^\frac{1}{p} u\|_p^{p-1}.
\]
Thus, $\| |b|^\frac{1}{p} u\|_p \leqslant c_p\delta \|f\|_p$, so we arrive at the estimate
\begin{equation}
\label{b__}
\||b|^{\frac{1}{p}}A^{-1}|b|^\frac{1}{p'}|f|\|_p \leq c_p\delta \|f\|_p. 
\end{equation}
To end the proof of (\textit{i}), it remains to apply  in the definition of $T_p$ the 
pointwise inequality (this is where the constant $m_d$ comes from)
\begin{equation}
\label{pointwise_ineq}
|\nabla_x (\mu-\Delta)^{-1}(x,y)| \leqslant m_d  (\kappa_d^{-1}\mu-\Delta)^{-\frac{1}{2}}(x,y), \quad x,y \in \mathbb R^d, \; x \neq y, 
\end{equation}
where $$\kappa_d:=\frac{d}{d-1},$$ and then apply \eqref{b__} to the result.

\begin{remark}
Similar estimates, without the gradient, were used earlier in \cite{BS,LS} in the study of Schr\"{o}dinger operators with form-bounded potentials.
\end{remark}

\begin{remark}
\label{m_d_rem}
By applying \eqref{pointwise_ineq} in the definition of $T_p$ we kill the gradient from the gradient term $b \cdot \nabla$. This allows us to apply to what is left the inequalities for symmetric Markov generators. By the way, this is why the interval $I_\delta$ of admissible $p$ in Theorem \ref{thm1_feller_weak} is symmetric, despite the fact that the operator $-\Delta + b \cdot \nabla$ is non-symmetric.
In the proof of Lemma \ref{lem_GPQ}, when dealing with the condition $b \in \mathbf{F}_\delta$, we  control the gradient in a more efficient way, which allows to impose less restrictive assumptions on $\delta$ than in Lemma \ref{lem_GPQ_weak}. 
It is not clear at the moment how to prove Lemma \ref{lem_GPQ_weak} without either resorting to \eqref{pointwise_ineq} or restricting the class $\mathbf{F}_\delta^{\scriptscriptstyle 1/2}$ to Morrey class $M_{1+\varepsilon}$, $\varepsilon>0$ (cf.\,the proof of Lemma \ref{lem_GPQ_Morrey} below).
\end{remark}

The interval $I_\delta$ expands to $]1,\infty[$ as $\delta \downarrow 0$. In particular, if $\delta$ is sufficiently small, $I_\delta$ contains $p>d-1$, which is what will be needed to construct the resolvent of a Feller generator in terms of $Q_p(q)$, $T_p$, $G_p(r)$ and some ``free'' Bessel potentials using the Sobolev embedding theorem. This is what is done in Theorem \ref{thm1_feller_weak} below.

\medskip

Set
\begin{equation}
\label{b_n1_weak}
b_n:= c_n \eta_{\varepsilon_n} \ast (\mathbf{1}_n b), 
\end{equation}
where $\mathbf{1}_n$ is the indicator of $\{x \mid |x| \leq n, |b(x)| \leq n\}$, $\eta_{\varepsilon_n}$ is the Friedrichs mollifier, and we choose $\varepsilon_n \downarrow 0$ (sufficiently rapidly) so that, for appropriate $c_n \uparrow 1$ (sufficiently slow), one has
\begin{equation}
\label{b_n2_weak}
b_n \rightarrow b \quad \text{ in } L^1_{\loc}(\mathbb R^d,\mathbb R^d)
\end{equation}
and
\begin{equation}
\label{b_n3_weak}
b_n \in \mathbf{F}^{\scriptscriptstyle 1/2}_{\delta} \quad \text{ with some $\lambda_\delta$ independent of $n=1,2,\dots$}  
\end{equation}
see Appendix \ref{approx_app}.1. The following theorem was proved in \cite{Ki}.

\begin{theorem}
\label{thm1_feller_weak}
Let $b \in \mathbf{F}^{\scriptscriptstyle \frac{1}{2}}_\delta$, $m_d\delta<1$. 
 The following is true for all $\mu \geq \kappa_d\lambda_\delta$.
\medskip

{\rm (\textit{i})} For every $p \in I_\delta=]\frac{2}{1+\sqrt{1-m_d\delta}},\frac{2}{1-\sqrt{1-m_d\delta}}[$, the function
\begin{equation}
\label{form_neu2}
u=(\mu - \Delta)^{-1}f - (\mu - \Delta)^{-\frac{1}{2}-\frac{1}{2q}} Q_{p}(q) (1 + T_p)^{-1} G_{p}(r) (\mu - \Delta)^{-\frac{1}{2}+\frac{1}{2r}}f, \quad f \in L^p
\end{equation}
is a weak solution to the elliptic equation 
\begin{equation}
\label{eq_e2}
(\lambda-\Delta + b \cdot \nabla)u=f, 
\end{equation}
i.e.\,$b \cdot \nabla u \in L^1_{\loc}$
$$
\mu \langle u,\psi\rangle + \langle \nabla u,\nabla \psi \rangle + \langle b \cdot \nabla u,\psi\rangle=\langle f,\psi\rangle \quad \text{ for all }\psi \in \mathcal S.
$$
Moreover, if $f \in L^p \cap L^2$, then $u$ is the unique in $\mathcal W^{\frac{3}{2},2}$ weak solution to \eqref{eq_e2}.

\smallskip

{\rm (\textit{ii})}  It follows from \eqref{form_neu2} that
$$u \in \mathcal W^{1+\frac{1}{q},p}, \quad q>p.$$
In particular,  if $m_d\delta<\frac{4(d-2)}{(d-1)^2}$, then in the interval $p \in I_\delta$ we can select $p>d-1$, and then select $q$ sufficiently close to $p$, so that by the Sobolev embedding theorem $u$ is H\"{o}lder continuous.

\smallskip

{\rm (\textit{iii})} The operator-valued function in \eqref{form_neu2} $$\Theta_p(\mu,b):=(\mu - \Delta)^{-1}- (\mu - \Delta)^{-\frac{1}{2}-\frac{1}{2q}} Q_{p}(q) (1 + T_p)^{-1} G_{p}(r) (\mu - \Delta)^{-\frac{1}{2}+\frac{1}{2r}}$$ on $\{\mu >\mu_0\}$ takes values in $\mathcal B(\mathcal W^{-1+\frac{1}{r},p},\mathcal W^{1+\frac{1}{q},p})$.

\smallskip

{\rm (\textit{iv})} Let $\delta$ satisfy $m_d\delta<\frac{4(d-2)}{(d-1)^2}$.
Fix $p \in I_\delta$ such that $p>d-1$. Then
\begin{equation*}
(\mu+\Lambda_{C_\infty}(b))^{-1}:=\bigl(\Theta_p(\mu,b) \upharpoonright L^p \cap C_\infty \bigr)_{C_\infty \rightarrow C_\infty}^{\clos}, \quad \mu \geq \kappa_d\lambda,
\end{equation*}
determines the resolvent of a Feller generator on $C_\infty$. This semigroup satisfies
\begin{equation*}
e^{-t\Lambda_{C_\infty}(b)}=\text{\small $s\text{-}C_\infty\text{-}$}\lim_n e^{-t\Lambda_{C_\infty}(b_n)} \quad \text{locally uniformly in $t \geq 0$},
\end{equation*} 
where $b_n$ are defined by \eqref{b_n1_weak}, 
and operators
$\Lambda_{C_\infty}(b_n):=-\Delta + b_n \cdot \nabla$ with domain $D(\Lambda_{C_\infty}(b_n)):=(1-\Delta)^{-1}C_\infty$ are Feller generators.

\smallskip

{\rm (\textit{v})} Feller semigroup $e^{-t\Lambda_{C_\infty}(b)}$ is conservative, i.e.\,its integral kernel $e^{-t\Lambda_{C_\infty}}(x,\cdot)$ satisfies 
\begin{equation*}
\langle e^{-t\Lambda_{C_\infty}(b)}(x,\cdot) \rangle=1 \quad \text{ for all } x \in \mathbb R^d,t>0.
\end{equation*}
\end{theorem}

For the proof, except for the part that concerns the weak solution to the elliptic equation, one can repeat the proof of Theorem \ref{thm1_feller} using Lemma \ref{lem_GPQ_weak} instead of Lemma \ref{lem_GPQ}.

\begin{remark}
The fact that $u$ is a weak solution was proved in \cite{Ki}. Moreover, since $f \in L^2 \cap L^p$, we have
$$
u=\Theta_2(\mu,b)f=\Phi_2(\mu,b)f,
$$
so $u \in \mathcal W^{\frac{3}{2},2}$.
(We could also prove that if in Theorem \ref{thm1_feller_weak} $p=2$, then one can take $q=r=2$.)
 The proof of the uniqueness of $u$ in $\mathcal W^{\frac{3}{2},2}$ goes as follows. 
Let $v \in \mathcal W^{\frac{3}{2}}$ be some weak solution of \eqref{eq_e2}. Then, selecting $\varphi=(\mu-\Delta)^{-\frac{1}{2}}\eta$, $\eta \in C_c^\infty$, we have
$$
\langle (\mu-\Delta)^{\frac{3}{4}}v, (\mu-\Delta)^{\frac{3}{4}}\eta\rangle+ \langle S (\mu-\Delta)^{\frac{3}{4}} u,H^* (\mu-\Delta)^{\frac{3}{4}}\eta\rangle=\langle f,(\mu-\Delta)^{-\frac{1}{2}}\eta\rangle,
$$
where $H= (\mu-\Delta)^{-\frac{1}{4}}|b|^{\frac{1}{2}}$, $S=b^{\frac{1}{2}} \cdot \nabla(\mu-\Delta)^{-\frac{3}{4}}$ are bounded on $L^2$ and $\|S\|_{2 \rightarrow 2}$, $\|H^\ast\|_{2 \rightarrow 2} \leq \sqrt{\delta}$ (for all $\mu \geq \lambda_\delta$ and thus for all $\mu \geq \kappa_d\lambda$), see \eqref{HS} and the discussion after that formula. Thus, the quadratic form
$$
\tau[v,\eta]:=\langle (\mu-\Delta)^{\frac{3}{4}}v, (\mu-\Delta)^{\frac{3}{4}}\eta\rangle+ \langle S (\mu-\Delta)^{\frac{3}{4}} v,H^* (\mu-\Delta)^{\frac{3}{4}}\eta\rangle
$$
is bounded and coercive on $\mathcal W^{\frac{3}{2},2}$ (endowed with the norm $\|\eta\|_{\mathcal W^{\frac{3}{2},2}}:=\|(\mu-\Delta)^{\frac{3}{2}}\eta\|_2$):
$$
|\tau[v,\eta]| \leq (1+\delta)\|v\|_{\mathcal W^{\frac{3}{2},2}}\|\eta\|_{\mathcal W^{\frac{3}{2},2}}
$$
and
$$
|\tau[v,\eta]| \geq (1-\delta)\|v\|_{\mathcal W^{\frac{3}{2},2}}\|\eta\|_{\mathcal W^{\frac{3}{2},2}},
$$
(the assumptions of Theorem \ref{thm1_feller_weak} imply, of course, that $\delta<1$).
 Extending $\tau$ to $\mathcal W^{\frac{3}{2},2}$, and returning to the discussion of the weak solution $u \in \mathcal W^{\frac{3}{2},2}$ constructed in \eqref{form_neu2}, we obtain that $\tau[u-v,\eta]=0$ for all $\eta \in \mathcal W^{\frac{3}{2},2}$, but on the other hand, by coercivity, $\tau[u-v,u-v] \geq (1-\delta)\|u-v\|_{\mathcal W^{\frac{3}{2},2}}^2$. Thus, $u=v$.
\end{remark}

\begin{remark} The operator-valued function $\mu \mapsto \Theta_p(\mu,b) \in \mathcal B(L^p)$ in Theorem \ref{thm1_feller_weak} is the resolvent of the generator of a \textit{quasi bounded} semigroup in $L^p$. This was proved in \cite{Ki}. 
\end{remark}

\begin{remark}
\label{measure_rem}
There is a non-trivial difference between the resolvent representations $\Phi_2$ and $\Theta_p$. For instance, $\Theta_p$ is nonlinear in $|b|$ even if $p=2$, however $\Phi_2$ is linear $|b|$. This circumstance was used in \cite{Ki_measure} to extend Theorem \ref{thm1_feller_weak}  to measure-valued drifts of the form
\begin{equation}
\label{b_measure}
b=\mathsf{f}dx+\mathsf{h},
\end{equation}
where $\mathsf{f}$ is a vector field in $\mathbf{F}_\delta^{\scriptscriptstyle 1/2}$ and $\mathsf{h}$ is a $\mathbb R^d$-valued measure  in the Kato class $\bar{\mathbf K}^{d+1}_\delta$, i.e.\,
\begin{equation}
\label{kato_measure}
\sup_{x \in \mathbb R^d} \int_{\mathbb R^d}(\lambda - \Delta)^{-\frac{1}{2}}(x,y) |\mathsf{h}|_1(dy) \leq \delta
\end{equation}
for some $\lambda=\lambda_\delta$. Here $|\mathsf{h}|_1$ denotes the sum of total variations of the components  of $\mathsf{h}$. This class contains, of course, the Kato class of vector fields $\mathbf{K}^{d+1}_\delta$ defined in the examples above. We could also define the class of weakly form-bounded measure-valued drifts. Indeed, Definition \ref{wfb} can be stated as
$$
\int_{\mathbb R^d} |b(x)|_1\bigl[(\lambda-\Delta)^{-\frac{1}{4}}f(x)\big]^2 dx \leq \delta\|f\|_2^2, \quad f \in \mathcal S,
$$ 
where $|b|_1$ is the sum of absolute values of the components of $b$. Replacing absolute values by total variations, we arrive at a more general condition
\begin{equation}
\label{b_measure2}
\int_{\mathbb R^d} |b(dx)|_1\bigl[(\lambda-\Delta)^{-\frac{1}{4}}f(x)\big]^2 \leq \delta\|f\|_2^2.
\end{equation}
 An example of a measure-valued $b$ satisfying \eqref{b_measure2} is \eqref{b_measure}. In \cite{Ki_measure}, the additional constraint that $b$ must be of the form \eqref{b_measure} comes from the construction of a regularization of $b$ that preserves weak form-bound $\delta$. Although from purely analytic point of view \eqref{b_measure2} is an $L^1$ condition on $|b|$, from the operator-theoretic point of view \eqref{b_measure2} is still an $L^2$ condition, i.e.\,is an $(L^2, L^2)$ operator norm inequality. Thus, \eqref{b_measure} may be viewed as a way that \eqref{b_measure2} filters out the singular measure component $\mathsf{h}$, because the latter satisfies, in the dual formulation of the Kato class, an $(L^1,L^1)$  operator norm inequality (cf.\,\eqref{kato_dual}).

Let us add that, at the level of SDEs, Bass-Chen \cite{BC} and Kim-Song \cite{KSo} considered Kato class measure-valued drifts.
\end{remark}

\textbf{2.~}We now state our result on weak well-posedness of SDEs with weakly form-bounded drifts.

\begin{theorem}
\label{thm1_weak} Let $b \in \mathbf{F}_\delta$ with $\delta<\frac{4(d-2)}{(d-1)^2}$.
Let $e^{-t\Lambda_{C_\infty}(b)}$ be the Feller semigroup constructed in Theorem \ref{thm1_feller_weak}. Fix $T>0$.
The following is true:

\smallskip

{\rm (\textit{i})} There exist probability measures $\{\mathbb P_x\}_{x \in \mathbb R^d}$ on the canonical space $(C[0,T],\mathcal B_t)$ such that 
\begin{equation*}
\mathbb E_{\mathbb P_x}[f(X_t)]=(e^{-t\Lambda_{C_\infty}(b)} f)(x), \quad f \in C_\infty, \quad x \in \mathbb R^d.
\end{equation*}
For every $x \in \mathbb R^d$ the measure $\mathbb P_x$ is a weak solution to SDE 
\begin{equation}
\label{sde1__}
X_t=x-\int_0^t b(X_r)dr + \sqrt{2}W_t, \quad 0 \leq t \leq T.
\end{equation}

{\rm (\textit{ii})} 
If $\{\mathbb Q_x\}_{x \in \mathbb R^d}$ is another weak solution to \eqref{sde1_} such that
$$
\mathbb Q_x=w{\mbox-}\lim_n \mathbb P_x(\tilde{b}_n) \quad \text{for every $x \in \mathbb R^d$},
$$
for some $\{\tilde{b}_n\} \subset \mathbf{F}_{\delta_1}$ with $\delta<\frac{4(d-2)}{(d-1)^2}$ and $\lambda_{\delta}$ independent of $n$, then $\{\mathbb Q_x\}_{x \in \mathbb R^d}=\{\mathbb P_x\}_{x \in \mathbb R^d}.$  
\end{theorem}

This result was proved in \cite{KiS_brownian}. In fact, in the series of papers on well-posedness of SDEs with form-bounded and form-bounded-type drifts that are discussed in this work \cite{KiS_brownian} appeared first.
In turn, \cite{KiS_brownian} was born out of the attempts to obtain a more detailed description of the corresponding Feller semigroup (i.e.\,the one constructed earlier in \cite{Ki}).

\begin{remark}
\label{a_rem}
There are two reasons why one might want to assume form-boundedness of  $b=b(x)$ and not its weak form-boundedness: the possibility to include discontinuous diffusion coefficients as in Section \ref{a_sect} and less restrictive assumptions on $\delta$. (One can compare, using $\mathbf{F}_\delta \subset \mathbf{F}_{\sqrt{\delta}}^{\scriptscriptstyle 1/2}$,
 the assumptions on $\delta$ in Theorem \ref{thm1} and in Theorem \ref{thm1_weak})
\end{remark}

The pointwise estimate \eqref{pointwise_ineq} is also valid for the resolvent of $-\nabla \cdot a \cdot \nabla$ provided that the uniformly elliptic matrix $a$ is H\"{o}lder continuous. If we were to extend Theorem \ref{thm1_weak} to non-constant diffusion coefficients in the spirit of Section \ref{a_sect}, then we could require that $a$ has H\"{o}lder continuous entries whose derivatives are weakly form-bounded.

\medskip

For a form-bounded drift $b \in \mathbf{F}_\delta$, we had two types of gradient bounds on solution $u$ to $(\mu-\Delta + b \cdot \nabla)u=f$, $f \in C_c^\infty$ (let us assume here, for simplicity, that $b$ is bounded and smooth, so we discuss gradient bounds with constants that do not depend on boundedness or smoothness of $b$ but depend only on $\delta$ and $\lambda_\delta$). That is, we had 
\begin{equation}
\label{grad_}
\|(\mu-\Delta)^{\frac{1}{2}+\frac{1}{q}}u\|^p_p \leq K_1\|f\|^p_p, \quad q>p,
\end{equation}
(Theorem \ref{thm1_feller}(\textit{ii}))
and
\begin{equation}
\label{grad_KS_}
\|\nabla |\nabla u|^{\frac{p}{2}}\|_2^2 \leq K_2\|f\|_p^p,
\end{equation}
proved in \cite{KS} using test function $\varphi=-\nabla \cdot (\nabla u|\nabla u|^{p-2})$. In both estimates $p \in [2,\frac{2}{\sqrt{\delta}}[$. These estimates were discussed in Remark \ref{grad_rem}.
 
Theorem \ref{thm1_weak}(\textit{ii}) provides an analogue of \eqref{grad_} for weakly form-bouned $b \in \mathbf{F}_\delta^{\scriptscriptstyle 1/2}$:
\begin{equation}
\label{grad__}
\|(\mu-\Delta)^{\frac{1}{2}+\frac{1}{2q}}u\|^p_p \leq K_1\|f\|^p_p, \quad q>p,
\end{equation}
for $p \in I_\delta$. 
Does there exist an analogue of \eqref{grad_KS_} for weakly form-bounded $b \in \mathbf{F}_\delta^{\scriptscriptstyle 1/2}$? The answer is  ``yes''. Let us note first that, by the solution representation \eqref{form_neu2},
\begin{equation}
\label{b_bd}
(1-c_pm_d\delta)\|||b|_1^{-\frac{1}{p'}}(\mu-\Delta)u\|_p \leq \||b|_1^{-\frac{1}{p'}}f\|_p, \quad |b|_1:=|b|+1
\end{equation}
(cf.\,Remark \ref{second_deriv_rem}), where  $c_p m_d\delta<1$ since $p \in I_\delta$.
Without loss of generality, $p$ is rational with odd denominator, so that we can raise functions taking negative values to power $p$. Also, since all our assumptions on $\delta$ are strict inequalities, we may assume, without loss of generality, that $\|T_p(\mu,|b|_1)\|_{p \rightarrow p} \leq m_dc_p\delta$ for $\mu$ sufficiently large. (Here $T_p(\mu,|b|_1)=|b|_1^{1/p}(\mu-\Delta)^{-\frac{1}{2}}|b|_1^{1/p'}$).
Now, we multiply equation $(\mu - \Delta + b \cdot \nabla)u=f$ by test function $$\varphi:=[(\mu-\Delta)^{\frac{1}{2}}u]^{p-1}$$ and integrate:
$$
\langle (\mu-\Delta)u, [(\mu-\Delta)^{\frac{1}{2}}u]^{p-1} \rangle + \langle b \cdot \nabla, [(\mu-\Delta)^{\frac{1}{2}}u]^{p-1}\rangle = \langle f,[(\mu-\Delta)^{\frac{1}{2}}u]^{p-1}\rangle.
$$
We treat each term separately:

1.~We have
\begin{align*}
& \langle b \cdot \nabla u, [(\mu-\Delta)^{\frac{1}{2}}u]^{p-1}\rangle = \langle b^{\frac{1}{p}}\cdot \nabla (\mu-\Delta)^{-1}|b|_1^{\frac{1}{p'}}|b|_1^{-\frac{1}{p'}}(\mu-\Delta)u, \bigl[|b|^{\frac{1}{p}}(\mu-\Delta)^{\frac{1}{2}}u\bigr]^{p-1}\rangle \\
& = 
\langle b^{\frac{1}{p}}\cdot \nabla (\mu-\Delta)^{-1}|b|_1^{\frac{1}{p'}}|b|_1^{-\frac{1}{p'}}(\mu-\Delta)u, \bigl[|b|^{\frac{1}{p}}(\mu-\Delta)^{-\frac{1}{2}}|b|_1^{\frac{1}{p'}}|b|_1^{-\frac{1}{p'}}(\mu-\Delta)u\bigr]^{p-1}\rangle \\
& (\text{we use \eqref{pointwise_ineq} and apply H\"{o}lder's inequality}) \\
& \leq m_d  \||b|^{\frac{1}{p}}(\kappa_d^{-1}\mu-\Delta)^{-\frac{1}{2}}|b|_1^{\frac{1}{p'}}\|_{p \rightarrow p} \||b|^{\frac{1}{p}}(\mu-\Delta)^{-\frac{1}{2}}|b|_1^{\frac{1}{p'}}\|_{p \rightarrow p}^{p-1} \||b|_1^{-\frac{1}{p'}}(\mu-\Delta)u\|_p^p \\
& \leq m_d \|T_p(\kappa_d^{-1}\mu,|b|_1)\|_{p \rightarrow p}  \|T_p(\mu,|b|_1)\|^{p-1}_{p \rightarrow p} \||b|_1^{-\frac{1}{p'}}(\mu-\Delta)u\|_p^p\\
&\leq m_d c_p^p\delta^p \||b|_1^{-\frac{1}{p'}}(\mu-\Delta)u\|_p^p \\
& (\text{we apply \eqref{b_bd}}) \\
& \leq m_d c_p^p\delta^p (1-c_pm_d\delta)^{-p}\|(|b|+1)^{-\frac{1}{p'}}f\|_p^p \leq m_d c_p^p\delta^p (1-c_p m_d\delta)^{-p}\|f\|_p^p.
\end{align*}

2.~Next,
\begin{align*}
 \langle (\mu-\Delta)u,[(\mu-\Delta)^{\frac{1}{2}}u]^{p-1}\rangle & = \langle (\mu-\Delta)^{\frac{1}{2}}\,(\mu-\Delta)^{\frac{1}{2}}u,[(\mu-\Delta)^{\frac{1}{2}}u]^{p-1}\rangle \\
& (\text{$(\mu-\Delta)^{\frac{1}{2}}$ is a symmetric Markov generator, so we apply \eqref{markov_ineq}}) \\
& \geq \frac{4}{pp'} \|(\mu-\Delta)^{\frac{1}{4}}[(\mu-\Delta)^{\frac{1}{2}}u]^{\frac{p}{2}}\|_2^2 \\
& \geq \frac{2}{pp'} \|(\lambda-\Delta)^{\frac{1}{4}}[(\mu-\Delta)^{\frac{1}{2}}u]^{\frac{p}{2}}\|_2^2 + \frac{C}{pp'} \|(\lambda-\Delta)^{\frac{1}{2}}u\|_p^p,
\end{align*}
where constant $C$ is from the fractional Sobolev embedding theorem.

3.~Also,
$$
\langle f,[(\lambda-\Delta)^{\frac{1}{2}}u]^{p-1}\rangle \leq \|f\|_p \|(\lambda-\Delta)^{\frac{1}{2}}u\|_p^{p-1}.
$$

Combining 1-3, we obtain the following result.
Let $b \in \mathbf{F}_\delta^{\scriptsize 1/2}$ with $m_d\delta<1$. Let $p \in I_\delta$. Then
\begin{equation}
\label{an_est}
\|(\lambda-\Delta)^{\frac{1}{4}}[(\lambda-\Delta)^{\frac{1}{2}}u]^{\frac{p}{2}}\|_2^2 \leq K\|f\|_p^p
\end{equation}
 for $\lambda$ sufficiently large. The estimate \eqref{an_est} is the analogue of \eqref{grad_KS_} for $b \in \mathbf{F}_\delta^{\scriptscriptstyle 1/2}$. Note that it gives the same H\"{o}lder continuity of $u$ as \eqref{grad__}, cf.\,Remark \ref{grad_rem}.

\medskip

The proof of Lemma \ref{lem_GPQ_weak}, and hence the proof of Theorem \ref{thm1_feller_weak}, use inequalities for \textit{symmetric} Markov generators, see Remark \ref{m_d_rem}, and thus depend in an essential manner on the fact that we are working in the elliptic setting. Below we will treat time-inhomogeneous drifts at expense of restricting the class $\mathbf{F}_\delta^{\scriptscriptstyle 1/2}$, but, from some points of view, not by much. At the same time, we will substantially strengthen all aspects of Theorems \ref{feller_parab} and \ref{sde_parab} except for their assumptions on $\delta$.

\medskip

\section{Time-inhomogeneous drifts in Morrey class}

\label{morrey_sect}

In this section we consider  drifts in the parabolic Morrey class $E_q$ with integrability parameter $q>1$ that can be chosen arbitrarily close to $1$. 
Define parabolic cylinder
$$
C_r(t,x):=\{(s,y) \in \mathbb R^{d+1} \mid t \leq s \leq t+r^2, |x-y| \leq r\}
$$
and, given a vector field $b:\mathbb R^{d+1} \rightarrow \mathbb R^d$ with components in $L^q_{\loc}(\mathbb R^{d+1})$,
 $q \in [1,d+2]$, set
\begin{align*}
\|b\|_{E_q}&:=\sup_{r>0, z \in \mathbb R^{d+1}} r\biggl(\frac{1}{|C_r|}\int_{C_r(z)}|b(t,x)|^q dtdx \biggr)^{\frac{1}{q}}\\
&=\sup_{r>0, z \in \mathbb R^{d+1}} r\biggl(\frac{1}{|C_r|}\int_{C_r(z)}|b(-t,x)|^q dtdx \biggr)^{\frac{1}{q}}.
\end{align*}

\begin{definition*}
We say that a vector field $b$ belongs to the parabolic Morrey class $E_{q}$ if
$
\|b\|_{E_q}<\infty.
$
\end{definition*}

One has $$\|b\|_{E_{q}} \leq \|b\|_{E_{q_1}}  \quad \text{ if } q<q_1.$$
So, the smaller is $q$ the larger is Morrey class $E_q$.

If above $b=b(x)$, then one obtains the usual elliptic Morrey class $M_q$ defined earlier.

\begin{examples}
1. The critical Ladyzhenskaya-Prodi-Serrin class
$$
|b| \in L^l(\mathbb R,L^p(\mathbb R^d)), \quad p \geq d, l \geq 2, \quad \frac{d}{p}+\frac{2}{l} \leq 1
$$
is contained in $E_q$.
To prove the inclusion it suffices to consider only the cases $l=2$, $p=\infty$ and $l=\infty$, $p=d$ (see the argument in Appendix \ref{example_sect}(3)).
In the former case the inclusion is trivial, in the latter case the inclusion follows using H\"{o}lder's inequality.

\smallskip

This example is strengthened in the next two examples.

\smallskip

2. Let $|b| \in L^{2,w}(\mathbb R,L^\infty(\mathbb R^d))$. Here and below, $L^{p,w}$ denotes weak Lebesgue spaces (Appendix \ref{example_sect}(4)). Then $b \in E_q$, $1<q<2$. Indeed, by a well known characterization of weak Lebesgue spaces, we have, setting $\tilde{b}(t):=\|b(t,\cdot)\|_{L^\infty(\mathbb R^d)}$,
\begin{align*}
r \biggl(\frac{1}{|C_r|}\int_{C_r}|b|^q dz \biggr)^{\frac{1}{q}}  \leq C r\biggl(\frac{1}{r^{2}}\int_t^{t+r^2}|\tilde{b}|^q ds \biggr)^{\frac{1}{q}} \leq C \|\tilde{b}\|_{L^{2,w}(\mathbb R)}.
\end{align*}
Hence, for example, a vector field $b$ that satisfies
\begin{equation}
\label{b_t}
\|b(t,\cdot)\|_{L^\infty(\mathbb R^d)} \leq \frac{C}{\sqrt{t}}, \quad t>0
\end{equation}
 (and defined to be zero for $t \leq 0$)
is in $E_q$ with $1<q<2$.

This example shows that the parabolic Morrey class  $E_q$ with $1<q<2$ contains vector fields that can have stronger singularities in the time variable than the vector fields in $L^\infty \mathbf{F}_\delta + L^2_{\loc}(\mathbb R)$ considered in the previous section. Namely, if, for simplicity, $b$ depends only on time, that it will be in $E_q$, $1<q<2$ if e.g.\,$|b(t)| \in L^{2,w}(\mathbb R)$ (as \eqref{b_t} above). However, to be in $L^\infty \mathbf{F}_\delta + L^2_{\loc}(\mathbb R)$ it would have to satisfy more restrictive condition $|b(t)| \in L^2_{\loc}(\mathbb R)$.

\smallskip

3. By the well known inclusion of the weak Lebesgue space $L^{d,w}(\mathbb R^d)$ in $M_q$, $$|b| \in L^\infty(\mathbb R,L^{d,w}(\mathbb R^d)) \quad \Rightarrow \quad b \in E_q\;\text{ with }1<q \leq d.$$

\smallskip

4. For every $\varepsilon>0$, one can find $b \in E_q$ such that $|b|$ is not in $ L^{q+\varepsilon}_{\loc}(\mathbb R^{d+1})$. So, selecting $q>1$ close to $1$, one obtains $b \in E_q$ that are not in $L^{1+\epsilon}_{\loc}(\mathbb R^{d+1})$, $\epsilon>0$.

\smallskip

5. In view of the inclusion $\mathbf{F}_\delta$ (with $c_\delta=0$) in $M_2$, we obtain that $\mathbf{F}_\delta \subset E_q$ with $1<q \leq 2$. Furthermore, combining this with example 2, we obtain that $$L^\infty \mathbf{F}_\delta + L^2_{\loc}(\mathbb R) \subset E_q, \quad 1<q \leq 2.$$
\end{examples}

\medskip

We now state our results for drifts $b$ in $E_q$, $1<q<2$. Set for $0<\alpha \leq 2$
\begin{equation}
\label{par_op_def0}
(\lambda-\partial_t-\Delta)^{-\frac{\alpha}{2}}h(t,x):=\int_{t}^\infty\int_{\mathbb R^d}e^{-\lambda(s-t)}\frac{1}{(4\pi (s-t))^{\frac{d}{2}}}\frac{1}{(s-t)^{\frac{2-\alpha}{2}}}e^{-\frac{|x-y|^2}{4(s-t)}}h(s,y)dsdy,
\end{equation}
\begin{equation}
\label{par_op_def}
(\lambda+\partial_t-\Delta)^{-\frac{\alpha}{2}}h(t,x):=\int_{-\infty}^t\int_{\mathbb R^d} e^{-\lambda(t-s)}\frac{1}{(4\pi (t-s))^{\frac{d}{2}}}\frac{1}{(t-s)^{\frac{2-\alpha}{2}}}e^{-\frac{|x-y|^2}{4(t-s)}}h(s,y)dsdy,
\end{equation}
where $\lambda \geq 0$. By a standard result, if $\lambda>0$, then these operators are bounded on $L^p(\mathbb R^{d+1})$, $1 \leq p \leq \infty$, with operator norm $\lambda^{-\frac{\alpha}{2}}$. If $\lambda>0$, then 
$(\lambda \pm \partial_t-\Delta)^{-1}$ is the resolvent of a Markov generator on $L^p(\mathbb R^{d+1})$, $1 \leq p<\infty$, which we will denote by $\lambda \pm \partial_t-\Delta$, respectively. In particular, one has well defined fractional powers $(\lambda \pm \partial_t-\Delta)^\frac{\alpha}{2}$. We refer to \cite{B,G} for the properties of these operators.

\medskip

Define for $p \in ]1,\infty[$
$$
G_p:=b^{\frac{1}{p}} (\lambda+\partial_t-\Delta)^{-\frac{1}{2p}}, \qquad R_p :=b^{\frac{1}{p}}\cdot \nabla (\lambda+\partial_t-\Delta)^{-\frac{1}{2}-\frac{1}{2p}},
$$
$$
Q_p \upharpoonright \mathcal E:=(\lambda+\partial_t-\Delta)^{-\frac{1}{2p'}}|b|^{\frac{1}{p'}}
$$
and
$$
T_p \upharpoonright \mathcal E:=R_p Q_p,
$$
where $\mathcal E:=\cup_{\varepsilon>0}e^{-\varepsilon |b|}L^p(\mathbb R^{d+1})$, a dense subspace of $L^p(\mathbb R^{d+1})$.

The following is an analogue of Lemmas \ref{lem_GPQ} and \ref{lem_GPQ_weak}.

\begin{lemma}
\label{lem_GPQ_Morrey}
Let
$
b=b_\mathfrak{s}+b_{\mathfrak b},
$
satisfy \eqref{H}.
Then, for every $p \in ]1,\infty[$, for all $\lambda>0$, the operators $G_p$, $R_p$, $Q_p$ admit extension to $L^p$ by continuity, and thus so does $T_p$. Moreover,
\begin{equation}
\label{e1}
\|G_p\|_{L^p(\mathbb R^{d+1}) \rightarrow L^p(\mathbb R^{d+1})} , \|R_p\|_{L^p(\mathbb R^{d+1}) \rightarrow L^p(\mathbb R^{d+1})} 
\leq C_{d,p,q}\|b_{\mathfrak{s}}\|^{\frac{1}{p}}_{E_q} + c \lambda^{-\frac{1}{2p}}\|b_{\mathfrak b}\|^{\frac{1}{p}}_{L^\infty(\mathbb R^{d+1})}
\end{equation}
\begin{equation}
\label{e2}
\|Q_p \|_{L^p(\mathbb R^{d+1}) \rightarrow L^p(\mathbb R^{d+1})} \leq C'_{p,q}\|b_{\mathfrak s}\|^{\frac{1}{p'}}_{E_q}  + c'\lambda^{-\frac{1}{2p'}}\|b_{\mathfrak b}\|^{\frac{1}{p'}}_{L^\infty(\mathbb R^{d+1})}.
\end{equation}
\end{lemma}

\begin{remark}
\label{morrey_vs_fbd}
The operators $G_p$ and $R_p$ in Lemma \ref{lem_GPQ_Morrey} correspond operators $G_p(q)$, $R_p(r)$ in Lemmas \ref{lem_GPQ} and \ref{lem_GPQ_weak} with ``$q=r=p$'', which is impossible in Lemmas \ref{lem_GPQ} and \ref{lem_GPQ_weak} if one is dealing with form-bounded and weakly form-bounded drifts (they require $r<p<q$). As a consequence, Lemma \ref{lem_GPQ_Morrey} deals with the operator $T_p$ very easily, since $T_p$ is now a product of two bounded operators in $L^p$. Such decomposition of $T_p$ is impossible for larger classes of form-bounded and weakly form-bounded drifts, see proofs of Lemmas \ref{lem_GPQ} and \ref{lem_GPQ_weak}.
\end{remark}

The first estimate \eqref{e1} follows from the  boundedness of parabolic Riesz transforms (see \cite{G}) and the following result: let $|b| \in E_q$ for some $q>1$ close to $1$,
then, for every $p \in ]1,\infty[$, there exists a constant $c_{p,q}$ such that
\begin{equation}
\label{lem_b}
\||b|^{\frac{1}{p}}(\pm\, \partial_t-\Delta)^{-\frac{1}{2p}}\|_{L^p(\mathbb R^{d+1}) \rightarrow L^p(\mathbb R^{d+1})} \leq c_{p,q}\|b\|^{\frac{1}{p}}_{E_q}.
\end{equation}
In the time homogeneous case $b=b(x)$, the estimate on $\||b|^{\frac{1}{p}}(\lambda-\Delta)^{-\frac{1}{2p}}\|_{L^p(\mathbb R^d) \rightarrow L^p(\mathbb R^d)}$ in terms of the elliptic Morrey norm of $|b|$ is due to \cite[Theorem 7.3]{A}. Similar estimates in the parabolic case were obtained in \cite{Kr3}. Lemma \ref{lem_GPQ_Morrey} is proved in \cite{Ki_Morrey} by adapting the arguments from \cite[proof of Prop.\,4.1]{Kr3}.

The  estimate \eqref{e2} follows from \eqref{lem_b} by duality.

Define
\begin{equation}
\label{b_n}
b_n:=\mathbf{1}_nb, 
\end{equation}
where $\mathbf{1}_n$ is the indicator of the set $\{(t,x) \in \mathbb R^{d+1} \mid |(t,x)| \leq n, |b(t,x)| \leq n\}$. We can  additionally mollify $b_n$ to obtain a $C^\infty$ smooth approximation of $b$ such that the Morrey norm of the approximating vector field does not exceed $(1+\varepsilon)\|b\|_{E_q}$ for any fixed $\varepsilon>0$. However, regularization \eqref{b_n} of $b$ will suffice. In particular, we will be able to apply It\^{o}'s formula to solutions of parabolic equations with drift $b_n$.

Armed with Lemma \ref{lem_GPQ_Morrey}, one obtains the following result \cite{Ki_Morrey}. For every $p \in ]1,\infty[$, there exist constants $c_{d,p,q}$ and $\lambda_{d,p,q}$ such that if
$$\|b_{\mathfrak s}\|_{E_q} < c_{d,p,q},$$ then,
for every $\lambda \geq \lambda_{d,p,q}$, solutions $u_n \in L^p(\mathbb R^{d+1})$ to the approximating parabolic equations
\begin{equation*}
(\lambda + \partial_t - \Delta + b_n \cdot \nabla)u_n=f, \quad f \in L^p(\mathbb R^{d+1})
\end{equation*}
converge in\footnote{Here $
\mathbb W^{\alpha,p}(\mathbb R^{d+1}):=(\lambda+\partial_t-\Delta)^{-\frac{\alpha}{2}}L^p(\mathbb R^{d+1})
$
endowed with the norm $\|h\|_{\mathbb W^{\alpha,p}}:=\|(\lambda+\partial_t-\Delta)^{-\frac{\alpha}{2}}h\|_p$.
} $\mathbb{W}^{1+\frac{1}{p},p}(\mathbb R^{d+1})$ to 
\begin{equation}
\label{u_repr}
u:=(\lambda+\partial_t-\Delta)^{-1}f - (\lambda+\partial_t-\Delta)^{-\frac{1}{2}-\frac{1}{2p}}Q_p (1+T_p)^{-1}R_p (\lambda+\partial_t-\Delta)^{-\frac{1}{2p'}} f.
\end{equation}
Moreover, this $u$ is a unique weak solution to $(\lambda + \partial_t - \Delta + b \cdot \nabla)u=f$, appropriately defined, see \cite{Ki2}.
If above $p>d+1$, then, by \eqref{u_repr} and by the parabolic Sobolev embedding, the convergence is uniform on $\mathbb R^{d+1}$ and $u \in C_\infty(\mathbb R^{d+1})$. 

Let us now construct a Feller evolution family.
Let $\delta_{s=r}$ denote the delta-function in the time variable $s$. Put
$$
(\lambda+\partial_t-\Delta)^{-1}\delta_{s=r}g(t,x):=\mathbf{1}_{t \geq r}e^{-\lambda (t-r)}(4\pi (t-r))^{-\frac{d}{2}}\int_{\mathbb R^d} e^{-\frac{|x-y|^2}{4(t-r)}} g(y)dy,
$$
$$
\nabla (\lambda+\partial_t-\Delta)^{-\frac{1}{2}-\frac{1}{2p'}}\delta_{s=r}g:=\mathbf{1}_{t \geq r}e^{-\lambda (t-r)} (t-r)^{-\frac{1}{2}+\frac{1}{2p'}} (4\pi (t-r))^{-\frac{d}{2}} \int_{\mathbb R^d} \nabla_x e^{-\frac{|x-y|^2}{4(t-r)}}g(y) dy.
$$
Fix $T>0$. For given $n=1,2,\dots$ and $0 \leq r<T$, let $v_n$ denote the classical solution to  Cauchy problem
\begin{equation}
\label{cauchy}
\left\{
\begin{array}{l}
(\lambda + \partial_t-\Delta + b_n(t,x) \cdot \nabla)v_n  =0 \quad (t,x) \in ]r,T] \times \mathbb R^d, \\[2mm]
 v_n(r,\cdot)=g(\cdot) \in C_\infty,
\end{array}
\right.
\end{equation}
where $b_n$'s are defined by \eqref{b_n}. 
By a standard result, for every $n$, the operators $$U_n^{t,r}g:=v_n(t), \quad 0 \leq r \leq t \leq T$$ constitute a Feller evolution family on $C_\infty$.
Recall $D_T=\{0 \leq r \leq t \leq T\}$.

\begin{theorem}
\label{thm2}
Let
$
b=b_{\mathfrak s}+b_{\mathfrak b},
$
where 
\begin{equation}
\label{H}
\text{$|b_\mathfrak{s}| \in E_q$ for some $q>1$ close to $1$, and $|b_{\mathfrak b}| \in L^\infty(\mathbb R^{d+1})$} 
\end{equation}
(indices $\mathfrak s$ and $\mathfrak b$ stand for ``singular'' and ``bounded'', respectively).
Fix $p>d+1$. There exist constants $c_{d,p,q}$ and $\lambda_{d,p,q}$ such that if $\|b_{\mathfrak s}\|_{E_q} < c_{d,p,q}$, then the following are true:

\smallskip

{\rm(\textit{i})} The limit
$$
U^{t,r}:=s\mbox{-}C_\infty(\mathbb R^d)\mbox{-}\lim_n U_n^{t,r}  \quad \text{uniformly in $(r,t) \in D_T$}
$$
exists and determines a Feller evolution family on $C_\infty(\mathbb R^d)$.

\smallskip

{\rm(\textit{ii})} For every initial function $g \in C_\infty(\mathbb R^d) \cap W^{1,p}(\mathbb R^d)$,  $v(t):=U^{t,r}g$, where $(r,t) \in D_T$, has representation
\begin{equation}
\label{v_repr}
v=(\lambda+\partial_t-\Delta)^{-1}\delta_{s=r}g - (\lambda+\partial_t-\Delta)^{-\frac{1}{2}-\frac{1}{2p}}Q_p (1+T_p)^{-1}G_p  S_p g,
\end{equation}
where
$
S_p g:=\nabla (\lambda+\partial_t-\Delta)^{-\frac{1}{2}-\frac{1}{2p'}}\delta_{s=r}g
$
satisfies 
$$\|S_pg\|_{L^p(\mathbb R^{d+1})} \leq C_{p,d} \|\nabla g\|_{L^p(\mathbb R^d)}.$$

\medskip

{\rm(\textit{iii})} As a consequence of \eqref{v_repr} and the parabolic Sobolev embedding, we obtain
$$
\sup_{(r,t) \in D_T, x \in \mathbb R^d}|v(t,x;r)| \leq C\|g\|_{W^{1,p}(\mathbb R^d)}.
$$
\end{theorem}

Theorem \ref{thm2} was proved in \cite{Ki_Morrey}.

\medskip

Define backward Feller evolution family ($0 \leq t \leq r \leq T$)
\begin{equation*}
P^{t,r}(b)=U^{T-t,T-r}(\tilde{b}), \quad \tilde{b}(t,x)=b(T-t,x),
\end{equation*}
where $U^{t,s}$ is the Feller evolution family from Theorem \ref{thm2}. 
As was explained in the previous section, weak well-posedness of SDE
\begin{equation}
\label{sde_morrey}
X_t=x - \int_0^b b(s,X_s)ds + \sqrt{2}W_t, \quad t \geq 0,
\end{equation}
follows from appropriate regularity results for the corresponding inhomogeneous parabolic equation \eqref{veq}.
Indeed, the solution representations \eqref{u_repr} and \eqref{v_repr} can be combined and, furthermore, localized, which yields an analogue of gradient estimates \eqref{rho_grad_est} and thus allows to prove (see \cite{Ki_Morrey}) the following  result:

\begin{theorem}
\label{thm2_morrey}

Under the assumptions of Theorem \ref{thm2}, the following are true:

\smallskip

{\rm (\textit{i})} The backward Feller evolution family $\{P^{t,r}\}_{0 \leq t \leq r \leq T}$ is conservative, i.e.\,the density $P^{t,r}(x,\cdot)$ satisfies
$$ 
\langle P^{t,r}(x,\cdot)\rangle=1 \quad \text{ for all } x \in \mathbb R^d,
$$ 
and determines probability measures $\mathbb P_x$, $x \in \mathbb R^d$ on $(C([0,T],\mathbb R^d),\mathcal B_t)$, such that
$$
\mathbb E_x[f(\omega_r)]=P^{0,r}f(x), \quad 0 \leq r \leq T, \quad f \in C_\infty(\mathbb R^d).
$$

\smallskip

{\rm (\textit{ii})} For every $x \in \mathbb R^d$, the probability measure $\mathbb P_x$ is a weak solution to \eqref{sde_morrey}.

{\rm (\textit{iii})} For every $x \in \mathbb R^d$ and $\mathsf{f}$ satisfying \eqref{H}, given a $p>d+1$ as in Theorem \ref{thm2} (generally speaking, the larger $p$ is the smaller $\|b_{\mathfrak s}\|_{E_q}$ has to be), there exists constant $c$ such that for all $ h \in C_c(\mathbb R^{d+1})$ 
\begin{equation}
\label{kr_est0_}
\mathbb E_{x}\int_0^T |\mathsf{f}(r,\omega_r)h(r,\omega_r)|dr \leq c\|\mathbf{1}_{[0,T]}|\mathsf{f}|^{\frac{1}{p}} h\|_{L^p(\mathbb R^{d+1})}
\end{equation}
(in particular, one can take $\mathsf{f}=b$). On the other hand, any martingale solution $\mathbb P'_x$ to \eqref{sde_morrey} that satisfies for some $p>d+1$ as in Theorem \ref{thm2} the estimate \eqref{kr_est0_} for $\mathsf{h}:=b$ coincides with $\mathbb P_x$.

{\rm ($\textit{iv}$)} For every $x \in \mathbb R^d$, given a $\nu>\frac{d+2}{2}$, there exists a constant $c$ such that for all $ h \in C_c(\mathbb R^{d+1})$ the following Krylov-type bound is true:
\begin{equation}
\label{krylov_type_}
\mathbb E_{x}\int_0^T |h(r,\omega_r)|dr \leq c\|\mathbf{1}_{[0,T]}h\|_{L^\nu(\mathbb R^{d+1})}.
\end{equation}
On the other hand, if additionally $|b| \in L_{\loc}^{\frac{d+2}{2}+\varepsilon}(\mathbb R^{d+1})$ for some $\varepsilon>0$ and $\|b_{\mathfrak s}\|_{E_q}$ is sufficiently small, then any martingale solution $\mathbb P'_x$ to \eqref{sde_morrey} that satisfies \eqref{krylov_type_} for some $\nu>\frac{d+2}{2}$ sufficiently close to $\frac{d+2}{2}$ (depending on how small $\varepsilon$ is)
coincides with $\mathbb P_x$.

\end{theorem}

We compared the uniqueness results of type (\textit{iii}), ($\textit{iv}$) in Remark \ref{unique_comp_rem}.

\medskip

Perhaps, the closest to our Theorems \ref{thm2}-\ref{thm2_morrey} results are contained in recent papers by Krylov 
\cite{Kr_weak}, \cite{Kr_Morrey_sdes}, \cite{Kr_Morrey_parab}, which also allow to deal with discontinuous (VMO) diffusion coefficients, see literature review in the introduction.

\medskip

The estimates \eqref{kr_est0_}, \eqref{krylov_type_} follow from the same kind of solution representations as \eqref{u_repr}, \eqref{v_repr} above, see \cite{Ki_Morrey}. The proofs of the uniqueness results in (\textit{iii}), ($\textit{iv}$) are similar. Let us prove the uniqueness result in ($\textit{iv}$).
Suppose that we have two martingale solutions $\mathbb P^1_x$, $\mathbb P^2_x$ of \eqref{sde_morrey} that satisfy, for $\nu>\frac{d+2}{2}$ close to $\frac{d+2}{2}$,
\begin{align}
\label{kr_est2}
\mathbb E_{x}^i\int_0^T |h(t,\omega_t)|dt  \leq c\|\mathbf{1}_{[0,T]}h\|_\nu, \quad  h \in C_c(\mathbb R^{d+1})
\end{align}
with constant $c$ independent of $h$ ($i=1,2$).
Here and below, $\mathbb E_{x}^1:=\mathbb E_{\mathbb P_x^1}$, $\mathbb E_{x}^2:=\mathbb E_{\mathbb P_x^2}$.
Our goal is to show: for every $f \in C_c(\mathbb R^{d+1})$,
\begin{equation}
\label{id4}
\mathbb E_x^1[\int_0^T f(t,\omega_t)dt]=\mathbb E_x^2[\int_0^T f(t,\omega_t)dt],
\end{equation}
which implies $\mathbb P_x^1=\mathbb P_x^2$.

So, let us prove \eqref{id4}. Let $u_n \in C([0,T],C_\infty(\mathbb R^d))$ be the solution to 
\begin{equation}
\label{eq_F}
(\partial_t + \Delta + b_n \cdot \nabla)u_n=f, \quad u_n(T,\cdot)=0,
\end{equation}
where, recall, $b_n=\mathbf{1}_n b$, and $\mathbf{1}_n$ is the indicator of $\{|b| \leq n\}$.
Set $\tau_R:=\inf\{t \geq 0 \mid |\omega_t| \geq R\}$, $R>0$.
By It\^{o}'s formula 
\begin{align}
\mathbb E_x^i u_n(T \wedge \tau_R,\omega_{T \wedge \tau_R})  & = u_n(0,x)+\mathbb E_x^i \int_0^{T \wedge \tau_R} f(t,\omega_t)dt \notag \\
& +  \mathbb E_x^i \int_0^{T \wedge \tau_R} \big[(b-b_n)\cdot \nabla u_n\big](t,\omega_t)dt \label{e_i}
\end{align}
($i=1,2$).
We have
\begin{align*}
\biggl| \mathbb E_x^i \int_0^{T \wedge \tau_R} \big[(b-b_n)\cdot \nabla u_n\big](t,\omega_t)dt \biggr| & \leq \mathbb E_x^i \int_0^{T \wedge \tau_R} \big[|b|(1-\mathbf{1}_n) |\nabla u_n|\big](t,\omega_t)dt \\
& (\text{we are applying \eqref{kr_est2}}) \\
& \leq c\|\mathbf{1}_{[0,T] \times B_R(0)}|b|(1-\mathbf{1}_n)|\nabla u_n|\|_\nu \\
& \leq
 c\|\mathbf{1}_{[0,T] \times B_R(0)}|b|(1-\mathbf{1}_n)\|_{s'}\|\nabla u_n\|_s, \qquad \frac{1}{s}+\frac{1}{s'}=\frac{1}{q}.
\end{align*}
At this point we note that $\tilde{u}_n(t):=e^{\lambda (T-t)}u_n(t)$ satisfies 
$$
(\lambda + \partial_t + \Delta + b_n \cdot \nabla)u_n=\mathbf{1}_{[0,T]}e^{\lambda (T-t)}f,
$$
so a solution representation of type \eqref{u_repr} (i.e.\,additionally taking into account the terminal value condition), see \cite{Ki_Morrey}, and the parabolic Sobolev embedding theorem, yield
$$
\|\nabla u_n\|_s \leq C\|f\|_{p} \quad \text{ for } s<\frac{d+2}{d+1}p\;\; \text{ close to } \frac{d+2}{d+1}p.
$$
Assuming that the Morrey norm $\|b\|_{E_q}$ is sufficiently small, we can select $p$ sufficiently large to make $s'$ close to $\nu$ and hence close to $\frac{d+2}{2}$. To be more precise, we have by our assumption $|b| \in L^{\frac{d+2}{2}+\varepsilon}$ for some $\varepsilon>0$, so we need $s' \geq \frac{d+2}{2}+\varepsilon$. Now, since $1-\mathbf{1}_n \rightarrow 0$ a.e.\,on $\mathbb R^{d+1}$ as $n \rightarrow \infty$, we have  $\|\mathbf{1}_{[0,T] \times B_R(0)}|b|(1-\mathbf{1}_n)\|_{s'} \rightarrow 0$ as $n \rightarrow \infty$.
Therefore,
$$
\mathbb E_x^i \int_0^{T \wedge \tau_R} \big[(b-b_n)\cdot \nabla u_n\big](t,\omega_t)dt \rightarrow 0 \quad (n \rightarrow \infty).
$$
We are left to note, using the convergence result in \cite[Cor.\,2]{Ki_Morrey} (of the same type as Theorem \ref{thm2}(\textit{i})),  that solutions $u_n$ converge to a function $u \in C([0,T], C_\infty(\mathbb R^d))$. Therefore, we can pass to the limit in \eqref{e_i}, first in $n$ and then in $R \rightarrow \infty$, to obtain
$$
0 =u(0,x)+\mathbb E_x^i \int_0^{T} f(t,\omega_t)dt \quad i=1,2,
$$
which gives \eqref{id4}. \hfill \qed

\bigskip

\section{SDEs driven by $\alpha$-stable process}

\label{stable_sect}

In this section we deal with the SDE
\begin{equation}
\label{eq0}
X_t=x-\int_0^t b(X_s)ds + Z_t-Z_0, \quad t \geq 0, \quad x \in \mathbb R^d,
\end{equation} 
where $Z_t$ be a rotationally symmetric $\alpha$-stable process, $1<\alpha<2$, i.e.\,a L\'{e}vy process with characteristic function 
\begin{displaymath}
\mathbb{E}[\exp (i\varkappa \cdot (Z_{t}-Z_0)]=\exp  (-t|\varkappa |^{\alpha}) \quad \text{ for every } \varkappa \in \mathbb{R}^{d}.
\end{displaymath}
The drift $b:\mathbb R^d \rightarrow \mathbb R^d$ is in general locally unbounded.

Recall that a weak solution to \eqref{eq0} is a process $X_t$ defined on some probability space having a.s.\,right continuous trajectories with left limits, such that $\int_0^t |b(X_s)|ds<\infty$ a.s.\,for every $t>0$, and such that $X_t$ satisfies \eqref{eq0} a.s.\,for a symmetric $\alpha$-stable process $Z_t$. 

A weak solution to \eqref{eq0}, when it exists (e.g.\,if $|b| \in L^\infty$, see \cite{Ko}), is called $\alpha$-stable process with drift $b$. It plays a central role in the study of stochastic processes which, in contrast to the Brownian motion, can have long range interactions.

The operator behind SDE \eqref{eq0} is the non-local operator
$(-\Delta)^{\frac{\alpha}{2}} + b \cdot \nabla$, i.e.\,one expects that the transition density of $X_t$ solves the corresponding parabolic equation for $(-\Delta)^{\frac{\alpha}{2}} + b \cdot \nabla$.

We are interested in the same question as in the previous sections: what are the minimal assumptions on the local singularities of the vector field $b$, not assuming additional structure such as symmetry or existence of the divergence, such that, for an arbitrary initial point, there exists a unique (in appropriate sense) weak solution to \eqref{eq0}? This question has been extensively studied in the literature. By the results of Portenko \cite{P2} and Podolynny-Portenko \cite{PP}, if 
\begin{equation}
\label{b_p}
|b| \in L^p+L^\infty, \quad \text{ for some } p>\frac{d}{\alpha-1}, 
\end{equation}
then there exists a unique in law weak solution to \eqref{eq0}. Although the exponent $\frac{d}{\alpha-1}$ is the best possible on the Lebesgue scale, the class \eqref{b_p} is far from being the maximal admissible: this result has been strengthened in \cite{CKS, CW, KSo} where the authors consider $b$ in the Kato class of vector fields $\mathbf{K}^{d,\alpha-1}_\delta$ (with $\delta$ arbitrarily small), i.e.\,
$|b| \in L^1_{\loc}$ and
$$\big\| \big(\lambda+(-\Delta)^{\frac{\alpha}{2}}\big)^{-\frac{\alpha-1}{\alpha}}|b|\big\|_{\infty} \leq \delta$$
for some $\lambda = \lambda_{\delta} \geq 0$.
This class contains some vector fields $b$ with $|b| \not \in L^{1+\varepsilon}_{\loc}$, $\varepsilon>0$, however, it does not contains the class $|b| \in L^{\frac{d}{\alpha-1}}+L^\infty$.

The Kato class $\mathbf{K}^{d,\alpha-1}_\delta$ with $\delta$ sufficiently small
provides the standard bounds on the heat kernel of the fractional Kolmogorov operator $$\Lambda(b) \supset (-\Delta)^{\frac{\alpha}{2}} + b \cdot \nabla,$$
i.e.
\begin{equation}
\label{bds1}
C^{-1}e^{-t(-\Delta)^{\frac{\alpha}{2}}}(x,y) \leq e^{-t\Lambda(b)}(x,y) \leq Ce^{-t(-\Delta)^{\frac{\alpha}{2}}}(x,y),  
\end{equation}
for all $x,y \in \mathbb R^d$ and $0<t<T$ for a constant $C=C_T$ \cite{BJ}.
It was established in \cite{CKS}, among many other results, that the probability measures $\{\mathbb P_x\}_{x \in \mathbb R^d}$ determined by $e^{-t\Lambda(b)}$ solve the martingale problem for $(-\Delta)^{\frac{\alpha}{2}} + b \cdot \nabla$ with test functions in $C_c^\infty$. The uniqueness in law of the weak solution to SDE \eqref{eq0} with $b \in \mathbf{K}_\delta^{d,\alpha-1}$ (with $\delta$ arbitrarily small) was established in \cite{CW}. Let us also mention that \cite{KSo} considered SDE \eqref{eq0} with Kato class measure-valued drift and established the corresponding heat kernel bounds.

We consider a larger class of \textit{weakly form-bounded vector fields}:

\begin{definition}
A vector field $b:\mathbb R^d \rightarrow \mathbb R^d$ with entries in $L^1_{\loc}$ is said to be weakly form-bounded if there exist $\delta>0$ such that
$$
\| |b|^{\frac{1}{2}}\big(\lambda+(-\Delta)^{\frac{\alpha}{2}}\big)^{-\frac{\alpha-1}{2\alpha}}\|_{2 \rightarrow 2} \leq \sqrt{\delta}
$$
for some $\lambda=\lambda_\delta>0$.
\end{definition}

This will be written as $b \in \mathbf{F}_{\delta}^{\scriptscriptstyle \frac{\alpha-1}{2}}$. This definition extends Definition \ref{wfb} from the previous section corresponding to the case $\alpha=2$.

Our assumptions concerning $\delta$ will involve only strict inequalities, so (using e.g.\,the Spectral theorem) we can re-state our hypothesis on the drift as
$$\| |b|^{\frac{1}{2}}(\lambda-\Delta)^{-\frac{\alpha-1}{4}}\|_{2 \rightarrow 2} \leq \sqrt{\delta}$$
for some $\lambda=\lambda_\delta>0$.

\begin{examples}

1.~Using the fractional Sobolev inequality, it is not difficult to show that
$$|b| \in L^{\frac{d}{\alpha-1}}+L^{\infty} \quad \Rightarrow \quad b \in \mathbf{F}_{\delta}^{\frac{\alpha-1}{2}},$$ where $\delta>0$ can be chosen arbitrarily small. 
More generally, vector fields  with entries in $L^{\frac{d}{\alpha-1},w}$ (the weak $L^{\frac{d}{\alpha-1}}$ class) are weakly form-bounded: 
$$
|b| \in L^{\frac{d}{\alpha-1},\infty} +  L^{\infty} \quad \Rightarrow \quad b \in \mathbf{F}_{\delta}^{\frac{\alpha-1}{2}},
$$
with
$$
\sqrt{\delta}=\Omega_d^{-\frac{\alpha-1}{2d}}\frac{2^{-\frac{\alpha-1}{2}}\Gamma\big(\frac{d-\alpha+1}{4}\big)}{\Gamma\big(\frac{d+\alpha-1}{4}\big)}\|b_1\|_{\frac{d}{\alpha-1},w}^{\frac{1}{2}},
$$
where $\Omega_d$ is the volume of the unit ball $B(0,1) \subset \mathbb R^d$. The proof is  obtained easily using \cite[Corollary 2.9]{KPS}.

2.~In particular, by the fractional Hardy inequality, 
\begin{equation}
\label{hardy_frac}
b(x)=\pm \sqrt{\delta} \kappa_{\alpha,d}|x|^{-\alpha}x,
\end{equation}
where
$$
\kappa_{\alpha,d}:=2^{\frac{\alpha-1}{2}-}\frac{\Gamma(\frac{d+\alpha-1}{4})}{\Gamma(\frac{d-\alpha+1}{4})},
$$
is in $\mathbf{F}_{\delta}^{\frac{\alpha-1}{2}}$ (with $\lambda=0$).

The drift \eqref{hardy_frac} destroys the standard heat kernel bounds \eqref{bds1} (and so it is not in the Kato class). However, for such $b$
sharp heat kernel bounds on $e^{-t\Lambda(b)}(x,y)$ exist but they depend explicitly on $\delta$ via an additional factor $\varphi_t(y)$,
$$
C^{-1}e^{-t(-\Delta)^{\frac{\alpha}{2}}}(x,y)\varphi_t(y) \leq e^{-t\Lambda(b)}(x,y) \leq Ce^{-t(-\Delta)^{\frac{\alpha}{2}}}(x,y)\varphi_t(y), \quad x,y \in \mathbb R^d, \quad y \neq 0.
$$
The factor $\varphi_t(y)$ either explodes at the origin or vanishes, depending on the sign of $\delta$ \cite{KiSSz}, \cite{KiS_RIMS}.

\smallskip

3.~The Kato class  vector fields are weakly form-bounded: $$b \in \mathbf{K}_{\delta}^{d,\alpha-1} \quad \Rightarrow \quad b \in  \mathbf{F}_{\delta}^{\frac{\alpha-1}{2}}.$$
Indeed, if $b \in \mathbf{K}_{\delta}^{d,\alpha-1}$, then by duality $\| |b|\big(\lambda+(-\Delta)^{\frac{\alpha}{2}}\big)^{-\frac{\alpha-1}{\alpha}}\|_{1 \rightarrow 1} \leq \delta$, and so by interpolation $\| |b|^{\frac{1}{2}}\big(\lambda+(-\Delta)^{\frac{\alpha}{2}}\big)^{-\frac{\alpha-1}{\alpha}}|b|^{\frac{1}{2}}\|_{2 \rightarrow 2} \leq \delta,$ hence $b \in \mathbf{F}_{\delta}^{\frac{\alpha-1}{2}}$.

\smallskip

4.~The elliptic Morrey class: $$|b|^{\frac{1}{\alpha-1}} \in M_{1+\varepsilon} \quad \Rightarrow \quad b \in \mathbf{F}_{\delta}^{\frac{\alpha-1}{2}}$$
with $\delta$ depending on the Morrey norm of $|b|^{\frac{1}{\alpha-1}}$ (see definition \eqref{elliptic_morrey2}).
Indeed, by \cite[Theorem 7.3]{A}, one has $\||b|^{\frac{1}{2(\alpha-1)}}(\lambda-\Delta)^{-\frac{1}{4}}\|_{2 \rightarrow 2} \leq \delta^{\frac{1}{2(\alpha-1)}}.$ Then, by the Heinz-Kato inequality (i.e.\,raising operators $|b|^{\frac{1}{2(\alpha-1)}}$ and $(\lambda-\Delta)^{-\frac{1}{4}}$ to power $\alpha-1<1$), we obtain 
$
\| |b|^{\frac{1}{2}}(\lambda-\Delta)^{-\frac{\alpha-1}{4}}\|_{2 \rightarrow 2} \leq \sqrt{\delta}$,  i.e. $b \in \mathbf{F}_\delta^{\scriptscriptstyle \frac{\alpha-1}{2}}.
$
This examples contains examples 1 and 2.

\end{examples}

\begin{remark}
\label{alpha_small_rem}
There is a rich literature on weak and strong well-posedness of SDE \eqref{eq0} (and its generalizations) in the case $0<\alpha \leq 1$, in which case $|b|$ is assumed to be (locally) H\"{o}lder continuous, say, with exponent $\beta$, and satisfy the balance condition  $\alpha+\beta>1$ (sub-critical) or $\alpha+\beta = 1$ (critical). See Zhao \cite{Zh2}, Song-Xie \cite{SX} who considered the case $\alpha+\beta = 1$. Regarding the corresponding heat kernel bounds, we refer to Xie-Zhang \cite{XZ} and Menozzi-Zhang \cite{MeZ} who proved the two-sided bound on the heat kernel of $(-\Delta)^{\frac{\alpha}{2}} + b \cdot \nabla$ in the case $\alpha+\beta>1$. Let us add that in the case $\alpha+\beta=1$ the behaviour of the heat kernel changes drastically, for instance, it can vanish, see \cite{KiMSe} who considered the heat kernel of operator $\Lambda=(-\Delta)^{\frac{\alpha}{2}}-\kappa |x|^{-\alpha}x \cdot \nabla$, $\kappa>0$ (Hardy-type drift) and proved upper bound of the form
$$
0 \leq e^{-t\Lambda}(x,y) \leq Ct^{-\frac{d}{\alpha}}[1 \wedge t^{-\frac{\gamma}{\alpha}}|y|^\gamma], \quad t \in ]0,1],
$$
where the order of vanishing $\gamma \in ]0,\alpha[$ is an explicit function of $\kappa$.
\end{remark}

For a given vector field $b \in \mathbf{F}_{\delta}^{\scriptscriptstyle \frac{\alpha-1}{2}}$, we fix a $C^\infty$ smooth approximation
\begin{equation*}
b_{n}:=c_n\eta_{\varepsilon_{n}} \ast (\mathbf{1}_{n}b), \quad \varepsilon_{n}\downarrow 0, \quad n=1,2,\dots,
\end{equation*}
where $\mathbf{1}_{n}$ is the indicator of  $\{x \in \mathbb{R}^{d} \mid |x| \leq n, |b(x)| \leq n \}$, $\eta_{\varepsilon}$ is the Friedrichs mollifier.
Selecting $\varepsilon_n \downarrow $ sufficiently rapidly and $c_n \uparrow 1$ sufficiently slow, one obtains that $$b_{n}  \in \mathbf{F}_{\delta}^{\scriptscriptstyle \frac{\alpha-1}{2}}, \quad n=1,2,\dots$$ with $\lambda$ independent of $n$.

Fix constant $m_{d,\alpha}$ by the pointwise inequality
\begin{equation}
\label{A}
\big|\nabla_x \big(\mu+(-\Delta)^{\frac{\alpha}{2}}\big)^{-1}(x,y)| \leq m_{d,\alpha}\big(\kappa^{-1}\mu+(-\Delta)^{\frac{\alpha}{2}}\big)^{-\frac{\alpha-1}{\alpha}}(x,y)
\end{equation}
for all $x,y \in \mathbb R^d$, $x \neq y$, and all $\mu>0$, for some $\kappa = \kappa_{d,\alpha}>0$. 
The following result was proved in \cite{KiM_stable} (one can find there an elementary estimate on $m_{d,\alpha}$ from above).

\begin{theorem}
\label{thm_stable}
Let $b \in \mathbf{F}_\delta^{\frac{\alpha-1}{2}}$ with  $\delta<m_{d,\alpha}^{-1}4\bigl[\frac{d-\alpha}{(d-\alpha+1)^2} \wedge \frac{\alpha(d+\alpha)}{(d+2\alpha)^2}\bigr]$. 
The following is true.

{\rm (\textit{i})} The limit 
\begin{equation*}
s{\mbox-}C_{\infty}\mbox{-}\lim_{n}e^{-t\Lambda_{C_{\infty}}(b_{n})} \quad \text{{\rm(}loc.\,uniformly in $t \geq 0${\rm)}},
\end{equation*}
where $$\Lambda_{C_{\infty}}(b_{n}):=(-\Delta)^{\frac{\alpha}{2}}+b_{n} \cdot\nabla \text{ with domain }  D(\Lambda_{C_{\infty}}(b_{n}))=\big(1+(-\Delta)^{\frac{\alpha}{2}}\big)^{-1}C_\infty,$$
exists and determines a Feller semigroup $T^t=:e^{-t\Lambda_{C_{\infty}}(b)}$. Its generator $\Lambda_{C_\infty}$ is an operator realization of the formal operator $(-\Delta)^{\frac{\alpha}{2}} + b \cdot \nabla$ in $C_\infty$.

\smallskip

{\rm (\textit{ii})} There exists $\mu_0 \geq 0$ such that for all $\mu \geq \mu_0$, for every $p \in [2,p_+[$, $p_+=\frac{2}{1 - \sqrt{1-m_{d,\alpha}\delta}}$, and all $1<r<p<q<\infty$,
$$
(\mu + \Lambda_{C_\infty}(b))^{-1} \upharpoonright C_\infty \cap L^p \text{ extends by continuity to $\mathcal B\big(\mathcal W^{-\frac{\alpha-1}{r'},p},\mathcal W^{1+\frac{\alpha-1}{q},p}\big)$}.
$$
In particular, if $p>d-\alpha+1$, then
$$\big(\mu+\Lambda_{C_\infty}(b)\big)^{-1}[C_\infty \cap L^p] \subset C^{0,\gamma}, \quad \gamma < 1-\frac{d-\alpha+1}{p}.$$
Also,
\begin{equation}
\label{L2_stable}
\big(\mu+\Lambda_{C_\infty}(b)\big)^{-1} \upharpoonright C_\infty \cap L^2 \text{ extends by continuity to $\mathcal B(\mathcal W^{-\frac{\alpha-1}{2},2},\mathcal W^{\frac{\alpha+1}{2},2})$}.
\end{equation}

\smallskip

{\rm (\textit{iii})} $e^{-t\Lambda_{C_{\infty}}(b)}$ is conservative, i.e.\,$\int_{\mathbb R^d }e^{-t\Lambda_{C_{\infty}}(b)}(x,y)dy=1$\; $\forall\,x \in \mathbb R^d$.

\smallskip

Let $\{\mathbb{P}_x\}_{x \in \mathbb R^d}$ be the probability measures on the canonical space $(D[0,T],\mathcal B'_t)$ determined by $e^{-t\Lambda_{C_{\infty}}(b)}$, i.e.
$$
\mathbb E_{\mathbb P_x}[f(X_t)]=(e^{-t\Lambda_{C_\infty}(b)} f)(x), \quad f \in C_\infty, \quad x \in \mathbb R^d.
$$

\smallskip

{\rm (\textit{iv})} For every $x \in \mathbb R^d$ and $t>0$, $\mathbb E_{\mathbb P_x}\int_0^t |b(X_s)|ds<\infty$
and there exists a process $Z_t$ with trajectories in $D(\mathbb R_+,\mathbb R^d)$, which is a symmetric $\alpha$-stable process under each $\mathbb P_x$, such that
$$
X_t=x-\int_0^t b(X_s)ds+Z_t-Z_0, \quad t \geq 0.
$$

\smallskip

{\rm (\textit{v})} The Feller property and property \eqref{L2_stable} determine $\{\mathbb P_x\}_{x \in \mathbb R^d}$ uniquely.
That is, suppose that for every $x \in \mathbb R^d$ we are given a weak solution $\mathbb Q_x$ to SDE \eqref{eq0}.
Define for every $f \in C_c^\infty$
$$
R^Q_\mu f(x):=\mathbb E_{\mathbb Q_x} \int_0^\infty e^{-\mu s} f(X_s)ds, \quad X_s \in D(\mathbb R_+,\mathbb R^d), \quad x \in \mathbb R^d, \quad \mu>\lambda_\delta.
$$ 
If $R^Q_\mu C_c^\infty \subset C_b$ and $R^Q_\mu \upharpoonright C_c^\infty$ admits extension by continuity to $\mathcal B(\mathcal W^{-\frac{\alpha-1}{2},2}, L^2)$, then $$\{\mathbb Q_x\}_{x \in \mathbb R^d}=\{\mathbb P_x\}_{x \in \mathbb R^d}.$$ 

\end{theorem}

We stated assertions (\textit{i}), (\textit{ii}) in a form that is somewhat different from Theorems \ref{thm1_feller} or \ref{thm1_feller_weak}, but we could have stated it in the form of these theorems as well. The construction of the Feller semigroup (\textit{i}) goes as in Theorems \ref{thm1_feller} (or rather as in Theorem \ref{thm1_feller_weak} since we use pointwise bound \eqref{A}). Same for the embedding properties (\textit{ii}), i.e.\,we can write an explicit operator-valued function representation for the resolvent as in the cited two theorems. 

The proof of the probability conservation property in (\textit{iii}) uses weighted estimates, similarly to the proof of Theorem \ref{thm1_feller} dealing with the local case $\alpha=2$ (see Section \ref{thm1_feller_proof_sect}, estimate \eqref{j_1_w} there). 
Set $$\eta(x):=(1+|x|^2)^{\nu}, \quad 0<\nu<\frac{\alpha}{2}.$$
Denote $L^p_\eta:= L^p(\mathbb R^d, \eta^2 dx)$, $\|\cdot\|_{p,\eta}^p:=\langle |\cdot|^p \eta ^2\rangle $. 

\begin{proposition}
\label{lem_weights}
\label{weight_lem}
Let $d \geq 3$, $b \in \mathbf{F}_\delta^{\frac{\alpha-1}{2}}$ with $\delta<m_{d,\alpha}^{-1}4\bigl[\frac{d-\alpha}{(d-\alpha+1)^2} \wedge \frac{\alpha(d+\alpha)}{(d+2\alpha)^2}\bigr]$. There exist $0<\nu<\alpha/2$, $p > (d-\alpha+1) \vee (\frac{d}{2\nu}+2)$ and $\mu_0>0$ such that
for every $h \in C_c$, $\mu \geq \mu_0$
\begin{align}
\label{j_1_w_}
\tag{$E_1$}
\|\eta^{-1} (\mu+\Lambda_{C_\infty}(b))^{-1} \eta h\|_\infty & \leq K_1\|h\|_{p,\eta}, \\
\label{j_2_w_}
\tag{$E_2$}
\|\eta^{-1} (\mu+\Lambda_{C_\infty}(b))^{-1}\eta |b_m| h\|_\infty & \leq K_2\| |b_m|^\frac{1}{p}  h\|_{p,\eta},  \\
\label{j_3_w}
\tag{$E_3$}
\|\eta^{-1} |b_m|^\frac{1}{p}(\mu+\Lambda_{C_\infty}(b))^{-1}\eta |b_m| h\|_{p,\eta} & \leq K_3 \| |b_m|^\frac{1}{p}  h\|_{p,\eta}, 
\end{align}
where $K_i>0$, $i=1,2,3$, do not depend on $m=1,2,\dots$ The constant $K_3$ can be chosen arbitrarily small at expense of increasing $\mu_0$.
\end{proposition}

The proof of these weighted estimates in \cite{KiM_stable} is, however, quite different from the proof of Theorem \ref{thm1_feller} where one can control  easily the commutator of the weight and the Laplacian.  \cite{KiM_stable} provides a different, rather interesting approach to the proof of Proposition \ref{weight_lem}, but we will not discuss it here. Note that if $b$ has compact support then we can take $\eta \equiv 1$.

Let us describe the proof of (\textit{iv}) in \cite{KiM_stable}, which uses the approach of \cite{PP,P2}  but in a  weighted $L^p$ space. Set
$$
Z_t:=X_t-X_0-\int_0^t b(X_s)ds, \quad t \geq 0.
$$
Our goal is to prove that under $\mathbb P_x$ the process $Z_t$ is a symmetric $\alpha$-stable process starting at $0$. We use notation introduced in the beginning of the previous section. 
For brevity, write $e^{-t\Lambda(b)} = e^{-t\Lambda_{C_\infty}(b)}$.

\smallskip

1. \textit{Define
\begin{equation*}
w(t,x,\varkappa)=\mathbb E_x \biggl[e^{i \varkappa \cdot \bigl (X_t-\int_0^t b(X_s)ds\bigr) } \biggr], \quad t \geq 0, \quad \varkappa \in \mathbb R^d.
\end{equation*} 
Then $w$ is a bounded solution to integral equation}
\begin{equation}
\label{eq_w}
w(t,x,\varkappa)=\int_{\mathbb R^d} e^{i \varkappa \cdot y}e^{-t\Lambda(b)}(x,y) dy - i\int_0^t \int_{\mathbb R^d} e^{-(t-s)\Lambda(b)}(x,z)(\varkappa \cdot b(z))w(s,z,\varkappa)dzds.
\end{equation}
Indeed, in view of
\begin{equation*}
e^{-i\cdot\varkappa\int_0^t b(X_\tau)d\tau}=1-i\int_{0}^{t}(\varkappa \cdot b(X_s))e^{-i\cdot\varkappa\int_s^t b(X_\tau)d\tau},
\end{equation*}
one has
\begin{align*}
w(t,x,\varkappa)&
=\mathbb E_x \biggl[e^{i \varkappa \cdot X_t}\biggr]-i\int_{0}^{t}\mathbb E_x \biggl[e^{i \varkappa \cdot X_t}(\varkappa \cdot b(X_s))e^{-i\cdot\varkappa\int_s^t b(X_\tau)d\tau}\biggr]ds\\
&=\mathbb E_x \biggl[e^{i \varkappa \cdot X_t}\biggr]-i\int_{0}^{t}\mathbb E_x\biggl[(\varkappa \cdot b(X_s))w(t-s,X_s,\varkappa)\biggr]ds\\
&=\int_{\mathbb R^d} e^{i \varkappa \cdot y}e^{-t\Lambda(b)}(x,y) dy - i\int_0^t \int_{\mathbb R^d} e^{-s\Lambda(b)}(x,z)(\varkappa \cdot b(z))w(t-s,z,\varkappa)dzds.
\end{align*}

\smallskip

2.\,\textit{Set $\tilde{w}(t,x,\varkappa):=e^{i \varkappa \cdot x - t|\varkappa|^\alpha}$. This is another bounded solution to \eqref{eq_w}.}
Indeed, multiplying  the Duhamel formula
$$
e^{-t\Lambda}(x,y)=e^{-t(-\Delta)^{\frac{\alpha}{2}}}(x,y) + \int_0^t \langle e^{-(t-s)\Lambda}(x,\cdot) b(\cdot)\cdot \nabla_\cdot e^{-s(-\Delta)^{\frac{\alpha}{2}}}(\cdot,y)\rangle ds
$$
 (which is proved in \cite[Corollary 1(\textit{iv})]{KiM_stable})
by $e^{i\varkappa \cdot y}$ and then integrating in $y$, we obtain the required.

\smallskip

Next, let us show that a bounded solution to \eqref{eq_w} is unique. We will need

\smallskip

\smallskip

3.\,\textit{For every $\varkappa \in \mathbb R^d$ there exists $T=T(\varkappa)>0$ such that the mapping
$$
(Hv)(t,x):=-i\int_0^t \int_{\mathbb R^d} e^{-(t-s)\Lambda(b)}(x,z)(\varkappa \cdot b(z))v(s,z)dsdz, \quad (t,x) \in [0,T] \times \mathbb R^d,
$$
is a contraction on $L^p(\mathbb R^d, |b|\eta^{-p+2}dx; L^\infty[0,T])$ (i.e.\,functions taking values in $L^\infty[0,T]$) for $p$ as in Proposition \ref{weight_lem}.}

Indeed, we have
\begin{align*}
|Hv(t,x)| & \leq \left| \int_0^t \langle e^{-(t-s)\Lambda(b)}(x,\cdot)(\varkappa \cdot b(\cdot))v(s,\cdot)\rangle ds \right|  \\
& \leq |\varkappa|\int_0^t \langle e^{-(t-s)\Lambda(b)}(x,\cdot)|b(\cdot)|^{\frac{1}{p'}}|b(\cdot)|^{\frac{1}{p}}|v(s,\cdot)|\rangle ds \\
& \leq |\varkappa|\int_0^t \langle e^{-(t-s)\Lambda(b)}(x,\cdot)|b(\cdot)|^{\frac{1}{p'}}|b(\cdot)|^{\frac{1}{p}}\sup_{\tau \in [0,T]}|v(\tau,\cdot)|\rangle ds \label{H_est} \tag{$\ast$}
\end{align*}
Let us note that, for every $x \in \mathbb R^d$,
\begin{align*}
&|b(x)|^{\frac{1}{p}}\eta^{-1}(x)\sup_{t \in [0,T]}\int_0^t \langle e^{-(t-s)\Lambda(b)}(x,\cdot)|b(\cdot)|^{\frac{1}{p'}}\eta(\cdot)\rangle ds \\
& (\text{we are applying the Dominated Convergence Theorem})\\
& |b(x)|^{\frac{1}{p}}\eta^{-1}(x)\sup_{t \in [0,T]}\lim_m\int_0^t \langle e^{-(t-s)\Lambda(b)}(x,\cdot)|b_m(\cdot)|^{\frac{1}{p'}}\eta(\cdot)\rangle ds,
\end{align*}
where, in turn, the last term
\begin{align*}
& |b|^{\frac{1}{p}}\eta^{-1}\sup_{t \in [0,T]}\lim_m\int_0^t e^{-(t-s)\Lambda(b)}|b_m|^{\frac{1}{p'}}\eta ds \\
& \leq |b|^{\frac{1}{p}}\eta^{-1}e^{\mu T}\lim_m(\mu+\Lambda_{C_\infty}(b))^{-1}|b_m|^{\frac{1}{p'}}\eta \in \mathcal B(L^p_{\eta}) \quad \text{by Proposition \ref{weight_lem}\eqref{j_3_w}}.
\end{align*}
Also by Proposition \ref{weight_lem}\eqref{j_3_w},
selecting $\mu$ sufficiently large, and then selecting $T$ sufficiently small, the $L^p_\eta \rightarrow L^p_\eta$ norm of the last operator can be made arbitrarily small. Applying this in \eqref{H_est}, we obtain that 
$H$ is indeed a contraction on $L^p(\mathbb R^d, |b|\eta^{-p+2}dx; L^\infty[0,T])$.

\medskip

We have $L^\infty([0,T] \times \mathbb R^d) \subset L^p(\mathbb R^d, |b|\eta^{-p+2}dx; L^\infty[0,T])$ since (see \cite[Lemma 5.1]{KiM_stable}) $|b|\eta^{-p+2} \in L^1(\mathbb R^d)$. Combining the assertions of Steps 1-3, we obtain that for every $\varkappa \in \mathbb R^d$
$$
w(t,x,\varkappa)=\tilde{w}(t,x,\varkappa) \quad \text{ in } L^p(\mathbb R^d, |b|\eta^{-p+2}dx; L^\infty[0,T]),
$$
and thus
$$
w(t,x,\varkappa)=\tilde{w}(t,x,\varkappa) \quad \text{ for a.e.\,$x \in \mathbb R^d$}
$$
(although $t<T(\varkappa)$, one can get rid of this constraint using the reproduction property of $e^{-t\Lambda(b)}$, so without loss of generality $T \neq T(\varkappa)$). Now, applying the continuity of $\int_0^t e^{-s\Lambda_{C_\infty}(b)}b\cdot w ds$  on $\mathbb R_+ \times \mathbb R^d$
(this is \cite[Corollary 1(\textit{iii})]{KiM_stable}) in the RHS of \eqref{eq_w}, we obtain that for every $\varkappa \in \mathbb R^d$ $w(t,x,\varkappa)$ is continuous in $t$ and $x$, and so $w=\tilde{w}$ everywhere. Thus, for all $t \leq T$, $x \in \mathbb R^d$
$$
\mathbb E_x \biggl[e^{i \varkappa \cdot \bigl (X_t-X_0-\int_0^t b(X_s)ds\bigr) } \biggr] = e^{-\varkappa \cdot x}w(t,x,\varkappa)=e^{-t|\varkappa|^\alpha}.
$$
By a standard result, $Z_t$ is a symmetric $\alpha$-stable process.
The proof of assertion (\textit{iv}) is completed.

\bigskip

\appendix

\section{Proof of Lemma \ref{lem_GPQ}}

\label{app_lem_GPQ}

The following proof was given in \cite{Ki2}.

Step 1. Let us show that 
 $$
\|T_p\|_{p \rightarrow p} \leq c_{\delta,p}, \quad \mu \geq \mu_0,
$$
which will give us assertion (\textit{i}).
We will also prove that operators
$$
G_p=b^{\frac{2}{p}} \cdot \nabla (\mu -\Delta )^{-1}, \quad 
Q_p=(\mu -\Delta )^{-1}|b|^{1-\frac{2}{p}}$$
satisfy
\begin{equation*}
\|G_p\|_{p \rightarrow p} \leq C_1 \mu^{-\frac{1}{2}+\frac{1}{p}}, \quad \|Q_p\|_{p \rightarrow p} \leq C_2 \mu^{-\frac{1}{2}-\frac{1}{p}}
\end{equation*}
The latter will be needed to prove assertions (\textit{ii}) and (\textit{iii}).

We will be using the operator-norm formulation of the form-boundedness condition:
\begin{equation*}
\|b (\lambda-\Delta)^{-\frac{1}{2}}\|_{2 \rightarrow 2} \leq \delta 
\end{equation*}
for some $\lambda=\lambda_\delta$, see \eqref{op_fb}.

\smallskip

 (\textbf{a}) Set $v:=(\mu-\Delta)^{-1}|b|^{1-\frac{2}{p}}f$, $0 \leq f \in L^p$. Then
\begin{align*}
&\|T_p f\|_p^p =\|b^{\frac{2}{p}} \nabla v\|_p^p = \langle |b|^2 |\nabla v|^p\rangle \notag \\
& = \||b|(\lambda-\Delta)^{-\frac{1}{2}} (\lambda-\Delta)^{\frac{1}{2}}|\nabla v|^{\frac{p}{2}}\|_2^2 \notag  \qquad (\lambda=\lambda_\delta) \notag  \\
& \leq \||b|(\lambda-\Delta)^{-\frac{1}{2}}\|_{2 \rightarrow 2}^2 \|(\lambda-\Delta)^{\frac{1}{2}}|\nabla v|^{\frac{p}{2}}\|_2^2 \notag  \\
& = \delta \|(\lambda-\Delta)^{\frac{1}{2}}|\nabla v|^{\frac{p}{2}}\|_2^2= \delta \bigl( \lambda\|\nabla v\|_p^p + \|\nabla |\nabla v|^{\frac{p}{2}}\|_2^2 \bigr). \notag
\end{align*}
It remains to prove the principal inequality
\[
\delta \bigl( \lambda\|\nabla v\|_p^p + \|\nabla |\nabla v|^{\frac{p}{2}}\|_2^2 \bigr) \notag  
\leq c^p_{\delta,p} \|f\|_p^p \tag{$\ast$} \label{ineq1}, 
\]
and conclude that $\|T_p\|_{p \rightarrow p} \leq c_{\delta, p}$. 

First, we prove an a priori variant of \eqref{ineq1}, i.e.\,for $v:=(\mu-\Delta)^{-1}|b|^{1-\frac{2}{p}}f$ with $b=b_n$. Since our assumptions on $\delta$ involve only strict inequalities, we may assume, upon selecting appropriate $\varepsilon_n \downarrow 0$, that $b_n \in \mathbf{F}_\delta$ with the same $\lambda=\lambda_\delta$ for all $n$.

 Set 
$$
w:=\nabla v, \quad I_q:=\sum_{r=1}^d\langle (\nabla_r w)^2 |w|^{p-2} \rangle, \quad
J_q:=\langle (\nabla |w|)^2 |w|^{p-2}\rangle.
$$
We multiply $(\mu-\Delta)v=|b|^{1-\frac{2}{p}}f$ by $\phi:=- \nabla \cdot (w|w|^{p-2})$ and integrate by parts to obtain
\begin{equation}
\label{main_id}
\mu \|w\|_p^p + I_p + (p-2)J_p = \langle |b|^{1-\frac{2}{p}}f, - \nabla \cdot (w|w|^{p-2})\rangle,
\end{equation}
where
\begin{align*}
&\langle |b|^{1-\frac{2}{p}}f, - \nabla \cdot (w|w|^{p-2})\rangle = \langle |b|^{1-\frac{2}{p}}f, (- \Delta v) |w|^{p-2} - (p-2) |w|^{p-3} w \cdot \nabla |w|\rangle \\
& (\text{use the equation } -\Delta v = - \mu v + |b|^{1-\frac{2}{p}}f) \\
&=\langle |b|^{1-\frac{2}{p}}f, \bigl( - \mu v + |b|^{1-\frac{2}{p}}f\bigr) |w|^{p-2}\rangle  - (p-2) \langle |b|^{1-\frac{2}{p}}f, |w|^{p-3} w \cdot \nabla |w|\rangle.
\end{align*}

\begin{remark}
Here we work with the same test function $\phi =- \nabla \cdot (w|w|^{p-2})$ as in \cite{KS}. 
\end{remark}

We have

1) $\langle |b|^{1-\frac{2}{p}}f, (- \mu v)|w|^{p-2} \rangle \leq 0$,

2) $|\langle |b|^{1-\frac{2}{p}}f, |w|^{p-3} w \cdot \nabla |w|\rangle| \leq \alpha J_p+\frac{1}{4\alpha} N_p\;\;( \alpha>0$), where $N_p:=\langle |b|^{1-\frac{2}{p}}f, |b|^{1-\frac{2}{p}}f |w|^{p-2} \rangle$,

\noindent so, the RHS of \eqref{main_id} $\leq  (p-2) \alpha J_p +  \big(1+\frac{p-2}{4\alpha}\big)N_p,
$
where, in turn,
\begin{align*}
N_p& \leq \langle |b|^{2}|w|^p\rangle^{\frac{p-2}{p}} \langle f^{p}\rangle^{\frac{2}{p}} \\
&\leq \frac{p-2}{p} \langle |b|^{2}|w|^p\rangle + \frac{2}{p} \|f\|_p^p \qquad \text{(use $b \in \mathbf{F}_\delta$ $\Leftrightarrow$ $\|b\varphi\|_2^2 \leq \delta \|\nabla \varphi\|_2^2 + \lambda\delta\|\varphi\|_2^2$, $\varphi \in W^{1,2}$)} \\
& \leq \frac{p-2}{p} \bigg(\frac{p^2}{4}\delta J_q + \lambda\delta \|w\|_p^p \bigg) + \frac{2}{p} \|f\|_p^p.
\end{align*}
Thus, applying $I_q \geq J_q$ in the LHS of \eqref{main_id}, we obtain
$$
\bigl(\mu - c_0 \bigr)\|w\|_p^p + \biggl[p-1-(p-2)\left(\alpha + \frac{1}{4\alpha}\frac{p(p-2)}{4}\delta\right) - \frac{p(p-2)}{4}\delta\biggr]\frac{4}{p^2}\|\nabla |\nabla v|^{\frac{p}{2}}\|_2^2 \leq \left(1+\frac{p-2}{4\alpha}\right)\frac{2}{p}\|f\|_p^p,
$$
where $c_0=\frac{p-2}{p}\lambda\delta \big(1+\frac{p-2}{4\alpha}\big)$. It is now clear that one can find a sufficiently large $\mu_0 = \mu_0(d,p,\delta)>0$ so that, for all $\mu > \mu_0$, \eqref{ineq1} (with $b=b_n$) holds  with
\begin{align*}
 c_{\delta,p}^p&=\delta \frac{p^2}{4}\frac{\left(1+\frac{p-2}{4\alpha}\right)\frac{2}{p}}{p-1-(p-2)\left(\alpha + \frac{1}{4\alpha}\frac{p(p-2)}{4}\delta\right) - \frac{p(p-2)}{4}\delta} \qquad \text{$\bigl($we  select $\alpha=\frac{p}{4}\sqrt{\delta}\bigr)$} \\
& = \frac{\frac{p}{2}\delta + \frac{p-2}{2}\sqrt{\delta}}{p-1-(p-1)\frac{p-2}{2}\sqrt{\delta} - \frac{p(p-2)}{4}\delta},
\end{align*}
as claimed.
Finally, we pass to the limit $n \rightarrow \infty$ using Fatou's Lemma. The proof of \eqref{ineq1} is completed.

\begin{remark}
It is seen that $\sqrt{\delta}<\frac{2}{p}\Rightarrow c_{\delta, p}<1$. We also note that the above choice of $\alpha$ is the best possible.
\end{remark}

\smallskip

(\textbf{b}) Set $v=(\mu-\Delta)^{-1}f$, $0 \leq f \in L^p$. Then
\begin{align*}
&\|G_p f\|_p^p=\|b^{\frac{2}{p}} \cdot \nabla v\|_p^p \\
& \text{(we argue as in (\textbf{a}))} \\
&\leq \delta \big(\lambda\|\nabla v\|_p^p + \|\nabla |\nabla v|^{\frac{p}{2}}\|_2^2 \big),
\end{align*}
where, clearly, $\|\nabla v\|^p_p \leq \mu^{-\frac{p}{2}}\|f\|_p^p$. In turn, arguing as in (\textbf{a}), we arrive at $\mu\|w\|_p^p + I_p + (p-2)J_p = \langle f,-\nabla \cdot (w|w|^{p-2})$ ($w=\nabla v$),
$$
\mu\|w\|_p^p + (p-1)J_p \leq \langle f^2,|w|^{p-2}\rangle + (p-2)\langle f, |w|^{p-3}w \cdot \nabla |w|\rangle),
$$
$$
\mu\|w\|_p^p + (p-1)J_p \leq \langle f^2,|w|^{p-2}\rangle + (p-2)\bigl(\varepsilon J_p + \frac{1}{4\varepsilon}\langle f^2,|w|^{p-2}\rangle  \bigr), \quad \varepsilon>0.
$$
Selecting $\varepsilon$ sufficiently small, we obtain
$$
J_p \leq C_0\|w\|_p^{p-2}\|f\|_p^2.
$$
Now,  applying $\|w\|_p \leq \mu^{-\frac{1}{2}}\|f\|_p$, we arrive at $\|\nabla |\nabla v|^{\frac{p}{2}}\|_2^2 \leq C\mu^{-\frac{p}{2}+1}\|f\|_p^p$. 
Hence,
$\|G_pf\|_{p} \leq C_1 \mu^{-\frac{1}{2}+\frac{1}{p}}\|f\|_p$ for all $\mu \geq \mu_0$.

\smallskip

(\textbf{c}) Set $v=(\mu-\Delta)^{-1}|b|^{1-\frac{2}{p}}f\;(=Q_p f)$, $0 \leq f \in L^p$. Then, multiplying  $(\mu  - \Delta )v = |b|^{1-\frac{2}{p}}f$ by $v^{p-1}$, we obtain
$$
\mu \|v\|_p^p + \frac{4(p-1)}{p^2} \|\nabla v^{\frac{p}{2}}\|_2^2 = \langle |b|^{1-\frac{2}{p}}f, v^{p-1} \rangle,
$$
where we estimate the RHS using Young's inequality:
\begin{align*}
\langle |b|^{1-\frac{2}{p}} v^{\frac{p}{2}-1}, f v^{\frac{p}{2}} \rangle \leq \varepsilon^{\frac{2p}{p-2}}\frac{p-2}{2p} \langle |b|^2 v^p \rangle + \varepsilon^{-\frac{2p}{p+2}}\frac{p+2}{2p} \langle f^{\frac{2p}{p+2}} v^{\frac{p^2}{p+2}}\rangle \quad \varepsilon>0.
\end{align*}
Using $b \in \mathbf{F}_\delta$ and selecting $\varepsilon>0$ sufficiently small, we obtain that for any $\mu_1>0$ there exists $C>0$ such that
$$
(\mu-\mu_1) \|v\|_p^p \leq C\langle f^{\frac{2p}{p+2}} v^{\frac{p^2}{p+2}}\rangle, \qquad \mu>\mu_1.
$$
Therefore,
$
(\mu-\mu_1) \|v\|_p^p \leq C\langle f^p \rangle^{\frac{2}{p+2}}\langle  v^{p}\rangle^{\frac{p}{p+2}}
$, so
$\|v\|_p \leq C_2 \mu^{-\frac{1}{2}-\frac{1}{p}}\|f\|_p
$. 

\medskip

Step 2: We now use the results of Step 1 to prove assertions (\textit{ii}) and (\textit{iii}). Below we use the following formula: For every $0<\alpha<1$, $\mu > 0$,
\begin{equation*}
(\mu-\Delta)^{-\alpha}=\frac{\sin \pi \alpha}{\pi}\int_0^\infty t^{-\alpha}(t+\mu-\Delta)^{-1}dt.
\end{equation*}

We have 
\begin{align*}
\|Q_{p}(q)f\|_p & \leq ~\|(\mu-\Delta)^{-\frac{1}{2}+\frac{1}{q}}|b|^{1-\frac{2}{p}}|f|\|_p \\
& \leq 
k_{q}\int_0^\infty t^{-\frac{1}{2}+\frac{1}{q}}\|(t+\mu-\Delta)^{-1}|b|^{1-\frac{2}{p}}|f|\|_p dt  \\
& \text{(we use (\textbf{c}))} \\
& \leq ~
k_{q}C_{2} \int_0^\infty t^{-\frac{1}{2}+\frac{1}{q}}(t+\mu)^{-\frac{1}{2}-\frac{1}{p}}dt \;\|f\|_p =K_{2,q}\|f\|_p,\quad f \in \mathcal E,
\end{align*}
where, clearly, $K_{2,q}<\infty$ due to $q>p$.

\smallskip

It suffices to consider the case $r>2$. We have
\begin{align*}
\|G_{p}(r)f\|_p &\leq 
k_{r}\int_0^\infty t^{-\frac{1}{2}-\frac{1}{r}}\|b^{\frac{2}{p}}\cdot \nabla (t+\mu-\Delta)^{-1} f\|_pdt \;\;
\\ 
&(\text{we use (\textbf{b})})\\ 
&\leq 
k_{r}C_{1} \int_0^\infty t^{-\frac{1}{2}-\frac{1}{r}}(t+\mu)^{-\frac{1}{2}+\frac{1}{p}}dt\;\|f\|_p = K_{1,r}\|f\|_p,  \quad f \in \mathcal E,
\end{align*}
where, clearly, $K_{1,r}<\infty$ due to $r<p$. 
The proof of Lemma \ref{lem_GPQ} is completed. \hfill \qed

\bigskip

\section{Some examples of form-bounded vector fields}
\label{example_sect}

Below we list some sub-classes of the class of form-bounded vector fields, defined in elementary terms.

\medskip

 1.~Let us prove that $$
b \in L^\infty(\mathbb R_+,L^d+L^\infty) \quad \Rightarrow \quad b \in L^\infty\mathbf{F}_{\delta}+L^2_{\loc}(\mathbb R_+)
$$
for appropriate $\delta$ and $g$ (see Definition \ref{parab_fb_def}).
Here we have, by definition, $b=b_1+b_2$, where $b_1 \in L^\infty(\mathbb R_+,L^d)$, $b_2 \in L^\infty(\mathbb R_+,L^\infty)$. By H\"{o}lder's inequality, for a.e.\,$t \in \mathbb R_+$ and all $\varphi \in C_c^\infty$,
\begin{align*}
\|b(t)\varphi\|_2^2 
& \leq (1+\varepsilon)\|b_1(t)\|_d^2 \|\varphi\|_{\frac{2d}{d-2}}^2 + (1+\varepsilon^{-1})\|b_2(t)\|_\infty^2 \|\varphi\|_2^2 \qquad (\varepsilon>0)\\
& (\text{apply the Sobolev embedding theorem}) \\
& \leq C_S (1+\varepsilon) \|b_1(t)\|_d^2 \|\nabla \varphi\|_2^2 + (1+\varepsilon^{-1})\|b_2(t)\|_\infty^2 \|\varphi\|_2^2.
\end{align*}
Thus,
$b \in L^\infty\mathbf{F}_{\delta}+L^2_{\loc}(\mathbb R_+)$ with $$\delta:=C_S(1+\varepsilon)\|b_1\|^2_{L^\infty(\mathbb R_+,L^d)}, \quad g(t):=(1+\varepsilon^{-1})\|b_2(t,\cdot)\|_\infty^2$$
(in this paper we mostly care about the value of $\delta$, so  $\varepsilon$ should be chosen sufficiently small).

\smallskip

2.~Next, let us show that $$b \in C(\mathbb R_+,L^d+L^\infty) \quad \Rightarrow \quad b \in \mathbf{F}_\delta \quad \text{ with $\delta$ that can be chosen arbitrarily small.}
$$
Without loss of generality, let us carry out the proof for $b \in C(\mathbb R_+,L^d)$.

First, let $b=b(x)$. Since $|b| \in L^d$, one can represent for every 
$\varepsilon>0$  $b=b_1+b_2$, where $\|b_1\|_d<\varepsilon$ and $\|b_2\|_\infty<\infty$. 
(For instance, $b_2=b\mathbf{1}_{|b| \leq m}$ and $b_1=b-b_2$, 
so by the Dominated convergence theorem $\|b_1\|_d$ 
can be made arbitrarily small by selecting $m$ sufficiently large.) Now the previous example applies and yields the required. 

In the general case $b \in C(\mathbb R_+,L^d)$, the continuity of $b$ in time allows us to represent
$b(t,\cdot)=b_1(t,\cdot)+b_2(t,\cdot)$, where $\|b_1(t,\cdot)\|_d <\varepsilon$ for all $t \in [0,1]$ and 
$b_2$ is bounded on $[0,1] \times \mathbb R^d$. Repeating this on every interval $[n,n+1]$ ($n \geq 1$), one obtains $\|b_1\|_{L^\infty(\mathbb R_+,L^d)}<\varepsilon$ and $b_2 \in L^\infty_{\loc}(\mathbb R_+,L^\infty)$.
In fact, the continuity in time is not necessary for the smallness of $\delta$, e.g.\,consider $b(t,x)=c(t)b_0(x)$ where $c \in L^\infty(\mathbb R_+)$ is discontinuous and $|b_0| \in L^d$.

\smallskip

3.~Any vector field $$
b \in L^p(\mathbb R_+,L^q), \quad \frac{d}{q}+\frac{2}{p} \leq 1, \quad p \geq 2, \quad q \geq d
$$
is in $L^\infty\mathbf{F}_{\delta}+L^2_{\loc}(\mathbb R_+)$ with appropriate $\delta$. 
Indeed, e.g.\,in the more difficult case $\frac{d}{q}+\frac{2}{p} = 1$, we have by Young's inequality
\begin{align*}
|b(t,x)|  =\frac{|b(t,x)|}{\langle |b(t,\cdot)|^q\rangle^{\frac{1}{q}}}\langle |b(t,\cdot)|^q\rangle^{\frac{1}{q}} \leq \frac{d}{q} \biggl(\frac{|b(t,x)|^q}{\langle |b(t,\cdot)|^q\rangle} \biggr)^{\frac{1}{d}} + \frac{2}{p}\bigl( \langle |b(t,\cdot)|^q\rangle^{\frac{1}{q}} \bigr)^{\frac{p}{2}},
\end{align*}
where the first term is in $ L^\infty(\mathbb R_+,L^d)$ (and so by the first example it is form-bounded) and the second term is in $L^2(\mathbb R_+,L^\infty)$ (the second term squared is to be absorbed by the function $g$). 
If $p<\infty$, $q>d$, then one can argue as in the previous example to show that $\delta$ can be chosen arbitrarily small.

\smallskip

4.~The class $\mathbf{F}_\delta$ contains vector fields $b=b(x)$ with $|b|$ in $L^{d,w}$, i.e.\,the weak $L^d$ class  (Section \ref{notation_sect}). Indeed, by  \cite[Prop.~2.5, 2.6, Cor.~2.9]{KPS}, if $|b| \in L^{d,w}$, then 
$b \in \mathbf{F}_{\delta_1}$ with 
\begin{align*}
\sqrt{\delta_1}&=\||b| (\lambda - \Delta)^{-\frac{1}{2}} \|_{2 \rightarrow 2} \\ & \leq
\|b\|_{d,w} \Omega_d^{-\frac{1}{d}} \||x|^{-1} (\lambda - \Delta)^{-\frac{1}{2}} \|_{2 \rightarrow 2}  \leq \|b\|_{d,w} \Omega_d^{-\frac{1}{d}} \frac{2}{d-2},
\end{align*}
where $\Omega_d=\pi^{\frac{d}{2}}\Gamma(\frac{d}{2}+1)$ is the volume of  $B_1(0) \subset \mathbb R^d$.

\smallskip

5.~The Chang-Wilson-Wolff class ${\rm CW}^2$ consists of vector fields $b=b(x)$ such that $$|b| \in L^2_{\loc}(\mathbb R^d)$$ and
$$
\|b\|_{{\rm CW}^2}:=\biggl(\sup_Q \frac{1}{|Q|}\int_Q |b(x)|^2\, l(Q)^2 \xi\big(|b(x)|^2\,l(Q)^2 \big) dx<\infty \biggr)^{\frac{1}{2}}<\infty,
$$
where $|Q|$ and $l(Q)$ are the volume and the side length of cube $Q \subset \mathbb R^d$, respectively,
$\xi:\mathbb R_+ \rightarrow [1,\infty[$ is an increasing function such that
$$
\int_1^\infty \frac{dx}{x\xi(x)}<\infty.
$$

One has, for every $\varepsilon>0$,
$$
M_{2+\varepsilon} \quad \subsetneq \quad {\rm CW}^2 \quad \subsetneq \quad \mathbf{F}_\delta
$$
with 
$\delta=\delta\big(\|b\|_{\rm CW^2}\big)$, see \cite{CWW}.

5. More generally, vector fields in $L^\infty(\mathbb R_+,M_q)$, $q>2$, or $L^\infty(\mathbb R_+,{\rm CW}^2)$ are form-bounded.

\bigskip

\section{Smooth approximations of form-bounded (-type) vector fields}

\label{approx_app}

We construct two kinds of smooth approximations (or regularizations) of a vector field $b$: defined by mollifying a cutoff of $b$, or by mollifying $b$ directly. 

\medskip

\textbf{1.~}If $b=b(x)$ is either in the class $\mathbf{F}_\delta$ (form-bounded) or in the class $\mathbf{F}_\delta^{\scriptscriptstyle 1/2}$ (weakly form-bounded), it is clear that multiplying $b$ by 
the indicator $\mathbf 1_m$ of $\{x \in \mathbb R^d \mid |b(x)| \leq m, |x| \leq m\}$ leaves us in the corresponding class with the same $\delta$ and $\lambda_\delta$. This gives us a bounded, compact-support approximation $\{\mathbf{1}_mb\}_{m=1}^\infty$ of $b$ (in $L^2_{\loc}$ or in $L^1_{\loc}$, respectively).

We can also go one step further and apply to $\mathbf{1}_mb$ a mollifier, which will still give us a uniformly form-bounded (or uniformly weakly form-bounded) approximation of $b$ by vector fields whose components are in the Schwartz space $\mathcal S$.
Below we provide the details of this simple construction for a time-inhomogeneous form-bounded vector field  $b \in L^\infty\mathbf{F}_\delta+L^2_{\loc}(\mathbb R_+)$:
\begin{equation}
\label{fb9}
\tag{$\bullet$}
\|b(t,\cdot)\varphi\|_2^2 \leq \delta \|\nabla \varphi\|_2^2+g(t)\|\varphi\|_2^2, \quad \varphi \in W^{1,2}
\end{equation}
for a.e.\,$t \in \mathbb R$, for a $0 \leq g \in L^1_{\loc}(\mathbb R)$ (Definition \ref{parab_fb_def}).

Let us extend $b$ to $\{t<0\}$ by $0$ and set
\begin{equation}
\label{b_m_cutoff}
\tag{$\star$}
b_m:=c_m E_\varepsilon (\mathbf 1_m b),
\end{equation}
where $E_\varepsilon \equiv E^{d+1}_\varepsilon $ is the De Giorgi or Friedrichs mollifier on $\mathbb R \times \mathbb R^{d}$ (Section \ref{notation_sect}), $\mathbf 1_m$ is the indicator of $\{(t,x) \in \mathbb R^{1+d} \mid |b(t,x)| \leq m, |x| \leq m, |t| \leq m\}$, and $\varepsilon_m \downarrow 0$ and $c_m \uparrow 1$ are to be chosen.
Clearly, 
\begin{equation}
\label{b_m_0}
b_m \in L^\infty \cap C^\infty (\mathbb R^{1+d},\mathbb R^d),
\end{equation}
also, provided that $\varepsilon_m \downarrow 0$ sufficiently rapidly,
\begin{equation}
\label{b_m_1}
b_m \rightarrow b \quad \text{ in $L^2_{\loc}(\mathbb R_+ \times \mathbb R^d, \mathbb R^d$}).
\end{equation}
and, provided that  $c_m \uparrow 1$ sufficiently slow,
\begin{equation}
\label{b_m_2}
\|b_m(t)\varphi\|_2^2   \leq \delta \|\nabla \varphi\|_2^2 +g(t)\|\varphi\|_2^2, \quad t \geq 0, \quad \varphi \in W^{1,2},
\end{equation}
i.e.\,$\{b_m\}$ are uniformly form-bounded.

\begin{proof}[Proof of \eqref{b_m_2}]  First, define
$
\tilde{b}_m=E_\varepsilon (\mathbf 1_m b)
$
and write
\begin{equation}
\label{st}
\tag{$\ast$}
\tilde{b}_m=\mathbf{1}_m b + \big(\tilde{b}_m - \mathbf{1}_m b\big).
\end{equation}
Clearly, the first term satisfies
$$
\|\mathbf{1}_m b(t)\varphi\|_2^2   \leq \delta \|\nabla \varphi\|_2^2 +g(t)\|\varphi\|_2^2,
$$
In turn, since  $\mathbf{1}_mb$ has compact support and is in $L^\infty(\mathbb R_+,L^r)$ for any $r>d$, given any $\gamma_m \downarrow 0$ we can select $\varepsilon_m \downarrow 0$ so that
$\|\tilde{b}_m - \mathbf{1}_m b\|^2_{L^r(\mathbb R_+ \times \mathbb R^d)} \leq \gamma_m$,
so, in view of example 3 in Appendix \ref{example_sect},  the second term in \eqref{st} satisfies
$$
\|\big(\tilde{b}_m(t) - \mathbf{1}_m b(t)\big)\varphi\|_2^2   \leq C \gamma_m \|\nabla \varphi\|_2^2.
$$
Therefore, 
$$
\|\tilde{b}_m(t)\varphi\|_2^2   \leq \delta_m \|\nabla \varphi\|_2^2 +g(t)\|\varphi\|_2^2
$$
with $\delta_m=(\sqrt{\delta}+\sqrt{C\gamma_m})^2$.
Now, multiplying $\tilde{b}_m$ by $c_m=\frac{\delta}{\delta_m}$  (clearly, $c_m \uparrow 1$) and recalling that $b_m=c_m\tilde{b}_m$, we obtain \eqref{b_m_2}. 
\end{proof}

\medskip

\textbf{2.~}In fact, do not need the cutoff function in \eqref{b_m_cutoff} to construct a smooth approximation of a $b=b(x)$, $b \in \mathbf{F}_\delta$:
$$
\|b\varphi\|_2^2 \leq \delta \|\nabla \varphi\|_2^2 + c_\delta\|\varphi\|_2^2, \quad \varphi \in W^{1,2}.
$$  This observation is important e.g.\,if one needs to control the divergence of the approximating vector fields. In \cite{KiS_heat} the authors defined
\begin{equation}
\label{E}
\tag{$\star\star$}
b_\varepsilon:=E_{\varepsilon} b,
\end{equation}
where $E_{\varepsilon} \equiv E_\varepsilon^d:=e^{\varepsilon\Delta}$ is the De Giorgi's mollifier on $\mathbb R^d$ and $\varepsilon \downarrow 0$ (at any rate). (We can also use Friedrichs' mollifier, see remark below.)
We have
\begin{equation}
\label{C_1}
b_\varepsilon \in L^\infty\cap C^\infty(\mathbb R^d,\mathbb R^d),
\end{equation}
\begin{equation}
\label{C_2}
\|b_\varepsilon \varphi\|_2^2 \leq \delta \|\nabla \varphi\|_2^2 + c_\delta\|\varphi\|_2^2, \quad \text{ i.e. } b_\varepsilon \in \mathbf{F}_\delta
\text{ with the same $c_\delta$},
\end{equation}
\begin{equation}
\label{C_3}
b_\varepsilon \rightarrow b \quad \text{ in } L^2_{\loc}(\mathbb R^d,\mathbb R^d).
\end{equation}

\begin{proof}[Proof of \eqref{C_1}-\eqref{C_3}] We repeat the argument from \cite{KiS_heat}. 

To prove \eqref{C_1}, we represent $b_\varepsilon=E_{\varepsilon/2}E_{\varepsilon/2}b$, so it suffices to only prove that $|b_\varepsilon| \in L^\infty$. Indeed, we have, using Fatou's lemma,
\begin{align*}
|b_\varepsilon(x)| & \leq \liminf_n\big\langle e^{\varepsilon\Delta}(x,\cdot)\mathbf{1}_{B_n(0)}(\cdot) |b(\cdot)|\big\rangle \\
& \leq \liminf_n\big\langle e^{\varepsilon\Delta}(x,\cdot)\mathbf{1}_{B_n(0)}(\cdot)|b(\cdot)|^2\big\rangle^{\frac{1}{2}} \leq \big(\delta \big\langle \big|\nabla \sqrt{e^{\varepsilon\Delta}(x,\cdot)}\big|^2\big\rangle + c_\delta\big)^{\frac{1}{2}},
\end{align*}
where $\big|\nabla_y \sqrt{ e^{\varepsilon\Delta}(x,y)}\big|=(4\pi\varepsilon)^{-\frac{d}{4}} \frac{|x-y|}{4\varepsilon}e^{-\frac{|x-y|^2}{8\varepsilon}} \leq C\varepsilon^{-\frac{d}{4}-\frac{1}{2}} e^{-\frac{c|x-y|^2}{\varepsilon}}$, and so $|b_\varepsilon| \in L^\infty$ for each $\varepsilon>0$.

Let us prove \eqref{C_2}. Indeed, $|b_\varepsilon|\leq\sqrt{E_\varepsilon |b|^2}$, and so
\begin{align*}
\|b_\varepsilon \varphi\|_2^2 &\leq \langle E_\varepsilon |b|^2,\varphi^2\rangle=\|b\sqrt{E_\varepsilon \varphi^2}\|^2_2 \\
&\leq \delta\|\nabla\sqrt{E_\varepsilon \varphi^2}\|_2^2+c_\delta\|\varphi\|_2^2, \quad \varphi \in W^{1,2},
\end{align*}
where
\begin{align}
\|\nabla\sqrt{E_\varepsilon \varphi^2}\|_2 & =\|\frac{E_\varepsilon(|\varphi||\nabla|\varphi|)}{\sqrt{E_\varepsilon \varphi^2}}\|_2 \label{E2} \tag{$\ast\ast$}\\
&\leq \|\sqrt{E_\varepsilon|\nabla |\varphi||^2}\|_2=\|E_\varepsilon|\nabla |\varphi||^2\|_1^\frac{1}{2} \notag \\
& \leq\|\nabla|\varphi|\|_2\leq \|\nabla \varphi\|_2 \notag,
\end{align}
i.e. $b_\varepsilon\in\mathbf{F}_{\delta}$. (The fact that $\|b\sqrt{E_\varepsilon \varphi^2}\|_2<\infty$ follows from $\mathbf{1}_{\{|b|\leq n\}}b\in \mathbf{F}_\delta$ and the inequality $\|\mathbf{1}_{\{|b|\leq n\}}b \sqrt{E_\varepsilon \varphi^2}\|_2^2\leq \delta\|\nabla \varphi\|_2^2+c_\delta\|\varphi\|_2^2$, using Fatou's lemma).

Regarding the proof of \eqref{C_3}, let us only demonstrate that $b_\varepsilon-b \rightarrow 0$ in $L^2(B_1)$, $B_1 \equiv B_1(0)$. To that end, we fix some $R>1$ and represent on $B_1$:
\begin{align*}
b_\varepsilon-b & = I_1 + I_2, \qquad I_1:=E_{\varepsilon} (1-\mathbf{1}_{B_R})b, \quad I_2:=E_{\varepsilon} (\mathbf{1}_{B_R}b) - \mathbf{1}_{B_R}b
\end{align*}
(of course, on $B_1$, one has $b=\mathbf{1}_{B_R}b$ since $R>1$).
Then $I_2 \rightarrow 0$ in $L^2(B_1)$  since $\mathbf{1}_{B_R}b$ has compact support. In turn, $I_1 \rightarrow 0$ in $L^2(B_1)$ by the separation property of the Gaussian density, i.e.\,$e^{\varepsilon\Delta}(x,y) \rightarrow 0$ uniformly in $x \in B_1$ if $y \in \mathbb R^d - B_R$ (here we have used $R>1$).
\end{proof}

\begin{remark}
If we were to use the Friedrichs mollifier in \eqref{E}, then we would get smoothness of $b_\varepsilon$ and convergence \eqref{C_3} from the usual properties of Friedrichs mollifiers. But we would have to be slightly more careful in \eqref{E2} in order to avoid division by zero, e.g.\,replace $|f|^2$ by $|f|^2+e^{-k|x|^2}$, carry out the estimates and then take $k \rightarrow \infty$.
\end{remark}

\textbf{3.~}The regularization \eqref{E} can also be used to handle time-inhomogeneous form-bounded drifts $b \in L^\infty \mathbf{F}_\delta + L^2_{\loc}(\mathbb R)$. That is, we can put
$$
b_\varepsilon:=E_\varepsilon^1 E_\varepsilon^d b, \quad \varepsilon \downarrow 0,
$$
where $E_\varepsilon^d$ is the De Giorgi or Friedrichs mollifier on $\mathbb R^d$ (in the spatial variables) and $E_\varepsilon^{1}$  is the Friedrichs mollifier on $\mathbb R$ (in the time variable; we use 
the Friedrichs mollifier here since, in the time variable, $b$ is in only locally in $L^1(\mathbb R)$, as is determined by our assumption on $g$, cf.\,\eqref{fb9}). 

This regularization $\{b_\varepsilon\}$ satisfies 
\begin{gather}
\label{D_1}
|b_\varepsilon| \in L^\infty_{\loc}(\mathbb R,L^\infty(\mathbb R^d)), \quad b_\varepsilon \text{ are } C^\infty \text{ smooth},\\[2mm]
\label{D_2}
\|b_\varepsilon(t,\cdot)\varphi\|_2^2 \leq \delta \|\nabla \varphi\|_2^2+g_\varepsilon (t)\|\varphi\|_2^2 \quad \text{ for all } t \in \mathbb R, \\
\label{D_22}
\sup_{\varepsilon>0}\int_{t_0}^{t_1}g_\varepsilon (s)ds<\infty \text{ for all finite } t_0,t_1, \text{ where } g_\varepsilon:=E_\varepsilon^{1 }g,\\
\label{D_3}
b_\varepsilon \rightarrow b \quad \text{ in } L^2_{\loc}(\mathbb R^{1+d},\mathbb R^d).
\end{gather}

\begin{proof}[Proof of \eqref{D_2}] First, we regularize $b$ only in the spatial variable: for every $t \in \mathbb R$, put $$\tilde{b}_\varepsilon(t,\cdot):=E_\varepsilon^d b(t,\cdot).$$
Then, since for a.e.\,$t$ $b(t,\cdot)$ is form-bounded on $\mathbb R^d$, we have by \eqref{C_2} 
\begin{equation}
\label{tilde_b}
\tag{$\ast\ast\ast$}
\|\tilde{b}_\varepsilon(t,\cdot)\varphi\|_2^2 \leq \delta \|\nabla \varphi\|_2^2+g(t)\|\varphi\|_2^2
\end{equation}
for a.e.\,$t \in \mathbb R$.
Next, recalling that $b_\varepsilon=E^{1}_\varepsilon \tilde{b}_\varepsilon$ and noting a pointwise inequality
$$
|b_\varepsilon(t)|^2 \leq (E_\varepsilon^{1}|\tilde{b}_\varepsilon|^2)(t), \quad t \in \mathbb R,
$$
we estimate, for every $\varphi \in W^{1,2}(\mathbb R^d)$, 
\begin{align*}
\langle |b_\varepsilon(t)|^2 \varphi^2 \rangle & \leq E_\varepsilon^{1}\langle|\tilde{b}_\varepsilon|^2 \varphi^2\rangle(t) \\
& (\text{we are applying $E_\varepsilon^{1}$ to both sides of the inequality \eqref{tilde_b}}) \\
& \leq \delta \|\nabla \varphi\|_2^2+g_\varepsilon(t)\|\varphi\|_2^2.
\end{align*}
Thus, we arrive at \eqref{D_2}. 
\end{proof}

A similar construction was considered in \cite{KiS_note}.

\bigskip

\section{Trotter's approximation theorem}

\label{trotter_sect}

Consider a sequence $\{e^{-tA_k}\}_{k=1}^\infty$ of $C_0$ semigroups on a (complex) Banach space $Y$.

\begin{theorem}[{H.F.\,Trotter \cite[Ch.\,IX]{Ka}}]
Let
$\sup_k\|(\mu+A_k)^{-m}\|_{Y \rightarrow Y} \leq M(\mu-\omega)^{-m}$, $m=1,2,\dots$, $\mu>\omega$, and 
$$s\mbox{-}\lim_{\mu \rightarrow \infty}\mu(\mu+ A_k)^{-1}=1 \quad \text{ uniformly in $k$,}
$$ and let $s\mbox{-}\lim_{k}(\zeta+ A_k)^{-1}$
exist for some $\zeta$ with ${\rm Re\,} \zeta>\omega$. Then there is a $C_0$ semigroup $e^{-tA}$ such that
$$
(z+ A_k)^{-1} \overset{s}{\rightarrow} (z+ A)^{-1} \quad \text{ for every } {\rm Re}\, z>\omega,
$$
and
$$
e^{-tA_k} \overset{s}{\rightarrow} e^{-tA}
$$
uniformly in any finite interval of $t \geq 0$.
\end{theorem}

The first condition of the theorem is satisfied if  e.g.
$$
\sup_k\|(\mu+A_k)^{-1}\|_{Y \rightarrow Y} \leq (\mu-\omega)^{-1}, \quad \mu>\omega
$$
(obviously)
of if
$$\sup_k\|(z+A_k)^{-1}\|_{Y \rightarrow Y} \leq C|z-\omega|^{-1}, \quad {\rm Re\,}z>\omega,$$ see \cite[IX.6.1]{Ka}. The second condition is what can be verified in practice when one is dealing with quasi bounded semigroups.


\bigskip

\begin{thebibliography}{99}

\bibitem[A]{A} D. Adams, Weighted nonlinear potential theory, {\em Trans. Amer. Math. Soc.}
\textbf{297} (1986), 73-94.

\bibitem[AF]{AF} R. Adams and J. Fournier, ``Sobolev Spaces'', Second Edition, Elsevier, 2003.


\bibitem[AD]{AD} D.\,Albrighton and H.\,Dong, Regularity properties of passive scalars with rough divergence-free drifts, arXiv:2107.12511.

\bibitem[B]{B} R. J. Bagby, Lebesgue spaces of parabolic potentials, {\em Illinois J. Math.} \textbf{15} (1971), 610-634.

\bibitem[BC]{BC} R.~Bass and Z.-Q.~Chen, \newblock Brownian motion with singular drift.
\newblock {\em Ann. Probab.}, 31 (2003), 791-817.

\bibitem[BG]{BG} P.\,Baras and J.\,A.\,Goldstein, {\newblock The heat equation with a singular potential.} {\em Trans. Amer. Math. Soc.},
\textbf{284} (1984), 121-139.


\bibitem[BFGM]{BFGM} L.~Beck, F.~Flandoli, M.~Gubinelli and M.~Maurelli, 
\newblock Stochastic ODEs and stochastic linear PDEs with critical
drift: regularity, duality and uniqueness. 
\newblock{\em Electr. J. Probab.}, \textbf{24} (2019),  Paper No. 136, 72 pp (arXiv:1401.1530).

\bibitem[BS]{BS}
A.G. Belyi and Yu.A. Semenov. 
\newblock On the $L^p$-theory of Schr\"{o}dinger semigroups. II. 
\newblock{ \em Sibirsk.
Math.~J.}, 31 (1990), p.~16-26; English transl. in \newblock {\em Siberian Math.~J.}, 31 (1991), 540-549.


\bibitem[BO]{BO} \'{A}. B\'{e}nyi and T.\,Oh, The Sobolev inequality on the torus revisited, {\em Publ. Math. Debrecen} \textbf{83}, no. 3 (2013), 359-374.


\bibitem[Bi]{Bi} M.\,S.\,Birman, ``On the spectrum of singular boundary-value problems'' (in Russian), \textit{Mat. Sbornik} \textbf{55}(97) (1961), 125-174.


\bibitem[BJ]{BJ}  K.\,Bogdan and T.\,Jakubowski, Estimates of heat kernel of fractional Laplacian perturbed by gradient operators, {\em Comm. Math. Phys.}, \textbf{271} (2007) 179-198.

\bibitem[BGe]{BGe} R.M. Blumenthal and R.K. Getoor, Markov Processes and Potential Theory.
Pure and Applied Mathematics 29. Academic Press, New York.

\bibitem[BJW]{BJW} D.\,Bresch, P.-E.\,Jabin and Z.\,Wang, Mean field limit and quantitative estimates with singular attractive kernels, {\em Duke Math. J.} \textbf{172} (2023), no. 13, 2591-2641 (arXiv:2011.08022).


\bibitem[CWW]{CWW} S.Y.A.\;Chang, J.M.\;Wilson and T.H.\;Wolff, Some weighted norm inequalities concerning
the Schr\"{o}dinger operator, \newblock { \em Comment.\;Math.\;Helvetici}, \textbf{60} (1985), 217-246.


\bibitem[CM]{CM} P.-\'{E}. Chaudru de Raynal and S.\,Menozzi, On multi-dimensional stable driven stochastic differential equations with Besov drift, {\em Electron. J. Probab.} \textbf{27} (2022), Paper No. 163, 52 pp. (arXiv:1907.12263).

\bibitem[CJM]{CJM} P.-\'{E}. Chaudru de Raynal, J.-F.\,Jabir and S.\,Menozzi, Multidimensional stable driven McKean-Vlasov SDEs with distributional interaction kernel: a regularization by noise perspective,  arXiv:2205.11866.




\bibitem[CKS]{CKS} Z.-Q.\,Chen, P.\,Kim and R.\,Song, Dirichlet heat kernel estimates for fractional Laplacian with gradient perturbation, {\em Ann.\,Prob.}, \textbf{40} (2012), 2483-2538.


\bibitem[CFKZ]{CFKZ} Z.-Q.\,Chen, P.\,J.\,Fitzsimmons, K.\,Kuwae and T.-S.\,Zhang, Perturbation of symmetric Markov processes, {\em Probab. Theory Related Fields} \textbf{140} (2008), 239-275.

\bibitem[CW]{CW} Z.-Q.\,Chen and L.\,Wang, Uniqueness of stable processes with drift, {\em Proc.\,Amer.\,Math.\,Soc}, \textbf{144} (2017), p.\,2661-2675 (arXiv:1309.6414).


\bibitem[C]{C} A.\,S.\,Cherny, On the uniqueness in law and the pathwise uniqueness for stochastic
differential equations, {\em Theory Probab. Appl.}, \textbf{46}(3) (2002), 406-419.


\bibitem[CE]{CE}
 A.\,S.\,Cherny and H.-J.\,Engelbert. 
 \newblock Singular Stochastic Differential Equations. 
 \newblock{\em LNM 1858}. Springer-Verlag, 2005.


\bibitem[CFr]{CF}  F. Chiarenza and M. Frasca, A remark on a paper by C. Fefferman, {\em Proc. Amer. Math. Soc.}, \textbf{108} (1990),
407-409. 


\bibitem[D]{D} H.\,Dong, Recent progress in the $L^p$ theory for elliptic and parabolic equations with discontinuous coefficients, {\em Anal. Theory Appl.} 36 (2020), no. 2, 161-199 (arXiv:2006.03966).



\bibitem[F]{F}  C. Fefferman, The uncertainty principle, {\em Bull. Amer. Math. Soc.} \textbf{9} (1983), 129-206.



\bibitem[FK]{FK} P.\,J.\,Fitzsimmons and K.\,Kuwae, Non-symmetric perturbations of symmetric
Dirichlet forms, {\em J.\,Funct.\,Anal.} \textbf{208} (2004), 140-162.

\bibitem[FIR]{FIR} F.\,Flandoli, E.\,Issoglio and F.\,Russo, Multidimensional stochastic differential equations with distributional drift, {\em Trans. Amerc. Math. Soc.}, \textbf{369} (2017), 1665-1688 (arXiv:1401.6010).

\bibitem[FGP]{FGP} F.\,Flandoli, M.\,Gubinelli and E.\,Priola, Well-posedness of the transport equation
by stochastic perturbation, {\em Invent.\,Math.} \textbf{180} (2010), 1-53.

\bibitem[FJ]{FJ} N.\,Fournier and B.\,Jourdain, Stochastic particle approximation of the Keller-Segel and two-dimensional generalization of Bessel process, {\em Ann.\,Appl. Probab.} \textbf{27} (2017), 2807-2861.


\bibitem[GM]{GM} B.\,Gess and M.\,Maurelli, Well-posedness by noise for scalar conservation laws, \textit{Comm. Partial
Differential Equations} \textbf{43} (2018), no. 12, 1702-1736.


\bibitem[GZa]{GZa} J.\,A.\,Goldstein and Qi.\,S.\,Zhang, {\newblock Linear parabolic equation with strongly singular potentials}, {\em Trans.
Amer. Math. Soc.} \textbf{355} (2003), 197-211.


\bibitem[Go]{G} V.R.\,Gopala Rao, A characterization of parabolic function spaces, {\em Amer.\,J.\,Math.}, \textbf{99} (1977), 985-993.



\bibitem[GO]{GO}
L.~Grafakos and  S.~Oh.
\newblock The Kato-Ponce inequality.
\newblock {\em Comm.~Partial~Diff.~Equ.}, 39 (2014), 1128-1157.



\bibitem[GC]{GvC} A.\,Gulisashvili and J.A.\,van Casteren, Non-autonomous Kato Classes and Feynman-Kac Propagators, \textit{World Scientific}, 2006.

\bibitem[G]{Gu} S.\,Gupta, Hardy and Rellich inequality on lattices, {\em Calc.\,Var.\,Partial Diff.\,Equations},  62, article no. 81 (2023).


\bibitem[H]{H} T.\,Hara, A refined subsolution estimate of weak subsolutions to second order
linear elliptic equations with a singular vector field, {\em Tokyo J.\,Math.}, \textbf{38}(1) (2015), 75-98.

\bibitem[He]{He}
E.\,Heinz, 
\newblock \textit{Beitr\"{a}ge zur St\"{o}rungstheorie der Spektralzerlegung},
\newblock {Math. Ann.}, \textbf{123} (1951) 415-438.


\bibitem[J]{J} P. Jin, Brownian motion with singular time-dependent drift. {\em J. Theoret. Probab.},
\textbf{30} (2017), 1499-1538 (arXiv:1710.05227).


\bibitem[Ka]{Ka}  T.\,Kato, Perturbation Theory for Linear Operators, Springer-Verlag, Berlin, Heidelberg,
1995.

\bibitem[KSo]{KSo}  P.\,Kim and R.\,Song, {Stable process with singular drift}, \newblock{\em 
Stoc.\,Proc.\,Appl.} \textbf{124} (2014), 2479-2516.




\bibitem[Ki1]{Ki_a_new_approach} D.\,Kinzebulatov, A new approach to the $L^p$-theory of $-\Delta + b\cdot\nabla$, and its applications to Feller processes with general drifts,
\newblock {\em Ann.~Sc.~Norm.~Sup.~Pisa (5)}, \textbf{17} (2017), 685-711 (arXiv:1502.07286). 


\bibitem[Ki2]{Ki2} D.\,Kinzebulatov, Regularity theory of Kolmogorov operator revisited, {\em Canadian Bull. Math.} \textbf{64} (2021), 725-736 (arXiv:1807.07597).

\bibitem[Ki3]{Ki} D.\,Kinzebulatov,
\newblock Feller evolution families and parabolic equations with form-bounded vector fields,
\newblock {\em Osaka J.\,Math.}, \textbf{54} (2017), 499-516 (arXiv:1407.4861).


\bibitem[Ki4]{Ki_measure} D. Kinzebulatov, Feller generators with measure-valued drifts, {\em Potential Anal.}, \textbf{48} (2018), 207-222.

\bibitem[Ki5]{Ki_Morrey} D.\,Kinzebulatov, Parabolic equations and SDEs with time-inhomogeneous Morrey drift, arXiv:2301.13805.

\bibitem[Ki6]{Ki_Orlicz} D.\,Kinzebulatov, Laplacian with singular drift in a critical borderline case, arXiv:2309.04436.

\bibitem[KiM1]{KiM} D.Kinzebulatov and K.R.Madou, Stochastic equations with time-dependent singular drift, {\em J.\,Differential Equations}, \textbf{337} (2022), 255-293 (arXiv:2105.07312).

\bibitem[KiM2]{KiM_stable} D. Kinzebulatov and K.R. Madou, On admissible singular drifts of symmetric $\alpha$-stable process, {\em Math. Nachr.}, \textbf{295}(10) (2022), 2036-2064 (arXiv:2002.07001).

\bibitem[KiM3]{KiM_strong} D. Kinzebulatov and K.R. Madou, Strong solutions of SDEs with singular (form-bounded) drift via Roeckner-Zhao approach, arXiv:2306.04825.

\bibitem[KiS1]{KiS_brownian} D.\,Kinzebulatov and Yu.A.\,Sem\"{e}nov, Brownian motion with general drift, \newblock{\em Stoch. Proc. Appl.}, \textbf{130} (2020), 2737-2750 (arXiv:1710.06729).


\bibitem[KiS2]{KiS_theory} D.\,Kinzebulatov and Yu.\,A.\,Sem\"{e}nov, On the theory of the Kolmogorov operator in the spaces $L^p$ and $C_\infty$, {\em Ann. Sc. Norm. Sup. Pisa (5)} \textbf{21} (2020), 1573-1647 (arXiv:1709.08598).


\bibitem[KiS3]{KiS_Osaka} D.\,Kinzebulatov and Yu.\,A.\,Sem\"{e}nov, Feller generators and stochastic differential equations with singular (form-bounded) drift,  \newblock{\em Osaka J.\,Math.}, 58 (2021), 855-883 (arXiv:1904.01268).

\bibitem[KiS4]{KiS_sharp} D.\,Kinzebulatov and Yu.\,A.\,Sem\"{e}nov, Sharp solvability for singular SDEs, {\em Electr.\,J.\,Probab.}, \textbf{28} (2023), article no. 69, 1-15. (arXiv:2110.11232).

\bibitem[KiS5]{KiS_RIMS} D.\,Kinzebulatov and Yu.\,A.\,Sem\"{e}nov, {Fractional Kolmogorov operator and desingularizing weights,}  {\em Publ.\,Res.\,Inst.\,Math.\,Sci.\,Kyoto}, to appear (arXiv:2005.11199).


\bibitem[KiS6]{KiS_heat} D.\,Kinzebulatov and Yu.\,A.\,Sem\"{e}nov, Heat kernel bounds for parabolic equations with singular (form-bounded) vector fields, {\em Math.\,Ann.}, \textbf{384} (2022), 1883-1929.

\bibitem[KiS7]{KiS_Nash} D.\,Kinzebulatov and Yu.\,A.\,Sem\"{e}nov, Kolmogorov operator with the vector field in Nash class, {\em Tohoku Math. J.}, \textbf{74}(4) (2022), 569-596 (arXiv:2012.02843).

\bibitem[KiS8]{KiS_BMO} D.\,Kinzebulatov and Yu.\,A.\,Sem\"{e}nov, Regularity for parabolic equations with singular non-zero divergence vector fields, arXiv:2205.05169.


\bibitem[KiS9]{KiS_note} D.\,Kinzebulatov and Yu.A.\,Sem\"{e}nov, Remarks on parabolic Kolmogorov operator, arXiv:2303.03993.


\bibitem[KMS]{KiMSe} D. Kinzebulatov, K.R. Madou and Yu.A. Sem\"{e}nov, On the supercritical fractional diffusion equation with Hardy-type drift, {\em J. d'Analyse Math\'{e}matique}, to appear (arXiv:2112.06329).

\bibitem[KSS]{KiSS_transport} D.\,Kinzebulatov, Yu.\,A.\,Sem\"{e}nov and R.\,Song, Stochastic transport equation with singular drift, {\em Ann.\,Inst.\,Henri Poincar\'{e} (B) Probab. Stat.}, to appear (arXiv:2102.10610).

\bibitem[KSSz]{KiSSz} D.\,Kinzebulatov, Yu.\,A.\,Sem\"{e}nov and K.\,Szczypkowski, { Heat kernel of fractional Laplacian with Hardy drift via desingularizing weights,} {\em J. London Math. Soc.}, 104 (2021), 1861-1900 (arXiv:1904.07368).

\bibitem[KiV]{KiV} D.\,Kinzebulatov and R.\,Vafadar, On divergence-free (form-bounded type) drifts, {\em Discrete Contin. Dyn.\,Syst.\,Ser.\,S.}, to appear (arXiv:2209.04537).




\bibitem[Ko]{Ko} T.\,Komatsu, {On the martingale problem for generators of stable processes with perturbations,} \newblock{\em Osaka J. Math. } \textbf{21} (1984), 113-132.


\bibitem[KS]{KS} V.\,F.\,Kovalenko and Yu.\,A.\,Sem\"{e}nov,
{\newblock $C_0$-semigroups in $L^p(\mathbb R^d)$ and $C_\infty(\mathbb R^d)$ spaces generated by differential expression $\Delta+b\cdot\nabla$.} 
(Russian) {\em Teor. Veroyatnost. i Primenen.}, \textbf{35} (1990), 449-458; translation in {\em Theory Probab. Appl.} \textbf{35} (1990), 443-453.

\bibitem[KPS]{KPS} V.\,F.\,Kovalenko, M.\,A.\,Perelmuter and Yu.\,A.\,Sem\"{e}nov, Schr\"{o}dinger operators with $L^{\frac{1}{2}}_{w}$($R^{l}$)-potentials, \newblock{\em J.\,Math.\,Phys.}, \textbf{22}  (1981), 1033-1044. 



\bibitem[Kr1]{Kr1} N.\,V.\,Krylov,  \newblock  On diffusion processes with drift in $L_d$, {\em Probab. Theory Related Fields} \textbf{179} (2021), no. 1-2, 165-199 (arXiv:2001.04950). 

\bibitem[Kr2]{Kr2} N.\,V.\,Krylov, \newblock On strong solutions of It\^{o}'s equations with $A \in W^{1,d}$ and $B \in L^d$, {\em Ann. Probab.} \textbf{49} (2021), no. 6, 3142-3167 (arXiv:2007.06040).



\bibitem[Kr3]{Kr3} N.V.\,Krylov, On strong solutions of It\^{o}'s equations with $D\sigma$ and $b$ in Morrey classes containing $L^d$, {\em Ann. Probab.} \textbf{51} (2023), no. 5, 1729-1751 (arXiv:2111.13795).


\bibitem[Kr4]{Kr_Adams} N.V.\,Krylov, On parabolic Adams's, the Chiarenza-Frasca theorems, and some other
results related to parabolic Morrey spaces, {\em Mathematics in Engineering}, \textbf{5}(2) (2022), 1-20 (arXiv:2110.09555).


\bibitem[Kr5]{Kr_weak} N.V.\,Krylov, On weak solutions of time-inhomogeneous It\^{o}'s equations with VMO diffusion and Morrey drift,  arXiv:2303.11238.

\bibitem[Kr6]{Kr_Morrey_sdes} N.V.\,Krylov, Once again on weak solutions of time-inhomogeneous It\^{o}'s equations with VMO diffusion and Morrey drift, arXiv:2304.04634.

\bibitem[Kr7]{Kr_Morrey_parab} N.V.\,Krylov, On parabolic equations in Morrey spaces with VMO $a$ and Morrey $b$, $c$, arXiv:2304.03736.

\bibitem[KrR]{KR} N.\,V.\,Krylov and M.\,R\"{o}ckner. 
\newblock Strong solutions of stochastic equations with singular time dependent drift. 
\newblock {\em Probab. Theory Related Fields}, 131 (2005), 154-196.

\bibitem[Ku]{Ku} H.~Kunita, Stochastic Flows and Stochastic Differential Equations, Cambridge Studies in Advanced Mathematics, vol. 24, Cambridge University Press, Cambridge, 1990.


\bibitem[LS]{LS}
V.~A.~Liskevich and Yu.~A.~Sem\"{e}nov, 
\newblock Some problems on Markov semigroups,
\newblock{\em 
``Schr\"{o}dinger Operators, Markov Semigroups, Wavelet Analysis, Operator Algebras''
(M. Demuth et al.,, Eds.), Mathematical Topics: Advances in Partial Differential
Equations, Vol. 11, Akademie Verlag, Berlin} (1996), 163-217.

\bibitem[LZ]{LZ} V.\,Liskevich and Q.\,S.\,Zhang, Extra regularity for parabolic equations with drift
terms, \textit{Manuscripta Math.} \textbf{113} (2004), 191-209.



\bibitem[MK]{MK}
A.\,J.\,Majda and P.\,R.\,Kramer, Simplified models for turbulent diffusion:
theory, numerical modelling, and physical phenomena, {\em Physics Reports} \textbf{314} (1999), 237-574. 

\bibitem[MV]{MV}
V.\,G.\,Mazya and I.\,E.\,Verbitsky, Form boundedness of the general
second-order differential operator, {\em Comm.\,Pure Appl.\,Math.} \textbf{59} (2006), 1286-1329. 


\bibitem[MeZ]{MeZ} S.\,Menozzi and X.\,Zhang, Heat kernel of supercritical non-local operators with unbounded drift, {\em J. \'{E}c. polytech. Math.} \textbf{9} (2022), 537-579 (arXiv:2012.14475).


\bibitem[MNS]{MNS} G.\,Metafune, L.\,Negro and C.\,Spina, ``{Sharp kernel estimates for elliptic operators with second-order discontinuous coefficients}'', {\em J. Evol. Equ.} \textbf{18} (2018), 467-514.


\bibitem[MP]{MP} T.\,Meyer-Brandis and F.\,Proske, Construction of strong solutions of SDE's via Malliavin calculus, {\em J.\,Funct.\,Anal.}, \textbf{258} (2010)(11), 3922-3953.

\bibitem[MNP]{MNP} A.\,Mohammed, T.\,Nilssen and F.\,Proske, Sobolev differentiable stochastic flows
for SDEs with singular coefficients: applications to the transport equation, {\em Ann.\,Probab.},
\textbf{43}(3) (2015), 1535-1576.


\bibitem[N]{N} K.\,Nam, Stochastic differential equations with critical drifts. {\em Stoch. Proc. Appl.}, \textbf{130} (2020), 5366-5393 (arXiv:1802.00074).

\bibitem[NU]{NU} A.\,I.\,Nazarov and N.\,N.\,Uraltseva, The Harnack inequality and related properties for solutions to elliptic and parabolic equations with divergence-free lower order coefficients, {\em Algebra i Analiz}, \textbf{23} (2011), 136-168.


\bibitem[O]{O} E.-M.\,Ouhabaz, Analysis of Heat Equations on Domains, Princeton Univ.\,Press, 2005.


\bibitem[OSSV]{OSSV} E.-M.\,Ouhabaz, P.\,Stollmann, K.-Th.\,Sturm and J.\,Voigt, \newblock The Feller property for absorption semigroups, {\em  J.\,Funct.\,Anal.}\, \textbf{138} (1996), 351-378.


\bibitem[PZ]{PZ} N.\,Perkowski and W. van Zuiljen, Quantitative heat kernel estimates for diffusions with
distributional drift, {\em Potential Anal.}, https://doi.org/10.1007/s11118-021-09984-3 (2022)

\bibitem[Ph]{Ph} T.\,Phan, Local $W^{1,p}$-regularity estimates for weak solutions of parabolic equations with singular divergence-free drifts, {\em Electr. J. Differential Equations} (2017), Paper No.\,75, 22 pp. 

\bibitem[P1]{P} N.\,I.\,Portenko, \newblock Generalized Diffusion Processes. \newblock {\em AMS}, 1990.


\bibitem[P2]{P2} N.\,I.\,Portenko, \newblock Some perturbations of drift-type for symmetric stable processes, {\em Random Oper. Stochastic Equations}, \textbf{2} (1994), 211-224.



\bibitem[PP]{PP}
S.\,I.\,Podolynny and N. I. Portenko, \newblock On multidimensional stable processes with locally
unbounded drift, {\em Random Oper. Stochastic Equations}, \textbf{3} (1995), 113-124.

\bibitem[Pr]{Pr}
E.\,Priola, Pathwise uniqueness for singular SDEs driven by stable processes. {\em Osaka J. Math.}
\textbf{49} (2012), 421-447.

\bibitem[R]{R} F.\,Rezakhanlou, {Regular flows for diffusions with rough drifts}, arXiv:1405.5856.




\bibitem[RZh1]{RZ} M.\,R\"{o}ckner and G.\,Zhao, {SDEs with critical time dependent drifts: weak solutions}, {\em Bernoulli}, \textbf{29} (2023), 757-784 (arXiv:2012.04161).

\bibitem[RZh2]{RZ2} M.\,R\"{o}ckner and G.\,Zhao,  SDEs with critical time dependent drifts: strong solutions,
arXiv:2103.05803.

\bibitem[S1]{S} Yu.\,A.\,Sem\"{e}nov, \newblock Regularity theorems for parabolic equations, \newblock {\em J.\,Funct.\,Anal.}, \textbf{231} (2006), 375-417.


\bibitem[S2]{S2} Yu.\,A.\,Sem\"{e}nov, On perturbation theory for linear elliptic and parabolic operators; the method of
Nash, Proceedings of the Conference on Applied Analysis, April 19-21 (1996), B\^{a}ton-Rouge, Louisiana, \textit{Contemp. Math.}, \textbf{221} (1999), 217-284.

\bibitem[Si]{Si} B.\,Simon, Schr\"{o}dinger semigroups, {\em Bull.\,Amer.\,Math.\,Soc.}, \textbf{7}(3) (1982), 447-526.

\bibitem[SX]{SX} R.\,Song and L.\,Xie, Weak and strong well-posedness of critical and supercritical SDEs with singular coefficients, {\em J. Differential Equations} \textbf{362} (2023), 266-313 (arXiv:1806.09033).


\bibitem[V]{V} A.\,Yu.\,Veretennikov, Strong solutions and explicit formulas for solutions of stochastic integral equations, {\em Matematicheski Sbornik} (in
Russian), \textbf{111}(3) (1980), 434-452, English translation in {\em Math. USSR-Sbornik}, \textbf{39}(3) (1981), 387-403.



\bibitem[WLW]{WLW} J.\,Wei, G.\,Lv and J.-L.\,Wu, On weak solutions of stochastic differential equations
with sharp drift coefficients, arXiv:1711.05058.

\bibitem[W]{W} R.\,J.\,Williams, Brownian motion with polar drift, {\em Trans.\,Amer.\,Math.\,Soc.}, \textbf{292} (1985), 225-246.


\bibitem[XXZZh]{XXZZ} P.\,Xia, L.\,Xie, X.\,Zhang and G.\,Zhao, $L^q(L^p)$-theory of stochastic differential equations, {\em Stoch.\,Proc.\,Appl.} \textbf{130} (2020), 5188-5211 (arXiv:1908.01255).


\bibitem[XZ]{XZ} L.\,Xie and X.\,Zhang, Heat kernel estimates for critical fractional diffusion operators, {\em Studia Math.} \textbf{224} (2014), no. 3, 221-263. 

\bibitem[YZ]{YZ} S.\,Yang and T.\,Zhang, Strong existence and uniqueness of solutions of SDEs
with time dependent Kato class coefficients, arXiv:2010.11467.

\bibitem[Za1]{Z_supercritical} Q.\,S.\,Zhang, A strong regularity result for parabolic equations, {\em Comm.\,Math.\,Phys.} \textbf{244} (2004) 245-260.

\bibitem[Za2]{Z_Gaussian} Q.\,S.\,Zhang, Gaussian bounds for the fundamental solutions of $\nabla (A \nabla u) + B \nabla u -u_t=0$, {\em Manuscripta Math.} 93 (1997), 381-390.


\bibitem[Z1]{Z} X.\,Zhang, Stochastic homeomorphism flows of SDEs with singular drifts and Sobolev diffusion coefficients, \newblock {\em Electr.\,J.\,Prob.}, \textbf{16} (2011), 1096-1116.


\bibitem[Z2]{Z0} X.\,Zhang, Strong solutions of SDEs with singular drift and Sobolev diffusion
coefficients, {\em Stoch.\,Proc. Appl.}, \textbf{115}(11) (2005), 1805-1818.

\bibitem[Z3]{Z4} X.\,Zhang, Stochastic homeomorphism flows of SDEs with singular drifts and Sobolev diffusion coefficients, {\em Electron. J. Probab.} \textbf{16} (2011), no. 38, 1096-1116.

\bibitem[Z4]{Z2} X.\,Zhang, Stochastic differential equations with Sobolev diffusion and singular
drift and applications, {\em Ann. Appl. Probab.}, \textbf{26}(5) (2016), 2697-2732.


\bibitem[Z5]{Z5} X.\,Zhang, Stochastic differential equations with Sobolev drifts and driven by $\alpha$-stable processes. {\em Ann.\,Inst.\,Henri Poincar\'{e} (B) Probab. Stat.} \textbf{49}(4) (2013), 1057-1079.

\bibitem[ZZh]{ZZ} X.\,Zhang, G.\,Zhao, Stochastic Lagrangian path for Leray solutions of $3D$ Naiver-Stokes equations, {\em  Comm.\,Math.\,Phys.}, \textbf{381}(2) (2021), 491-525.

\bibitem[ZZh2]{ZZ2} X.\,Zhang and G.\,Zhao, Heat kernel and ergodicity of SDEs with distributional drifts, arXiv:1710:10537.

\bibitem[Zh]{Zh} G.\,Zhao, Stochastic Lagrangian flows for SDEs with rough coefficients, arXiv:1911.05562.

\bibitem[Zh2]{Zh2} G.\,Zhao, Weak uniqueness for SDEs driven by supercritical stable processes with H\"{o}lder drifts,
{\em Proc. Amer. Math. Soc.} \textbf{147} (2019), 849-860.

\bibitem[Zv]{Zv} A.\,K.\,Zvonkin, A transformation of the phase space of a diffusion process that
removes the drift, {\em Math. USSR Sbornik} \textbf{22} (1974), 129-149.



\end{thebibliography}
\end{document}